\documentclass[doublespacing]{utdthesis}
\usepackage{microtype}
\usepackage{url}
\usepackage[numbers]{natbib}
\let\cite=\citep
\usepackage[pagewise]{lineno}

\usepackage{amsmath}
\usepackage{tabularx}
\usepackage{enumitem}
\usepackage[misc]{ifsym} 
\usepackage{mathrsfs}
\usepackage[colorlinks=true]{hyperref}
\hypersetup{urlcolor=blue, citecolor=red}
\usepackage{float}
\allowdisplaybreaks
\usepackage[T1]{fontenc}

\usepackage{amsmath,amsfonts,amssymb,amscd,amsthm,bm}
\usepackage{cases}
\usepackage{latexsym}
\usepackage{pstricks}
\usepackage{amssymb}
\usepackage{amsmath}
\usepackage{amsfonts}
\usepackage{amsthm,array}
\usepackage{array}
\usepackage{epsfig}
\usepackage{verbatim}

\usepackage{tikz}
\usepackage{tikz-cd}
\usepackage{color}
\usetikzlibrary{positioning}
\usetikzlibrary{arrows.meta}
\usepackage{xcolor}

\usepackage{amsmath}
\usepackage{mathrsfs}
\usepackage{float}
\usepackage[colorlinks=true]{hyperref}
\hypersetup{urlcolor=blue, citecolor=red}
\usepackage{xcolor}
\usepackage{cases}
\usepackage{pstricks}
\usepackage{graphicx}
\usepackage{listings}

\usepackage{subcaption}

\usepackage{placeins}

\usetikzlibrary{external}
\graphicspath{{./figures/}}

\usepackage{mathrsfs}
\usepackage{graphicx,pstricks}
\usepackage{epsfig} 
\usepackage{multirow}
\usepackage{extarrows}

\newcommand\old[1]{}
\newcommand{\pend}{\hfill \thicklines \framebox(6.6,6.6)[l]{}}
\newcommand{\be}{\begin{eqnarray}}
	\newcommand{\ba}{\begin{array}}
		\newcommand{\ben}{\begin{eqnarray*}}
			\newcommand{\ee}{\end{eqnarray}}
		\newcommand{\ea}{\end{array}}
	\newcommand{\een}{\end{eqnarray*}}

\def\beq{\begin{equation}}
	\def\eeq{\end{equation}}

\numberwithin{equation}{section}

\renewenvironment{proof}{\noindent {\bf  Proof.} \rm}{\pend}

\newrgbcolor{violet}{.6 .1 .8}
  \newrgbcolor{lightyellow}{1 1 .8}
  \newrgbcolor{lightblue}{.80 1 1}
  \newrgbcolor{mygreen}{0 .66 .05}
  \definecolor{mygreen}{rgb}{0,.66,.05}
  \definecolor{lightyellow}{rgb}{1,1,.80}
  \newrgbcolor{orange}{1 .6 0}
  \newrgbcolor{GreenYellow}{.85 1 .31}
  \newrgbcolor{Yellow}{1  1  0}
  \newrgbcolor{Goldenrod}{1  .90  .16}
  \newrgbcolor{Dandelion}{1  .71  .16}
  \newrgbcolor{Apricot}{1  .68  .48}
  \newrgbcolor{Peach}{1  .50  .30}
  \newrgbcolor{Melon}{1  .54  .50}
  \newrgbcolor{YellowOrange}{1  .58  0}
  \newrgbcolor{Orange}{1  .39  .13}
  \newrgbcolor{BurntOrange}{1  .49  0}
  \newrgbcolor{Bittersweet}{1.  .4300  .24}
  \newrgbcolor{RedOrange}{1  .23  .13}
  \newrgbcolor{Mahogany}{1.  .4475  .4345}
  \newrgbcolor{Maroon}{1.  .4084  .5376}
  \newrgbcolor{BrickRed}{1.  .3592  .3232}
  \newrgbcolor{Red}{1  0  0}
  \newrgbcolor{OrangeRed}{1  0  .50}
  \newrgbcolor{RubineRed}{1  0  .87}
  \newrgbcolor{WildStrawberry}{1  .04  .61}
  \newrgbcolor{CarnationPink}{1  .37  1}
  \newrgbcolor{Salmon}{1  .47  .62}
  \newrgbcolor{Magenta}{1  0  1}
  \newrgbcolor{VioletRed}{1  .19  1}
  \newrgbcolor{Rhodamine}{1  .18  1}
  \newrgbcolor{Mulberry}{.6668  .1180  1.}
  \newrgbcolor{RedViolet}{.9538  .4060  1.}
  \newrgbcolor{Fuchsia}{.5676  .1628  1.}
  \newrgbcolor{Lavender}{1  .52  1}
  \newrgbcolor{Thistle}{.88  .41  1}
  \newrgbcolor{Orchid}{.68  .36  1}
  \newrgbcolor{DarkOrchid}{.60  .20  .80}
  \newrgbcolor{Purple}{.55  .14  1}
  \newrgbcolor{Plum}{.50  0  1}
  \newrgbcolor{Violet}{.98 .15 .95}
  \newrgbcolor{RoyalPurple}{.25  .10  1}
  \newrgbcolor{BlueViolet}{.84  .38  .98}
  \newrgbcolor{Periwinkle}{.43  .45  1}
  \newrgbcolor{CadetBlue}{.38  .43  .77}
  \newrgbcolor{CornflowerBlue}{.35  .87  1}
  \newrgbcolor{MidnightBlue}{.4414  .9259  1.}
  \newrgbcolor{NavyBlue}{.06  .46  1}
  \newrgbcolor{RoyalBlue}{0  .50  1}
  \newrgbcolor{Blue}{0  0  1}
  \newrgbcolor{Cerulean}{.06  .89  1}
  \newrgbcolor{Cyan}{0  1  1}
  \newrgbcolor{ProcessBlue}{.04  1  1}
  \newrgbcolor{SkyBlue}{.38  1  .88}
  \newrgbcolor{Turquoise}{.15  1  .80}
  \newrgbcolor{TealBlue}{.1572  1.  .6668}
  \newrgbcolor{Aquamarine}{.18  1  .70}
  \newrgbcolor{BlueGreen}{.15  1  .67}
  \newrgbcolor{Emerald}{0  1  .50}
  \newrgbcolor{JungleGreen}{.01  1  .48}
  \newrgbcolor{SeaGreen}{.31  1  .50}
  \newrgbcolor{Green}{0  1  0}
  \newrgbcolor{ForestGreen}{.1992  1.  .2256}
  \newrgbcolor{PineGreen}{.3100  1.  .5575}
  \newrgbcolor{LimeGreen}{.50  1  0}
  \newrgbcolor{YellowGreen}{.56  1  .26}
  \newrgbcolor{SpringGreen}{.74  1  .24}
  \newrgbcolor{OliveGreen}{.6160  1.  .4300}
  \newrgbcolor{RawSienna}{.53  .28  .16}
  \newrgbcolor{Sepia}{1.  .7510  .70}
  \newrgbcolor{Brown}{.41  .25  .18}
  \newrgbcolor{TAN}{.86  .58  .44}
  \newrgbcolor{Gray}{1.  1.  1.}
  \newrgbcolor{Black}{1  1  1}
  \newrgbcolor{White}{1  1  1}

\def\itemc{\itemindent=10pt\labelsep=8pt\labelwidth10pt\itemsep=4pt}
\def\id{{\text{\rm Id}}}

\def\re{\text{\rm Re\,}}

\def\vs{\vskip.2cm}

\def\vp{\varphi}
\def\ve{\varepsilon}
\def\id{\text{\rm Id\,}}
\def\sign{\text{\rm sign\,}}
\def\dim{\text{\rm dim\,}}

\def\ker{\text{\rm Ker\,}}
\def\br{\mathbb R}

\def\bn{{\mathbb{N}}}

\def\br{{\mathbb{R}}}

\def\bz{{\mathbb{Z}}}
\def\be{{\mathbb{E}}}

\def\br{\mathbb R}

\def\id{{\text{\rm Id}}}

\def\re{\text{\rm Re\,}}

\def\vs{\vskip.3cm}

\def\vp{\varphi}
\def\ve{\varepsilon}
\def\id{\text{\rm Id\,}}
\def\sign{\text{\rm sign\,}}
\def\dim{\text{\rm dim\,}}

\def\ker{\text{\rm Ker\,}}

\def\gdeg{G\text{\rm -deg\,}}

\newcommand\cV{\ensuremath{\mathcal V}}
\newcommand\cW{\ensuremath{\mathcal W}}

\newcommand\bbN{\ensuremath{\mathbb N}}

\newcommand\bbR{\ensuremath{\mathbb R}}

\newcommand\bbZ{\ensuremath{\mathbb Z}}

  \newcommand{\abs}[1]{\left\lvert#1\right\rvert}
  \newcommand{\norm}[1]{\left\lVert#1\right\rVert}

\newtheorem{theorem}{Theorem}[section]
\newtheorem{proposition}{Proposition}[section]
\newtheorem{lemma}{Lemma}[section]
\newtheorem{corollary}{Corollary}[section]

\newtheorem{definition}{Definition}[section]

\newtheorem{remark}{Remark}[section]

\newtheorem{remark-definition}{Remark and Definition}[section]
\newtheorem{rem-not}{Remark and Notation}[section]

\lstdefinelanguage{GAP}{
  basicstyle=\ttfamily,
  keywords={true, false, function, return, fail, if, in, while, do, od, else, elif, fi, break, continue},
  keywordstyle=\color{blue}\bfseries,
  otherkeywords={>, <, ==},
  identifierstyle=\color{black},
  sensitive=True,
  comment=[l]{\#},
  commentstyle=\color{cyan},
  stringstyle=\color{red},
  morestring=[b]',
  morestring=[b]"
}

\def\id{{\text{\rm Id}}}
\def\sign{\text{\rm sign\,}}
\def\dim{\text{\rm dim\,}}
\def\ker{\text{\rm Ker\,}}

\def\re{\text{\rm Re\,}}

\def\gdeg{G\text{\rm -deg\,}}
\def\sdeg{S^1\text{\rm -deg\,}}

\author{Casey Maikalani Crane}
\title{Consensus-Breaking Global Hopf Bifurcation in Memory-Based Multi-Agent Systems}
\thesistype{Dissertation}  % or "Thesis"
\degreefull{Doctor of Philosophy}
\degreeabbr{PhD}
\subject{Mathematics}
\graduationmonth{August}
\graduationyear{2026}
\prevdegrees{BS, MS} % comma-separated list of PREVIOUS degrees

\committeemember*{Wieslaw Krawcewicz}
\begin{document}

\frontmatter

% \signaturepage

% \copyrightpage{2026} % optional

% \begin{dedication} % optional
% This thesis class file \\
% is dedicated to my students, \\
% who suffered without a proper one \\
% until the present time.
% \end{dedication}

\maketitle
\copyrightpage{2026} % optional
\begin{acks}{June 2026} % date when thesis first submitted to committee
I would like to thank my advisor, Professor Wieslaw Krawcewicz, for his guidance, patience, and support throughout this research, as well as my committee members Dmitry Rachinskiy, Mieczyslaw Dabkowski, and Carlos Garcia Azpeitia, for their useful feedback and comments. I am also grateful to Chaoquan Chen and Travis Hensley, my collaborators on the work which became the foundation of Chapter \ref{chapter:distributed}. I would also like to thank Huafeng Xiao, Ziad Ghanem, Jinghzou Liu, and Yameng Duan for enjoyable and productive collaborations across many other projects which helped inspire this work both directly and indirectly. 

I would also like to thank my partner Kyle for her continuous support, my mother, my father, my sister Hannah, along with Garrett, Greg, Savannah, Dylan, Aaron, and all my other friends and family who have sat politely and listened to me try to explain my ideas over all these years.
\end{acks}

\begin{abstract}
This dissertation provides the first systematic study of symmetric consensus-breaking bifurcation to periodic multiconsensus in multi-agent systems. It analyzes this for three classes of multi-agent systems based on three different types of memory, whose closed-loop dynamics equations form delay differential equations of retarded type, neutral type, and pseudoneutral type---a subclassification of retarded type equations introduced in this dissertation which bridges retarded and neutral type delay equations. Equivariant twisted degree is used to analyze the symmetric global Hopf bifurcation problem in these systems, i.e. bifurcation from a stable consensus to periodic multiconsensus. This shows how the effects of memory allow self-organizing agents to move beyond mere stationary consensus: They can experience disagreement, form factions with antagonistic or competitive relationships; or split up to form higher-order patterns of synchronization and cooperation beyond mere stationary consensus, mirroring long-term cooperation and divisions of labor among humans and other intelligent beings. Theoretical results for the global Hopf bifurcation and symmetric classification of periodic multiconsensus solutions across all three systems are provided, and numerical results are conducted to both validate and enhance the theoretical predictions by providing stability information on the branches which is not obtainable by the degree alone. These principles are demonstrated in three real-world applications: one involving the control of formations of UAVs, allowing them to maintain their overall spatial relationships while dancing in complex selectable oscillations; and two more in networked asset markets featuring different traders with different memory-based strategies, showing how similar mechanisms can be responsible for economic cycles of bubbles and crashes. Finally, we also numerically investigate resonant double Hopf bifurcations in the neutral delay system, showing strong evidence of a breakdown to chaos via the Ruelle-Takens-Newhouse scenario and the existence of riddled basins. Through the study of these three systems and their applications, this dissertation showcases a comprehensive toolset for studying consensus-breaking bifurcation in multi-agent systems, and illustrates the ways that equivariant degree theory can work hand-in-hand with numerical analysis to enhance each theory's conclusions and buttress their limitations.
\end{abstract}

\tableofcontents
\listoffigures % required if you have any figures
\listoftables % required if you have any tables

\mainmatter
% !TEX root = main_new.tex

\def\sdeg{S^1\text{\rm -deg\,}}

\chapter{Introduction}
Human beings did not gather together and collectively decide to build great cities, establish complex economies and social hierarchies, or organize themselves into tribes, corporations, and nations whose collective goals supersede those of their constituents. When flocks of birds dance in the sky in dazzling, folding, interpenetrating patterns, they are not following a choreography rehearsed in advance or the lead of a single individual. From complex organisms such as humans and birds, to the most atomic models of decision such as neurons or cellular automata, to new frontiers in networked artificial intelligences, to autonomous self-organization and coordination in robots and unmanned vehicles, there is a unifying paradigm in the incredible properties of emergent synchronization and self-guided pattern formation exhibited across all these systems. 
\vs
In each case, individuals make decisions based on independent strategies and processes, informed only by localized information on their own state and their surroundings. These individuals have no knowledge of the global state of the vast systems to which they belong, and no ability to meaningfully change that system. Nor do they have any ability to agree or conspire on any unified course of action which could effect large-scale change. In spite of this, by pursuing their own independent goals within their own limited scope and means, and mutually influencing only themselves and their neighbors, they are able to create incredibly intricate structures and patterns on a global scale (quite literally in the case of human beings), transform their environments, and send cascading ripples of influence across networks of massive size. 
\vs
Systems like these can be described as multi-agent systems. Such systems consist of individual agents who measure their own state and the states of their neighbors, defined according to some interaction topology, and use a protocol to transform this state information into action or decision which can change their state or affect their local environment. These systems can be very robust and resilient because they are often leaderless, decentralized, or distributed.
\vs
Researchers are often interested in situations where agents self-organize and synchronize to have the same state, a phenomenon known as consensus. A common question in multi-agent systems is how to design a protocol which achieves consensus, and under what parameter conditions this consensus will be stable. This question naturally begs another: what happens beyond these conditions---when consensus breaks down? This dissertation answers this question in a systematic fashion for symmetric multi-agent systems.

This breakdown of consensus and what happens afterward is both highly relevant and all too familiar to humans. Our individual and universal experiences, along with history both ancient and recent, offer some suggestions as to what might happen. Agents who were once in cooperation and agreement will separate and form factions, lines will be drawn, and different groups of agents will have very different ideas about what their goals are and how to achieve them. This can be envisaged as conflict and competition where there was once cooperation and synergy. However, there are other less depressing interpretations. Two dancers might make symmetrically opposing steps and move in opposite directions, but this is not conflict---it is simply a higher form of organization, one which may be more effective and elegant for some applications. In this case, the breakdown of consensus is really just the emergence of higher order cooperation. 
\vs
As this dissertation will show, all of these outcomes are possible. The clustering of agents into different groups, where each cluster achieves consensus internally but differs from other clusters, is called multiconsensus.  We will consider an application in formation control for symmetric formations of UAVs, where the breakdown of common consensus into periodic multiconsensus states corresponds to the emergence of tunable and selectable symmetric patterns of oscillation around their formation set points. In this context, much like the dancers, this can be viewed as a higher order of cooperation and synchronization, useful for sensing, reconnaisance, evasion, and other distributed tasks. We will also study a networked system of asset markets, where each market is populated by two types of traders following different trading strategies. In this case, the breakdown of consensus into periodic multiconsensus states corresponds to boom-bust cycles, where one market experiences a speculative bubble while another seemingly unrelated market experiences a crash.
\vs
In particular, we are interested in studying multi-agent systems with three types of memory: a continuous moving average of past states, over some finite time horizon up to the present; the instantaneous rate of change of a continuous moving average over past states; and a memory of an agent's own past decisions. As our applications illustrate, these memory terms are very natural modeling choices for a wide variety of phenomena, and also very advantageous design choices for designing agent protocols with certain goals in mind. For example, the continuous averaging over past states naturally smooths out noisy input or state data, which is useful for agents who must make decisions in environments where their own state information or history may be noisy or unreliable. Our application in UAV formation control illustrates their design utility, and our economic applications illustrate their usefulness in modeling the decision-making processes of highly intelligent actors. Through these applications, we address the viewpoint of protocol design as well as the viewpoint of modeling emergent phenomena in nature.
\vs
The closed-loop dynamics of multi-agent systems with these types of memory are governed by three classes of delay differential equations with distributed delays. As such, the results in this dissertation also pertain to such systems of delay equations abstractly, and do not depend on any particular feature of multi-agent systems to be used or interpreted. For this reason, we will take a ``bilingual'' approach to presenting and proving this dissertation's results. On one hand. we will show consensus-breaking bifurcation to symmetrically classified orbits of periodic multiconsensus states in three novel classes of multi-agent systems with three types of continuous memory. On the other hand, we will show symmetric global Hopf bifurcation of non-constant periodic solutions in three classes of delay differential equations featuring distributed delays along with pseudoneutral\footnote{This is a new term describing a certain class of retarded delay equations which are closely related to neutral delay equations, which we define formally in Section \ref{sec:prelim:pseudoneutral}} and neutral delays.
\vs
Accordingly, this work contains several aspects of novelty which may be of interest to mathematicians across several different fields. It systematically studies the relatively unexplored frontier of consensus-breaking bifurcation, offers a new framework to study periodic multiconsensus scenarios that emerge post-bifurcation, and comprehensively classifies them according to their spatio-temporal symmetries. From a dynamical systems point of view, it studies the global bifurcation problem of non-constant periodic solutions in three classes of symmetric systems of delay equations with three types of delays. Finally, it shows how equivariant degree techniques can be used in conjunction with numerical simulations to answer questions of great interest in applied mathematics, such as showing the stability of the symmetric branches of solutions whose existence and symmetries are topologically guaranteed by the degree.
\vs
We also intend to show that multi-agent systems are a very natural arena for the application of equivariant degree techniques, and hope to inspire further work in this area. When multi-agent systems have homogeneous protocols and symmetric interaction topologies, the closed-loop dynamics equation of the system naturally form highly symmetrically coupled systems of differential equations. Bifurcation analysis in such systems using classical techniques is extremely challenging, since symmetries force high multiplicities of eigenvalues, and high-dimensional center manifolds complicate Lyapunov-Schmidt methods. Equivariant degree methods thrive in such environments, where these high multiplicities of critical eigenvalues are not an obstacle but rather an asset which provides an abundance of symmetry information on bifurcating solutions. This analysis lets us show the fascinating and beautiful range of patterns, self-organization, and synchrony which are possible in multi-agent systems.
\vs
The first two chapters of this dissertation are intended to provide historical background and technical preliminaries on the topics of delay differential equations, multi-agent systems, and equivariant degree theory. Chapter \ref{chapter:background} is intended both to provide a brief overview and historical summary, as well as to contextualize this work as a natural development of the current state-of-the-art across all three fields. Chapter \ref{chapter:preliminaries} provides technical background on these three fields, intended to make this dissertation as self-contained as is reasonably possible, especially with respect to the equivariant degree techniques used, which are built up from the local Brouwer degree to the equivariant twisted Nussbaum-Sadovskii degree used to show Hopf bifurcation in symmetrically coupled neutral equations. 
\vs
The following chapters provide our main results across three classes of multi-agent systems which increase in complexity. All systems considered are homogeneous multi-agent systems featuring a symmetric interaction topology with state-dependent effective weights. All involve certain types of memory. For each system, its protocol depends on a continuous memory term which averages past states over a finite time horizon. We believe this type of memory, which can be viewed as a distributed delay, is a natural modeling choice for a variety of reasons. On one hand, it has natural smoothing properties and copes well with noise. On the other hand, most intelligent decision-making agents do not merely consider a single discrete point in their past, but rather view past events over a wide variety of time scales, and averaging represents a way that an intelligent being might turn a series of specific events into a general impression which is more useful for making decisions. We also consider systems where agents additionally base their decisions on the instantaneous rate of change of a continuous moving average, possibly transformed nonlinearly. This represents a sensitivity to trends, and approximates a derivative in a certain sense. These types of agents not only remember an average of previous states, but are also estimating how their state changed in the past. Finally, we consider systems where agents' protocols depend on their own derivative at a certain point in the past. Since all the systems studied here have first order dynamics, the derivative of an agent can be interpreted as its decision. Therefore, this memory of past derivatives, which we call ``momentum memory,'' represents an exact memory of its past decisions at some fixed point in the past.
\vs
In Chapter \ref{chapter:distributed}, we show the local asymptotic stability of a trivial consensus equilibrium in a multi-agent system of homogeneous agents with symmetric interaction topology and continuous memory, whose closed-loop dynamics are given by a highly symmetrically coupled delay differential equation with distributed delays. We then show the local bifurcation and global continuation of symmetrically classified branches of periodic multiconsensus solutions after consensus is broken. We offer an application of such a system as a formation control scheme for an octahedral formation of six UAVs. We show how the abstract results can be applied in this setting to create 27 distinct conjugacy classes of three-dimensional synchronized patterns of motion around the UAVs' formation set points. These ``dances'' are selectable and tunable in amplitude by varying a single parameter, and can be quieted back to the stationary formation state by adjusting the parameter back below the critical threshold. We confirm the predictions of the theoretical results numerically, showing that for each spatial axis of motion, there are three stable branches of solutions which can be selected by forcing an appropriate perturbation of the UAVs in the direction of the subspace fixed by the desired orbit type. In this context, consensus-breaking bifurcation can be viewed as a controllable and useful form of higher-order coordination and synchronization, useful for imaging, surveillance, evasion, and many other purposes.
\vs
In Chapter \ref{chapter:pseudo} we consider a multi-agent system of homogeneous agents with a symmetric interaction topology whose protocols depend not only on a continuous memory of past states, but also on the instantaneous rate of change of a nonlinearly transformed memory of past states, which we term ``trend memory''. In this case, the closed-loop dynamics form a delay equation of retarded type, but which is closely related in both form and properties to delay equations of neutral type. We again show that a trivial consensus exists and conditions for its local asymptotic stability, and show the local bifurcation, global continuation, and symmetric classification of non-constant periodic multiconsensus solutions bifurcating from the trivial consensus branch, using the framework of symmetric Hopf bifurcation.
\vs
Our application in this chapter relates this system to the popular economic modeling framework of heterogeneous agent models, which study the evolution of markets filled with heterogeneous populations of traders following different trading strategies. In our model, each asset market is an agent (in the multi-agent system sense) whose ``goal'' is price discovery, and whose protocol represents the competing influences of fundamentalist traders (represented by the continuous memory term) and chartists (represented by the trend memory term). In this sense, each agent in our multi-agent system is a separate but identical heterogeneous agent model. Our multi-agent system then comprises eight coupled markets, where each market has identical traders and is characterized by three binary properties, such as sector, size, and geographic location. Markets which share two properties and differ in one property are assumed to be more closely linked in the ambient economy and thus have a higher weight of coupling. Markets that share one property and differ by two have a medium weight of coupling, and markets that differ in all three properties have the lowest weight of coupling.
\vs
Here, we show that at a critical parameter threshold of the relative aggressiveness of fundamentalist traders (or the relative size of their population), periodic oscillations of overvaluation and undervaluation emerge, which can be interpreted as periodic cycles of pricing bubbles followed by crashes. We investigate these patterns numerically, and use these numerical results to infer the stability of the periodic multiconsensus branch. We also observe the emergence of arbitrages between markets, where bubbles in one market correspond with crashes in another. Interestingly, it is the markets which share the most properties (and thus are most strongly coupled) that experience these anti-phase relationships in their bubble-crash cycles. The formation of these cycles can be understood as a consequence of the ``forgetfulness'' of the memory of the traders, who look at average deviations of asset prices from their underlying fundamental value over a finite time horizon. The length of this memory can also be seen to have a direct influence on the time scale of these economic cycles.
\vs
In Chapter \ref{chapter:neutral} we consider a homogeneous multi-agent system where protocols depend on both continuous memory and momentum memory---that is, a memory of its own past decisions. The closed-loop dynamics of this system form a neutral functional differential equation with distributed retarded delays. We likewise show the existence and local stability of a trivial equilibrium in this system, along with the local bifurcation, global continuation, and symmetric classification of non-constant periodic multiconsensus solutions. We again use the framework of Hopf bifurcation to do this, building on a framework established in the previous chapter to handle the special functional-analytic difficulties related to neutral equations, which are detailed in Chapter \ref{chapter:preliminaries}. This system can be viewed in certain ways as a limiting case of the type studied in Chapter \ref{chapter:pseudo}, or, conversely, that system could be viewed as an approximation of the one studied in this chapter.
\vs
We consider a similar economic application in this chapter, but replace the population of trend-chasing chartist traders from the previous chapter with momentum traders. These can be viewed as representing algorithmic traders or high-frequency traders. Rather than analyzing trends over a time scale of days or weeks, these traders make decisions based on the instantaneous derivative of stock prices, delayed only by the intrinsic lags and processing delays in obtaining and acting on this information. As such, their delay term is on a far shorter time scale than that of their chartist counterparts in the previous application, on the order of minutes or hours, and this large proportional difference in time scales between the momentum traders and fundamentalists characterize this type of system. We again use numerical simulations to validate the conclusions of the theory, and investigate the stable symmetric periodic multiconsensus solution which emerges after the breakdown of the consensus state of market equilibrium. 
\vs
We also explore the parameter space of the non-bifurcation parameters for this system numerically, showing the existence of smooth curves corresponding to double Hopf bifurcation thresholds, and curves corresponding to resonances between limit frequencies. Where these curves coincide, we numerically investigate the double Hopf bifurcation scenario which exhibits quasiperiodic consensus and chaotic multiconsensus solutions. We give strong numerical evidence for the existence of chaos, suspected to occur via the Ruelle-Takens-Newhouse scenario. We show quasiperiodic and chaotic trajectories where periodic waves (whose frequencies are governed by the predicted limit frequencies of the Hopf bifurcation) are amplitude-modulated by quasiperiodic or chaotic signals on a much slower timescale. Our numerical basis of evidence for chaos is obtained by constructing first-return maps comparing peak-to-peak amplitudes of these waves which demonstrate fractal structure (strange attractors). We also extract the chaotic envelope function via numerical Hilbert transform and performing Fourier analysis, showing a broadband power spectrum characteristic of chaos, and we numerically compute the largest Lyapunov exponent of both the signal itself and the envelope function, showing the positive largest Lyapunov exponents also characteristic of chaos. We also note that the boundaries of the basins of attraction of the chaotic attractors appear to be intermingled with those of quasiperiodic and steady-state solutions at every scale. This, combined with simulated trajectories showing \emph{chaotic saddles}, where trajectories obey chaotic dynamics for extremely long timescales before returning to non-chaotic dynamics, are strongly suggestive of the presence of riddled basins, in the sense of \cite{alexander1992riddled}, where basins of attraction of chaotic and nonchaotic attractors are intermingled in a measure-theoretic sense at every scale.
\vs
In Chapter \ref{chapter:conclusion}, we summarize the conclusions of this dissertation and discuss several promising directions of future research. We also provide appendices for supplemental information on certain topics. Specifically, Appendix \ref{appendix:notation} motivates and describes the notation used for twisted orbit types throughout this dissertation. Appendix \ref{appendix:software} details the software used in this dissertation, specifically the methods used for the numerical simulations, as well as certain problems and considerations necessary for numerically simulating neutral delay equations, and details on the steps taken to address these issues. Finally, Appendix \ref{appendix:supplementary-figs} provides some supplementary figures for the system studied in Chapter \ref{chapter:neutral}.

\chapter{Background and Motivation}\label{chapter:background}
\section{Introduction to delay differential equations}
The purpose of this section is to provide some preliminaries on delay differential equations (DDEs), their historical background and development, and the specific challenges posed by their analysis. Our goal is to give some motivation and reasoning as to why mathematicians have studied delay equations, both historically and today; to show why delay equations are paradoxically very natural and simplifying from the point of view of modeling real world phenomena, but very difficult and complicated from the point of view of analysis; and to show why equivariant degree methods are a natural fit for addressing bifurcation problems in these types of systems. 

We will begin with a brief historical overview of the development of delay equations. We will then formulate an illustrative and contrastive example of a relatively simple ODE, the famous Verhulst-Pearl or logistic equation, and show how the addition of a discrete delay term greatly enriches the dynamics of the model, while also adding fundamental challenges to its analysis. We will show how many of the problematic features of this delayed logistic equation (also known as the Hutchinson equation) are representative of the challenges posed by DDEs as a whole. We will also give a brief survey of methods often used to analyze such systems, and discuss the many unique advantages of equivariant topological methods of the type used in subsequent chapters. 
\subsection{Historical background}
Delay differential equations (DDEs) are a class of differential equations which depend explicitly on their own past state. Equations of this type have been studied for more than two and a half centuries, and their study has always been strongly motivated by applications. The earliest known appearance of a differential equation with a delay is often attributed to the Marquis de Condorcet in 1767, who encountered such an equation in the context of the astronomical three-body problem \cite{Keller2011}. Throughout the 19th century, similar equations appeared sporadically, typically arising from specific problems in geometry, number theory, and physics \cite{Hale2006History}. For example, the work of E. Schmitt in 1911 provides an early survey of the properties of linear equations with deviating arguments, indicating that the study of such systems had already been ongoing for over a century \cite{Hale2006History}.

The study of these problems was initially sporadic and motivated by specific applied problems. Engineers and theorists of the early 20th century were well aware of the impact of ``hereditary'' effects and memory, but it was not yet understood that these equations, then termed ``equations with deviating arguments'', were not merely particularly odd or difficult ODEs, but actually constituted a new and different class of differential equation altogether. 

This revolution began with the work of Myshkis in 1949 and 1950, who developed the first systematic framework for analyzing delay equations of retarded type. This systematization continued with the work of Krasovskii in 1959. Krasovskii's key realization was that the state of a DDE is a function segment in an infinite dimensional space. This led to the development of the wider theory of functional differential equations (FDEs), which was cemented in 1977 with the publication of Hale's \textit{Theory of Functional Differential Equations} \cite{Hale1977}. These early insights allowed for a unification and systematic analysis of DDEs for the first time, replacing the various tricks and ad-hoc techniques developed in the early 20th century, and catalyzing a flurry of theoretical development and results. 

Much of this development was spurred by practical problems in control theory, a field which was also coalescing in the mid-20th century. The problem of time lags in feedback loops in control systems, as well as designing effective controls for systems with large intrinsic time lags, directly led to the formulation of many delay equations, which in turn necessitated the development of the functional-analytic framework.

\subsection{A motivating example: the Hutchinson equation}
Here we will show concretely how the addition of a discrete delay to a relatively simple first order ODE transforms it into a much more complicated problem with much more interesting dynamics.
Consider the Verhulst-Pearl equation, a first order ODE designed to model the growth of a population with a limited carrying capacity:
\begin{equation}\label{eq:intro:dde:logistic}
\frac{dN}{dt} = rN\left(1-\frac{N}{K}\right),
\end{equation}
where $N(t)$ represents the population at time $t$, and $r>0$ and $K>0$ are constants specifying the growth rate and carrying capacity, respectively. This forms a Bernoulli type equation with the well-known explicit solution
\[
N(t) = \frac{KN_0e^{rt}}{K+N_0(e^{rt}-1)},
\]
where $N_0 := N(0)$. Since this is a first order ODE, one can conclude---even in the absence of an explicit solution---that it cannot admit any periodic or oscillatory solutions, and the relatively simple dynamics at each of its two equilibria can be easily deduced from its linearization. In 1948, G.E. Hutchinson \cite{Hutchinson1948} introduced his famous modification of the Verhulst-Pearl equation by adding a discrete delay term, motivated by the obvious fact that in most living organisms of sufficient complexity, there is a delay between birth and reproductive maturity. This leads to the Hutchinson equation:
\begin{equation}\label{eq:intro:dde:hutchinson}
\frac{dN}{dt} = rN(t)\left(1-\frac{N(t-\tau)}{K}\right).
\end{equation}
In contrast to the Verhulst-Pearl equation, this equation has no closed-form solution. Unlike \eqref{eq:intro:dde:logistic}, which is a first order ODE of Bernoulli type, the Hutchinson equation is a delay differential equation of retarded type, and thus belongs to the broader class of FDEs. In order to specify its initial conditions, it is now necessary to provide a state history in the form of a continuous function $N_t : [-\tau,0]\to \mathbb R$. Therefore, the phase space of the Hutchinson equation is the infinite-dimensional Banach space $C([-\tau,0];\mathbb R)$.

As one might reasonably expect, this also means that \eqref{eq:intro:dde:hutchinson} can exhibit more interesting dynamics and many more types of solution than its first order ODE counterpart. In particular, it exhibits the same types of stable periodic solutions that are often seen in real-world population data. 

To show this fact, and to show the impact of the delay on analysis of this system, let us perform a short parallel study of the dynamics of the two systems near the carrying capacity equilibrium. We first note that $N(t) \equiv K$ is an equilibrium for both equations. Setting $\tilde N := N-K$ and taking the linearization of both systems around this equilibrium, for the Verhulst-Pearl equation, we obtain
\[
\frac{d\tilde N}{dt} = -r\tilde N(t),
\]
which has the corresponding characteristic equation
\[
\lambda = -r.
\]
This indicates a single eigenvalue which is negative for all $r>0$, and therefore the $N(t) \equiv K$ solution is asymptotically stable. Performing the same process for Hutchinson's equation, we obtain the linearization
\[
\frac{d\tilde N}{dt} = -r\tilde N(t-\tau),
\]
which has the characteristic equation
\[
\lambda +re^{-\lambda\tau}=0.
\]
In contrast to the characteristic equation of the \eqref{eq:intro:dde:logistic} linearized around the same equilibrium, this is a transcendental equation admitting infinitely many complex roots, and the $N\equiv K$ solution is asymptotically stable if and only if all of these roots have negative real part. Characteristic equations of this type are termed \textit{characteristic quasipolynomials}, and here, as in subsequent chapters, we will denote such equations as
\[
P(\lambda):= \lambda +re^{-\lambda\tau}.
\]
We would like to make a few brief observations about $P(\lambda)$ and its roots which are very typical of the characteristic equations of many DDEs, including all those which will be studied in detail in subsequent chapters:
\begin{itemize}
    \item $P(\lambda)$ is a transcendental equation which has no closed-form solution.
    \item $P(\lambda)$ has infinitely many roots which are generically complex.
    \item The roots of $P(\lambda)$ can be viewed as depending continuously on certain parameters.
    \item The presence of complex roots---which could be purely imaginary for certain critical parameter choices---suggests the possibility of periodic solutions.
\end{itemize}
These factors make Hopf bifurcation a very natural framework for showing the existence of periodic solutions---at least, small-amplitude periodic solutions near an equilibrium. Choose a bifurcation parameter $\alpha \in \{\tau, r\}$ and denote by $P_\alpha(\lambda)$ the parametrized characteristic quasipolynomial. The classical local Hopf bifurcation theorem for DDEs as formulated in \cite{Hale1977} requires that a single simple purely imaginary root $\lambda = i\beta_0$ crosses the imaginary axis transversally (i.e. with non-zero velocity, $\text{Re}(\tfrac{d\lambda}{d\alpha})|_{\alpha=\alpha_0} \neq 0$), and that no other root of $P_{\alpha_0}(\lambda)$ is an integer multiple of $i\beta_0$.

Let us attempt to apply this theorem to \eqref{eq:intro:dde:hutchinson}, setting $\alpha = \tau$. We will then discuss some issues with this approach when applied to other, more general DDEs. Putting $\lambda := u+iv$ and separating real and imaginary parts yields the system
\begin{align*}
    u + re^{-u\tau}\cos(v\tau) &= 0 \\
    v - re^{-u\tau}\sin(v\tau) &= 0
\end{align*}
By inspection, we can see that $\tau=0$ reduces to the characteristic equation of the Verhulst-Pearl equation as expected. Some elementary algebra on these equations shows that purely imaginary roots occur when $r\tau = \pi/2$. Noting that these roots occur when $v=r$ and performing some simple calculus, we can also see that the transversality condition is satisfied for all $r>0$:
\[
    \frac{du}{d\tau}\Big\vert_{\tau=\tau_0}=r^2>0.
\]
However, now it must be shown that $\lambda = ir$ is a simple root, and that the non-resonance condition is satisfied. Both of these facts are easy to show for the Hutchinson equation, but for more general DDEs this analysis is often far from trivial. The non-resonance condition can be particularly difficult, as there are no general algorithms for showing that only one pair of conjugate roots is on the imaginary axis. Moreover, this only guarantees the existence of a small amplitude periodic solution near the equilibrium, and in practical terms can often only be used for the first root to cross the imaginary axis. If we are considering a coupled system, or an equation with other spatial symmetries, then this presents an additional obstacle to analyzing Hopf bifurcation, as symmetries often force non-simplicity of eigenvalues. 

Even if these issues are resolved, gaining more information about the global fate of these branches of periodic solutions requires us to look beyond the linear picture, and here we find a very difficult barrier. The standard route to obtaining a more detailed local description—determining the stability of the bifurcating periodic orbit, its direction of branching, or the effects of symmetry—is to perform a reduction to a finite-dimensional center manifold or to apply Lyapunov--Schmidt reduction \cite{FariaMagalhaes1995, Hale1977}. While the theoretical foundations for these techniques are well-established for functional differential equations, their practical implementation for DDEs is often intractable.

The primary difficulties are threefold. First, the phase space is the infinite-dimensional Banach space \(C([-\tau,0], \mathbb{R}^n)\). Projecting onto the center eigenspace requires the computation of eigenfunctions and their adjoints, which involve solving transcendental equations and evaluating bilinear forms that are rarely expressible in closed form. Second, the normal form computations are combinatorially explosive. Even for scalar DDEs with a single delay, deriving the cubic normal form coefficients can be an arduous task, and the complexity grows exponentially with the dimension of the system and the number of delays. Third, when the system possesses symmetry---a central theme of this dissertation---the center manifold is higher-dimensional due to degenerate eigenvalues, and the normal form must be equivariant with respect to the symmetry group. Deriving such equivariant normal forms for DDEs has been successfully accomplished for only a small number of canonical examples \cite{2013guo}.

Finally, even when such a reduction is feasible, the resulting dynamical description is inherently local. The normal form provides an accurate portrait of the flow only in a sufficiently small neighborhood of the bifurcation point. It offers no insight into the global behavior of the bifurcating branch of periodic solutions---whether it terminates, becomes unbounded, or connects to other bifurcating branches. We will show how equivariant degree theory provides an elegant solution to this problem in section \ref{sec:intro:equideg}.

\subsection{Distributed delays}\label{sec:intro:dist-delay}
In many systems, the influence of the past on the present does not act through a single, discrete time lag, but is spread over an interval or a distribution of delays. In neural field equations, signal propagation delays depend on the distance between neurons \cite{Hutt2005}. In epidemiological models, the latent period between exposure and infection is often variable \cite{Wearing2005}. In models of network traffic and congestion, the actual latency of packets is stochastic and can be modeled by a probability distribution \cite{Hespanha2007,wang2021global}. The proper mechanism for incorporating these effects into delay equation models is distributed delay. A distributed delay is a delay term of the form
\[
\int_0^\infty h(s)x(t-s)ds,
\]
where $h(s)$ is called the \emph{delay kernel} and is often chosen so as to satisfy
\[
\int_0^\infty h(s)ds =1.
\]
As this suggests, distributed delays are named by analogy to probability distributions. In subsequent chapters, we will consider distributed delays over a finite time horizon, often of the form
\[
\int_0^\tau x(t-s)ds.
\]
This choice of a delay kernel (implicitly a step function of length $\tau$) is often called a ``boxcar kernel'' in literature. 
\vs
On the other hand, when the delay kernel takes the form of an Erlang distribution, there is a method known as the ``linear chain trick'' which can reduce the system to a system of ODEs by using the Laplace transform. This is often convenient when the statistical processes influencing the delay can be assumed to follow a memoryless distribution. However, this trick cannot be used for general memory kernels or the boxcar kernels of the type used throughout this dissertation, and so we are forced to use more technical methods to show the existence of periodic solutions.
\vs
There are two major interpretations of distributed delay terms, depending on the type of phenomenon being modeled. The first is that the distributed delay models the expected value of a variable delay modeled by a probability distribution in deterministic mean-field models. This is common in the types of applications mentioned above, and this motivates the normalization requirement often described in literature (but not enforced in this dissertation) that $\int_0^\infty h(s)ds=1$.
\vs
The second interpretation is that distributed delays, rather than averaging over a stochastic process, are reflecting a weighted and continuous dependence on past states. This type of continuous memory has been used to model reaction times in drivers \cite{Sipahi2015}, the dependence of economic markets on past levels of investment \cite{Yu2016,Tarasov2021}, and fading memory in neural models \cite{Sharma2013}. These latter interpretations often relax the normalization condition on the delay kernel, as it is no longer interpreted as a probability distribution. It is in this spirit that the distributed delays in this dissertation are formulated. 
\vs
Distributed delays also have a number of interesting properties. They tend to have stronger ``smoothing'' properties than discrete delays, in terms of how discontinuities propagate through higher derivatives. This is because the distributed delay operator acts by convolution. Recent research by Cassidy \cite{Cassidy2025} has shown that distributed delays with compact support can be modeled equivalently by a suitable quadrature of discrete delay terms. This gives theoretical support to the observations of Tavakoli and Longtin \cite{tavakoli2020multi} that increasing the number of discrete delays in several prototypical delay equations increases complexity up to a certain point, where complexity suddenly sharply decreases. This suggests that this ``complexity collapse'' could occur when a multidelay equation sufficiently well approximates a more well-behaved distributed delay equation. 
\vs
These results are also important for establishing error bounds on numerical simulations of distributed DDEs, since numerical simulations inherently require the discretization of distributed delay terms.

\subsection{Pseudoneutral equations}\label{sec:prelim:pseudoneutral}
We also delineate a special class of retarded type delay equations in this dissertation. Motivated by distributed delays, we call an equation \emph{pseudoneutral} if it takes the form
\begin{equation}\label{eq:prelim:pseudoneutral}
\frac{d}{dt}\left[x - \int_0^\tau \phi(s)x(t-s)ds\right] = F(x_t)
\end{equation}
where $\phi\in L^1[0,\tau]$ is an ordinary integrable function (i.e. it contains no Dirac-delta atoms). By the fundamental theorem of calculus, the distributed delay inside the derivative can be expressed as a finite difference of purely retarded type, i.e.
\[
\frac{d}{dt}\left[\int_0^\tau \phi(s)x(t-s)ds\right]=x(t)\phi(0) - x(t-\tau)\phi(\tau) - \int_0^\tau\phi'(s)x(t-s)ds.
\]
This definition is motivated by the relationship of such equations to neutral equations, which can be viewed in two complementary ways: On one hand, if we take a sequence of delay kernels $(\phi_n)_{n=1}^\infty$ which converges weakly to a Dirac delta funciton $\delta_{\tau_0}(s)$, then the limit of the corresponding sequence of delay equations will be the neutral equation
\[
\dot x = F(x_t) + \dot x(t-\tau_0).
\]
On the other hand, if we view a distributed delay as representing a mean field approximation based on the expected value of a discrete delay whose delay length follows a distribution, then the pseudoneutral equation \eqref{eq:prelim:pseudoneutral} is the purely retarded equation representing a mean field approximation of neutral equations whose neutral delay terms are distributed according to $\phi(s)$, where $\phi$ is a distribution with compact support. 
\vs
pseudoneutral equations with distributed retarded delays can thus be viewed as occupying a middle continuum between neutral equations with distributed retarded delays, and distributed delay equations of purely retarded type, both of which can be realized as limiting cases of the pseudoneutral family in the above manner. We also note that a pseudoneutral equation with distributed retarded delays will appear as a retarded equation containing both discrete and distributed delays.

\section{Introduction to multi-agent systems}
Multi-agent systems have been a topic of rapidly increasing for the past several decades. They are a natural framework for modeling cooperation, competition, flocking and swarming behaviors, and synchronization across a wide variety of domains, from biology and economics to artificial intelligence, differential game theory, and control theory. In this section, we will give some historical background and motivation on the development of this relatively young and highly promising field. For the technical definitions and preliminaries needed for this dissertation, see Section \ref{sec:prelim:mas}. 
\vs
One can intuitively view the study of multi-agent systems as a branch or subfield of distributed control theory for networked systems. If control theory can be thought of as a conductor leading an orchestra playing sheet music, multi-agent systems are more like a jazz ensemble---leaderless and self-organizing, but capable of producing equally rich and complex structures. 
\vs
Early work in this field was inspired by observations of self-organizing behaviors in the natural world. The synchronous motions and flocking behaviors of birds inspired Craig Reynold's 1986 development of the Boids artificial life algorithm \cite{reynolds1987boid}, which is often considered one of the earliest works in the field. This revolutionary model prefigured the development of multi-agent systems as an independent field of study concerned with distributed self-organizing behaviors in agents of all types. In 1995, Vicsek et al. \cite{vicsek1995novel} introduced a minimal physics model where particles align their velocities with neighbors and showed a phase transition from disordered motion to collective flocking behavior.
\vs
Near the turn of the millennium, in 2003, Jadbabaie, Lin, and Morse \cite{jadbabaie2003coordination} gave the first rigorous Lyapunov-based proof that the discrete-time Vicsek model converges to a state where all agents move in the same direction, provided the sequence of neighbor graphs is periodically connected. This was followed in 2004 by the foundational work of Olfati-Saber and Murray \cite{saber2004consensus}, which formalized the synchronization behaviors shown by Jadbabaie et al. in the Vicsek model into the notion of \emph{consensus}, introduced the now-canonical continuous-time consensus protocol (cf. Section \ref{sec:prelim:consensus}).
\vs
This established one of the central problems of multi-agent systems theory: creating protocols such that agents will eventually align their decisions over time. This formalizes the search for self-organizing behaviors and synchrony as a problem which can be expressed and understood in dynamical systems terms as the search for agent protocols (i.e. their decision-making process) such that the closed-loop dynamics admit stable solutions where all agents agree. This was followed by works of Ren and Beard in 2005 \cite{ren2005consensus}, which showed that consensus could be achieved even when agent interaction topologies change dynamically, provided the union of the graphs over bounded time intervals contains a spanning tree.
\vs
Time delays are ubiquitous in real networked systems. There may be time lags in the exchange of state information between agents, or internal delays in agents' processing of information to make decisions. Such transmission and processing delays are very natural in both artificial systems (e.g. modeling network latency) or in natural systems (e.g. modeling the time humans take to think and make decisions). Olfati-Saber and Murray \cite{saber2004consensus} already studied the destabilizing effect of constant communication delays, and gave a formula for the maximum allowable delays before consensus is lost. The following decade produced a large body of work aimed at finding delay margins for various agent dynamics, network topologies, and delay distributions (see \cite{munz2011consensus} and further references therein). The focus of virtually all of these studies is showing the stability of the consensus manifold, and establishing conditions on the delay terms guaranteeing this stability up to some critical margin. What happens \emph{after} the delay exceeds this critical margin, when the consensus equilibrium loses stability, has received far less attention.
\vs
Moreover, these works predominately view delays in discrete terms, and as representing delays in receiving or processing information. There are two other possibilities, however, which have very promising applications (several of which are explored in this dissertation). One is to extent from discrete delays to continuous delays, where a delay function in convolved with a delay kernel, often chosen as a probability distribution. These are called distributed delays, and are often used to represent the expected value of a variable delay which follows a certain probability distribution, taken in a mean field approximation. The other aspect is to interpret delays as memory. These two ideas especially flourish when taken together: intelligent beings base their current decisions on their experiences of the past, but we do not look backwards to a fixed, discrete distance from the present. Rather, such agents view a spectrum of moments of their past, possibly weighted by their immediacy or otherwise transformed, and use these accumulations of experiences to guide decision-making. In this context, a distributed delay represents an agent's continuous, weighted memory of the past. We believe this underexplored model can be used to realistically represent many types of behavior, both in artificial and natural systems, and we advance three types of continuous memory for agents in this dissertation.
\vs
In dynamical systems, when a stable equilibrium loses stability because one of its eigenvalues crosses the imaginary axis as a parameter moves through a critical value, this typically results in a Hopf bifurcation, and gives rise to a branch of periodic solution. Several papers have used this framework to study delayed multi-agent systems \cite{chen2014multiconsensus,xie2015second,wang2021global,Wang2024}. These papers, however, are concerned with using Hopf bifurcation as the threshold where a consensus loses stability. While some consider the direction and stability of the Hopf bifurcation, they do not investigate the resulting spatio-temporal oscillations which can emerge after the stable consensus breaks down, nor do they consider the global continuation of these oscillating branches. One reason for this is because in homogeneous multi-agent systems with highly symmetric interaction topologies, these high degrees of symmetry often cause the eigenvalues crossing the imaginary axis in Hopf bifurcation to be non-simple. This is a major degeneracy from the point of view of classical Hopf bifurcation theories, and using those theories it is very difficult to say much more about the patterns of oscillations which might emerge.
\vs
However, elsewhere in dynamical systems, equivariant degree theory has been used very successfully to show global Hopf bifurcation in many different dynamical systems, including systems with time delays. In the context of equivariant degree theory, these high levels of symmetry and non-simplicity of critical eigenvalues are not degeneracies but rather an essential and useful property of the system which enable a rich spatio-temporal classification of branches of solutions, and allow us to prove their global unbounded continuation away from the bifurcation point. We will provide an overview and some historical background on equivariant degree theory in the next section.

\section{Introduction to equivariant degree}\label{sec:intro:equideg}
Equivariant degree theory is a synergy of two distinct and powerful ideas. A map is equivariant if it commutes with a group action. More specifically, if $G$ is a group, $V$ is a $G$-representation, and $f:V\to V$ a map, then $f$ is said to be $G$-equivariant if, for all $g\in G$, $f(gv) = gf(v)$. This notion is both complementary and closely related to the concept of invariance, where a subset $\Omega \subseteq V$ is said to be $G$-invariant if for all $v \in \Omega$ and all $g \in G$, $gv \in \Omega$. On an intuitive level, equivariant theories are based on the idea that if a problem is equivariant, then its solutions will inherit some of these symmetries.

In this dissertation, we use the equivariant degree to show the bifurcation of non-constant periodic multiconsensus solutions from stable consensus in three novel classes of multi-agent systems featuring three types of memory. The closed-loop dynamics equations governing each of these systems forms three different classes of delay differential equations. This allows us to view the goals and results of this thesis from a complementary perspective as showing global Hopf bifurcation of spatio-temporally symmetric branches of non-constant periodic solutions from a locally asymptotically stable trivial equilibrium in three different families of symmetrically coupled delay differential equations. Based on the properties of these systems as DDEs (in particular whether they are of retarded or neutral type), different degrees must ultimately be used.

periodic solutions in three separate classes of distributed delay differential equations. Based on the properties of these systems, three different degrees are used. We will detail these degrees and the reasons why they are the appropriate choices for their respective problems in this section. First, similar to the exposition given for delay equations in the previous section, we will give a brief overview of topological degree methods, as well as some background on the historical development of the degree and its equivariant extension.
\vs
The core idea of topological degree methods is to define certain topological invariants on homotopy classes of maps, where these homotopy classes are defined by homotopies satisfying appropriate admissibility conditions. These invariants---if they can be effectively computed---can guarantee the existence of solutions for all maps in the same homotopy class. Therefore, if a map belongs to the same homotopy class as its own linearization, then one can compute this invariant on the linearized map and apply these results to the nonlinear system. Every type of degree therefore faces the dual challenges of both defining these invariants abstractly and rigorously establishing the admissibility criteria for homotopies, and also finding ways that these invariants can actually be computed for sufficiently simple maps in a given homotopy class (e.g. linear maps). 
\vs
For linear maps, these invariants are frequently computable using the spectrum. This means that, assuming the map in question is well-behaved enough to be admissibly homotopic to its own linearization, we can obtain very strong results from spectral data which is already necessary to compute in order to perform many other types of analysis---such as analyzing the asymptotic stability of an equilibrium, for example.
\vs
Additionally, the conditions under which it is possible to construct an admissible homotopy from a map to its linearization on a sufficiently small open neighborhood are fairly relaxed, requiring neither smoothness nor analyticity, only continuity and some very modest regularity assumptions. 
\vs
However, the classical Brouwer degree for finite-dimensional maps and its extension into Banach spaces, the Leray-Schauder degree, have certain shortcomings when used to study complicated symmetric systems, or when one is interested in particular types of solutions (e.g. when one is interested in non-constant periodic solutions). This is because the degree must be taken within a certain region, and if that region contains several solutions, then it is possible their individual degree coefficients may sum to zero.
\vs
To use the Leray-Schauder degree to show the existence of solutions to differential equations, one must first reformulate the equation as an operator between suitably chosen functional spaces. Let $\mathscr E$ denote this functional space and $\mathscr F:\mathscr E\to \mathscr E$ the operator such that $\mathscr F(x) = 0$ implies that $x(t)$ is a solution to our differential equation.
\vs
If the differential equation is autonomous and $\mathscr F(0)=0$ is a solution, and one can show the existence of \emph{a priori} bounds on the norm of solutions in our functional space, then one standard approach (which we will call the ``big ball/small ball'' approach) is to show that for sufficiently large $R>0$, there is an admissible homotopy defined on the ball $B_R(0)$ from $\mathscr F$ to the identity map, and therefore the degree is trivial. If certain regularity assumptions are satisfied (or $\mathscr F$ admits a regular normal approximation), then for sufficiently small $r>0$ there exists an admissible homotopy on the ball $B_r(0)$ from $\mathscr F$ to its linearization around 0, $\mathscr A$. By computing the degree of this linearization using its spectrum, and using the additivity property of the degree, we can obtain the degree on the annulus $\Omega :=B_R(0)\setminus \overline{B_r(0)}$. If this degree is nonzero, this guarantees the existence of a nontrivial solution somewhere in $\Omega$.
\vs
However, we do not know where this solution is or anything about its nature other than its nontriviality. Additionally, it is entirely possible that although there are many nontrivial solutions to $\mathscr F$ in $\Omega$, they are such that their individual degrees all sum to zero, and in this case the degree tells us nothing about their existence.
\vs
The equivariant degree elegantly addresses this problem and greatly extends the capabilities of topological degree methods by noting that, if $\mathscr E$ is a $G$-space, and $\mathscr F$ is $G$-equivariant for some compact Lie group $G$ (i.e. if for all $g \in G, \mathscr F(gx) = g\mathscr F(x)$), then solutions to $\mathscr F$ must inherit some of these symmetries. In particular, each solution must lie on exactly one $G$-orbit (or equivalently, be fixed by exactly one subgroup of $G$), and solutions on the same orbit (or on conjugate orbits) will share the same local Leray-Schauder degree. Accordingly, the integer value of the Leray-Schauder degree can be decomposed into a formal sum of conjugacy classes of subgroups of $G$. 
\vs
It should be noted that there are certain properties of solutions which equivariant degree cannot tell us. In particular, the stability of detected solutions as well as their minimal period are not detectable by the degree. However, stability is also not detectable through the classical Hopf bifurcation theorem from the linear picture alone, requiring the computation of nonlinear terms up to cubic order which, as mentioned above, is a very difficult and often computationally intractable problem. 
\vs
These properties can be explored numerically, although the numerical analysis of DDEs is not without its own considerable challenges. We perform some numerical simulations on the example systems in chapters \ref{chapter:distributed} and \ref{chapter:neutral}. In section \ref{appendix:software} of this chapter, we discuss the numerical methods used in this dissertation and their limitations, some of the challenges of numerical analysis of DDEs, and potential ways that deeper analysis could be performed. Numerical analysis of DDEs is a deep subject and is not the focus of this dissertation, and the use of these methods here is merely intended to illustrate how numerical techniques can complement the use of equivariant degree and fill in some of its blind spots.
\subsection{Historical background}
It could be argued that topological methods of this type have their earliest foundations in the intermediate value theorem, which itself has its roots in deep antiquity. From the very beginnings of the understanding of functions as objects which can be analyzed, the notion of continuity (and the implicitly topological character of such a notion) was present in the thinking of many mathematicians. The first tool which can be called a topological degree in the modern sense was developed in 1912 by L. E. J. Brouwer \cite{Brouwer1912}. 
\vs
The Brouwer degree, defined for continuous maps between finite-dimensional spaces, is still the anchor of many subsequent degree theories which ultimately reduce to it in the final computation. In 1934, Leray and Schauder \cite{LeraySchauder1934}, motivated by the search for solutions to nonlinear PDEs, generalized the Brouwer degree to continuous maps between Banach spaces for compact perturbations of identity (cf. Section \ref{sec:prelim:leray-schauder:compact-properties}).
\vs
The use of topological degree techniques to study local bifurcation was pioneered by M.A. Krasnosel'skii \cite{Krasnoselskii1956} in the 1950s. Rabinowitz extended this to the foundational global bifurcation result in 1985 \cite{rabinowitz1971,rabinowitz1985}. Both of these results used the Leray-Schauder degree, and were extended to the equivariant Leray-Schauder degree by Krawcewicz et al. in \cite{KrawcewiczWu1997}. On the other hand, following Darbo's introduction of condensing maps in 1955 \cite{darbo1955}, Nussbaum and Sadovskii independently developed an extension of the Leray-Schauder degree to condensing maps in the 1970s \cite{nussbaum1972,sadovskii1971}. These ideas were combined in the creation of the composite coincidence degree by Erbe, Krawcewicz, and Wu in 1993 \cite{Erbe1993}, which was directly applied to study boundary problems in neutral equations. This was extended to the symmetric case by Krawcewicz and Wu in 1997 \cite{KrawcewiczWu1997}. 
\vs
The first equivariant extensions of the Leray-Schauder degree occurred in the late 1980s. In 1988, Dylawerski \cite{Dylawerski1988} developed an $S^1$-equivariant degree for maps between representation spheres. This was extended by Dylawerski, G\k{e}ba, Jodel, and Marzantowicz \cite{dylawerski1991s1} in 1991 to connect with the Fuller index for periodic orbits. Independently and simultaneously, the group of Ize, Massab\'o, and Vignoli also constructed an $S^1$-equivariant degree and extended it to general compact Lie groups \cite{Ize1989,Ize1992,izevignoli2003}. However, this degree was very difficult to compute for general nonlinear problems, and no algorithmic approach yet existed.
\vs
The use of this degree to study Hopf bifurcation required an extension of the $S^1$ degree which could detect spatio-temporal symmetries. Gęba, Krawcewicz, and Wu \cite{geba1994equivariant} introduced the \emph{twisted} equivariant degree. This tool was further developed by Balanov and Krawcewicz \cite{balanov2008symmetric} and by Balanov, Krawcewicz, and Steinlein \cite{BalanovKrawcewiczSteinlein2006}; a comprehensive treatment of equivariant degree theory can be found in \cite{KrawcewiczWu1997,izevignoli2003}. 
\vs
The other great theoretical achievement which allowed for the practical computation of the general $G$-equivariant degree was the recurrence formulas for computing degree coefficients, and the splitting lemma allowing for the $G$-degree of a map to be computed as a Burnside ring product of basic degrees. This was developed in 2006 by Balanov, Krawcewicz, and Ruan \cite{BalanovKrawcewiczRuan2006a}, along with an axiomatization of the degree. These recurrence formulas provided an algorithmic approach for computing the coefficients of orbit types having finite Weyl group and, by considering a slice normal to an $S^1$ orbit, can also be applied to the $S^1$-degree and twisted degree.

\chapter{Preliminaries}\label{chapter:preliminaries}
Here we will give some technical preliminaries on equivariant degree methods and multi-agent systems. The first section on the topological degree was designed to give enough background for a reader relatively unfamiliar with topological degree methods and applications to follow the main results in Chapters \ref{chapter:distributed}, \ref{chapter:pseudo}, and \ref{chapter:neutral}, and understand how they were proved and obtained. As such, it is much longer and more thorough than Section \ref{sec:prelim:mas}, which provides needed definitions for multi-agent systems. 

\section{Topological degree}

In this section we will give the necessary technical preliminaries on topological degree methods. These include the Brouwer, Leray-Schauder, and Nussbaum-Sadovskii degrees. Once these have been developed in a non-equivariant setting, we will then define the two relevant equivariant degrees used in subsequent chapters: the $G$-equivariant degree, and the twisted degree. As will be shown, these equivariant degrees can be defined (or rather, computed) in terms of any of the above non-equivariant degrees as dictated by the setting of the application. 
\vs
The definitions and results presented in this section are standard in the literature, and are presented here in the interest of completeness and to allow this work to be as self-contained as possible. Most of the statements and proofs in this section are based on those in \cite{BalanovEtAl2025}, tailored to our notational conventions, but can also be found in other modern texts concerning equivariant degree. 
\vs
In order to avoid excessive bloat, we will not regurgitate all of the proofs used in constructing the degree, some of which are rather technical. Rather, we will focus on the theorems and proofs which are judged to be most important and useful in actually applying the degree to real problems. In particular, we will omit the existence proofs of the degree, and the proofs of some intermediate lemmas used to prove important results will be streamlined, sketched, or omitted depending on their immediacy to equivariant degree computations.
\vs
Additionally, we define these degrees in a slightly less general setting than found in other sources, tailored to our particular set of applications. For example, in other sources the degree is defined for maps $F:X\to Y$ where $X$ and $Y$ are, depending on the degree, two finite-dimensional vector spaces, two Banach spaces,  or two metric spaces. Since our applications concern maps $F:X\to X$, in order to avoid unnecessary notation which would not be reused, we define all the degrees in these terms, i.e. a map $F:X\to X$ from a space to itself. However, all definitions and results presented here apply to that more general setting in the obvious way.
\vs
There are two distinct challenges in establishing and using any topological degree method. The first is to define the degree and its associated classes of admissible maps, and show that it actually is a well-defined topological invariant over homotopy classes of these maps. The second challenge is to show that this degree can actually be computed, and establish some sort of general procedure for computing it over a reasonably large and interesting subset of those homotopy classes of admissible maps. This second challenge is typically at least as difficult as the first, and this is holds even more true for the equivariant degree. For this reason, we will take the approach of first defining each degree axiomatically, and then showing how they are computed. The existence proofs showing that a degree function actually exists satisfying these axioms are typically the most technical, and will be omitted here.
\vs
Finally, we will make one last preemptive remark. We will reuse notation throughout this section, for example using $\deg(f,\Omega)$ to refer to the Brouwer, Leray-Schauder, or Nussbaum-Sadovskii degrees. Due to the different contexts in which these degrees are used, their axioms, and the ways by in which they are computed, this use does not create any inconsistency and greatly economizes notation.
\section{Brouwer degree}
Let $V$ be a finite-dimensional vector space, $\Omega \subset V$ an open bounded set, and $f:V\to V$ a continuous map. The map $f$ is said to be \textit{$\Omega$-admissible} (equivalently, $(f,\Omega)$ is called an \textit{admissible pair}) if $\forall\; x \in \partial \Omega, f(x) \neq 0$. Denote by $\mathcal M$ the set of all admissible pairs. Then the Brouwer degree is the unique map $\text{deg}:\mathcal M\to \mathbb Z$ which satisfies the following properties:
\begin{enumerate}[label=\textbf{(B.\arabic*)}]

\item \textbf{(Additivity)}\label{prelim:brouwer:p1} Let $\Omega_{1}$ and $\Omega_{2}$
be two disjoint open subsets of $\Omega$ such that
$f^{-1}(0)\cap\Omega\subset\Omega_{1}\cup\Omega_{2}$. Then,
\begin{align*}
\deg(f,\Omega)=\deg(f,\Omega_{1})+\deg
(f,\Omega_{2}).
\end{align*}

    \item \textbf{(Homotopy)}\label{prelim:brouwer:p2} If $h:[0,1]\times V\to V$ is an
$\Omega$-admissible homotopy, then
\begin{align*}
\deg(h_{t},\Omega)=\mathrm{constant}.
\end{align*}

    \item \textbf{(Normalization)}\label{prelim:brouwer:p3} Let $\Omega$ be an
open bounded neighborhood of $0$ in $V$. Then,
\begin{align*}
\deg(\id,\Omega)=1.
\end{align*}

\item \textbf{(Existence)}\label{prelim:brouwer:p4} If $\deg(f,\Omega)\neq0$, then there exists $x_0 \in \Omega$ such that $f(x_0)=0$.

\item \textbf{(Excision)}\label{prelim:brouwer:p5} Suppose that $(f,\Omega)\in \mathcal M$ and $\Omega_0$ is an open subset of $\Omega$ such that $f^{-1}(0) \cap \Omega \subset \Omega_0$. Then $\deg(f,\Omega) = \deg(f,\Omega_0)$.

\item \textbf{(Multiplicativity)}\label{prelim:brouwer:p6} For any $(f_{1},\Omega
_{1}),(f_{2},\Omega_{2})\in\mathcal{M}$,
\begin{align*}
\deg(f_{1}\times f_{2},\Omega_{1}\times\Omega_{2})=
\deg(f_{1},\Omega_{1})\cdot \deg(f_{2},\Omega_{2}).
\end{align*}
\end{enumerate}
The full list of prescribed properties of the Brouwer degree varies somewhat in the literature. This is because only properties \ref{prelim:brouwer:p1}---\ref{prelim:brouwer:p3} are necessary for the axiomatic definition of the Brouwer degree, and thsee can be used to derive properties \ref{prelim:brouwer:p4}---\ref{prelim:brouwer:p6}, and others. We use this list on the basis that these properties are the most useful for our applications.
\vs
We note that although the above definition assumes that $f$ is defined and continuous on the whole space $V$, clearly it need only be continuous and defined on $\overline \Omega$, as the Tietze-Dugundji theorem provides the existence of a continuous extension to all of $V$. This fact will be useful in the construction of the equivariant Brouwer degree.
\vs
\subsection{Computing the Brouwer degree}\label{sec:prelim:brouwer:computing}
The computation of all the degrees discussed in this dissertation can be defined in terms of computing the degrees of linear isomorphisms, and by the conditions under which nonlinear maps are guaranteed to be locally homotopic to linear isomorphisms near their zeroes. The following proposition is the foundation of all Brouwer degree computations and, by extension, computations for the Leray-Schauder and Nussbaum-Sadovskii degrees. The proof of this is standard and can be obtained in many references (cf. \cite{KrawcewiczWu1997,BalanovKrawcewiczRuan2006a,BalanovEtAl2025}), but we present a brief sketch of it here for completeness.
\vs
\begin{proposition}\label{prop:prelim:brouwer-compute}
Let $A:V\to V$ be a continuous linear isomorphism and $B(V)$ the open unit ball on $V$. Then $\deg(A,B(V)) = \sign\det A$.
\end{proposition}
\begin{proof}
Notice that $\mathrm{GL}(V)$ (which, given a choice of basis, can be identified with $\mathrm{GL}(n,\mathbb R)$) can be written as the disjoint union of two path-connected components, i.e. $\mathrm{GL}(n,\mathbb R) = \mathrm{GL}^+(n,\mathbb R) \sqcup \mathrm{GL}^-(n,\mathbb R)$, where
\begin{align*}
\mathrm{GL}^+(n,\mathbb R) &= \{A \in \mathrm{GL}(n,\mathbb R):\det A>0\}\\
\mathrm{GL}^-(n,\mathbb R) &= \{A \in \mathrm{GL}(n,\mathbb R):\det A<0\}.
\end{align*}
This implies that there are only two homotopy classes of linear isomorphisms between finite-dimensional real vector spaces. The normalization property \ref{prelim:brouwer:p3} implies that if $A \in \mathrm{GL}^+(n,\mathbb R)$, then $\deg(A,B(V)) = 1$. If $A\in \mathrm{GL}^-(n,\mathbb R)$, then homotopy \ref{prelim:brouwer:p2} and multiplicativity \ref{prelim:brouwer:p6} imply $\deg(A,B(V)) = \deg(-x,(-1,1))$, i.e. the degree of the one-dimensional minus identity map. Let $h_t:[0,1]\times \mathbb R\to \mathbb R$ be defined $h_t(x) = (1-t) + t(\vert x \vert -1)$. Let $\Omega = (-2,2)$ and note that $h_t(x)$ is $\Omega$-admissible. 
\vs
Since $h_0(x)\equiv1$ is a constant map with no zeroes, $\deg(h_0,(-2,2)) = 0=\deg(h_t,(-2,2))$. By additivity \ref{prelim:brouwer:p1}, we also have $\deg(h_1,(-2,2)) = \deg(-x-1,(-2,0)) +\deg(x-1,(0,2))$. Since $\deg(-x-1,(-2,0))= \deg(-x,(-1,1))$ and $\deg(x-1,(0,2)) = \deg(x,(-1,1))$, this implies that $0 = \deg(-x,(-1,1)) + \deg(x,(-1,1)) = \deg(-x,(-1,1)) + 1$ and so $\deg(-x,(-1,1))=-1$. 
\end{proof}
\vs
The excision property \ref{prelim:brouwer:p5} suggests that we can use this approach to compute the degree of a general admissible pair $(f,\Omega)$ only if $f^{-1}(0)$ is a finite set, and if for every $x_j \in f^{-1}(0)$, there is a linear isomorphism admissibly homotopic to $f$ on a sufficiently small ball. If 0 is a regular value of $f\vert_\Omega$, then $f$ is called an \emph{$\Omega$-admissible regular map} and $(f,\Omega)$ an \emph{$\Omega$-admissible regular pair}. 
\vs
If $f$ is not an $\Omega$-admissible regular map, Sard's lemma guarantees that for any $\varepsilon>0$, there is an $\Omega$-admissible regular map $\tilde f$ such that $\sup_{x\in\Omega} |f(x)-\tilde f(x)|<\varepsilon$. These ideas can be summarized by the following proposition:
\vs
% Let $x_0 \in f^{-1}(0)$, and denote by $Df(x_0)$ the Jacobian of $f$ at $x_0$. By the definition of the derivative, $Df(x_0)$ is admissibly homotopic to $f(x)-x_0$ on a sufficiently small ball. Therefore, the local Brouwer degree of $f$ on an isolating neighborhood of $x_0$ is given by:
% \[
% \deg(f,B_\varepsilon(x_0)) = \text{sign}(\det Df(x_0)),
% \]
% where $B_\varepsilon(x_0)$ is a sufficiently small open ball that $f^{-1}(0)\cap B_\varepsilon = \{x_0\}$, and $f$ is $B_\varepsilon$-admissible.
% It follows that if $f^{-1}(0) \cap \Omega$ is a finite set, then
% \[
% \deg(f,\Omega) = \sum_{x_0\in f^{-1}(0)\cap \Omega}\text{sign}(\det Df(x_0)).
% \]
\begin{proposition}
Let $V:= \mathbb R^n$ and $\Omega \subset V$ an open bounded set. Let $f:V\to V$ be a continuous $\Omega$-admissible regular map. Then
\[
\deg(f,\Omega) = \sum_{x_0\in f^{-1}(0)\cap \Omega}\text{sign}(\det Df(x_0)).
\]
\end{proposition}
\begin{proof}
Since $f$ is an $\Omega$-admissible regular map, the set $f^{-1}(0)\cap \Omega$ is finite, and for each $x_j \in f^{-1}(0)\cap\Omega$, $Df(x_j)$ exists and is a linear isomorphism, and there exists $\varepsilon>0$ sufficiently small such that $f$ is $B_\varepsilon(x_j)$-admissibly homotopic to $Df(x_j)$. The conclusion then follows directly from the additivity, homotopy, and excision properties of the Brouwer degree.
\end{proof}
\section{Leray-Schauder degree}\label{sec:prelim:leray-schauder}
The Leray-Schauder degree extends the Brouwer degree to completely continuous fields on Banach spaces. Let $\mathscr E$ be a Banach space, and $F:\mathscr E\to \mathscr E$ a continuous map. If $\overline{F(X)}$ is compact for any subset $X\subset \mathscr E$, then $F$ is called a \emph{compact map}. A map $\mathscr F:\mathscr E\to \mathscr E$ is called a \emph{compact perturbation of identity}
\footnote{In older literature, $F$ is also called a \emph{completely continuous map} if $\overline{F(A)}$ is compact for any \emph{bounded} $A\subset \mathscr E$. $\mathscr F:= \id - F$ is called a \emph{compact field} or a \emph{completely continuous field} if $F$ is a compact map or completely continuous map, respectively. However, these definitions coincide when $F$ is restricted to an open bounded subset $\Omega \subset \mathscr E$, so we use the term ``compact perturbation of identity'' for consistency and clarity.} 
if it can be written $\mathscr F = \id - F$, where $F$ is a compact map. Given a compact perturbation of identity $\mathscr F$ and an open bounded subset $\Omega\subset \mathscr E$, $(\mathscr F,\Omega)$ is called an admissible pair if $\mathscr F$ is a compact perturbation of identity. We denote the set of all admissible pairs over $\mathscr E$ by $\mathscr M(\mathscr E)$, and the set of all admissible pairs as
\[
\mathscr M : = \bigcup_{\mathscr E}\mathscr M(\mathscr E).
\]
Then the Leray-Schauder degree is the unique function $\deg:\mathscr M\to \mathbb Z$ satisfying the following properties:
\begin{enumerate}[label=\textbf{(L.\arabic*)}]

    \item \textbf{(Additivity)}\label{prelim:leray-schauder:p1} Let $\Omega_{1}$ and $\Omega_{2}$
be two disjoint open subsets of $\Omega$ such that
$\mathscr F^{-1}(0)\cap\Omega\subset\Omega_{1}\cup\Omega_{2}$. Then,
\begin{align*}
\deg(\mathscr F,\Omega)=\deg(\mathscr F,\Omega_{1})+\deg
(\mathscr F,\Omega_{2}).
\end{align*}

    \item \textbf{(Homotopy)} If $H_t:[0,1]\times \mathscr E\to \mathscr E$ is an
$\Omega$-admissible homotopy, then
\begin{align*}
\deg(H_{t},\Omega)=\mathrm{constant}.
\end{align*}

    \item \textbf{(Normalization)}\label{prelim:leray-schauder:p3} For any open bounded subset $\Omega \subset \mathscr E$ where $\mathscr E$ is a Banach space, and any $x_0 \in \mathscr E$ such that $x_0 \not\in \partial \Omega$, one has
\begin{align*}
\deg(\id - x_0,\Omega) = \begin{cases}
    1\quad &\text{if } x_0 \in \Omega\\
    0\quad &\text{if } x_0 \not\in \Omega\\
\end{cases}
\end{align*}

    \item \textbf{(Existence)} If $\deg(\mathscr F,\Omega)\ne
0$, then there
exists $x\in\Omega$ such that $\mathscr F(x)=0$.

\item \textbf{(Excision)}
Suppose that $(\mathscr F,\Omega)\in \mathscr M$ and $\Omega_0$ is an open subset of $\Omega$ such that $\mathscr F^{-1}(0) \cap \Omega \subset \Omega_0$. Then $\deg(\mathscr F,\Omega) = \deg(\mathscr F,\Omega_0)$.

    \item \textbf{(Multiplicativity)}\label{prelim:leray-schauder:p6} For any $(\mathscr F_{1},\Omega
_{1}),(\mathscr F_{2},\Omega_{2})\in\mathscr{M}$,
\begin{align*}
\deg(\mathscr F_{1}\times \mathscr F_{2},\Omega_{1}\times\Omega_{2})=
\deg(\mathscr F_{1},\Omega_{1})\cdot \deg(\mathscr F_{2},\Omega_{2}).
\end{align*}
\end{enumerate}

\subsection{Properties of compact maps and compact perturbations of identity}\label{sec:prelim:leray-schauder:compact-properties}
The distinction between a compact map $F$ and a compact perturbation of identity $\mathscr F:=\id - F$ can be viewed in two ways. On one hand, it transforms the problem of finding fixed points of $F$ into the problem of finding zeroes of $\mathscr F$, which can be treated with the degree. 
\vs
From a spectral point of view, we note that if $F$ is compact and its Fr\'echet derivative exists at some $x_j$, then $A:=DF(x_j)$ is a compact linear operator. By the Riesz-Schauder theorem, this means that eigenvalues of $A$ can only accumulate at 0, and in particular, $\mu > 1$ for only finitely many eigenvalues $\mu\in \sigma(A)$. 
\vs
If $A$ exists and is a compact map, then $\mathscr A:= D\mathscr F(x_j)$ exists and is a compact perturbation of identity. If the map $F$ is assumed to be regular, then $A$ is additionally a linear isomorphism. Therefore $\mathscr A$ is also a Fredholm operator of index 0, and is invertible. This allows us to write $\mathscr E = E^+ \oplus E^-$, where
\[
E^+ := \bigcup_{\mu \in \sigma_+(\mathscr A)} E(\mu)\quad \text{ and } \quad E^- := \bigcup_{\mu \in \sigma_-(\mathscr A)} E(\mu),
\]
where $\sigma_+(\mathscr A):=\{\mu \in \sigma(\mathscr A):\re(\mu)>0\}, \sigma_-(\mathscr A):=\{\mu \in \sigma(\mathscr A):\re(\mu)<0\}$, and $E(\mu)$ represents the generalized eigenspace of $\mu\in \sigma(\mathscr A)$. Since $A$ is a compact map, $\sigma_-(\mathscr A)$ is finite. 

\subsection{Computing the Leray-Schauder degree}
Bearing in mind that we can only really compute the Brouwer degree of linear isomorphisms on finite dimensional spaces, the key to computing the Leray-Schauder degree lies in the fact that if $\mathscr F$ is a completely continuous field, then $\sigma_-(\mathscr F)$ is finite. 
\begin{proposition}\label{prop:prelim:leray-schauder}
Let $\mathscr F:\mathscr E\to \mathscr E$ be a compact perturbation of identity, $(\mathscr F,\Omega)$ be an $\Omega$-admissible regular pair, and let $x_j \in \mathscr F^{-1}(0)\cap \Omega$ be an isolated regular point of $\mathscr F$ such that $\mathscr F$ is Fr\'echet differentiable at $x_j$. Put $\mathscr A:= D\mathscr F(x_j)$. Then
\[
\deg(\mathscr F,B_\varepsilon(x_j)) = \deg(\mathscr A,B(\mathscr E)) = (-1)^{\vert\sigma_-(\mathscr A)\vert},
\]
where $\deg$ is the Leray-Schauder degree and $B_\varepsilon(x_j)$ is a sufficiently small isolating neighborhood of $x_j$.
\end{proposition}
\begin{proof}
First we note that by shifting $\mathscr F$, because it is a regular map, it is admissibly homotopic to $\mathscr A$ on a sufficiently small ball, and so
\[
\deg(\mathscr F,B_\varepsilon(x_j)) = \deg(\mathscr F(x-x_j),B_\varepsilon(0))=\deg(\mathscr A,B(\mathscr E)).
\]
Note that $E^-$ is a finite-dimensional vector space. Put $\mathscr A_- := \mathscr A\vert_{E^-}, \mathscr A_+ := \mathscr A\vert_{E^+}$.Since $\mathscr A(E^-)\subseteq E^-$, the map $\mathscr A_-$ is a linear isomorphism of the finite dimensional vector space $E^-$. On the other hand, since $\mathscr E = E^- \oplus E^+$ the multiplicativity property \ref{prelim:leray-schauder:p6} implies that $\deg(\mathscr A,B(\mathscr E)) = \deg(\mathscr A_-,B(E^-))\cdot \deg(\mathscr A_+,B(E^+))$. The term $\deg(\mathscr A_-,B(E^-))$ reduces to the Brouwer degree. 
\vs
To compute $\deg(\mathscr A_+,B(E^+))$, we construct the homotopy $h_t:[0,1]\times E^+ \to E^+$ such that $h_t=t \id\vert_{E^+} + (1-t)\mathscr A_+$. This homotopy is $B(E^+)$-admissible unless $0\in \sigma(h_t)$ for some $t\in[0,1]$. If $\mu_+ \in \sigma(\mathscr A_+)$, then $1+t(\mu_+-1)\in \sigma(h_t)$. Since $\mu_+>0$ for all $\mu_+\in \sigma(\mathscr A_+)$, $1+t(\mu_+-1)>0$ for all $\mu_+\in \sigma(\mathscr A_+)$ and all $t\in [0,1]$. Finally, by the normalization property \ref{prelim:leray-schauder:p3}, $\deg(\mathscr A_+,B(E^+)) = 1$, and so $\deg (\mathscr A,B(\mathscr E)) = \deg(\mathscr A_-,B(E^-))$, and the result follows directly from the computational formula for the Brouwer degree.  
\end{proof}

\section{Nussbaum-Sadovskii degree}\label{sec:prelim:n-s}
The nonlinear operators which arise from neutral equations typically fail to be compact operators. This is because the neutral operator, as already described, can contribute to the essential spectrum. As such, we require a degree which is computable on a broader class of maps than compact perturbations of identity. This degree turns out to be the Nussbaum-Sadovskii degree, defined on condensing maps. To define these terms, we must first define the related notion of a measure of noncompactness.

\subsection{Measures of noncompactness}
\vs
\begin{definition}
\normalfont
Let $\mathscr E$ be a metric space (for our applications, it will also be a Banach space), and denote by $\mathcal B$ the class of all bounded subsets of $\mathscr E$. Then a function $\mu: \mathcal B \to [0,\infty)$ is called a \textit{measure of noncompactness} on $\mathscr E$ if, for all $A,B\in \mathcal B$, the following conditions are satisfied:
\begin{enumerate}
    \item $\mu(A) = 0 \iff \bar{A} \text{ is compact}$
    \item $\mu(A) = \mu(\bar{A})$
    \item $\mu(\text{conv}(A)) = \mu(A)$
    \item $\mu(A \cup B) = \max\{\mu(A),\mu(B)\}$
    \item $\mu(rA) = |r|\cdot\mu(A),\quad r\in\mathbb R$
    \item $\mu(A+B) \leq \mu(A) + \mu(B)$
\end{enumerate}
\end{definition}
\vs
There are many classical measures of noncompactness, but the one of greatest utility to us is the \textit{Kurotowski measure of noncompactness} $\alpha:\mathcal B \to [0,\infty)$ defined for $A\in\mathcal B$
\begin{equation}\label{eq:kurotowski}
\alpha(A) = \inf\left\{\varepsilon > 0: \bigcup_{i=1}^N A_i = A,\quad\text{diam}(A_i) < \varepsilon\quad \forall i\in 1,\dots,N\right\}
\end{equation}
\subsection{Condensing maps}
Given a metric space $\mathscr E$, a bounded subset $A\subset \mathscr E$, and a measure of noncompactness $\mu$; let $F:\mathscr E\to \mathscr E$ be a continuous map such that $F(A)$ is bounded. Then $F$ is
\normalfont Let $V,W$ be metric spaces, $\mu$ a measure of noncompactness on $V$, $X\subset V$, $Y \subset W$, and $F:X \to Y$ a continuous map taking bounded subsets of $X$ to bounded subsets of $Y$. Then $F$ is
\begin{enumerate}
    \item a $\mu$-\textit{Lipschitzian map} with constant $\kappa \geq 0$ if $\mu(F(A)) \leq \kappa\mu(A)$ for all bounded subsets $A \subset \mathscr E$.
    \item a \textit{compact map} if it is $\mu$-Lipschitzian with $L=0$
    \item a \textit{Darbo map} (or \textit{$\mu$-set contraction}) if it is $\mu$-Lipschitzian with $L<1$
    \item a \textit{condensing map} if it is $\mu$-Lipschitzian with $L = 1$ and $\mu(F(A)) < \mu(A)$ for every bounded subset $A \subset \mathscr E$ such that $\mu(A) > 0$
\end{enumerate}
Likewise, $\mathscr F:=\id - F$ is called a \emph{condensing perturbation of identity} if $F$ is a condensing map. One can immediately see that any compact perturbation of identity is also a condensing perturbation of identity. Additionally, any contraction map is also condensing. Indeed, if there exists $0<q<1$ such that for all $x,y\in \mathscr E$, $\norm{F(x)-F(y)} \leq \norm{x-y}$, then $F$ is a condensing map with respect to the Kuratowski measure of non-compactness.
\vs
The same definitions of admissibility apply as with the Leray-Schauder degree, with the relaxed condition that $\mathscr F$ is a condensing perturbation of identity. To ease the notational burden, we will again use $\mathscr M$ to refer to the set of all admissible pairs when it is clear from the context that we are interested in condensing perturbations of identity. Two condensing perturbations of identity $\mathscr F_0$ and $\mathscr F_1$ are called $\Omega$-admissibly homotopic if there exists an $\Omega$-admissible condensing map $h_t:[0,1]\times \mathscr E\to \mathscr E$ such that $h_0=\mathscr F_0$ and $h_1=\mathscr F_1$, and $h_t$ is called an $\Omega$-admissible homotopy. 
\vs
If we allow ourselves the slight reuse (or abuse) of notation to interpret $\mathscr F$ as a condensing perturbation of identity instead of a compact perturbation of identity, and $\deg$ as referencing the Nussbuam-Sadovskii degree instead of the Leray-Schauder degree, then the Nussbaum-Sadovskii degree is the unique function $\deg:\mathscr M\to \mathbb Z$ satisfying the properties \ref{prelim:leray-schauder:p1}---\ref{prelim:leray-schauder:p6}. As such, we will not reprint the identical properties here. 
\subsection{Properties of condensing maps and condensing perturbations of identity}
Similarly to Section \ref{sec:prelim:leray-schauder:compact-properties}, we will give a few basic properties of condensing maps which will aid us in the computation of the Nussbaum-Sadovskii degree, and also show why these maps are the appropriate framework for studying certain types of differential equations, particularly neutral functional differential equations.
\vs
A simple intuition for condensing maps can be given by considering the family of bounded linear operators $F_\kappa:=\kappa\id$, i.e. homotheties. As is well-known, the unit ball $B(\mathscr E)$ in an infinite-dimensional Banach space is not compact, and so $F_\kappa$ is not a compact map unless $\kappa=0$. For $0<\kappa<1$, $F_\kappa$ is clearly Darbo (and hence condensing) because 
\[
\alpha(F_\kappa(B(\mathscr E))) = \kappa\alpha(B(\mathscr E)) = 2\kappa,
\]
where $\alpha$ is the Kuratowski measure of noncompactness and $\alpha(B(\mathscr E))$ follows from the fact that $\mathrm{diam}(B(\mathbb R^n)) = 2$ for any $n$. Since $DF_\kappa(0) = F_\kappa$, we can see that the Fr\'echet derivative of a condensing operator can also be condensing. 
\vs
Recall that the \emph{essential spectrum} of a bounded linear operator $A:\mathscr E\to \mathscr E$, denoted $\sigma_{\mathrm{ess}}(A)$, is the set of all $\lambda \in \mathbb C$ such that $A-\lambda\id$ is not Fredholm, i.e. such that $A-\lambda\id$ is unbounded, or its kernel or cokernel is infinite-dimensional. The situation where $\dim\ker(A-\lambda\id)=\infty$ can also be interpreted as $\lambda$ having a generalized eigenspace of infinite dimension. 
\vs
The essential spectrum is invariant under compact perturbations, and for a compact linear operator $C$, $\sigma_{\mathrm{ess}}(C)=\{0\}$. Consequently, discrete eigenvalues of a compact linear operator can accumulate only at 0.
The \emph{essential spectral radius} is denoted $\rho_{\mathrm{ess}}(T)$, and is defined
\[
\rho_{\mathrm{ess}}(A) := \sup\{\vert \lambda \vert : \lambda \in \sigma_{\mathrm{ess}}(A)\}.
\]
\begin{proposition}\label{prop:prelim:condensing-linear}
    $A$ is a condensing linear operator if and only if $\rho_{\mathrm{ess}}(A)<1$. 
\end{proposition}
\begin{proof}
    Let $B\subset \mathscr E$ be a bounded subset, $E(\lambda)$ the generalized eigenspace corresponding to $\lambda \in \sigma(A)$, and put $B_\lambda := B \cap E(\lambda)$. Then clearly $B = \bigcup_{\lambda \in \sigma(A)} B_\lambda$. Note that if $\lambda$ is a discrete eigenvalue having finite multiplicity, then $A\vert_{E(\lambda)}$ is a compact linear operator, and $\alpha(A(B_\lambda))=0$. Since $\alpha(X\cup Y) = \max\{\alpha(X),\alpha(Y)\}$, this implies that $\alpha(A(B))=\alpha(A(B_\mathrm{ess}))$, where $B_\mathrm{ess} := \bigcup_{\lambda \in \sigma_{\mathrm{ess}}(A)} B_\lambda$. Moreover, $\alpha(A(B_\mathrm{ess})) \leq \rho_\mathrm{ess}(A)\alpha(B)$, and the conclusion follows.
\end{proof}
\subsection{Computing the Nussbaum-Sadovskii degree}
\begin{proposition}
Let $\mathscr F:\mathscr E\to \mathscr E$ be a condensing perturbation of identity, $(\mathscr F,\Omega)$ be an $\Omega$-admissible regular pair, and let $x_j \in \mathscr F^{-1}(0) \cap \Omega$ be an isolated regular point of $\mathscr F$ such that $\mathscr F$ is Fr\'echet differentiable at $x_j$, and put $\mathscr A:= D\mathscr F(x_j)$. If $\mathscr A$ is an isomorphism, then
\[
\deg_{\mathrm {NS}}(\mathscr F,B_\varepsilon(x_j)) = \deg_{\mathrm{NS}}(\mathscr A,B(\mathscr E)) = (-1)^{\vert \sigma_-(\mathscr A)\vert},
\]
where $\deg_{\mathrm{NS}}$ stands for the Nussbaum-Sadovskii degree, and $B_\varepsilon(x_j)$ is a sufficiently small isolating neighborhood of $x_j$.  
\end{proposition}
Note the similarity to Proposition \ref{prop:prelim:leray-schauder}, with the difference being that we interpret the degree differently, and that $\mathscr F$ is taken to be a condensing perturbation of identity.

\begin{proof}
We first show that $\mathscr F$ is admissibly homotopic to $\mathscr A$ on a sufficiently small ball. Since a compact map (or compact linear operator) is also condensing, and the convex combination of condensing maps is condensing, this implies that the linear homotopy $h_t = (1-t)\mathscr F(x-x_j) + t\mathscr A$ is condensing and admissible, and therefore remains within the domain of the Nussbaum-Sadovskii degree for all $t\in[0,1]$.
\vs
Now we take a similar approach as the computation of the Leray-Schauder degree. We note that $\mathscr E = E^- \oplus E^+$, where
\[
E^+ := \bigcup_{\mu \in \sigma_+(\mathscr A)} E(\mu)\quad \text{ and } \quad E^- := \bigcup_{\mu \in \sigma_-(\mathscr A)} E(\mu).
\]
Since $\mathscr A := \id - A$, where $A$ is a condensing linear operator, $\rho_\mathrm{ess}(A)<1$, and so $\rho_\mathrm{ess}(\mathscr A) > 0$. This means that $E^-$ is a union over finitely many discrete eigenvalues having finite multiplicity, and hence $\dim E^- <\infty$, and $\mathscr A_- := \mathscr A\vert_{E^-}:E^-\to E^-$ is a linear isomorphism whose Brouwer degree is well defined on the unit ball $B(E^-)$. On the other hand, on the infinite dimensional space $E^+$, a homotopy can be constructed between $\mathscr A_+:= \mathscr A\vert_{E^+}$ and $\id\vert_{E^+}$, and just as with the Leray-Schauder degree, the normalization and multiplicativity properties imply
\begin{align*}
\deg_{\mathrm{NS}}(\mathscr A,B(\mathscr E)) &= \deg_{\mathrm{NS}}(\mathscr A_-,B(E^-)) \cdot \deg_{\mathrm{NS}}(\mathscr A_+,B(E^+))
\\&= \deg_{\mathrm{NS}}(\mathscr A_-,B(E^-)) \cdot 1 = (-1)^{\vert \sigma_-(\mathscr A)\vert}
\end{align*}
\end{proof}

\section{Equivariant degree}
\label{subsec:G-degree}
We will now define the equivariant extensions of the above degrees. There are two types of equivariant degree which will be used in subsequent chapters: the \emph{$G$-degree}, and the \emph{twisted degree}. Each depends on certain characteristics of the group action of a compact Lie group $G$ and is capable of detecting different types of orbits of solutions. As in the previous section, for each degree we will first define it, and then show how it can be computed. The definition and computation aspects each require some background in group theory and representation theory, which will be provided in this section. Additionally, these equivariant constructions are somewhat agnostic of the exact degree used, as long as certain equivariant principles are satisfied by the maps, bounded sets, and homotopies in question. We will note when extra assumptions are required for the Leray-Schauder and Nussbaum-Sadovskii degrees, but we will take the equivariant Brouwer degree and twisted Brouwer degree as our primary setting of interest.
\subsection{Isotropy groups, orbit types, and the Burnside ring}
One common characteristic of the equivariant versions of these degrees is that they no longer take integer values, but rather count solutions according to their isotropy groups. To make this notion precise, and to make it compatible with the additivity and multiplicativity properties of the Brouwer degree, it is necessary to define these groups and a compatible algebraic structure where their sums and products can be taken---the Burnside ring.
\vs
Let $G$ be a compact Lie group and $\mathscr E$ a Banach $G$-space, which we take without loss of generality to be an isometric Banach $G$-representation.
For $u \in \mathscr E$, we denote by $G_u$ the \emph{isotropy group} of $u$, defined
\[
G_u := \{g \in G: g\cdot u = u\}.
\]
We denote by $G(u)$ the \emph{$G$-orbit} of $u$, 
\[
G(u) := \{gu:g\in G\}.
\]
For a subgroup $H \leq G$, we denote by $\mathscr E^H$ the \emph{$H$-fixed point subspace} of $\mathscr E$, defined
\[
\mathscr E^H := \{u \in \mathscr E:\forall\; h \in H, \;h\cdot u = u\},
\]
and by $V_H$ the set
\[
V_H := \{u\in \mathscr E:G_u = H\}
\]
We denote by $(H)$ the conjugacy class of $H$ in $G$. We denote by $\mathscr E_{(H)} := G(\mathscr E_H)$. We call $(H)$ an \emph{orbit type} if there exists some $u \in \mathscr E\setminus\{0\}$ such that $G_u = H$ (which implies $\dim \mathscr E^H >0$). We will denote by $\Phi(G)$ the set of all conjugacy classes of subgroups of $G$, and by $\Phi(G;\mathscr E)$ the set of all orbit types of $(G)$ on $\mathscr E$. There is a natural partial order on $\Phi(G)$ (and $\Phi(G;\mathscr E)$) given by subconjugacy, i.e. $(H)\leq (K)$ if and only if there exists $H' \in (H)$ and $K' \in (K)$ such that $H'\leq K'$. If $(H) \in \Phi(G;\mathscr E)$ is such that $(H) < (G)$, and $(H)\leq (K)$ implies that $(K)=(H)$ or $(K)=(G)$, then we call $(H)$ a \emph{maximal orbit type}. Maximal orbit types have special significance in degree calculations, as will be shown. 
% We will denote by $\mathfrak M$ the set of all maximal orbit types in $\Phi(G;\mathscr E)$.

% Note that when $G$ is an infinite group, there may not exist any maximal orbit types in $\Phi(G;\mathscr E)$. For example, if $G = S^1$ and $\mathscr E = C_{2\pi}(\mathbb R;\mathbb R)$, then $(\mathbb Z_k)<(\mathbb Z_{2k})$ for all $k$, and clearly each of these fixes a subspace of $\mathscr E$. However, by restricting to a finite-dimensional subspace $\mathscr E_{k,j} \subset \mathscr E$, we guarantee that $\Phi(G;\mathscr E_{k,j})$ contains only finitely many orbit types, and therefore there must exist a maximal orbit type in $\Phi(G;\mathscr E_{k,j})$, even if this orbit type is \emph{not} maximal in the full set $\Phi(G;\mathscr E)$. We call such orbit types \emph{orbit types of maximal kind}, and $\mathfrak M_{k,j}$ will refer to all orbit types of maximal kind which are maximal in $\Phi(G;\mathscr E_{k,j})$. 
\vs
We denote by $N(H)$ the normalizer of $H$, and by $W(H) := N(H)/H$ the Weyl group of $H$. The Weyl group is significant because $N(H)$ is, by definition, the largest subgroup of $G$ whose action leaves $\mathscr E^H$ invariant. The Weyl group, therefore, is the largest subquotient of $G$ which acts faithfully on $\mathscr E^H$, and the local (i.e. Brouwer, Leray-Schauder, etc.) degree of a solution $u \in \mathscr E^H$ will be the same for every solution on that $W(H)$-orbit in $\mathscr E^H$. Thus, the size of the Weyl group determines the count of individual solutions in $\mathscr E^H$ which have this degree. Since we can only take the local Brouwer degree of a linear isomorphism, we require that the orbit of a solution $x_0$ consists of isolated points. Since $G$ is a compact Lie group, this means that the orbit must be finite, and therefore that $W(H)$ must be a finite group. We denote the set of all orbit types having finite Weyl group as
\[
\Phi_0(G;\mathscr E) := \{(H) \in \Phi(G;\mathscr E): \dim(W(H)) = 0\}.
\]
We also denote by $\mathfrak M(G;\mathscr E)$ the set of all maximal orbit types in $\Phi_0(G;\mathscr E)$.

% We also put $\mathscr E_H := \{u\in \mathscr E:G_u = H\}$. Since every $u \in \mathscr E$ has some isotropy group, it is clear that these sets partition $\mathscr E$, i.e.
% \[
% \mathscr E = \bigcup_{\substack{H'\in (H)\\(H) \in \Phi(G;\mathscr E)}}\mathscr E_{H'}.
% \]
\subsection{Equivariant and invariant maps and sets}
Consider a group $G$ and a $G$-space $V$. A map $f:V\to V$ is called \emph{$G$-equivariant} if, for all $g \in G$, $gf(x) = f(gx)$. A subset $\Omega \subset V$ is called \emph{$G$-invariant} if, for all $g \in G$ and all $x \in \Omega$, $gx \in \Omega$. It is called invariant if for all $g \in G$, $f(gx)=f(x)$. 
\vs
\begin{proposition}
If $f$ is $G$-equivariant and Fr\'echet differentiable at $x_0\in V$, then $A:= Df(x_0)$ is $G_{x_0}$-equivariant. 
\end{proposition}
\begin{proof}
Recall that $f$ is Fr\'echet differentiable at $x_0$ with $A:=Df(x_0)$ if, for all $h\in V$ with $\norm h$ sufficiently small
\[
f(x_0+h) = f(x_0)+Ah + o(\norm h).
\]
Fix some $g\in G_{x_0}$. By $G$-equivariance of $f$, we have
\[
gf(x_0+h) = f(g(x_0+h)) = f(gx_0 + gh) = f(x_0+gh)
\]
Taking the Fr\'echet derivative expansion of the right hand side, we get
\[
f(x_0 + gh) = f(x_0) + Agh + o(\norm{h}).
\]
On the other hand, we also have
\[
gf(x_0+h) = gf(x_0) + gAh + o(\norm h) = f(x_0) + gAh + o(\norm h).
\]
Equating these and subtracting $f(x_0)$ from both sides yields
\[
gAh = Agh.
\]
Note that if we pick $g' \not\in G_{x_0}$, then this inequality does not hold, since $f(g'x_0)\neq f(x_0)$ and so the Fr\'echet derivative expansion would yield a different liearization at the point $g'x_0$ which does not allow us to draw conclusions on the equivariance of $A$.
\end{proof}
\vs
Since we are often considering equations derived from autonomous systems, where $f(0)=0$, and faithful $G$-actions where $G_0 = G$, this allows us to consider $Df(0)$ as a $G$-equivariant map.
\vs
Additionally, if $f$ is $G$-equivariant, then the solution set $f^{-1}(0)$ is $G$-invariant. More generally, the $G$-orbit of any point is a $G$-invariant set.

\subsection{The Burnside ring}
The Burnside ring $A(G):= \mathbb Z[\Phi_0(G)]$ is a $\mathbb Z$-module over $\Phi_0(G)$. For orbit types $(H),(K)\in \Phi_0(G)$, we define $n(H,K) := \vert\{K'\leq G:K'\in (K),\; H\leq K'\}\vert$, where $\vert X\vert$ stands for the cardinality of the set $X$. Then the coefficient $m_L$ of the orbit type $(L)$ in the Burnside ring product $(H) \cdot (K)$ is defined through the following recurrence formula
\begin{equation}\label{eq:prelim:burnside-mult}
m_L := \frac{n(L,H)\vert W(H)\vert n(L,K)\vert W(K)\vert - \sum_{(\tilde L)>(L)}m_{\tilde L}n(L,\tilde L)\vert W(\tilde L)\vert}{\vert W(L)\vert}
\end{equation}
and the full product $(H)\cdot(K)$ is given by
\[
(H) \cdot (K) := \sum_{(L)\in \Phi_0(G)} m_L(L).
\]
Note that if $(L)>(H)$ or $(L)>(K)$, then $n(L,H) = 0$ or $n(L,K) = 0$, respectively. So in practice, the coefficients need only be computed for $(L)<(H\cap K)$. 
\vs
We also define the $\mathbb Z$-isomorphism $\text{coeff}^H:A(G)\to \mathbb Z$ by
\[
\text{coeff}^H(\alpha):= m_H,\quad \text{where }\alpha := \sum_{(L)\in\Phi_0(G)}m_L(L)\in A(G).
\]
We also note a few important properties of the Burnside ring product which are useful for degree computations:
\begin{enumerate}
    \item $(G)$ acts as the unit element of the commutative ring $A(G)$, i.e. $(G) \cdot \alpha = \alpha$ for all $\alpha \in A(G)$
    \item For any $(H)\in \Phi_0(G)$, one has $(H)\cdot (H) = \vert W(H)\vert(H) +c$, where $c\in A(G)$ is such that $\text{coeff}^L(c)\neq 0 \Rightarrow (L)<(H)$. 
    % \item As a corollary of the above, if $(H)\in\mathfrak M(G;\mathscr E)$, then $((G)-(H)) \cdot ((G)-(H)) = (G) + c$, where as above $c$ is a remainder term consisting of orbit types strictly less than $(H)$ under subconjugacy. This is because
    % \begin{align*}
    % ((G)-m_H(H)) \cdot ((G)-m_H(H)) &= (G)\cdot (G) - (G)\cdot m_H(H) - m_H(H)\cdot (G) + {m_H}^2(H)\cdot(H) \\
    % &= (G) -2m_H(H) + {m_H}^2\vert W(H) \vert(H)\\
    % &= (G)  +({m_H}^2-2m_H\vert W(H) \vert)(H)
    % \end{align*}
\end{enumerate}
The Burnside ring also has a physical, geometrical meaning. In order to ease computations of the degree, it is necessary to exploit the multiplicativity property of the degree in an equivariant setting. This allows us to decompose the overall degree computation into a product of degrees of so-called basic maps. Taking this product of degrees in a way that is compatible with the multiplicativity property requires us to consider the product space of orbits $G/H \times G/K$. This space, equipped with the diagonal $G$-action, may contain many other orbit types which emerge through interaction. The recurrence formula \eqref{eq:prelim:burnside-mult} captures this through a recursive counting argument based on the orbit types emerging from the diagonal $G$-action on this product space.
% \vs
% However, we also note that the coefficients of each orbit type in this product must be computed individually. This is a highly algorithmic process and can easily be performed by software such as GAP, but for arbitrary products of Burnside ring elements, it is often quite difficult to know if any non-maximal orbit type will appear with non-trivial coefficient.
\vs
It is also important to note that the Burnside ring is defined over $\Phi_0(G)$, consisting of conjugacy classes of subgroups of $G$ which have finite Weyl group. These conjugacy classes of subgroups need not be orbit types \emph{per se}, in the sense that they do not need to fix anything in any $G$-space, nor do they need to be the isotropy group of any vector in any $G$-space. However, we often refer to an element $(H)\in \Phi_0(G)$ (and in a Burnside ring element) as an ``orbit type'' regardless of this fact. This is because we are not particularly interested in arbitrary Burnside ring elements, but rather elements which are generated by equivariant degree computations. As we will show in the next section, these computations yield trivial coefficients for any $(H)\in \Phi_0(G)$ which is not an orbit type on a relevant $G$-space. The recurrence formula \eqref{eq:prelim:burnside-mult} also, by definition and design, only yields non-trivial coefficients for conjugacy classes of subgroups which are orbit types in the associated product space of $G$-orbits.

\subsection{Regular normal maps}
Much as the notion of a regular map gives us a guarantee that zeroes are isolated regular values, permitting us to construct a homotopy between a map and its own linearization, which is a linear isomorphism, we require an additional notion which will allow us to canonically separate orbits of zeroes of a regular map. Recall that for a submanifold $M \subset V$ and $x \in M$, we denote by $\tau_x(M)$ the tangent space to $M$ at $x$, and by $\nu_x(M)$ the normal space of $M$ at $x$, and $\tau_x(M) \oplus \nu_x(M) = V$. If $(H)$ is an orbit type of $G$ in $V$, then the space $V_{(H)}$ is a $G$-invariant submanifold of $V$. In particular, if $\Omega\subset V$ is an open $G$-invariant set, then $\Omega_{(H)} := \Omega \cap V_{(H)}$ is a $G$-invariant submanifold.
\vs
If $f:V\to V$ is a $G$-equivariant $\Omega$-admissible map, we say that $f$ is \emph{normal} in $\Omega$ if for every orbit type $(H) \in \Phi(G;\Omega)$ and every $x \in f^{-1}(0)\cap \Omega_H$, the following $(H)$-normality condition is satisfied: There exists some $\delta_x>0$ such that for all $w \in \nu_x(\Omega_H)$ with $\norm w < \delta_x$, $f(x+w)=w$.
\vs
We say that $f$ is a \emph{regular normal map} in $\Omega$ if $f$ is of class $C^1$, $f$ is normal in $\Omega$, and for every $(H) \in \Phi(G;f^{-1}(0)\cap\Omega),$ zero is a regular value of $f^H:= f\vert_{\Omega^H}$. 
\vs
By the \emph{regular normal approximation theorem} (cf. \cite{BalanovEtAl2025} Theorem D.26), if $f$ is a $G$-equivariant map and $\Omega$ is an open, bounded $G$-invariant subset such that for all $x \in \partial \Omega$, $f(x) \neq 0$, then for all $\varepsilon>0$, there exists a regular normal $G$-equivariant map $f_0$ such that
\[
\sup_{x\in \Omega}\norm{f(x)-f_0(x)} < \varepsilon.
\]
This also implies that $f$ is $G$-equivariantly $\Omega$-admissibly homotopic to $f_0$. Regular normality guarantees that orbits of solutions can be cleanly separated and thus allows orbits of zeroes to be counted. 
\vs
\begin{remark}\rm
We note that this decomposition of $V$ into orthogonal subspaces $\tau_x(M)$ and $\nu_x(M)$ of the submanifold $M$ implicitly relies on the existence of a $G$-invariant inner product. Such an inner product can always be defined on Euclidean or Hilbert spaces via averaging over the $G$-action, and the exact form of this inner product is not necessary for the above construction. 
\vs
The lack of such an inner product in the general Banach space setting of the Leray-Schauder degree and metric space setting of the Nussbaum-Sadovskii degree is not an obstacle becuase, as was shown in the non-equivariant case, the nontrivial part of the degree reduces to the Brouwer degree on the generalized eigenspace corresponding to the negative spectrum, which is finite for both compact perturbations of identity and condensing perturbations of identity, respectively.
\end{remark}
\subsection{The equivariant Brouwer degree}
We are now in a position to formulate the equivariant Brouwer degree. Let $G$ be a compact Lie group and $V$ a finite-dimensional $G$-representation, which we assume to be orthogonal without loss of generality. Let $f:V\to V$ be a continuous, $G$-equivariant map, and $\Omega\subset V$ an open, bounded $G$-invariant subset. $(f,\Omega)$ is called an \emph{admissible $G$-pair} if $f(x)\neq 0$ for all $x \in \partial \Omega$. Similar to before, we denote the set of all admissible pairs on $V$ by $\mathscr M^G(V)$, and the set of all admissible $G$-pairs by
\[
\mathscr M^G := \bigcup_V \mathscr M^G(V).
\]
Given $(f_0,\Omega),(f_1,\Omega)\in \mathscr M^G$, a continuous $G$-equivariant $\Omega$-admissible map $h_t:[0,1]\times V\to V$ is called a $G$-equivariant $\Omega$-admissible homotopy (or an $\Omega$-admissible $G$-homotopy) between $f_0$ and $f_1$ if $h_0 = f_0$ and $h_1 = f_1$.
\vs
The $G$-equivariant Brouwer degree is the unique function $\gdeg : \mathscr M^G\to A(G)$ satisfying the following properties:
\begin{enumerate}[label=\textbf{(G.\arabic*)}]

\item \textbf{(Additivity)} Let $\Omega_{1}$ and $\Omega_{2}$
be two disjoint open $G$-invariant subsets of $\Omega$ such that
$f^{-1}(0)\cap\Omega\subset\Omega_{1}\cup\Omega_{2}$. Then
\begin{align*}
\gdeg(f,\Omega)=\gdeg(f,\Omega_{1})+\gdeg
(f,\Omega_{2}).
\end{align*}

\item \textbf{(Homotopy)} If $h:[0,1]\times V\to V$ is an
$\Omega$-admissible $G$-homotopy, then
\begin{align*}
\gdeg(h_{t},\Omega)=\mathrm{constant}.
\end{align*}

\item \textbf{(Normalization)} Let $\Omega$ be a $G$-invariant
open bounded neighborhood of $0$ in $V$. Then,
\begin{align*}
\gdeg(\id,\Omega)=(G).
\end{align*}

\item \textbf{(Existence)} If $\gdeg(f,\Omega)\ne
0$, i.e. there is some $(H) \in \Phi_0(G)$ such that $\text{coeff}^H(\gdeg(f,\Omega))\neq 0$, then there
exists $x\in\Omega$ such that $f(x)=0$ and $(G_{x})\geq(H)$.

\item \textbf{(Excision)}
Suppose that $(f,\Omega)\in \mathscr M^G$ and $\Omega_0$ is an open $G$-invariant subset of $\Omega$ such that $f^{-1}(0) \cap \Omega \subset \Omega_0$. Then $\gdeg(f,\Omega) = \gdeg(f,\Omega_0)$.

\item \textbf{(Multiplicativity)}\label{prelim:equi-brouwer:p6} For any $(f_{1},\Omega
_{1}),(f_{2},\Omega_{2})\in\mathscr{M} ^{G}$,
\begin{align*}
\gdeg(f_{1}\times f_{2},\Omega_{1}\times\Omega_{2})=
\gdeg(f_{1},\Omega_{1})\cdot \gdeg(f_{2},\Omega_{2}),
\end{align*}
where the product is taken in the Burnside ring $A(G )$.

\item \textbf{(Recurrence Formula)}\label{prelim:equi-brouwer:p7} For an admissible $G$-pair
$(f,\Omega)$, put $n_H := \text{coeff}^H(\gdeg(f,\Omega))$. Then the $G$-degree can be computed using the
following recurrence formula:
\begin{equation}
\label{eq:RF-0}n_{H}=\frac{\deg(f^{H},\Omega^{H})- \sum_{(K)>(H)}{n_{K}\,
n(H,K)\, \left|  W(K)\right|  }}{\left|  W(H)\right|  },
\end{equation}
where $\left|  X\right|  $ stands for the number of elements in the set $X$
and $\deg(f^{H},\Omega^{H})$ is the Brouwer degree of the map $f^{H}%
:=f|_{V^{H}}$ on the set $\Omega^{H} := \Omega \cap V^{H}$.
\end{enumerate}
\vs
\subsection{Leray-Schauder and Nussbaum-Sadovskii extensions}
The equivariant Brouwer degree can be directly extended to the equivariant Leray-Schauder degree (resp. Nussbaum-Sadovskii degree) by considering a Banach space (resp. metric space) $\mathscr E$ instead of the finite-dimensional vector space $V$, a compact perturbation of identity (resp. condensing perturbation of identity) $\mathscr F:\mathscr E\to \mathscr E$ instead of the continuous map $f:V\to V$. In this case, the properties of the degree remain the same.

\section{Computing the equivariant Brouwer degree}
The recurrence formula \ref{prelim:equi-brouwer:p7} gives a way of computing the coefficients of orbit types in terms of the local Brouwer degree of $f^H:= f\vert_{V^H}$ restricted to the fixed point subspace $\Omega^H$. However, computing this seems to require us to know \emph{a priori} which fixed point subspaces contain zeroes and which do not. Even by passing to the Fr\'echet derivative $A:=Df(x_0)$ taken at a regular point $x_0\in f^{-1}(0)$, explicitly computing the fixed point subspaces and restricted maps for an arbitrary orbit type $(H)\in \Phi_0(G;V)$ is not very tractable. The solution is to leverage the multiplicativity property \ref{prelim:equi-brouwer:p6} to decompose $V$ into subspaces and $A$ into restricted maps where these computations are possible. Schur's lemma provides a natural choice for such subspaces: the irreducible representations of $G$, on which the restriction of any $G$-equivariant linear isomorphism $A$ acts as a homothety.

\subsection{\texorpdfstring{$G$-spaces}{G-spaces} and isotypic components}
Let $G$ be a compact Lie group and $V$ an isometric $G$-space (i.e. an isometric $G$-representation, possibly an infinite-dimensional Banach space). $V$ admits a decomposition into irreducible $G$-representations $\mathcal V_j$ indexed by $j$. Though this decomposition is not necessarily unique, it can be made unique if we group irreducible representations together into \emph{$G$-isotypic components} $V_j$, where 
\[
V_j := \underbrace{\mathcal V_j \oplus \mathcal V_j \oplus \dots \oplus \mathcal V_j}_{\text{$m_j$ times}},\quad \text{ and }\quad  m_j := \frac{\dim V_j}{\dim \mathcal V_j}
\]
where $m_j$ denotes the \emph{isotypic multiplicity} of $V_j$, and we say that $V_j$ is \emph{modeled on the irreducible representation $\mathcal V_j$}. This \emph{$G$-isotypic decomposition} is unique, and we can write
\[
V = \bigoplus_j V_j.
\]
Only real irreducible representations are of interest from the point of view of the Brouwer degree. To see why this is the case, consider $A\in \operatorname{GL}(n,\mathbb C)$. We can produce its realification $A_{\mathbb R} \in \operatorname{GL}(2n,\mathbb R)$ by taking each complex entry and expanding it to the $2\times2$ block $\left(\begin{smallmatrix}a &-b\\b &a\end{smallmatrix}\right)$. Note that for each such block, $\left\vert \begin{smallmatrix}a &-b\\b &a\end{smallmatrix}\right\vert >0$, and $\det_\mathbb R A_\mathbb R = \vert \det_\mathbb C A\vert^2 >0$. The Brouwer degree is defined in terms of the real determinant, and therefore this implies that $\operatorname{GL}(n,\mathbb C)\subset \operatorname{GL}^+(2n,\mathbb R)$. In other words, for real-linear maps, it is possible to distinguish orientation-preserving and orientation-reversing $G$-equivariant linear maps into two connected components (i.e. homotopy classes), but all complex-linear maps are orientation-preserving when viewed in terms of the real determinant. 
\vs
\subsection{Schur's lemma}
Schur's lemma is the key result that allows us to decompose a $G$-equivariantl linear isomorphism into a direct sum of linear homotheties over real irreducible representations, where equivariant degree computations are significantly more tractable. We will formulate it here in the following way:
\begin{lemma}
Let $G$ be a compact Lie group, and $\mathcal V$ an irreducible real $G$-representation. Then the set of all $G$-equivariant linear maps $\text{End}_G(\mathcal V) := \{A\in \mathrm{GL}(\mathcal V):\forall\,g\in G,\;gA=Ag\}$ is a finite-dimensional real division algebra, i.e. 
\[
\text{End}_G(\mathcal V) \cong \mathbb D,\quad \text{ where }\mathbb D\in\{\mathbb R,\mathbb C,\mathbb H\},
\]
which implies that for all $A \in \text{End}_G(\mathcal V)$, $A= \mu\id$, where $\mu \in \mathbb D$.
\end{lemma}
By the Frobenius theorem, the only finite-dimensional real division algebras are $\mathbb R$, $\mathbb C$, and $\mathbb H$.  If $\text{End}_G(\mathcal V) \cong \mathbb R$, then $\mathcal V$ is said to be of \emph{real type}; if $\text{End}_G(\mathcal V)\cong \mathbb C$, then $\mathcal V$ is said to be of \emph{complex type}; and if $\text{End}_G(\mathcal V)\cong \mathbb H$, then $\mathcal V$ is said to be of \emph{quaternionic type}. This also means that for any $G$-equivariant linear isomorphism $A:V\to V$, we can write
\[
A\vert_{\mathcal V} := \lambda \id\vert_{\mathcal V},
\]
where $\lambda \in \mathbb R, \mathbb C \text{ or } \mathbb H$, depending on the type of the irreducible representation $\mathcal V$. However, if $\mathcal V$ is of complex or quaternionic type, then the same issues apply as for complex representations, i.e. its real determinant is always positive, and so it contributes trivially to the degree by the multiplicativity property. 
\vs
\subsection{Representations of real, complex, and quaternionic type}
Another way of viewing the type of a representation is by examining its complexification. We define the complexification $\mathcal V^{\mathbb C} := \mathbb C\otimes_{\mathbb R} \mathcal V$. Then we can look at the decomposition of $\mathcal V^{\mathbb C}$ into irreducible $G$-representations. There are three possibilities:
\begin{enumerate}
    \item $\mathcal V^{\mathbb C}$ is an irreducible complex $G$-representation. In this case, $\mathcal V$ is of real type, and it is also called an \emph{absolutely irreducible real representation}.
    \item $\mathcal V^{\mathbb C}$ splits as $\mathcal V^{\mathbb C} = \mathcal U \oplus \overline{\mathcal U}$, where $\mathcal U$ is an irreducible complex representation and $\overline {\mathcal U}$ is its complex conjugate, and $\mathcal U \not\cong \overline{\mathcal U}$ as complex G-linear maps. In this case, $\mathcal V$ is of complex type.
    \item $\mathcal V^{\mathbb C}$ splits as $\mathcal V^{\mathbb C} = \mathcal W \oplus \mathcal W$, where $\mathcal W$ is self-conjugate, but $W$ does not come from a real representation (i.e. it has no $G$-invariant real structure). In this case, $\mathcal V$ is of quaternionic type.
\end{enumerate}
We are often interested in product groups of the form $G:= G_1 \times G_2$, where $G_1$ and $G_2$ are compact Lie groups. In this case, $G$ is also a compact Lie group, and if $\mathcal V_1$ and $\mathcal V_2$ are irreducible $G_1$ and $G_2$ representations, respectively, then $\mathcal V := \mathcal V_1 \otimes \mathcal V_2$ is an irreducible $G$ representation. $\mathcal V$ is of real type if and only if both $\mathcal V_1$ and $\mathcal V_2$ are also of real type.
\vs
Therefore, in light of these considerations, if $V$ has a $G$-isotypic decomposition into $V_0 \oplus \dots \oplus V_r$, where each $V_j$ is an isotypic component modeled on an absolutely irreducible real representation, then any $G$-equivariant linear isomorphism can be block diagonalized into the form
\[
A =\begin{bmatrix}
    A_0&&&\\
    &A_1&&\\
    &&\ddots\\
    &&&A_r
\end{bmatrix}\quad \text{ where } A_j := A\vert_{V_j}.
\]
Since $V_j = \mathcal V_j \otimes \mathbb R^{m_j}$, we have $A\vert_{V_j}=\id_{\mathcal V_j} \otimes B_j$, where $B_j \in \operatorname{GL}(m_j,\mathbb R)$ is called the \emph{multiplicity matrix} and describes the (generically different) real-scalar homothetic action of $A_j$ on each copy of the absolutely irreducible representation $\mathcal V_j$ in $V_j$. $B_j$ also describes the way that $A_j$ permutes the copies of $\mathcal V_j$ within $V_j$. However, we can always choose a basis of $\mathbb R^{m_j}$ in which $B_j$ is in real Jordan form or diagonalizable. Put $d_j := \dim \mathcal V_j$. Then $\det A_j = (\det B_j)^{d_j}$. This implies that only absolutely irreducible real representations of odd dimension can contribute nontrivially to the degree, and only when $\det B_j <0$
\vs
In subsequent chapters, we will commonly make the simplifying assumption that for $\mu_j \in \sigma(A)$, $E(\mu_j) = V_j$. In this case, we obtain $\det A = \prod_{j=0}^r (\det A\vert_{\mathcal V_j})^{m_j}$. This implies that the only irreducible representations which contribute nontrivially to the degree are those of real type, of odd dimension, having odd isotypic multiplicity, and whose associated eigenvalues belong to the negative spectrum of A. This assumption merely makes the form of some equivariant degree computations more convenient. If we do not make this assumption, then we must simply take the $G$-isotypic decomposition of each generalized eigenspace $E(\mu_j)$ where $\mu_j \in \sigma_-(A)$, count how many of each irreducible representations appear, and thereby find a modified isotypic multiplicity tracking how many of each irreducible representation in each isotypic component correspond to negative eigenvalues of $A$ and satisfy the oddness requirements for their dimensionality and isotypic multiplicity. Note that this can also be accomplished by considering the $G$-isotypic decomposition of $E^- := \bigcup_{\mu\in \sigma_-(A)} E(\mu)$ and computing isotypic multiplicities relative to this space instead of $V$.

\subsection{Basic maps and basic degrees}
Let us consider a single irreducible $G$-representation of real type denoted $\mathcal V_j$. Let $A:V\to V$ be a $G$-equivariant linear isomorphism, suppose that $V$ contains $\mathcal V_j$ as a subspace, and assume that $E(\mu) = \mathcal V_j$ for some $\mu \in \sigma_-(A)$. The multiplicativity property of the equivariant Brouwer degree allows us to consider the contribution of $\deg(A\vert_{\mathcal V_j},B(\mathcal V_j))$ to the overall degree $\deg(A,B(V))$ by multiplication in the Burnside ring. Since $A\vert_{\mathcal V_j} = \mu_j\id_{\mathcal V_j}$, and $\mu_j <0$, $\mu_j\id$ is $B(\mathcal V_j)$-admissibly homotopic to $-\id_{\mathcal V_j}$. We call $-\id_{\mathcal V_j}$ a \emph{basic map} and $\deg_{\mathcal V_j}:=\deg(-\id_{\mathcal V_j},B(\mathcal V_j))$ the \emph{$\mathcal V_j$-basic degree}.
\vs
The recurrence formula \ref{prelim:equi-brouwer:p7} gives us a simple way to calculate the coefficient of any orbit type in this basic degree. For any $(H) \in \Phi_0(G;\mathcal V_j)$, note that $\deg(-\id_{\mathcal V_j}^H,B(\mathcal V_j)^H) = (-1)^{\dim({\mathcal V_j}^H)}$. Putting $n_H := \text{coeff}^H(\deg_{\mathcal V_j})$, we obtain the recurrence formula for orbit type coefficients in basic degrees:
\[
n_H = \frac{(-1)^{\dim({\mathcal V_j}^H)} - \sum_{(K)>(H)}n_K n(H,K)\vert W(K)\vert}{\vert W(H) \vert}.
\]
This immediately implies that $n_G = 1$ for any basic degree. This logically corresponds to the fact that every $G$-invariant linear subspace contains one solution with isotropy exactly $G$, namely $0$. If $(H) \in \mathfrak M(G;\mathcal V_j)$ (i.e. $(H)$ is a maximal orbit type), then because $n(H,G) = 1$, $n_G = 1$, and $\vert W(G)\vert=1$, this simplifies to
\[
n_H = \frac{(-1)^{\dim({\mathcal V_j}^H)} - 1}{\vert W(H) \vert}.
\]
Note that the numerator $(-1)^{\dim({\mathcal V_j}^H)} - 1$ must equal 1 or 2. Since the degree coefficient must be an integer, this implies that for a maximal orbit type $(H)$, we must have $\vert W(H)\vert =1$ or $2$.
\vs
Therefore, the computation of the equivariant degree around an isolated zero $x_0$ of a regular normal $G$-equivariant map $f:V\to V$ which is Fr\'echet differentiable ar $x_0$ can be summarized as follows:
\begin{enumerate}
    \item Take the Fr\'echet derivative $A:=Df(x_0)$ and compute its negative spectrum $\sigma_-(A)$. 
    \item Take the $G$-isotypic decomposition of $E^- := \bigcup_{\mu\in \sigma_-(A)}E(\mu)$. Put $E^- := V_0\oplus V_1 \oplus \dots \oplus V_r$, where each $V_j$ is modeled on the absolutely irreducible real representation $\mathcal V_j$ and has isotypic multiplicity $m_j$.
    \item For each $\mathcal V_j$, compute $\deg_{\mathcal V_j}$.
    \item Compute $\gdeg(A,B(V))$ using the multiplicativity property:
    \[
        \gdeg(f,B_\varepsilon(x_0)) = \gdeg(A,B(V)) = \prod_{j=0}^r (\deg_{\mathcal V_j})^{m_j}.
    \]
\end{enumerate}
Step 2 can be modified depending on simplifying assumptions of how the eigenspaces of $\mu_j \in \sigma_-(A)$ coincide with the isotypic components, or can be computed in terms of the determinants of the multiplicity matrices $B_j$ for the isotypic component $V_j = \mathcal V_j \otimes B_j$. We also note that the above procedure applies equally to the equivariant Brouwer, Leray-Schauder, and Nussbaum-Sadovskii degrees, since it relies only on the properties of the negative spectrum of $A$, which is finite-dimensional for compact and condensing perturbations of identity.

\section{Using the equivariant degree}
Here we will discuss the two most common ways in which the equivariant degree is actually applied. Our discussions of the degree so far have consisted in finding the degree of a single isolated zero of a function satisfying certain assumptions. However, this obviously requires us to know that a zero exists, which seems to be the main point of the degree. The systems we study in these chapters are autonomous systems of differential equations of various types. Such systems, when formulated as compact or condensing perturbations of identity on functional spaces, typically admit a trivial equilibrium solution $x \equiv 0$ (or some other constant solution which can be shifted to the origin). Indeed, let $G$ be a compact Lie group, $\mathscr E$ a Banach space, and $\mathscr F:\mathscr E\to \mathscr E$ a compact or condensing perturbation of identity ($\mathscr F:= \id -F$ where $F:\mathscr E\to \mathscr E$ is a compact or condensing map), such that $\mathscr F(0) = 0$.
\vs
Taking the linearization around this zero, $\mathscr A := D\mathscr F(0)$, and analyzing its spectrum is already the basis of many types of useful analysis for such systems (e.g. showing local asymptotic stability of the trivial equilibrium). Using this same information, one can use the method of basic degrees outlined above to find the $G$-degree of this trivial solution. From here, two distinct possibilities emerge:
\vs
\subsection{Existence via \emph{a priori} bounds}
$\mathscr F$ is said to admit \emph{a priori} bounds if there is some $R>0$ such that for all $x_0 \in \mathscr F^{-1}(0)$, $\norm {x_o} < R$. The existence of such\emph{a priori} bounds can be obtained from, for exmaple, a Nagumo condition. Then $\mathscr F$ is $B_R(0)$-admissible, and  it is possible to construct a $B_R(0)$-admissible $G$-equivariant homotopy to a map reflecting the asymptotic behavior of $\mathscr F$, whose $G$-degree can be conveniently computed. 
\vs
For example, if $F$ is sublinear, i.e. $\norm {F(x)}/\norm{x} \to 0$ as $x \to \infty$, then $h_t := x-tF(x)$ is a $B_R(0)$-admissible $G$-homotopy to identity, and therefore
\[
\gdeg(\mathscr F,B_R(0)) = (G).
\]
If $0$ is a regular isolated zero of $\mathscr F$ and the corresponding linearization $\mathscr A$ exists and is a $G$-equivariant linear isomorphism, then there exists $r>0$ sufficiently small where $\mathscr F$ is $B_r(0)$-admissibly homotopic to $\mathscr A$, and
\[
\gdeg(\mathscr F,B_r(0)) = \gdeg(\mathscr A,B_r(0)),
\]
where $\gdeg(\mathscr A,B_r(0))$ is computed as discussed in the previous section.
Then, by the additivity property, the degree on the annulus $\Omega := B_R(0)\setminus B_r(0)$ can be computed
\[
\gdeg(\mathscr F,\Omega) = \gdeg(\mathscr F,B_R(0)) - \gdeg(\mathscr F,B_r(0)).
\]
The nontriviality of this degree then implies the existence of nontrivial solutions in $\Omega$. Moreover, for each orbit type $(H)\in \Phi_0(G;\mathscr E)$ such that $\text{coeff}^H(\gdeg(\mathscr F,\Omega)) \neq 0$, there exists a solution $x_0 \in \Omega$ satisfying $(G_{x_0}) \geq (H)$. This also implies that if $(H)$ is a maximal orbit type, then $(G_{x_0}) = (H)$.

\begin{figure}[h]
\centering
\begin{tikzpicture}[
  scale=1.2,
  annulus/.style={pattern=north east lines},
  big circle/.style={draw=black!70!black, thick, dashed},
  small circle/.style={draw=black!70!black, thick, dashed},
  label style/.style={font=\footnotesize\itshape}
]

% Coordinates
\coordinate (O) at (0,0);

% The annulus (shaded between big and small circle)
% \fill[annulus, even odd rule] (O) circle[radius=3] (O) circle[radius=1.2];

% Draw the big and small circles
\draw[big circle]   (O) circle[radius=3];
\draw[small circle] (O) circle[radius=0.6];

% Center point
\fill (O) circle[radius=1.2pt] node[below left] {$0$};

% Radius arrows
% R arrow
% \draw[->, shorten >=2pt, shorten <=2pt] (O) -- node[above] {$R$} +(45:3);
% r arrow
% \draw[->, shorten >=2pt, shorten <=2pt] (O) -- node[above] {$r$} +(-30:1.2);

% Labels for the sets
\node[label style] at (2.4,2.4) {$B_R(0)$};
\node[label style] at (0.6,0.75) {$B_r(0)$};
\node[label style, fill=white, inner sep=1pt] at (2,0) {$\Omega$};

\end{tikzpicture}
\end{figure}

\subsection{One-parameter bifurcation}\label{sec:prelim:one-parameter}
The other common approach is to use the degree to detect bifurcating branches of solutions. Let $\mathscr F(\alpha,x) : \mathbb R \times \mathscr E \to \mathscr E$ be a $G$-equivariant compact or condensing perturbation of identity parametrized by $\alpha$, and such that $\mathscr F(\alpha,0)=0$ for all $\alpha$. We call the set $\{(\alpha,0)\in \mathbb R \times \mathscr E\}$ the \emph{trivial branch}. The \emph{one-parameter bifurcation} problem is to find critical values $\alpha=\alpha_0$ for which there are continua of nontrivial solutions bifurcating out from the trivial branch. However, this approach requires some subtlety. First of all, if there is a bifurcating branch of solutions at $\alpha=\alpha_0$, then this implies that $\mathscr A_{\alpha_0} = D_x\mathscr F(\alpha_0,0)$ fails to be an isomorphism because it has a nontrivial kernel, the nonlinear continuation of which comprises the bifurcating branch. In such a case, $\mathscr A_{\alpha_0}$ is not an admissible map on any open bounded subset and so its degree is not well defined. 
\vs
We therefore require an approach which allows us to construct some admissible map whose degree equals the degree of this bifurcating branch. It perhaps seems intuitively appealing that if $\mathscr F(\alpha_0,0)$ is a bifurcation point (which will be formally defined shortly) where the linearization $\mathscr A$ fails to be an isomorphism and a branch of nontrivial solutions emerges, then the degree of the bifurcating branch ought to be somehow related to the difference 
\begin{equation}\label{eq:prelim:bifurcation:gdeg-diff}
\gdeg(\mathscr F(\alpha-\varepsilon,\cdot),\Omega) - \gdeg(\mathscr F(\alpha + \varepsilon,\cdot),\Omega)
\end{equation}
for $\varepsilon>0$ sufficiently small, and $\Omega$ sufficiently small. In other words, with an intuition of the $G$-degree as a sort of ``conserved quantity'' along the parameter range, it seems like if the negative spectrum $\sigma_-(\mathscr A)$ loses an eigenvalue as $\alpha$ sweeps across the critical point $\alpha_0$, then that missing eigenvalue (and its associated eigenspace) should emerge as a branch of bifurcating solutions locally transverse to the trivial branch, and the nonlinear continuation of that eigenspace comprises the invariant submanifold of the bifurcating branch of solutions. 
\vs
The local bifurcation theorem of M.A. Krasnosel'skii makes this intuitive notion precise. Crucially, it rigorously justifies the fact that \eqref{eq:prelim:bifurcation:gdeg-diff} actually corresponds \emph{exactly} to the $G$-degree of the bifurcating branch, and is sufficient to prove its existence. 
\subsection{Krasnosel'skii bifurcation theorem}\label{sec:prelim:krasnoselskii}
\begin{figure}[h]
\centering
\begin{tikzpicture}[
  scale=1,
  % annulus/.style={pattern=north east lines},
  % big circle/.style={draw=black!70!black, thick, dashed},
  % small circle/.style={draw=black!70!black, thick, dashed},
  label style/.style={font=\footnotesize\itshape}
]

% Coordinates
\coordinate (O) at (0,0);

% The annulus (shaded between big and small circle)
% \fill[annulus, even odd rule] (O) circle[radius=3] (O) circle[radius=1.2];

% Draw the big and small circles
\draw[gray, dashed] (3,2) arc(90:270:0.5cm and 2cm);
\draw[black, dashed] (3,-2) arc(270:360:0.5cm and 2cm) arc(0:90:0.5cm and 2cm);
\draw[dashed] (-3,0) ellipse (0.5cm and 2cm);
% \draw (3,0) ellipse (0.5cm and 2cm);

\draw (-3,2) -- (3,2);
\draw (-3,-2) -- (3,-2);

\draw[blue] (-6,0) -- (-3,0);
\draw[blue] (3,0) -- (6,0);
\draw[blue, thick] (-3,0) -- (3,0);

\draw[blue,thick] (0,0) .. controls (0,1) and (1,2) .. (1,2);

\draw[blue] (1,2) .. controls(1.4,2.4) and (4,2.5) .. (6,2.5);

\draw (-6,-3) -- (6,-3) node[right] {$\alpha$};
\draw (-3,-2.9) -- (-3,-3.1) node[below] {$a$};
\draw (3,-2.9) -- (3,-3.1) node[below] {$b$};

\filldraw[red] (0.67,1.6) circle (2pt);
\filldraw[black] (1,2) circle (2pt);

% \draw (-6.5,-2.5) -- (-6.5,2) node[left] {$V$};
\draw[dashed] (-5.5,2) -- (-6.5,2.5) -- (-6.5,-2) -- (-5.5,-2.5) -- (-5.5,2) node[below left] {$V$};

\filldraw[black] (-3,0) circle (2pt);
\filldraw[black] (3,0) circle (2pt);
\filldraw[blue] (-6,0) circle (1pt);

\node[label style] at (-3,0.6) {$\Omega$};
\node[label style] at (3,0.6) {$\Omega$};
\node[label style] at (0,2.2) {$\partial\Omega$};

\end{tikzpicture}
\caption{The horizontal blue line represents a branch of solutions in $\mathbb R\times V$ which persists for all parameter values. The cylinder represents the neighborhood $(a,b) \times \Omega$. The black points are values $(\alpha,x)\in \{a,b\}\times \Omega \cup (a,b)\times \partial \Omega$ where $\mathscr F(\alpha,x) = 0$, but $\eta(\alpha,x)\neq 0$. The red point is $(\alpha_0,x_0)\in (a,b)\times \Omega$ such that $\mathscr F(\alpha_0,x_0)=0$ and $\eta(\alpha_0,x_0)=0$.}\label{fig:prelim:eta}
\end{figure}
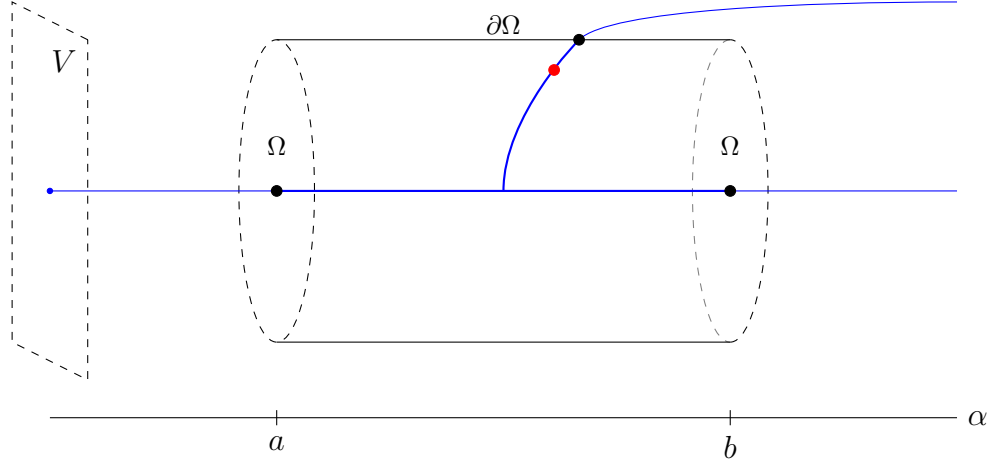

For the sake of convenience, we will formulate this theorem and its preliminaries in terms of the equivariant Leray-Schauder degree for compact perturbations of identity, but note that all these definitions and results can be trivially extended to the equivariant Nussbaum-Sadovskii degree for condensing perturbations of identity.
\vs
As usual, let $G$ be a compact Lie group and $\mathscr E$ an isometric Banach $G$-representation. Let $\Omega\subset \mathscr E$ be an open bounded $G$-invariant subset, and $\mathscr F:\mathbb R \times \mathscr E \to \mathscr E$ a compact perturbation of identity such that for $a<b$, $(\mathscr F_a,\Omega)$ and $(\mathscr F_b,\Omega)$ are admissible $G$-pairs, where $\mathscr F_a := \mathscr F(a,\cdot)$. A continuous $G$-equivariant function $\eta:\mathbb R \times \mathscr E \to \mathscr E$ is called an \emph{$\Omega$-auxiliary function} for $\mathscr F_a$ on the interval $[a,b]$ if
\[
\begin{cases}
    \eta(\alpha,x) < 0 &\text{ if } \alpha\in \{a,b\}, x\in \Omega, \mathscr F_\alpha(x) = 0\\
    \eta(\alpha,x) > 0 &\text{ if } \alpha\in (a,b), x\in \partial\Omega.
\end{cases}
\]
Given such an $\Omega$-auxiliary function, we define the \emph{complemented map} $\mathscr F_\eta:\mathbb R\times \mathscr E \to \mathbb R \times \mathscr E$ by $\mathscr F_\eta(\alpha,x) = (\eta(\alpha,x),\mathscr F(\alpha,x))$. Clearly $(\mathscr F_\eta,(a,b)\times \Omega)\in \mathscr M^G$, and so $\gdeg(\mathscr F_\eta,(a,b)\times \Omega)$ is well-defined. Given two $\Omega$-auxiliary functions $\eta_1$ and $\eta_2$ for $\mathscr F_\alpha$ with $\alpha \in[a,b]$, the linear homotopy $\eta_\tau(\alpha,x) = (1-\tau)\eta_1 + \tau\eta_2$ is itself an $\Omega$-auxiliary function, and therefore the value of $\gdeg(\mathscr F_\eta,(a,b)\times \Omega)$ is independent of the choice of auxiliary function.
\vs
We also have the following proposition (proved in \cite{BalanovEtAl2025}, pg. 163):
\begin{proposition}\label{prop:prelim:bif-invariant}
    Let $\mathscr F:\mathbb R \times \mathscr E\to \mathscr E$ be a $G$-equivariant compact perturbation of identity such that $(\mathscr F_a,\Omega),(\mathscr F_b,\Omega) \in \mathscr M^G(\mathscr E)$ and $\eta:\mathbb R \times \mathscr E\to \mathscr E$ is an $\Omega$-auxiliary function for $\mathscr F_\alpha$ with $\alpha \in[a,b]$. Then
    \[
    \gdeg(\mathscr F_\eta,(a,b)\times \Omega) = \gdeg(\mathscr F_a,\Omega) - \gdeg(\mathscr F_b,\Omega).
    \]
\end{proposition}
The proof uses the homotopy, additivity, and excision properties of the degree to find open sets inside $U:=(a,b)\times \Omega$ where $\mathscr F_\eta$ can be made admissibly homotopic to $\mathscr F_a$ and $\mathscr F_b$, and shows that $\mathscr F_\eta$ can have no other zeroes outside of these sets. Rather than reproducing the full proof, we will instead discuss a few properties of the complemented map and describe why it can accomplish what we need for bifurcation. First, notice that $\partial U = \{a,b\}\times \Omega \cup (a,b) \times \partial \Omega$, and so $\eta$ is $U$-admissible by construction. This also means that $\mathscr F_\eta(\alpha,x)\neq 0$ for all $\alpha,x$ on the ends of the cylinder as shown in Figure \ref{fig:prelim:eta}, so $\gdeg(\mathscr F_\eta,U)$ will not detect the existing branch of solutions from which new solutions bifurcate. On the other hand, the sign conditions on $\eta$ force a zero surface $\{(\alpha,x)\in \mathbb R \times \mathscr E: \eta(\alpha,x)=0\}$ to exist somewhere within the cylinder $U$ which must surround the branch which enters and exits $U$ at $(a,x_a) \in \{a\} \times \Omega$ and $(b,x_b)\in \{b\}\times \Omega$, respectively. The exact shape of this zero surface doesn't matter, but the zeroes of $\mathscr F\eta$ are precisely the points where a solution of $\mathscr F_\alpha$ intersects this zero surface.

\vs
Consider the $d$-parameter bifurcation problem $\mathscr F:\mathbb R^d \times \mathscr E \to \mathscr E$, where $\mathscr F$ is a compact perturbation of identity, and suppose there is a given $G$-invariant $d$-dimensional manifold $M \subset \mathscr F^{-1}(0)$, called the \emph{trivial solutions} of $\mathscr F(\alpha,x)=0$. We then denote by $\mathscr S$ the set of all \emph{nontrivial solutions} to $\mathscr F(\alpha,x)=0$, i.e.
\[
\mathscr S:= \{(\alpha,x)\in \mathbb R \times \mathscr E : \mathscr F(\alpha,x)=0 \text{ and }(\alpha,x)\not\in M\}.
\]
A nonempty set $\mathscr C \subset \mathscr S$ is called a \emph{branch} of nontrivial solutions if for some connected component $\mathscr D \subset \overline {\mathscr S}$, one has $\mathscr C = \mathscr S \cap \mathscr D$. If $(\alpha_0,x_0) \in \overline{\mathscr C}\cap M$, then we say $\mathscr C$ \emph{bifurcates} from $(\alpha_0,x_0)$ and $(\alpha_0,x_0)$ is called a \emph{branching point}. A branching point is called a \emph{bifurcation point} if for any neighborhood $U$ containing $(\alpha_0,x_0)$, $U \cap \mathscr S \neq \emptyset$. 
\vs
Given a group $H\leq G$, a branch $\mathscr C$ is called a \emph{branch with symmetries at least $(H)$} if $\mathscr C' = \mathscr D'\cap \mathscr S^H$, where $\mathscr D'$ is a connected component of $\overline{\mathscr S^H}$, and if $(\alpha_0,x_0)\in \overline{\mathscr C'}\cap M^H$, then we say that $\mathscr C'$ \emph{bifurcates} from $(\alpha_0,x_0)$. 
\vs
We consider $\mathbb R^d$ to be a parameter space, typically of either one or two parameters, depending on the application. The precise definition of which solutions are considered to be ``trivial'' obviously depends on the context and the application, and the above definitions could also be easily generalized to describe secondary bifurcations.
\vs
We also require the equivariant Kuratowski lemma:
\begin{lemma}
Let $\mathscr E$ be an isometric Banach $G$-representation, $B_0,B_1 \subset \mathscr E$ two $G$-invariant disjoint closed sets in $\mathscr E$, and $K$ a $G$-invariant compact set in $\mathscr E$ such that $K \cap B_0 \neq \emptyset \neq K \cap B_1$. If the set $K$ does not contain a connected component $C$ such that $C\cap B_0 \neq \emptyset \neq C\cap B_1$, then there exist two $G$-invariant disjoint open sets $V_0,V_1$ such that $B_0 \subset V_0, B_1 \subset V_1$, and $B_0\cup B_1 \cup K \subset V_0 \cup V_1$. 
\end{lemma}
We now have all necessary bifurcation-related definitions to provide the equivariant version of the Krasnosel'skii bifurcation theorem. Due to the foundational nature of this result in applying topological degree methods to bifurcation problems, we will supply the proof here.
\begin{theorem}\label{thm:prelim:krasnoselksii}
Let $\mathscr F:\mathbb R \times \mathscr E \to \mathscr E$ be a compact perturbation of identity, $a<\alpha_0<b$ and $x:[a,b]\to M\cap \mathscr E^G$ a continuous curve such that
\begin{enumerate}
    \item For all $\alpha \in [a,b]$ with $\alpha\neq \alpha_0$, the point $(\alpha,x(\alpha))$ is not a bifurcation point of $\mathscr F$.
    \item There exists $\delta>0$ such that if $\alpha=a$ or $\alpha =b$ and $0<\norm {x-x(\alpha)} \leq \delta$, then $\mathscr F(\alpha,x) \neq 0$.
    \item For $\Omega_\alpha  := B_\delta(x(\alpha))$, one has
    \[
    \omega_G(\alpha_0,x_0):= \gdeg(\mathscr F_a,\Omega_a) - \gdeg(\mathscr F_b,\Omega_b) \neq 0.
    \]
\end{enumerate}
Then there exists a branch $\mathscr C$ of nontrivial solutions to $\mathscr F(\alpha,x)=0$ bifurcating from $(\alpha_0,x_0)$. Moreover, if for some orbit type $(H) \in \Phi_0(G;\mathscr E)$ one has
\[
\text{coeff}^H(\omega_G(\alpha_0,x_0)) \neq 0,
\]
then there exists a branch $\mathscr C'$ of nontrivial solutions with symmetries at least $(H)$, bifurcating from $(\alpha_0,x_0)$.
\end{theorem}
\begin{proof}
To simplify the proof, we will assume that $x(\alpha)=0$ and $\Omega := \Omega_a = \Omega_b$. In this case, any $G$-equivariant continuous map $\theta:\mathbb R\times \mathscr E \to \mathscr E$ satisfying
\[
\begin{cases}
    \theta(\alpha,x) < 0  &\text{ if } \alpha \in \{a,b\} \text{ and }x(\alpha)=0\\
    \theta(\alpha,x)>0 &\text{ if }\alpha \in (a,b) \text{ and }\norm{x}=\delta
\end{cases}
\]
is an auxiliary function for $\mathscr F_\alpha$ on $[a,b]$, for $\delta>0$ sufficiently small. Pick any $\rho$ satisfying $0<\rho<\delta$ and define $\theta(\alpha,x) = \norm x - \rho$. Then $\theta$ is an auxiliary function for $\mathscr F_\alpha$ on $[a,b]$. 
\vs
Let $U:= (a,b)\times \Omega$. We claim that $\mathscr F^{-1}(0) \cap U$ contains a compact sonnected subset $K_0$ such that $K_0\cap \{(\alpha,x)\in \mathbb R\times \mathscr E:\alpha \in(a,b), \norm x = \rho\}\neq \emptyset$ and $(\alpha_0,x_0)\in K_0$, in which case there is obviously a branch $\mathscr C$ bifurcating from $(\alpha_0,x_0)$.
\vs
Suppose for contradiction that such a component does not exist. Put $K := \mathscr F^{-1}(0)\cap \overline U, A_0:= \{(\alpha_0,0)\}, A_1 := \{(\alpha,x)\in \mathbb R\times \mathscr E:\alpha\in (a,b), \norm x = \rho\}.$ Then by the equivariant Kuratowski lemma, there exist tow $G$-invariant disjoint open sets $V_0,V_1$ such that $A_0\subset V_0,A_1\subset V_1$, and $A_0\cup A_1 \cup K \subset V_0 \cup V_1$. Put $K_0 := K \cap V_0$ and $K_1 := K \cap V_1$. Clearly, $K_0$ and $K_1$ are $G$-invariant compact sets. Assume that $\mu:\mathbb R \times V\to [0,1]$ is a continuous $G$-invariant Urysohn function such that
\[
\mu(\alpha,x)=\begin{cases}
1 &\text{ if }(\alpha,x)\in K_0 \cup A_0\\
0 &\text{ if } (\alpha,x) \in K_1 \cup A_1,
\end{cases}
\]
and define the auxiliary function $\phi:\mathbb R \times \mathscr E \to \mathbb R$ by
\[
\phi(\alpha,x) := \norm x - \mu(\alpha,x)\rho,\quad (\alpha,x)\in \mathbb R\times \mathscr E.
\]
Then by Proposition \ref{prop:prelim:bif-invariant}, we have
\[
\gdeg(\mathscr F_\phi,\Omega) = \gdeg(\mathscr F_a,\Omega)-\gdeg(\mathscr F_b,\Omega) \neq 0.
\]
\end{proof}
\vs
Hence there exists an orbit type $(H)$ with $\text{coeff}^H(\gdeg(\mathscr F_\phi,\Omega))\neq 0$ and, by the existence property of the degree, a solution $(\alpha^*,x^*)\in U$ to $\mathscr F(\alpha,x)=0$ such that $\norm{x^*} - \mu(\alpha^*,x^*)\rho =0$ and $G_{(\alpha^*,x^*)}\geq H$. 
\vs
Clearly $(\alpha^*,x^*)\in K_0$ or $(\alpha^*,x^*)\in K_1$. If $(\alpha^*,x^*)\in K_1$, then $\norm{x^*} -\mu(\alpha^*,x^*)\rho = \norm{x^*} - \rho = 0$, which implies $(\alpha^*,x^*)\in K_1\cap A_1$, which is a contradiction. On the other hand, if $(\alpha^*,x^*)\in K_1$, then $\norm{x^*} - \mu(\alpha^*,x^*)\rho = \norm{x^*} = 0$, which implies $(\alpha^*,x^*)=(\alpha_9,x_0) \in K_0\cap A_0$, which is also a contradiction.
We call $\omega_G(\alpha_0,x_0)$ the \emph{local bifurcation invariant}. It is also useful to formulate some conditions under which the Krasnosel'skii bifurcation theorem can be applied and computed. We call $\Lambda$ the \emph{critical set}. It will be defined in several generally equivalent ways in subsequent chapters depending on the particular application, but a useful working definition is
\[
\Lambda := \{(\alpha,0)\in\mathbb R\times \mathscr E: D\mathscr F_\alpha(0) \text { is not an isomorphism}\}.
\]
Then clearly we need $\Lambda$ to be a discrete set in order that $\omega_G(\alpha_0,0)$ is computable for all $(\alpha_0,0)\in \Lambda$. Often we can obtain the discreteness of $\Lambda$ from the analyticity (or other properties) of the characteristic quasipolynomial of $\mathscr A$. 
\vs
Theorem \ref{thm:prelim:krasnoselksii} is often termed a ``local'' bifurcation theorem, and in non-topological formulations it only guarantees the existence of a sequence of nontrivial solutions. The topological versionn presented here, however, is actually quite global (or at least ``non-local'') in contrast with other non-topological bifurcation results. For example, it guarantees the existence of actual connected components of solutions, not just nontrivial solutions accumulating at the trivial branch. 
\vs
We may also consider how far these branches are guaranteed to extend away from the bifurcation point as parameters are varied. These questions belong to the domain of global bifurcation results, and there are two main theorems we use to show this. The first is the so-called \emph{continuation theorem}, which provides a process for extending branches.
\begin{theorem}
Let $\mathscr E$ be an isometric Banach $G$-representation and $\Omega\subset \mathscr E$ a $G$-invariant bounded open set. Assume that $\mathscr F:[a,b]\times \mathscr E\to \mathscr E, a<b$ is a compact perturbation of identity such that $(\mathscr F_\alpha,\Omega)\in \mathscr M^G$ for every $\alpha \in [a,b]$, and for some orbit type $(H) \in \Phi_0(G,\Omega)$, one has 
\[
\text{coeff}^H(\gdeg(\mathscr F_\alpha,\Omega))\neq 0, \quad \alpha\in[a,b],
\]
then there exists a nontrivial continuum $\mathscr K \subset \mathscr F^{-1}(0)\cap [a,b]\times\Omega^H$ such that
\[
\mathscr K\cap\left(\{a\}\times \Omega^H\right) \neq 0 \neq \mathscr K \cap \left(\{b\}\times \Omega^H\right).
\]
Moreover, if $(H)$ is a maximal orbit type in $\Phi_0(G;\Omega)$, points from $\mathscr K$ have symmetries exactly $(H)$. 
\end{theorem}
Note that, by constructing successive intervals $[a_1,a_2],[a_2,a_3],\dots$, one can extend branches indefinitely. This implies that the branches of solutions guaranteed by the Krasnosel'skii theorem cannot collapse or cease to exist away from the bifurcation point. They must either go to infinity (i.e. they are unbounded) or they must rejoin the trivial branch at a different bifurcation point. 
\vs
\subsection{Rabinowitz alternative}\label{sec:prelim:rabinowitz}
\begin{figure}[h]
\centering
\begin{tikzpicture}[
  scale=1,
  % annulus/.style={pattern=north east lines},
  % big circle/.style={draw=black!70!black, thick, dashed},
  % small circle/.style={draw=black!70!black, thick, dashed},
  label style/.style={font=\footnotesize\itshape},
  >=Latex
]

% Coordinates
\coordinate (O) at (0,0);

\draw[blue,thick,<->] (-6,0) -- (6,0);

\draw[blue,thick] (-4,0) .. controls (-3.5,1) and (-2.5,1) .. (-2,0);

\draw[blue,thick] (-2,0) .. controls (-1,-1.6) and (-0,-1.6) .. (1,0);

\draw[blue,thick,->] (1,0) .. controls (2,2) .. (6,2);

\draw[dashed] (-5,0) ellipse (0.5cm and 1.5cm);
\draw[dashed] (5,0) ellipse (0.5cm and 1.5cm);

\draw[dashed] (-5,1.5) -- (5,1.5);
\draw[dashed] (-5,-1.5) -- (5,-1.5);

\node[label style] at (-4,-0.3) {$(\alpha_1,0)$};
\node[label style] at (-2,-0.3) {$(\alpha_2,0)$};
\node[label style] at (1,-0.3) {$(\alpha_3,0)$};
\node[label style] at (-3,1.15) {$\mathscr C_1$};
\node[label style] at (-0.5,-0.8) {$\mathscr C_2$};
\node[label style] at (4,2.3) {$\mathscr C_3$};
\node[label style] at (6.3,0) {$M$};
\node[label style] at (0,1.8) 
{$U$};

\filldraw[blue] (-4,0) circle (2pt);
\filldraw[blue] (-2,0) circle (2pt);
\filldraw[blue] (1,0) circle (2pt);

\end{tikzpicture}
\caption{It is possible that $\mathscr C_1,\mathscr C_2,$ and $\mathscr C_3$ all have the same isotropies, and $\mathscr C_3$ is the unbounded branch guaranteed by the Rabinowitz alternative.}\label{fig:prelim:rabinowitz2}
\end{figure}
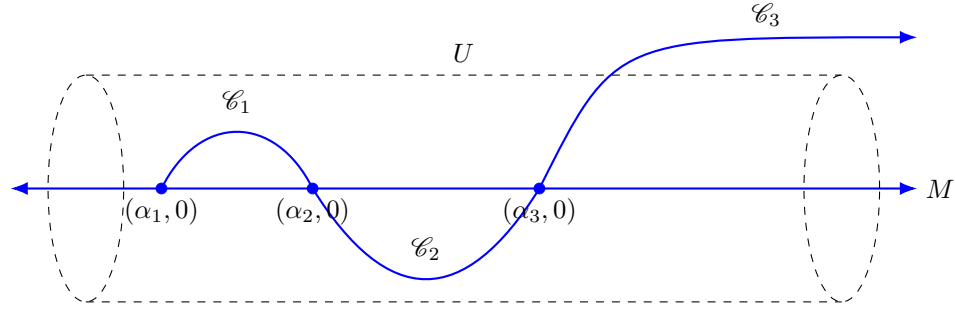

The Rabinowitz alternative formalizes these two possibilities. We present its equivariant version here, formulated in the language of a $d$-parameter bifurcation problem (where $d=1$ or $2$ for our purposes). 
\begin{theorem}\label{thm:prelim:rabinowitz}
Let $\mathscr F:\mathbb R^d\times \mathscr E \to \mathscr E$ be a compact perturbation of identity such that $\mathscr F(\alpha,0)=0$ for all $\alpha \in \mathbb R^d$ and $\mathscr F$ is Fr\'echet differentiable with respect to $x$ for all $\alpha\in \mathbb R^d$. Put $\mathscr A_\alpha := D_x\mathscr F(\alpha,0)$, and denote the critical set $\Lambda := \{(\alpha,0)\in \mathbb R^d \times \mathscr E:\mathscr A \text{ is not an isomorphism}\}$. Assume $\Lambda$ is discrete, and consider an open bounded $G$-invariant set $U \subset \mathbb R^d\times \mathscr E$ such that $\partial U \cap \Lambda = \emptyset$ and $U \cap \Lambda =\{(\alpha_1,0),(\alpha_2,0),\dots,(\alpha_{m},0)\}\neq \emptyset$ for some $m \in \mathbb N$. If $\mathscr C$ is a branch of nontrivial solutions to $\mathscr F(\alpha,x)=0$ bifurcating from $(\alpha_j,0)$, then one has the following alternative:
\begin{enumerate}
    \item either $\mathscr C\cap \partial U = \emptyset$
    \item or $\mathscr C\cap U =\{(\alpha_{j_1},0),\dots,(\alpha_{j_k},0)\}\subset \Lambda$ for some $k\in \mathbb N$, and 
    \[
    \sum_{l=1}^k \omega_G(\alpha_{j_k},0) =0.
    \]
\end{enumerate}
\end{theorem}
Henceforth, a branch $\mathscr C$ bifurcating from $(\alpha_0,0)$ will be called \emph{unbounded} if for every open bounded $G$-invariant set $U$ with $(\alpha_0,0)\in U$, $\mathscr C\cap U\neq \emptyset$.
\vs
This theorem immediately yields two useful corollaries. First of all, note that although Theorem \ref{thm:prelim:rabinowitz} is formulated for a bounded neighborhood $U$ containing only finitely many points from the critical set $\Lambda$, one can always inductively take successively larger neighborhoods. This leads to the following corollary, which we formulate here as an independent theorem guaranteeing the global unboundedness of branches:
\begin{theorem}\label{thm:prelim:rabinowitz2}
    Let $\mathscr F:\mathbb R^d \times \mathscr E \to \mathscr E$ be a compact perturbation of identity such that $\mathscr F(\alpha,0)=0$ for all $\alpha \in \mathbb R^d$, $\mathscr A_\alpha :=D_x\mathscr F(\alpha,0)$ exists for all $(\alpha,0)\in \mathbb R^d\times \mathscr E$, and the critical set $\Lambda := \{(\alpha,0)\in \mathbb R^d \times \mathscr E:\mathscr A \text{ is not an isomorphism}\}$ is discrete. Then if there is an orbit type $(H) \in \Phi_0(G;\mathscr E)$ such that
    \[
    \sum_{(\alpha,0)\in \Lambda} \text{coeff}^H(\omega_G(\alpha,0)) \neq 0,
    \]
    then there exists an unbounded branch of solutions $\mathscr C$ having symmetries at least $(H)$.
\end{theorem}
Secondly, neither Theorem \ref{thm:prelim:rabinowitz} nor Theorem \ref{thm:prelim:rabinowitz2} prevent a bifurcating branch from returning to the trivial branch (cf. Figure \ref{fig:prelim:rabinowitz2}). When we want to show that this cannot happen, we require the following stronger condition:
\begin{theorem}
    Let $\mathscr F:\mathbb R^d \times \mathscr E \to \mathscr E$ be a compact perturbation of identity such that $\mathscr F(\alpha,0)=0$ for all $\alpha \in \mathbb R^d$, $\mathscr A_\alpha :=D\mathscr F(\alpha,0)$ exists for all $(\alpha,0)\in \mathbb R^d\times \mathscr E$, and the critical set $\Lambda := \{(\alpha,0)\in \mathbb R^d \times \mathscr E:\mathscr A \text{ is not an isomorphism}\}$ is discrete. Then if there is an orbit type $(H) \in \Phi_0(G;\mathscr E)$ such that
    \[
    \sum_{(\alpha,0)\in \Lambda} \text{coeff}^H(\omega_G(\alpha,0)) \neq 0,
    \]
    then there exists an unbounded branch of solutions $\mathscr C$ having symmetries at least $(H)$.
\end{theorem}

We will often informally use the term ``Rabinowitz alternative'' to refer not only to Theorem \ref{thm:prelim:rabinowitz} but also to Theorem \ref{thm:prelim:rabinowitz2}, since the latter is more directly useful for applications. 

\subsection{Periodic solutions}
Our main interest throughout all the systems studied in this thesis is showing the existence or bifurcation of periodic solutions. The major consideration in doing this is the functional setting, and the proper temporal symmetry group. If the period of solutions can be assumed \emph{a priori}, then we can consider the functional space $\mathscr E:= C_{2\pi}(\mathbb R;\mathbb R^n)$, and if the system is time-reversible, then the temporal symmetries can be described by $O(2)$, where $\kappa\in O(2)$ acts on $x \in \mathscr E$ as $\kappa x(t) = x(-t)$. This allows the $G$-equivariant degree to be used, in either existence or bifurcation settings, because the subgroups of $O(2)$ containing $\kappa$ are of the form $D_m$, and thus have finite Weyl groups. However, if the system is not time-reversible, the relevant temporal symmetry group is $S^1 \cong SO(2)$, whose subgroups do \emph{not} have finite Weyl group. This precludes the use of the $G$-equivariant degree as defined above. To make matters worse, all nontrivial irreducible representations of $S^!$ are of complex type. These situations can be handled using the \emph{twisted degree}, which we will define in the next section.
\vs
Moreover, if the period is unknown (or equivalently, if we wish to be able to detect solutions of different periods), then we must take $u(t) := x(\beta t)$, where $\beta >0$ is now a parameter encoding the frequency of the solution. In this case, the one-parameter bifurcation approach discussed in the previous section is generalized to \emph{two-parameter bifurcation}, where symmetric Hopf bifurcation is the natural framework for detecting periodic solutions.

\section{Twisted degree}
Note that $S^1$ is the natural symmetry group of non-reversible periodic functions lying in $C_{2\pi}(\mathbb R;\mathbb R)$. However, the group $S^1$ has three intractable issues relative to the $G$-equivariant degree defined in the previous section:
\begin{enumerate}
    \item $S^1$-orbits of solutions are homeomorphic to circles, and therefore do not contain any isolated points.
    \item All nontrivial irreducible $S^1$-representations are of complex type.
    \item All subgroups $\mathbb Z_k < S^1$ have infinite Weyl group. In particular $W(\mathbb Z_k) = S^1$ and so $\dim W(\mathbb Z_k) = 1$.
\end{enumerate}
This section will give an overview of the degree tools needed to resolve these issues. The first is the local $S^1$-degree, which allows us to take the degree of $S^1$-orbits. The second is the twisted degree, which is to the local $S^1$-degree as the $G$-equivariant Leray-Schauder degree is to the local Brouwer degree. The construction of both of these degrees, and especially the twisted degree, is quite subtle. There are many delicate issues, such as the orientation of $S^1$-orbits relative to their ambient spaces, handling the complex structure of $S^1$-orbits, defining analogous algebraic structures counterpart to the Burnside ring $A(G)$ and defining a product in that new structure which counts orbits in the correct way, etc. Due to its central importance in many of our later applications, we will endeavor here to define this degree with an appropriate level of detail, addressing these concerns as much as possible and striving to give an intuitive picture of how the twisted degree is defined and computed. However, many of the more technical proofs, as before, will be omitted. For a modern reference on the construction of the twisted degree, see \cite{BalanovEtAl2025}, where all the detailed proofs abridged from this section can be found.
\subsection{Local \texorpdfstring{$S^1$-degree}{S¹-degree}}
As with the local Brouwer degree, let $V:=\mathbb R^n$ be a finite-dimensional real vector space, let $G:=S^1$, and let $f:\mathbb R\times V\to V$ be a $G$-equivariant continuous map parametrized by $\alpha\in \mathbb R$. Let $\Omega \subset \mathbb R \times V$ be an open bounded $S^1$-invariant set such that for all $(\alpha,x)\in \partial \Omega$, $f(\alpha,x)\neq 0$. Then $f$ is called $\Omega$-admissible and $(f,\Omega)$ an \emph{admissible $S^1$-pair}. We put $\mathscr M_1^{S^1}(\mathbb R \times V)$ the set of all admissible $S^1$-pairs and 
\[
\mathscr M_1^{S^1} := \bigcup_{V}\mathscr M_1^{S^1}(\mathbb R \times V),
\]
where the union is taken over all orthogonal $S^1$-representations $V$. We define admissible $S^1$-homotopies exactly as in previous sections, and call the function $\deg_{S^1}:\mathscr M_1^{S^1}\to A_1(S^1)$ the \emph{$S^1$-degree}, where $A_1(S^1)$ is the free $\mathbb Z$-module generated by the conjugacy classes $(\mathbb Z_k),\;k=1,2,\dots$. We also put $\Phi_1(S^1):=\{(H)\in \Phi(S^1):\dim(W(H)) = 1\}$. Then $A_1(S^1) = \mathbb Z[\Phi_1(S^1)]$. Then an element $\mathbf{a} \in A_1(S^1)$ can be written
\[
\mathbf{a} = \sum_{k=1}^\infty n_k(\mathbb Z_k),
\]
where the coefficient $n_k$ is nonzero for only finitely many $k \in \mathbb N$, and we write $\text{coeff}^k(\mathbf a) = n_k$. Note that as with the equivariant Brouwer degree, the assignment of such an element to an admissible pair $(f,\Omega)$ depends upon the existence of a regular normal approximation admissibly $S^1$-homotopic to $f$. 
\vs
$S^1$ has a one-dimensional real trivial irreducible representation, and infinitely many two-dimensional real irreducible representations, $\mathcal W_k$ for $k \in \mathbb N$, all of complex type, where $e^{i\theta}$ acts on each $z \in \mathcal W_k$ by $e^{i\theta}z = e^{ik\theta}\cdot z$, where here `$\cdot$' denotes complex multiplication.
\vs
Given an orthogonal $S^1$-representation $V$ and a point $v:=(\alpha,x)\in \mathbb R\times V$ having isotropy ${S^1}_v = \mathbb Z_k$ with $k\in \mathbb N$, consider the $S^1$-orbit $S^1(v)$. To define a degree, we need a concept of a positive orientation on this manifold. Notice that the $S^1$-action induces a natural orientation. Indeed, consider the $S^1$-homomorphism (which is also an $S^1$-equivariant diffeomorphism) $\rho:\mathbb R \to {S^1}(v)$ given by $\rho(\theta) = e^{i\theta}v$, and use the velocity field defined by $\dot \rho$ to define the orientation.
\vs
We define the \emph{slice} $S_v$ to $S^1(v)$ at $v \in S^1(v)$ as the orthogonal complement to the vector $\dot\rho(0)$,
\[
S_v := \{w \in \mathbb R \times V: w \perp \dot \rho(0)\}.
\]
A \emph{tube} around $S^1(v)$ is the set
\[
U:= S^1(e^{i\theta}(v+B_{\eta}(S_v)),\quad\text{ where } B_\eta(S_v):=\{u\in S_v:\vert u \vert < \eta \},
\]
for some $\eta >0$. Since $U$ is constructed as an $S^1$ orbit of a bounded open set, $U$ is also a bounded open $S^1$-invariant set. We denote the natural orientation of $S^1(v)$ by $\mathfrak o_v$, the standard orientation on $V$ by $\mathfrak o_V$, and the standard orientation on $\mathbb R \times V$ by $\mathfrak o_{\mathbb R \times V} := \mathfrak o_{\mathbb R} \oplus \mathfrak o_V$, where $\mathfrak o_{\mathbb R}$ is given by the vector $e_1 := (1,0,\dots,0) \in \mathbb R \times V$. Then the \emph{positive orientation} on the slice $S_v$ is the orientation $\mathfrak o_{S_v}$ such that $\mathfrak o_{S_v} \oplus \mathfrak o_v = \mathfrak o_{\mathbb R \oplus V}$. 
\vs
Note that for any $(f,\Omega) \in \mathscr M^{S^1}_1$, $f^{-1}(0) \cap \Omega$ can be written as a finite union of closed curves $C_1,\dots,C_m$, and each $C_j$ is either an $S^1$-orbit or $C_j \subset \Omega^{S^1}:=\Omega \cap (\mathbb R \times V)^{S^1}$. 
\vs
For any $k\in \mathbb N$, we define the \emph{$k$-folding homomorphism} $\Psi_k:S^1 \to S^1$ as $\Psi_k(e^{i\theta}) = {(e^{i\theta}})^k = e^{ik\theta}$, where $e^{i\theta} \in S^1$. Take any orthogonal $S^1$-representation $V$ with an $S^1$-action $\psi:S^1\times V \to V$. Then $\Psi_k$ induces the action $\psi_k:S^1\times V\to V$ given by $\psi_k(e^{i\theta},v) = \Psi_k(e^{ik\theta})v$. We will use the symbol ${}^kV$ for this representation where $S^1$ acts on $V$ through $\psi_k$, and we will call ${}^k V$ the \emph{$k$-folding of $V$}. Notice that if a map $f$ is $S^1$-equivariant under any $S^1$-action, then it is also $S^1$-equivariant under any $k$-folding of that action. Moreover, if $\mathcal V_m$ is the irreducible $S^1$-representation where $e^{i\theta} \in S^1$ acts on $z \in \mathcal V_m$ as $e^{i\theta}z = e^{im\theta}\cdot z$, then ${}^k\mathcal V_m$ acts as $e^{i\theta}z = e^{ikm\theta}\cdot z$, and so ${}^k\mathcal V_m = \mathcal V_{km}$. There is also a corresponding $\mathbb Z$-module homomorphism $\Psi_k:A_1(S^1)\to A_1(S^1)$ which acts on generators as $\Psi_k(\mathbb Z_m) = \mathbb (Z_{km})$ for all $m\in \mathbb N$.
\vs
Then the $S^1$-degree is the unique function $S^1\text{-deg}:\mathscr M_1^{S^1}\to A_1(S^1)$ satisfying the following properties:
\begin{enumerate}[label=\textbf{(S.\arabic*)}]

\item \textbf{(Additivity)} \label{prelim:s1deg:p1} Let $\Omega_{1}$ and $\Omega_{2}$
be two disjoint open $S^1$-invariant subsets of $\Omega$ such that
$f^{-1}(0)\cap\Omega\subset\Omega_{1}\cup\Omega_{2}$. Then
\begin{align*}
\sdeg(f,\Omega)=\sdeg(f,\Omega_{1})+\sdeg
(f,\Omega_{2}).
\end{align*}

\item \textbf{(Homotopy)} If $h:[0,1]\times V\to V$ is an
$\Omega$-admissible $S^1$-homotopy, then
\begin{align*}
\sdeg(h_{t},\Omega)=\mathrm{constant}.
\end{align*}

\item \textbf{(Normalization)} \label{prelim:s1deg:p3} For $(f,\Omega)\in \mathscr M_1^{S^1}$, let $f$ be a regular normal map in $\Omega$, $f^{-1}(0)\cap \Omega$ is connected, and $v_0 \in f^{-1}(0)\cap \Omega$. Then
\begin{align*}
\sdeg(f,\Omega)=\begin{cases}
        \rho_0(\mathbb Z_k) &\text{ if }S^1_{v_0} = \mathbb Z_k \text{ for some }k\in \mathbb N\\
        0 &\text{ if }S^1_{v_0} = S^1
    \end{cases}
\end{align*}
where $\rho_0 := \sign \det (Df(v_0)\vert_{S_{v_0}})$, and $S_{v_0}$ is the positively oriented slice to the orbit $S^1(v_0)$ at $v_0$.

\item \textbf{(Existence)} If for some $k\in \mathbb N$, $\text{coeff}^k(\sdeg(f,\Omega)) \neq 0$, then there
exists $v:=(\lambda,x)\in\Omega$ such that $S^1_{v}\geq\mathbb Z_k$.

\item \textbf{(Excision)}
Suppose that $(f,\Omega)\in \mathscr M_1^{S^1}$ and $\Omega_0$ is an open $S^1$-invariant subset of $\Omega$ such that $f^{-1}(0) \cap \Omega \subset \Omega_0$. Then $\sdeg(f,\Omega) = \sdeg(f,\Omega_0)$.

\item \textbf{(Product)} \label{prelim:s1deg:p6} Let $W$ be a finite-dimensional vector space with the trivial $S^1$-action, $U\subset W$ an open bounded neighborhood of zero and $(g,U)\in\mathscr M(W)$ an admissible pair in the sense of the Brouwer degree. Then
\[
\sdeg(g \times f,U\times \Omega)=\deg(g,U)\cdot \sdeg(f,\Omega)
\]
where $\deg(g,U)$ is the Brouwer degree.

\item \textbf{(Folding)} \label{prelim:s1deg:p7} If $k \in \mathbb N$ and $(f,\Omega)\in \mathscr M_1^{S^1}$, then
\[
\sdeg(f,{}^k\Omega) = \Psi_k(\sdeg(f,\Omega))
\]

\end{enumerate}
% \item \textbf{(Multiplicativity)}\label{prelim:s1deg:p6} For any $(f_{1},\Omega
% _{1}),(f_{2},\Omega_{2})\in\mathscr{M}_1 ^{S^1}$,
% \begin{align*}
% \sdeg(f_{1}\times f_{2},\Omega_{1}\times\Omega_{2})=
% \sdeg(f_{1},\Omega_{1})\cdot \sdeg(f_{2},\Omega_{2}),
% \end{align*}
% where the product is taken in $A_1(S^1)$.

\subsubsection{Cumulative $S^1$-degree}
Finally, we will define the \emph{cumulative $S^1$-degree}, $\deg_{S^1}:\mathscr M_1^{S^1}\to \mathbb Z$ as follows:
\[
\deg_{S^1}(f,\Omega):= \sum_{(\mathbb Z_k)\in \Phi_1(S^1;\Omega)}\text{coeff}^k(\sdeg(f,\Omega)).
\]
This provides a total algebraic count of nontrivial $S^1$-orbits in $f^{-1}(0)\cap \Omega$. One can immediately see that $\deg_{S^1}$ satisfies properties \ref{prelim:s1deg:p1}---\ref{prelim:s1deg:p6} if $\sdeg$ is replaced with $\deg_{S^1}$, $A_1(S^1)$ is replaced with $\mathbb Z$, and the normalization property \ref{prelim:s1deg:p3} is replaced with the following:
\begin{enumerate}[label=\textbf{(\~S.\arabic*)}, start=3]
    \item \textbf{(Normalization)} For $(f,\Omega)\in \mathscr M_1^{S^1}$, let $f$ be a regular normal map in $\Omega$, $f^{-1}(0)\cap \Omega$ is connected, and $v_0 \in f^{-1}(0)\cap \Omega$. Then
    \[
    \deg_{S^1}(f,\Omega) = \sign \det (Df(v_0)\vert_{S_{v_0}}),
    \]
    wehere $S_{v_0}$ is the positively oriented slice to the orbit $S^1(v_0)$ at $v_0$.
\end{enumerate}
\begin{remark}\rm
    Readers familiar with the concepts of Poincar\'e sections and Poincar\'e  map might see an intuitive connection between the Poincar\'e section through a point $v_0 \in \mathbb R \times V$ lying on a periodic orbit, and the normal slice $S_{v_0}$ to the orbit $S^1(v_0)$ at $v_0$. This connection can be made quite precise. Given a Poincar\'e section $\Sigma$ through $v_0$ and the Poincar\'e map $P:\Omega\subset\Sigma\to\Sigma$, the $S^1$-degree is exactly the Brouwer degree of the Poincare map taken on an admissible neighborhood on the slice, $\deg(P,\Omega)$, also called the Poincar\'e index. For more details on this connection, see \cite{Ize1992}.
\end{remark}

The twisted degree can be viewed as a counterpart to the equivariant Brouwer degree which uses the local $S^1$-degree as its basis, as opposed to the local Brouwer degree. Analogously to the development of the equivariant Leray-Schauder and Nussbaum-Sadovskii degrees, once this is defined for maps on finite-dimensional spaces, it can be extended to compact and condensing perturbations of identity on Banach spaces in precisely the same way. On an intuitive level, the twisted degree can be thought of as taking the \emph{equivariant} Brouwer degree on the normal slice $S_{v_0}$ to an $S^1$-orbit, and thereby obtaining orbit types corresponding to spatio-temporal symmetries, where the spatial symmetry group $\Gamma$ interacts with the temporal symmetry group $S^1$.
\subsection{Twisted subgroups}\label{sec:prelim:twisted-subgroups}
Let $\Gamma$ be a compact Lie group and put $G:= \Gamma \times S^1$. Subgroups of $G$ split into two classes: \emph{product subgroups} $K\times S^1$ where $K \leq \Gamma$, and \emph{twisted subgroups} of the form
\[
K^{\theta,l}:=\{(\gamma,z)\in K\times S^1: \theta(\gamma)=z^l\},
\]
where $\theta:K\to S^1$ is a group homomorphism and $l=0,1,2,\dots$, and we call $K^{\theta,l}$ a \emph{$\theta$-twisted $l$-folded subgroup}. Notice that product subgroups are the isotropy groups of constant solutions, and twisted subgroups are the isotropy groups of non-constant periodic solutions. The twisted degree can detect solutions whose isotropy groups $\mathscr H\leq G$ satisfy $\dim(W(\mathscr H)) = 1$. This immediately implies that $(K^{\theta,l})\in \Phi_1(G)$ if and only if $(K)\in \Phi_0(\Gamma)$, i.e. $W_\Gamma(K)$ is finite. 
\vs
Denote by $\Phi_1^t(G)$ the set of all conjugacy classes of $\theta$-twisted $l$-folded subgroups, and put $A_1^t(G) := \mathbb Z[\Phi_1^t(G)]$. Then elements $a \in A_1^t(G)$ are written
\[
a:= \sum_{(\mathscr H)\in \Phi_1^t(G)}n_{\mathscr H}(\mathscr H),
\]
where $n_{\mathscr H} = 0$ for all but finitely many $(\mathscr H)\in \Phi_1^t(G)$. The twisted degree takes values in $A_1^t(G)$. We also define $\text{coeff}^{\mathscr H}(a) := n_{\mathscr H}$ as before. We also carry over the definitions of the $k$-folding homomorphism, defining $\psi_k:G\to G$ as $\psi_k(\gamma,z) = (\gamma,z^k)$ and the induced $k$-folding homomorphism of $\mathbb Z$-modules $\Psi_k:A_1^t(G)\to A_1^t(G)$, which acts on generators $(\mathscr H)\in \Phi_1^t(G)$ as $\Psi_k(\mathscr H) = \psi_k^{-1}(\mathscr H)$.
\subsection{Multiplication in \texorpdfstring{$A_1^t(G)$}{A₁ᵗ(G)}}
In order to define the twisted degree and endow it with the multiplicativity property which is necessary for all practical computations, we need to define the multiplicative structure of $A_1^t(G)$. First we define $\cdot :A(\Gamma)\times A_1^t(G) \to A_1^t(G)$. For $(K)\in A(\Gamma)$ and $(H^{\theta,l})\in A_1^t(G)$, the coefficient $n_L$ of $(L^{\theta,l}) \in A_1^t(G)$ is given by
\[
n_L = \frac{n(L,K)\vert W(K)\vert n(L^{\theta,l},H^{\theta,l})\vert W(H^{\theta,l})/S^1\vert-\sum_{(\tilde L)>(L)}n_{\tilde L}n(L^{\theta,l},\tilde L^{\theta,l})\vert W(\tilde L^{\theta,l})/S^1\vert}{\vert W(L^{\theta,l})\vert}.
\]
Given $(K)\in A(\Gamma)$ and $(H^{\theta,l})\in A_1^t(G)$, the action of $\Psi_k$ is given by
\[
\Psi_k((K)\cdot (H^{\theta,l})) = (K)\cdot \Psi_k(H^{\theta,l}) = (K)\cdot (H^{\theta,kl}).
\]
\subsection{Definition of the twisted \texorpdfstring{$G$-equivariant}{G-equivariant} degree}
Let $V$ be an orthogonal $G$-representation, $\Omega \subset \mathbb R \times V$ an open bounded $G$-invariant set (where $G$ acts trivially on $\mathbb R$), and $f:\mathbb R \times V \to V$ a continuous $G$-equivariant $\Omega$-admissible map. Then $(f,\Omega)$ is called an \emph{admissble $G$-pair} in $\mathbb R \times V$, we denote by $\mathscr M_1^G(\mathbb R \times V)$ the set of all admissible pairs in $V$, and by $\mathscr M_1^G:=\bigcup_V \mathscr M^G(\mathbb R \times V)$ the set of all admissible $G$-pairs over all orthogonal $G$-representations. Take $(f,\Omega)\in \mathscr M_1^G(\mathbb R \times V)$ and $(\mathscr H)\in \Phi_1^t(G;\Omega)$. Then since $\mathbb R \times V^{\mathscr H}$ is a $W(\mathscr H)$-space and $S^1\leq W(\mathscr H)$, it is also an $S^1$-space, and so $(f^{\mathscr H},\Omega^{\mathscr H})$ is an admissible $S^1$-pair and its cumulative $S^1$-degree $\deg_{S^1}(f^{\mathscr H},\Omega^{\mathscr H})$ is well-defined. We equip the slice $S_v$ to the orbit $W(\mathscr H)(v)$ with the positive orientation induced by the $S^1$-action in the same way as described for the local $S^1$-degree. Then for $(\mathscr H)\in \Phi_1^t(G;\Omega)$, we define the numbers $d_{\mathscr H}$ as follows:
\[
d_{\mathscr H}(f,\Omega) := \frac{\deg_{S^1}(f^\mathscr H,\Omega^\mathscr H) - \sum_{(\mathscr K)>(\mathscr H)}d_{\mathscr K}(f,\Omega)n(\mathscr H,\mathscr K)\vert W(\mathscr K)/S^1\vert}{\vert W(\mathscr H)/S^1\vert},
\]
and for any $(\mathscr H)$ which is not an orbit type in $\Omega$, we put $d_{\mathscr H}(f,\Omega):= 0$. Then the numbers $d_\mathscr H(f,\Omega)$ are called the \emph{$\mathscr H$-twisted degrees} of $f$ in $\Omega$. We define the twisted $G$-equivariant degree $\gdeg:\mathscr M_1^G\to A_1^t(G)$ by
\[
\gdeg(f,\Omega)= \sum_{(\mathscr H)\in \Phi_1^t(G)} d_{\mathscr H}(f,\Omega)(\mathscr H),
\]
and it satisfies the following properties:
\begin{enumerate}[label=\textbf{(S.\arabic*)}]

\item \textbf{(Additivity)} \label{prelim:twistdeg:p1} Let $\Omega_{1}$ and $\Omega_{2}$
be two disjoint open $G$-invariant subsets of $\Omega$ such that
$f^{-1}(0)\cap\Omega\subset\Omega_{1}\cup\Omega_{2}$. Then
\begin{align*}
\gdeg(f,\Omega)=\gdeg(f,\Omega_{1})+\gdeg(f,\Omega_{2}).
\end{align*}

\item \textbf{(Homotopy)} If $h:[0,1]\times V\to V$ is an
$\Omega$-admissible $G$-homotopy, then
\begin{align*}
\gdeg(h_{t},\Omega)=\mathrm{constant}.
\end{align*}

\item \textbf{(Normalization)} \label{prelim:twistdeg:p3} For $(f,\Omega)\in \mathscr M_1^{G}$, let $f$ be a regular normal map in $\Omega$. If $f^{-1}(0)\cap \Omega = G(v_0)$ for some $v_0 \in\Omega$ with $G_{v_0} =\mathscr H$ and $(\mathscr H)\in \Phi_1^t(G;\Omega)$, then
\begin{align*}
    \gdeg(f,\Omega) = \sign \det (Df(w_0)\vert_{S_{v_0}}),
\end{align*}
where $S_{v_0}$ is the positively oriented slice to the orbit $W(\mathscr H)(v_0)$ at $v_0$ in $\Omega_{\mathscr H}$, and if $\Phi_1^t(G;f^{-1}(0)\cap \Omega) = \emptyset$, then $\gdeg(f,\Omega)=0$.

\item \textbf{(Existence)} If $(f,\Omega)\in \mathscr M_1^G$ and $\text{coeff}^{\mathscr H}(\gdeg(f,\Omega)) \neq 0$ for some $(\mathscr H)\in \Phi_1^t(G)$, then there
exists $v:=(\lambda,x)\in\Omega$ such that $G_{v}\geq\mathscr H$ and $f(t,x)=0$.

\item \textbf{(Excision)}
Suppose that $(f,\Omega)\in \mathscr M_1^{S^1}$ and $\Omega_0$ is an open $G$-invariant subset of $\Omega$ such that $f^{-1}(0) \cap \Omega \subset \Omega_0$. Then $\gdeg(f,\Omega) = \gdeg(f,\Omega_0)$.

\item \textbf{(Product)}\label{prelim:twistdeg:p6} Let $(f,\Omega)\in \mathscr M_1^G$, $W$ be an orthogonal $\Gamma$-representation, $U\subset W$ an open bounded neighborhood of zero and $(g,U)\in\mathscr M^\Gamma(W)$ an admissible $\Gamma$-pair. Then
\[
\gdeg(g \times f,U\times \Omega)=\operatorname{\text{$\Gamma$-deg}}(g,U)\cdot \gdeg(f,\Omega).
\]
where multiplication is taken in the $A(\Gamma)$-module $A_1^t(G)$. 

\item \textbf{(Folding)} \label{prelim:twistdeg:p7} If $k \in \mathbb N$ and $(f,\Omega)\in \mathscr M_1^{G}$, then
\[
\gdeg(f,{}^k\Omega) = \Psi_k(\sdeg(f,\Omega)).
\]
\end{enumerate}

Note that since this definition depends on the $G$-equivariant degree, it can be computed in terms of the Brouwer, Leray-Schauder, or Nussbaum-Sadovkii degrees, depending on the context. 
\subsection{Computing the twisted \texorpdfstring{$G$-equivariant}{G-equivariant} degree}
There is an analogous notion of basic maps and basic twisted degrees which allow the effective computation of the twisted $G$-equivariant degree by looking at irreducible $G$-representations and leveraging the product property. The definition and structure of these basic maps hinges on the fact that all nontrivial irreducible $S^1$-representations have complex structure, and this is therefore inherited by all irreducible $G$-representations whose expression as a tensor product of an irreducible $S^1$-representation and an irreducible $\Gamma$-represntation contains a nontrivial $S^1$ component (i.e. all the irreducible $G$-representations which admit twisted orbit types). This also means that all such irreducible $G$-representations will be even-dimensional. The exact construction of the corresponding basic maps and basic twisted degrees is rather technical, and we refer the interested reader to \cite{BalanovEtAl2025}, Proposition 9.5. 
\vs
We will simply provide the computational formula here: For $k>0$, let $\mathcal V_{k,j}$ be an irreducible $G$-representation. We call $\deg_{\mathcal V_{k,j}}$ the \emph{twisted $\mathcal V_{k,j}$-basic degree}, and $n_{\mathscr H}:=\text{coeff}^{\mathscr H}(\deg_{\mathcal V_{k,j}})$ for $(\mathscr H)\in \Phi_1^t(G)$ is given by
\[
n_\mathscr H = \frac{\frac{1}{2}\dim {\mathcal V_{k,j}}^{\mathscr H}-\sum_{(\mathscr K)>(\mathscr H)} n_{\mathscr K}n(\mathscr H,\mathscr K)\vert W(\mathscr K)/S^1\vert}{\vert W(\mathscr H) / S^1 \vert}.
\]
\vs
We will now show another useful formula for computing the degree, first generally, and then in a practical setting tailored to our later applications. Take an orthogonal $G$-representation $W$ with $G$-isotypic decomposition given by
\[
W = \bigoplus_{k,j} W_{k,j},
\]
where $W_{k,j}$ is the $G$-isotypic component of $W$ modeled on the irreducible $G$-representation $\mathcal W_{k,j} = \mathcal U_k \otimes \mathcal V_j$ with isotypic multiplicity $m_{k,j}$, where $\mathcal U_k$ is an irreducible $S^1$-representation and $\mathcal V_j$ an irreducible $\Gamma$-representation. Note that $\mathcal W_{k,j}$ is an absolutely irreducible complex $\Gamma$-representation if $k>0$. A linear $G$-equivariant operator $T:W_{k,k}\to W_{k,j}$ can thus be identified with a $m_{k,j}\times m_{k,j}$ complex matrix in $\operatorname{GL}(m_{k,j},\mathbb C)\cong \operatorname{GL}^G(W_{k,j})$. Given a map $a:S^1 \to \operatorname{GL}^G(W)$, denote by $a_{k,j}$ its restriction to the isotypic component $W_{k,j}$ and take the complex $m_{k,j}\times m_{k,j}$ matrix $\widetilde a_{k,j}(\gamma):\mathbb C^{m_{k,j}}\to \mathbb C^{m_{k,j}}$ with $\gamma \in S^1$ representing $a_{k,j}(\gamma)$. Then we put $d_{k,j} := \deg \det_{\mathbb C} \widetilde a_{k,j}$ for $j=0,1,\dots$ and $k>0$, where $\deg \det_{\mathbb C} \widetilde a_{k,j}$ denotes the algebraic count of zeroes of ${\det}_\mathbb C \widetilde a_{k,j}$\footnote{Formally, for $\sigma:S^n\to S^n$, $\deg(\sigma)$ is defined as the Brouwer degree $\deg(f_\sigma,B_1(0))$, where $f_\sigma:\mathbb R^{n+1}\to \mathbb R^{n+1}$ is a continuous extension given by the Tietze extension theorem. This can also be understood in the sense of the agrument principle from complex analysis. For more details, see \cite{BalanovEtAl2025}, Sec. 2.7}.
\vs
For $k=0$, $W_{0,j}$ is modeled on $\mathcal V_j$ with isotypic multiplicity $m_{0,j}$, and we define $\rho_j = \sign \det (a_{0,j}(\gamma))$ where $a_{0,j}(\gamma):=a(\gamma)\vert_{W_{0,j}}$ and $\gamma \in S^1$. Since $\mathcal W_{0,j}$ is the tensor product of the trivial $S^1$-representation with an irreducible $\Gamma$-representation $\mathcal V_j$, clearly $\rho_j$ is independent of the choice of $\gamma \in S^1$. 
\vs
Put $\mathfrak D:=\{(\lambda,v)\in \mathbb C\times W:\abs v < 2, \frac{1}{2}<\abs \lambda < 2$. Take a continuous map $a:S^1\to GL^G(W)$ and define $F_a:\mathfrak D\subset \mathbb C\times W\to W$ as $F_a(\lambda,v) = a \left(\frac{\lambda}{\abs \lambda}\right) v$. Let $\theta:\mathbb C\times W \to \mathbb R$ be a $G$-invariant continuous auxiliary function in $\mathfrak D$ for $F_a$ satisfying
\[
\begin{cases}
\theta(\lambda,0) > 0 &\text{ if }\abs \lambda = 2\\
\theta(\lambda,0)<0 &\text{ if } \abs \lambda = \frac{1}{2},
\end{cases} 
\]
and let $f_{\theta,a}:\mathbb C\times W \to \mathbb R \times W$ be given by
\[
f_{\theta,a}(\lambda,v) = (\theta(\lambda,v),F_a(\lambda,v)).
\]
Then
\[
\gdeg(f_{\theta,a},\mathfrak D) = \prod_j \delta_j \cdot \sum_{\substack{k,j\\k>0}}d_{k,j}\deg_{\mathcal V_{k,j}},
\]
where
\[
\delta_j := \begin{cases}
    (\deg_{\mathcal V_{j}})^{m_{k,j}} &\text{ if }\rho_j=-1,\\
    (\Gamma) &\text{ if } \rho_j = 1.
\end{cases}
\]
Since the $S^1$-isotypic components and associated eigenspaces of an $S^1$-equivariant linear operator correspond to the components of a Fourier decomposition, the index $k$ of irreducible $S^1$-representations is often idiomatically referred to as a ``Fourier mode'' or simply a ``mode''. Such a usage is intuitive and obvious when $S^1$ acts on the space of $2\pi$-periodic functions $C_{2\pi}(\mathbb R;\mathbb R^n)$. 
\vs
One obvious advantage of the twisted degree over the standard $G$-equivariant degree is that its computation mainly involves \emph{adding} basic degrees, with the exception of the multiplication of the $\Gamma$-degrees on the subspace of constant solutions. Since the orbit types which appear on the $\mathcal V_{0,j}$ mode are disjoint from those which appear on the higher modes (since $\mathcal V_{0,j}$ has the trivial $S^1$-action), this multiplication cannot destroy any orbit types on the higher modes. Therefore, it is only necessary to keep track of the indices $d_{k,j}$.
\vs
\section{Two parameter bifurcation}\label{sec:prelim:two-parameter-bifurcation}
It is most natural to use the twisted degree in the context of a system whose period is parametrized by an external parameter $\beta>0$. This can be accomplished by taking a system with solutions in any functional space, for which periodic solutions of any period are known or can be assumed to exist, performing the period-normalizing change of variables $u(t) := x(\beta t)$, and taking the functional setting of $C_{2\pi}(\mathbb R;\mathbb R^n)$ (or similar $2\pi$-periodic functional spaces, e.g. Sobolev spaces), where all solutions are $2\pi$-periodic and therefore are either constant, or are non-constant with period $\frac{2\pi}{k}$ for some $k>0$. This space is a natural $S^1$-representation, and this application is why the $S^1$-degree and twisted degree were formulated in terms of parametrized maps from, for example, $\mathbb R\times V$ to $V$.
\vs
However, although $\beta$ is a parameter, it is not exactly a natural choice for bifurcation in the one-parameter bifurcation sense described in Section \ref{sec:prelim:one-parameter}. This is because the relevant mode of bifurcation for the study of periodic solutions is Hopf bifurcation, and in this case, $\beta$ will be seen to correspond to the limit frequency of periodic bifurcating solutions. To study Hopf bifurcation, we need another parameter upon which eigenvalues of our linearized system continuously depend. We label this second parameter $\alpha$, and consider a system parametrized by $(\alpha,\beta)\in \mathbb R^2_+$, where $\mathbb R^2_+ := \mathbb R \times \mathbb R_+$ (since $\beta >0$). 
\vs
The idea is that for some isolated values of $\alpha$, the characteristic equation of a linear differential equation of interest will have purely imaginary roots, and the corresponding linear operator will fail to be an isomorphism. As long as these points $(\alpha,\beta)$ are isolated in a suitable sense, it is possible to define local bifurcation invariants analogous to those in Section \ref{sec:prelim:one-parameter}, and the Krasnosel'skii and Rabinowitz bifurcation theorems can be applied. This is the conceptual bridge between the one-parameter bifurcation regime described earlier and the framework used by the twisted degree to study Hopf bifurcation. We will consider a simple ODE system as a backdrop to formulate the required definitions.
\vs
Consider a first order system parametrized by $\alpha$ and given by
\[
\dot {\bm{x}}(t) = A(\alpha,\bm x) + F(\alpha,\bm x),\quad \bm x=(x_1(t),\dots,x_n(t))^T\in \mathbb R^n
\]
where $A_\alpha := A(\alpha,\cdot)$ is a linear operator, and $F_\alpha :=F(\alpha,\cdot)$ represents higher order terms. Then this system has a linearization at $\bm x\equiv 0$
\begin{equation}\label{eq:prelim:linear-hopf}
\dot{\bm x}(t) = A_\alpha \bm x,
\end{equation}
whose corresponding characteristic equation given by
\[
\triangle_{\alpha}(\lambda):=\lambda\id - A_\alpha.
\]
Let $\Gamma_0$ be the symmetry group of the state space and interactions, where $\Gamma_0\leq S_n$ acts by permuting components of $\bm x$. Put $\Gamma:= \Gamma_0 \times \mathbb Z_2$, where $\mathbb Z_2$ acts antipodally, i.e. as minus identity (this will be useful later for distinguishing nonconstant periodic solutions). Assume $A_\alpha$ and $F_\alpha$ are $\Gamma$-equivariant (note $\mathbb Z_2$-equivariance implies they are odd). Take $V:=\mathbb R^n$ as an isometric $\Gamma$-representation. Then $V$ admits a $\Gamma_0$-isotypic decomposition
\[
V = V_0\oplus V_1\oplus \dots \oplus V_r,
\]
where each $V_j$ is modeled on the irreducible $\Gamma$-representation $\mathcal V_j^-$ and has $\Gamma$-isotypic multiplicity $m_j$, where $\mathcal V_j^-$ is the irreducible $\Gamma_0$-representation $\mathcal V_j$ where $\mathbb Z_2$ acts antipodally. The characteristic equation then admits a decomposition by restriction to complexified irreducible $\Gamma$-representation, given by
\[
\triangle_{\alpha,j}(\lambda):=\triangle_\alpha\vert_{\mathcal V_j^{\mathbb C}}
\]
where $\mathcal V_j^{\mathbb C}:= \mathbb C \otimes_\mathbb R \mathcal V_j^-$. Then we have
\[
{\det}_{\mathbb C}(\triangle_{\alpha}(\lambda)) = \prod_{j=0}^r \big({\det}_{\mathbb C}\triangle_{\alpha}(\lambda)\big)^{m_j}.
\]
The $\bm x \equiv 0$ solution to \eqref{eq:prelim:linear-hopf} at $\alpha=\alpha_0$ is called a \emph{center} if there exists a corresponding $\beta_0>0$ (called the \textit{limit frequency}) such that $\det_{\mathbb C}\triangle_{\alpha_0}(i\beta_0) =0$. If, in addition, there exists $\varepsilon>0$ such that $0<|\alpha-\alpha_0|+|\beta-\beta_0|<\varepsilon$ implies
\[
{\det}_{\mathbb C}\triangle_{\alpha}(ik\beta) \not=0 
\]
for all $k\in\mathbb N$, then it is called an \emph{isolated center}.
\vs
Assuming that solutions of an arbitrary period $p>0$ exist, we can put $\beta = \tfrac{2\pi}{p}$ and perform the period normalizing transformation $\bm u(t):= \bm x(\beta t)$. We then obtain the period normalized system
\begin{equation}\label{eq:prelim:period-normalized-hopf}
\dot{\bm u}(t) := \frac{1}{\beta}F(\alpha,u(t)),
\end{equation}
which is still $\Gamma$-invariant. This can then be formulated as a compact perturbation of identity $\mathscr F:\mathbb R^2_+\times\mathscr E \to \mathscr E$, where $\mathscr E:= C_{2\pi}(\mathbb R;\mathbb R^n)$, such that $\mathscr F(\alpha,\beta,u)=0$ if $u$ if a $2\pi$-periodic solution to \eqref{eq:prelim:period-normalized-hopf}, and therefore a $\tfrac{2\pi}{\beta}$-periodic solution to our original system. Note that $\mathscr E$ is a natural $S^1$-representation, and the two-parameter equation $\mathscr F(\alpha,\beta,u)=0$ is clearly $S^1$-equivariant. We put $G:=\Gamma \times S^1$. Then $\mathscr E$ has a $G$-isotypic decomposition
\[
\mathscr E = \overline{\bigoplus_{k=0}^\infty \bigoplus_{j=0}^r \mathscr E_{k,j}}
\]
where each $G$-isotypic component $\mathscr E_{k,j}$ is modeled on the irreducible $G$-representation $\mathcal W_{k,j}:= \mathcal U_k \otimes \mathcal V_k^-$, where $\mathcal U_k$ is, for $k>0$, the two-dimensional irreducible $S^1$-representation where $S^1$ acts with the $k$-folded $S^1$-action, and for $k=0$ is the trivial representation. Note that each $\mathscr E_{k,j}$ is equivalent as a $G$-space to
\[
\mathscr E_{k,j} = \{a\cos(kt) + b\sin(kt):a,b\in V_j\}1.
\]
Then $\mathscr F$ has a linearization $\mathscr A(\alpha,\beta):= D_u\mathscr F(\alpha,\beta,0)$, and we can also restrict $\mathscr A(\alpha,\beta)$ to $G$-isotypic compoents, obtaining
\[
\mathscr A_{k,j}(\alpha,\beta) := \mathscr A(\alpha,\beta)\vert_{\mathscr E_{k,j}},
\]
and can be written in terms of the characteristic equation of the original system as
\begin{align*}
    \mathscr A_{k,j}(\alpha,\beta) &= \frac{1}{\beta(1+ik)}\triangle_{\alpha,j}(ik\beta) &\text{ if }k>0,\\
\mathscr A_{0,j}(\alpha,\beta) &= \frac{1}{\beta}A_\alpha &\text{ if }k=0.
\end{align*}
Note that purely imaginary roots to the characteristic equation $\triangle_\alpha(ik\beta)=0$ correspond to parameter values of $\alpha$ and $\beta$ where $\mathscr A_{k,j}(\alpha,\beta)$ fails to be an isomorphism. 
\subsection{Critical set}
This leads us to two definitions of the critical set. We define the \emph{critical set} $\Lambda := \{(\alpha,\beta,0)\in \mathbb R^2_+\times \mathscr E:\mathscr A_{k,j}(\alpha,\beta) \text{ is not an isomorphism for some $k>0$}\}$. The other definition consists of those values of $\alpha$ and $\beta$ which correspond to purely imaginary roots $\triangle_{\alpha}(i\beta) =0$ directly. This definition can be obtained before the functional space reformulation, and the two definitions are closely related. We often also refine this set slightly to exclude steady-state critical points where $\mathscr A_{0,j}$ is not an isomorphism, obtaining 
\[
\widetilde \Lambda := \{(\alpha,\beta,0)\in \mathbb R^2_+\times \mathscr E:\mathscr A_{k,j}(\alpha,\beta) \text{ is not an isomorphism for some $k>0$, and }\det\tfrac{1}{\beta}A_\alpha\neq 0 \}.
\]
The twisted equivariant degree can be used to show Hopf bifurcation in situations where the critical set $\Lambda$ consists of isolated centers, provided that crossing numbers can be computed.

\subsection{Transversality and crossing numbers}
\begin{remark}
Before defining the crossing numbers, we will first comment on the relationship between classical Hopf transversality conditions and those required by this method. Recall that for classical Hopf bifurcation, the transversality condition requires that, at a critical point $(\alpha_0,i\beta_0)$, the eigenvalues of the linear system cross the imaginary axis transversally, i.e. with nonzero velocity with respect to the parameter. In other words, $\tfrac{d}{d\alpha}\left[\re \lambda(\alpha_0)\right](i\beta)\neq 0$. The topological version of Hopf bifurcation used here does not actually have this explicit requirement. All that is necessary is that the net algebraic count of eigenvalues in a small neighborhood of $i\beta_0$ on the right-half plane changes as $\alpha$ sweeps through $\alpha_0$. This net change is captured by the \emph{crossing number}, which we will define shortly.
\vs
In principle, this is a more relaxed condition than the classical transversality condition, as it allows eigenvalues to cross the imaginary axis tangentially, as long as a net non-zero number of eigenvalues change the sign of their real part in any arbitrarily small $\alpha$-neighborhood of the critical point. Obviously, the case where such a valid crossing happens non-transversally is highly non-generic. In practice, it is actually often easier to compute normal Hopf bifurcation transversality conditions anyway, and use this stronger condition to obtain information on the signs of crossing numbers, which can be difficult to compute explicitly.
\end{remark}

\normalfont Take an isolated center $(\alpha_0,\beta_0,0) \in \tilde{\Lambda}$. Then for some sufficiently small $\varepsilon > 0$ and $\delta > 0$, there exists $U := (0,\varepsilon) \times (\beta_0 - \delta, \beta_0 + \delta) \subset \bbR^2_+$ such that for all $(p,\beta) \in U$, one has that
\[
{\det}_{\mathbb C} \triangle_\alpha(p + i\beta) \not= 0
\]
holds whenever $\alpha \in [\alpha_0-\delta,\alpha_0+\delta]$ and either
\begin{enumerate}
    \item $|\beta-\beta_0| = \delta$ and $p \in (0,\varepsilon)$
    \item $\beta \in (\beta_0 - \delta, \beta_0 + \delta)$ and $p = \varepsilon$
    \item $\beta \in (\beta_0 - \delta, \beta_0 + \delta)$, $p = 0$, and $\alpha \not= \alpha_0$
\end{enumerate}
\vs
In other words, at the critical point $(\alpha_0,\beta_0,0)$, roots $\lambda_{\alpha_0} \in {\mathbb C}_+$ of the characteristic equation can only enter or exit the region $U$ when $\lambda_{\alpha_0} = i\beta$ and when $\alpha = \alpha_0$. The crossing number is the net count of eigenvalues crossing the imaginary axis counted with their algebraic multiplicities. This value can be expressed using the Brouwer degree for values of $\alpha$ very near the critical point $\alpha_0$. Indeed, put $\alpha_\pm := \alpha_0 \pm \delta$ for $\delta>0$ sufficiently small. Then the crossing number is given by
\[
\mathfrak{t}(\alpha_0,\beta_0) := \mathfrak{t}_-(\alpha_0,\beta_0) - \mathfrak{t}_+(\alpha_0,\beta_0),
\]
where 
\[
\mathfrak{t}_\pm(\alpha_0,\beta_0) := \text{deg}(\triangle_{\alpha_\pm},U(\beta_0)) = \sign({\det}_{\mathbb C}\triangle_{\alpha_\pm}(i\beta_0)).
\]

We can further decompose these crossing numbers along the $\Gamma$-isotypic components of $\mathbb R^n$ to obtain the net count of eigenvalues crossing the imaginary axis on particular component, counted with their algebraic multiplicities. We denote the \emph{$\Gamma$-isotypic crossing numbers}
\[
\mathfrak{t}_{j}(\alpha_0,\beta_0) := \deg({\det}_{\mathbb C}\triangle_{\alpha_-,j}(i\beta_0)) - \deg({\det}_{\mathbb C}\triangle_{\alpha_+,j}(i\beta_0)).
\]

In order to obtain the crossing numbers corresponding to the period normalized problem, we consider the count with respect to the $G$-isotypic components, obtaining the \emph{$k$-resonant $\Gamma$-isotypic crossing numbers}
\[
\mathfrak{t}_{k,j}(\alpha_0,\beta_0) := \mathfrak{t}_j(\alpha_0,k \beta_0),\quad k \in {\mathbb N}.
\]

As mentioned before, we can also relate the signs of crossing numbers to the transversality condition as follows:
\begin{lemma}\label{lem:prelim:sign-crossing-num}
Let $(\alpha_0,\beta_0,0)\in \tilde\Lambda$ be an isolated center. If $\tfrac{\partial}{\partial \alpha}P_j(\alpha_0,i\beta_0) \neq 0$, then 
\[
\mathfrak t_{k,j}(\alpha_0,\beta_0) = -\sign\big(\tfrac{d}{d \alpha}\triangle_{\alpha,j}(i\beta_0)\big)m_j
\]
\end{lemma}
\begin{proof}
This follows directly from the smooth dependence of $\triangle_{\alpha,j}$ with respect to $\alpha$, and the definitions of the crossing numbers and $\Gamma$-isotypic multiplicity $m_j$.
\end{proof}

\subsection{Local bifurcation invariants}\label{sec:prelim:local-bif-twisted}
Now we need a way to extend the two fundamental topological bifurcation results of the Krasnosel'skii bifurcation theorem and the Rabinowitz alternative to this setting. It is interesting to note that, in this setting, Hopf bifurcation can be studied using the exact same theoretical framework as the one-parameter $G$-equivariant bifurcation given in Section \ref{sec:prelim:krasnoselskii}. As before, the method of auxiliary functions and complementing maps is used to allow us to take the degree of bifurcating branches by computing a difference of degrees at parameter values on either side of the bifurcation point. 
\vs
If $\widetilde \Lambda$ consists of isolated centers, then at each $(\alpha_0,\beta_0,0) \in \widetilde{\Lambda}$, for some sufficiently small $\varepsilon,\delta > 0$, one can find an \textit{isolated $G$-invariant neighborhood} $\Omega(\alpha_0,\beta_0,0)$ given by
\[
\Omega(\alpha_0,\beta_0,0) := \{(\alpha,\beta,u) \in \bbR^2_+ \times \mathscr E:\;\; |(\alpha,\beta)-(\alpha_0,\beta_0)| < \delta,\;\; ||u|| < \varepsilon\}
\]
such that
\[
\Omega(\alpha_0,\beta_0,0) \cap \overline{\mathscr S} = \emptyset.
\]
\vs
Now, following the approach of Ize et al. (cf. \cite{Ize1992}) and the equivariant generalizations of this approach (cf. \cite{krawcewiczwuxia1993}), by using an equivariant analog of the Tietze-Dugundji theorem, one can construct a continuous $G$-equivariant auxiliary function $\eta:\bbR^2_+ \times \mathscr E \to \bbR$ satisfying
\[
\begin{cases} \eta(\alpha,\beta,0)<0 &\text{ if } \; |(\alpha,\beta)-(\alpha_0,\beta_0)|  = \delta,\\
\eta(\alpha,\beta,v)>0 &\text{ if } \; |(\alpha,\beta)-(\alpha_0,\beta_0)|\le \delta \text{ and }\; \|v\|=\ve.\\
\end{cases}
\]
\vs
Now define $\mathscr F_\eta : \overline{\Omega(\alpha_0,\beta_0,0)} \to \mathscr \bbR \times \mathscr E$ by 
\[
\mathscr F_\eta(\alpha,\beta,u) = (\eta(\alpha,\beta,u),\mathscr F(\alpha,\beta,u)).
\]
Then $(\mathscr F_\eta,\Omega(\alpha_0,\beta_0,0))$ forms an admissible $G$-pair, and we can define the local bifurcation invariant $\omega_G(\alpha_0,\beta_0)\in A_1^t(G)$ by
\[
\omega_G(\alpha_0,\beta_0) := \gdeg(\mathscr F_\eta, \Omega(\alpha_0,\beta_0,0)),
\]
where $G$-deg is the twisted $G$-equivariant degree. It is clear from the homotopy invariance of the degree that $\omega_G(\alpha_0,\beta_0)$ is independent of the choice of $\varepsilon,\delta,$ and the auxiliary function $\eta$.
\vs
Recall that the computational formula for the twisted $G$-equivariant degree depended on the numbers $d_{k,j}:=\deg \det_{\mathbb C}(\widetilde a_{k,j})$, where $\deg$ is equivalent to the winding number. The relevant $S^1$-equivariant linear map in this case is 
$\widetilde \triangle_{\alpha,j}$, which is the $m_j\times m_j$ complex matrix equivalent to $\triangle_{\alpha,j}$ (in other words, a block matrix with $m_j$ diagonal blocks of $\triangle_{\alpha,j}$). Recall that we have
\[
    \mathscr A_{k,j}(\alpha,\beta) = \frac{1}{\beta(1+ik)}\triangle_{\alpha,j}(ik\beta),\quad \text{ if }k>0.
\]
Put $\varphi_{k,j}(\alpha,i\beta) := \det_\mathbb C( {\mathscr A}_{k,j}(\alpha,\beta))$. Since $\varphi_{k,j}(\alpha,i\beta)$ is analytic in $\beta$, by taking a Taylor expansion of $\varphi_{k,j}$ and applying the argument principle, $\deg(\varphi_{k,j})=\deg\det_{\mathbb C}({\mathscr A}_{k,j}(\alpha_0,\beta_0))$ is equivalent to the crossing number $\mathfrak t_{k,j}(\alpha_0,\beta_0)$. (For a full and detailed proof of this fact, see \cite{BalanovEtAl2025}, Chapter 9). Then the local bifurcation invariant is given by
\[
\omega_G(\alpha_0,\beta_0) =  \gdeg(\mathscr F_\eta, \Omega(\alpha_0,\beta_0,0)) = \gdeg(\mathscr F(\alpha_-,\beta_0,0),U) - \gdeg(\mathscr F(\alpha_+,\beta_0,0),U),
\]
where $\alpha_\pm$ are chosen sufficiently close to $\alpha_0$ and satisfy $\alpha_-<\alpha_0<\alpha_+$. Denote $d_{k,j}(\alpha,\beta) = \deg{\det}_{\mathbb C}\widetilde\triangle_{\alpha,j}(i\beta)$. Then we have
\begin{align*}
\omega_G(\alpha_0,\beta_0) &= \prod_{j=0}^r\Gamma\text{-deg}(\mathscr A_{0,j}({\alpha_-,\beta_0),B(V))}\cdot \sum_{k=1}^\infty\sum_{j=0}^r \mathfrak t^-_{k,j}(\alpha_0,\beta_0)\deg_{\mathcal W_{k,j}}\\
&- \prod_{j=0}^r\Gamma\text{-deg}(\mathscr A_{0,j}({\alpha_+,\beta_0),B(V))}\cdot \sum_{k=1}^\infty\sum_{j=0}^r
\mathfrak t^+_{k,j}(\alpha_0,\beta_0)\deg_{\mathcal W_{k,j}},
\end{align*}
where $\mathfrak t^\pm_{k,j}(\alpha_0,\beta_0) := \deg(\widetilde \triangle_{\alpha_\pm,j},U(k\beta_0)$, analogously to the definition of the non-isotypic crossing numbers. Since $\widetilde \Lambda$ excludes points of steady-state bifurcation which could cause a change in the $\Gamma$-degrees of the zero mode, this is equivalent to
\[
\omega_G(\alpha_0,\beta_0) = \prod_{j=0}^r\Gamma\text{-deg}(\mathscr A_{0,j}({\alpha_0,\beta_0),B(V))}\cdot\sum_{k=1}^\infty\sum_{j=0}^r \big(\mathfrak t^-_{k,j}(\alpha_0,\beta_0)-\mathfrak t^+_{k,j}(\alpha_0,\beta_0)\big)\deg_{\mathcal W_{k,j}}
\]
which, by applying the definition of the crossing numbers, yields the computational formula
\[
\omega_G(\alpha_0,\beta_0) = \prod_{j=0}^r\Gamma\text{-deg}(\mathscr A_{0,j}({\alpha_0,\beta_0),B(V))}\cdot\sum_{k=1}^\infty\sum_{j=0}^r \mathfrak t_{k,j}(\alpha_0,\beta_0)\deg_{\mathcal W_{k,j}}.
\]
Note that the above construction, like that of the twisted equivariant degree in general, does not depend on whether the Leray-Schauder or Nussbaum-Sadovskii degree is ultimately used, and will yield identical results for each degree provided the $G$-equivariant field $\mathscr F$ is respectively a compact or condensing perturbation of identity. 
\vs
We can now formulate the version of the Krasnosel'skii theorem adapted to the setting of two-parameter Hopf bifurcation with the twisted degree:
\begin{theorem}\label{thm:prelim:kransoselskii-hopf}
Put $G:=\Gamma \times S^1$, and let $F(\alpha,x):\mathbb R\times V \to V$ be a $\Gamma$-equivariant continuous map such that for all $\alpha\in \mathbb R$, $f(\alpha,0)=0$ and the derivative $A(\alpha):=D_xf(\alpha,0):V\to V$ exists, depends continuously on $\alpha$, and 
\[
\lim_{(\alpha',v)\to(0,0)} \frac{\vert f(\alpha',v) - \mathscr A(\alpha')v) \vert}{\vert v \vert}=0.
\]
Let $\mathscr E$ be an isometric Banach $G$-representation, and let $\mathscr F(\alpha,\beta,u):\mathbb R^2_+\times \mathscr E\to \mathscr E$ be the functional operator reformulation of the system $\dot x = F(\alpha,x)$ such that if $\mathscr F(\alpha,\beta,u)=0$ then $u(t) = x(\beta t)$ is a $\tfrac{2\pi}{\beta}$-periodic solution of $\dot x = F(\alpha,x)$. Assume the critical set $\Lambda$ for $A(\alpha)$ is isolated, and for some $(\alpha_0,\beta_0,0)\in \Lambda$ and some twisted orbit type $(H) \in \Phi_1^t(G;\mathscr E_{k,j})$ one has
\[
\text{coeff}^H(\omega_G(\alpha_0,\beta_0)) \neq 0.
\]
Then there exists a branch $\mathscr C$ of periodic solutions to $\dot x = F(\alpha,x)$ bifurcating from the trivial branch at $(\alpha_0,0)$ with symmetries at least $(H)$.
\end{theorem}
We will now make a few observations about this theorem and its implications. First of all, if the bifurcation occurs on a nonzero Fourier mode, then the bifurcating solutions near the bifurcation point will be $\tfrac{2\pi}{\beta}$-periodic. However, this does not let us say anything about the \emph{minimal} period of the bifurcating solutions. It is also possible that the frequencies of the periodic solutions can change due to nonlinear effects as they move away from the trivial branch. 
\vs
Moreover, given a $\theta$-twisted $k$-folded orbit type $(H^{\theta,k}) \in \Phi_1^t(\Gamma_0\times S^1)$, one can see that in fact $(H^{\theta,k})\leq (H\times S^1)$. This product orbit type $(H\times S^1)$ represents the same spatial symmetries in $\Gamma_0$, but corresponds to constant solutions. This is obviously problematic from the point of view of finding non-constant periodic solutions. This is why we put $\Gamma:=\Gamma_0 \times \mathbb Z_2$. The antipodal $\mathbb Z_2$ action allows us to find twisted symmetries which can only correspond to odd solutions, and thus there \emph{are} twisted orbit types $(H^{\theta,k})\in \Phi_1^t(G)$ such that $(H^{\theta,k})\not < (H\times S^1)$. We can filter our search for just these orbit types through reduction to an appropriate fixed point space. We will detail this process in Section \ref{sec:prelim:fixed-point-reduction}.
\vs
This approach has a number of advantages over classical Hopf bifurcation. First of all, the critical frequency $i\beta_0$ does not have to satisfy any non-resonance condition. Second of all, we do not require it to have simple multiplicity as a root of the characteristic equation. These two facts alone give this approach tremendous power in addressing Hopf bifurcation in highly symmetric contexts such as the multi-agent systems studied in this dissertation. It is quite common in highly symmetric systems for eigenvalues to be forced into high multiplicities by symmetries. This is a major degeneracy for classical Hopf bifurcation theorems which must be handled with extreme care, but poses no problem at all for the equivariant topological Hopf bifurcation approach described here.
\vs
\subsection{Global bifurcation}
The global bifurcation (or unbounded continuation) theorem of Rabinowitz applies here essentially unchanged. We will present it here, in the form of Theorem \ref{thm:prelim:rabinowitz2}, modified to the twisted degree two-parameter Hopf bifurcation context.
\begin{theorem}\label{thm:prelim:rabinowitz2-global}
    Let $\mathscr F:\mathbb R^2_+ \times \mathscr E \to \mathscr E$ be a compact (resp. condensing) perturbation of identity such that $\mathscr F(\alpha,\beta,0)=0$ for all $(\alpha,\beta) \in \mathbb R^2_+$, $\mathscr A(\alpha,\beta) :=D_u\mathscr F(\alpha,\beta,0)$ exists for all $(\alpha,\beta,0)\in \mathbb R^2_+\times \mathscr E$, and the critical set $\Lambda := \{(\alpha,\beta,0)\in \mathbb R^2_+ \times \mathscr E:\mathscr A \text{ is not an isomorphism}\}$ is discrete. Then if there is an orbit type $(H) \in \Phi_1^t(G;\mathscr E)$ such that
    \[
    \sum_{(\alpha,\beta,0)\in \Lambda} \text{coeff}^H(\omega_G(\alpha,\beta,0)) \neq 0,
    \]
    then there exists an unbounded branch of solutions $\mathscr C$ having symmetries at least $(H)$.
\end{theorem}
Suppose there is some $\mathfrak t_{k,j}(\alpha_0,\beta_0)\neq 0$ and therefore an associated orbit type $(H)\in \Phi_1^t(G;\mathcal W_{k,j})$ with $\text{coeff}^H(\deg_{\mathcal W_{k,j}}) \neq 0$. Then
\[
\sum_{(\alpha,\beta,0)\in \Lambda} \text{coeff}^H(\omega_G(\alpha,\beta,0)) \neq 0 \Leftrightarrow \sum_{\alpha,\beta,0)\in \Lambda} \mathfrak t_{k,j}(\alpha,\beta,0)\neq 0.
\]
Therefore, we mainly focus on showing $\sum_{\alpha,\beta,0)\in \Lambda} \mathfrak t_{k,j}(\alpha,\beta,0)\neq 0$. Since our critical sets are infinite and depend on an infinite number of discrete roots of a transcendental equation, the most practical way to achieve this is to give conditions that $\sign \mathfrak t_{k,j}(\alpha_{n,j},\beta_{n,j},0)$ is constant for fixed $j$ and all $n$ (or at least all $n$ sufficiently large). 

\subsection{Fixed point reduction}\label{sec:prelim:fixed-point-reduction}
The equivariant Krasnosel'skii theorem tells us about the (minimum) symmetries of a branch at bifurcation, and the Rabinowitz alternative can tell us about its global unbounded continuation, but knowing how the symmetries evolve as that branch continues globally is a more delicate question. It is always possible that the global branch will undergo secondary bifurcation, fold backwards, intersect with another branch non-locally, etc. 
\vs
One situation of particular interest arises in the search for non-constant periodic solutions, a problem which will be a major focus of later chapters. It is possible that a branch of non-constant periodic solutions with symmetries at least $(H)$ can bifurcate and then ``collapse'' non-locally to a branch of constant solutions. Although this branch may be unbounded and continue globally in the sense of the Rabinowitz alternative, it may only be periodic and non-constant for a short time before undergoing some secondary bifurcation, and neither the Krasnosel'skii nor the Rabinowitz theorems can tell us this.
\vs
The solution is a technique known as \emph{fixed point reduction}, where we restrict the problem to a specially chosen subspace which includes only solutions with symmetries of interest to us. By applying the Krasnosel'skii and Rabinowitz theorems to this restricted problem, we obtain guarantees that global branches must possess certain symmetries. Crucially, any solution to this restricted problem must also be a solution to the original problem. This also allows us to resolve resonances or other situations where multiplicities of eigenvalues at critical points are higher than our liking, by constructing appropriate fixed point subspaces where the eigenvalue crossings are $G$-isotypically simple, for example.
\vs
The most common application of fixed point reduction in subsequent chapters involves the group $G:=\Gamma_0\times \mathbb Z_2 \times S^1$, where $\Gamma_0$ is a compact Lie group, and is typically a finite group of spatial symmetries. We then take the cyclic subgroup $\bm K$ generated by the following element:
\[
\bm K := \langle(e_{\Gamma_0},-1,e^{i\frac{\pi}{\kappa}}) \rangle
\]
where $\kappa \in \mathbb N$. In many cases it suffices to take $\kappa =1$. The choice of $\kappa$ acts as a kind of filter to only take Fourier modes which are odd multiples of $\kappa$, which can be useful in eliminating degeneracies related to non-isotypically simple critical points, as discussed above. Suppose we have a compact perturbation of identity $\mathscr F(\alpha,\beta,u)$, which we consider as a two-parameter problem without loss of generality. Then we have the restricted problem $\mathscr F^{\bm K}:= \mathscr F\vert_{\mathscr E^{\bm K}}$, where $\mathscr E^{\bm K}$ can be understood by considering which $G$-isotypic components $\mathscr E_{k,j}$ of $\mathscr E$ will be fixed under the $\bm K$-action. Since $\bm K$ is a cyclic group, it is sufficient to consdier those which are fixed by the single generator $(e_{\Gamma_0},-1,e^{i\frac{\pi}{\kappa}})$. If we take
\[
\mathscr E_{k,j} := \{c\cos(kt)+d\sin(kt):c,d\in V_j\},
\]
where $(e_{\Gamma_0},-1,e^{i\frac{\pi}{\kappa}})$ acts with the $k$-folded $S^1$-action, then we have
\[
(e_{\Gamma_0},-1,e^{i\frac{\pi}{\kappa}})(c\cos(kt)+d\sin(kt)) = (-c\cos(kt+k\tfrac{\pi}{\kappa})-d\sin(kt+k\tfrac{\pi}{\kappa})).
\]
This is fixed if and only if $\tfrac{k}{\kappa} \in 2\mathbb N-1$. In this case, $-c\cos(kt+k\tfrac{\pi}{\kappa}) = -c\cos(kt + (2l-1)\pi)=c\cos(kt)$ and $-d\sin(kt+k\tfrac{\pi}{\kappa}) = -d\sin(kt + (2l-1)\pi)=d\sin(kt)$, for some $l \in \mathbb N$. Therefore, we have
\[
\mathscr E_{k,j}^{\bm K} = \begin{cases}
    \mathscr E_{k,j} &\text{ if }k=(2l-1)\kappa\text{ for some }l \in \mathbb N\\
    \{0\} &\text{ otherwise}.
\end{cases}
\]
and so we can write
\[
\mathscr E^{\bm K} = \overline{\bigoplus_{l=1}^\infty \bigoplus_{j=0}^r \mathscr E_{(2l-1)\kappa, j}}.
\]
Finally, notice that $\bm K$ is normal in $G$, and $G/\bm K \cong \Gamma_0 \times S^1$, and so $\mathscr E^{\bm K}$ is a natural $G/\bm K$-representation, and the above $G$-isotypic decomposition of $\mathscr E$ is isomorphic to an equivalent $G/\bm K$-isotypic decomposition of $\mathscr E^{\bm K}$. All of the above has two important consequences: first of all, solutions to $\mathscr F^{\bm K}(\alpha,\beta,u)=0$ must be solutions to $\mathscr F(\alpha,\beta,u)=0$, and such solutions must be odd (and therefore non-constant). By applying the Rabinowitz alternative to $\mathscr F^{\bm K}(\alpha,\beta,u)=0$, we obtain the global continuation of solutions in $\mathscr E^{\bm K}$, and therefore a guarantee that such solutions must remain periodic and non-constant across the entire global branch (although secondary bifurcation is still possible). To do this, we can refine the critical set $\Lambda$ to $\Lambda^{\bm K}$ which, as might be expected, is defined by considering only the critical points in $\Lambda$ which correspond to $G$-isotypic components $\mathscr E_{(2l-1)\kappa,j}$ for $l \in \mathbb N$, i.e.
\[
\Lambda^{\bm K} := \{(\alpha,\beta,0)\in \mathbb R^2_+\times \mathscr E:\mathscr A_{(2l-1)\kappa,j} \text{ is not an isomorphism for some }l\in \mathbb N\},
\]
where $\mathscr A_{(2l-1)\kappa,j}$ is the Fr\`echet derivative of $\mathscr F^{\bm K}$ at $u=0$, restricted to one of the $G/\bm K$-isotypic components. This leads to the following formula for the local bifurcation invariant taken at $(\alpha_0,\beta_0,0)\in \Lambda^{\bm K}$:
\[
\omega_G(\alpha_0,\beta_0) = \sum_{j=0}^r\sum_{l=1}^\infty \mathfrak t_{(2l-1)\kappa,j}(\alpha_0,\beta_0)\deg_{\mathcal W_{(2l-1)\kappa,j}},
\]
where $\mathcal W_{(2l-1)\kappa,j}$ is the irreeducible $G$ (or $G/\bm K$) representation on which the isotypic component $\mathscr E_{(2l-1)\kappa,j}$ is modeled. This formula is advantageous because it includes no multiplication at all. This means that the condition for global continuation of branches of non-constant periodic solutions is given by
\[
\sum_{(\alpha_0,\beta_0,0)\in \Lambda^{\bm K}}\sum_{j=0}^r\sum_{l=1}^\infty \mathfrak t_{(2l-1)\kappa,j}(\alpha_0,\beta_0)\deg_{\mathcal W_{(2l-1)\kappa,j}} \neq 0.
\]
This follows from simply applying the Rabinowitz alternative to $\mathscr F^{\bm K}(\alpha,\beta,u)=0$. In applications, most critical points are $G$-isotypically simple, and in this case, the local bifurcation invariant consists of a single $k$-resonant $\Gamma_0$-isotypic crossing number, and the above formula can be further simplified.

\section{Multi-agent systems}\label{sec:prelim:mas}
Much like control theory and other related subjects, the field of multi-agent systems has its own particular language and nomenclature for various concepts, to which direct parallels in dynamical systems and differential equations can often be drawn. Since the bulk of the existing applications of equivariant degree in the literature are to general families of differential equations or dynamical systems, we judge it worthwhile to adopt something of a dual perspective in the following chapters of this dissertation. 
\vs
On one hand, we are defining three novel classes of multi-agent sytems, and we are interested in using equivariant degree to explore their dynamics, the breakdown of consensus, the emergence of competition and mixed cooperative-competitive regimes, the formation of schismatic factions of subpopulations of agents, and periodic multiconsensus states.
\vs
On the other hand, setting aside this interpretation, we can frame many of the above questions in terms more familiar to other applications of equivariant degree, e.g. the search for nonstationary periodic solutions to symmetric systems of delay equations. We will adopt this dynamical systems terminology throughout most of the chapters, as we feel it makes it easier to evaluate both the novelty and validity of these results by way of easier comparison with the established literature and the cited references. However, when introducing these systems, motivating them, and interpreting our results, we will gravitate towards the multi-agent system description. 
\vs
This section will outline the most essential concepts and jargon from multi-agent systems which are relevant here, and provide the bridge between these concepts and other familiar notions in the study of dynamical systems. We will also show why multi-agent systems are a natural and promising arena for applying equivariant degree methods to obtain useful and novel results which cannot be obtained using other methods, opening new doors in both subjects.
 
\subsection{Agents}
An \emph{agent} is an autonomous entity capable of sensing its local environment, processing this information, and acting upon it in order to achieve individual or collective goals. A \emph{multi-agent system} (MAS) is a collection of such agents which interact with each other, either through communication, their shared environment, or both. This framework is a natural choice for modeling complex networks of intelligent decision-makers, such as humans, animals, or artificial intelligences. Other more abstract composite entities such as polities, nation-states, corporations, or political movements can also be modeled in this way. The study of MAS therefore sits at the intersection of control theory, dynamical systems, distributed computing, and artificial intelligence research.
\vs
An agent $i\in \{1,\dots,n\}$ possesses a \emph{state variable} $x_i(t) \in \mathbb R^d$ (we will consider $d=1$ for simplicity). The \emph{open loop dynamics} of an agent describe its intrinsic evolution without any interaction or internal decision-making. We call $u_i(t) \in \mathbb R$ a \emph{control input}, and consider it to be an arbitrary function of time $t$ which may be freely manipulated or designed, but which does not depend on the state $x_i(t)$. Then the open-loop dynamics are given by
\[
\dot x_i(t) = F(x_i,u_i)\quad \text{ or }\quad  x_i^{(k)} = F(x_i,\dot x_i,\ddot x_i,\dots,x_i^{(k-1)},u_i)
\]
Common open-loop dynamics are \emph{single-integrator dynamics}, given by
\[
\dot x_i = u_i,
\]
and \emph{double-integrator dynamics}, given by 
\[
\ddot x_i = u_i.
\]
If $u_i$ is defined as a function of the system state $x:=(x_1,\dots,x_n)^T\in \mathbb R^{n}$, it is called a \emph{protocol} (or equivalently \emph{control law} or \emph{feedback law}). This yields the \emph{closed-loop dynamics}. For single-integrator dynamics, for example, the closed-loop dynamics of a single agent take the form
\[
\dot x_i = u_i(x_1,\dots,x_n,t),
\]
where the protocol $u_i$ for each agent $i$ may only depend on certain ``nearby'' agents in some neighborhood $\mathcal N_i \subset \{1,\dots,n\}$.
\subsection{Interaction topology}
Let $\mathcal G=(\mathcal V,\mathcal E)$\footnote{We will consider the graph notation used here as confined to this section and used only for definitional purposes. Later chapters will describe the interaction topology implicitly and the letters chosen here (e.g. $\mathcal V$) will be used for other purposes.} be a directed or undirected graph with vertices $\mathcal V = \{1,\dots,n\}$ and edges 
\[
\begin{cases}
\mathcal E\subset \{\{v_1,v_2\}:v_1,v_2\in \mathcal V,\; v_1\neq v_2\} &\text{ if $\mathcal G$ is undirected}\\
\mathcal E\subset \{(v_1,v_2)\in \mathcal V\times \mathcal V:\; v_1\neq v_2\} &\text{ if $\mathcal G$ is directed}.
\end{cases}
\]
We also define the \emph{adjacency matrix} $A:= [a_{ij}]_{i,j=1}^n$ of non-negative \emph{weights}, where the weight $a_{i,j}>0$ if and only if $(j,i) \in \mathcal E$. Note that this convention assigns $a_{i,j}$ as the influence of agent $j$ on agent $i$, which will be notationally convenient later. We then call the pair $(\mathcal G,A)$ the \emph{interaction topology}, and we define the \emph{neighborhood} of agent $i$, denoted $\mathcal N_i$, as
\[
\mathcal N_i := \{j\in\mathcal V:a_{i,j}>0\}.
\]
Given an undirected graph, we will simply view it as its equivalent directed graph with symmetric weights, where $(i,j)\in \mathcal E$ iff $(j,i)\in \mathcal E$ and where $a_{ij}=a_{ji}$ for all $i,j,\in \mathcal V$. 
\vs
% We define the \emph{in-degree} of vertex $i$ as
% \[
% d_i^{\mathrm{in}}:= \sum_{j=1}^n a_{ij},
% \]
% the \emph{degree matrix} as the diagonal matrix $D:=\text{diag}(d_1^{\mathrm{in}},\dots,d_n^{\mathrm {in}})$, and the \emph{graph Laplacian} $L$ as $L:=D-A$. Then 
% \vs
The interaction topology can be state-dependent as well. We say that a MAS has \emph{state-dependent effective weights} if the protocol for agent $i$ can be written
\[
u_i = f(x_i,t) + h_i(x_1,\dots,x_n),
\]
where $f(x_i,t)$ is a self-regulation term, and $h_i$ is a nonlinear coupling function where the influence of other agents on agent $i$ is state-dependent, and where, given $h(x):=(h_1(x),\dots,h_n(x))^T \in \mathbb R^{n}$, its Jacobian $Dh(x_0)\in \mathbb R^{n\times n}$ is an adjacency matrix in the interaction topology, and represents the effective strength of the interaction at the state $x(t)\in \mathbb R^n$. 
\subsection{Consensus}\label{sec:prelim:consensus}
Following the standard definition of Olfati-Saber \& Murray \cite{saber2004consensus}, a point $x \in \mathbb R^n$ is called a \emph{consensus state} if $x_1=x_2=\dots=x_n$, i.e. all agents have an identical state. Note that this is a purely geometric subspace of the state space $\mathbb R^n$, and has no relationship whatsoever to the closed-loop dynamics. A consensus state which is also an equilibrium of the closed-loop dynamics is called a \emph{consensus equilibrium}. A system is said to \emph{achieve consensus} if, for all initual conditions in some specified set $\mathcal I \subseteq \mathbb R^n$,
\[
\lim_{t\to\infty} \norm{ x_i(t) - x_j(t)} = 0 \quad \text{ for all }i,j,\in\{1,\dots,n\}.
\]
If $\mathcal I = \mathbb R^n$, then the system is said to achieve \emph{global} consensus. If $\mathcal I$ is some neighborhood of the consensus manifold, or of a particular consensus equilibrium, then the system is said to achieve \emph{local} consensus. Much work in MAS literature is related to designing protocols which achieve local or global consensus. The classical global consensus protocol for $n$ single-integrator agents is given by
\[
u_i = \sum_{j\in \mathcal N_i}a_{ij}(x_j - x_i),
\]
which leads to the autonomous closed-loop dynamics
\[
\dot x_i = \sum_{j\in \mathcal N_i}a_{ij}(x_j-x_i),\quad i=1,\dots,n.
\]
It is easy to see that any point in the \emph{consensus manifold} $\mathcal C:=\{x\in \mathbb R^n:x_1 = x_2=\dots=x_n\}$ is a consensus equilibrium for this system, and so the single-integrator MAS with this protocol achieves global consensus. 
\vs
In the language of dynamical systems, whether a system achieves consensus is related to the existence of a consensus equilibrium having local or global Lyapunov stability. In this dissertation, we will mainly be concerned with the consensus equilibrium $x \equiv 0$, and conditions under which the system achieves local consensus around this point (i.e. under which $x\equiv 0$ is locally asymptotically stable).  
\vs
\subsection{Periodic consensus, multiconsensus, cooperation, and competition}
The concept of consensus has been generalized in several natural and important ways. A nonstationary $T$ periodic solution $x(t) = (x_1(t),\dots,x_n(t))$ of the closed loop dynamics is called a \emph{periodic consensus solution} if, for all $t\in \mathbb R$ and all $i,j\in \{1,\dots,n\}$, $x_i(t) = x_j(t)$. This corresponds to a $T$-periodic function on the consensus manifold $\mathcal C$. A system achieves periodic consensus if there exists a $T$-periodic consensus solution which is locally asymptotically stable. 
\vs
Most of the discussion of consensus thus far has implied cooperation and synchronization between agents. In order to study competitive behaviors, which will be described in more detail shortly, it is necessary to develop some broader ideas of consensus states. A \emph{multiconsensus} state involves splitting $\mathcal V$ into disjoint clusters of agents, $\mathcal V = C_1\sqcup C_2 \sqcup\dots \sqcup C_k$ such that within each cluster, agents achieve consensus, i.e.
\[
x_i = x_j\quad \text{for all }i,j\in C_p, p=1,\dots,k,
\]
where different clusters can attain different values. A system achieves multiconsensus if
\[
\lim_{t\to\infty} \norm{x_i(t) - x_j(t)} = 9 \quad \text{ for all }i,j\in C_p, \text{ for each }p=1,\dots,k.
\]
\vs
A \emph{bipartite} consensus state occurs when $\mathcal V$ can be partitioned into two populations $\mathcal V = C_1 \sqcup C_2$ such that for all $i,j\in C_1$, $x_i=x_j$, and for all $i,j\in C_2$, $x_i=x_j$, and for all $i\in C_1,j\in C_2$, $x_i=-x_j$. In other words, the two populations take equal and opposite values. A system achieves bipartite consensus if
\begin{align*}
    \lim_{t\to \infty} \norm{x_i(t) - x_j(t)} = 0 \quad &\text{ for all } i,j \in C_k, k=1,2,\\
    \lim_{t\to \infty} \norm{x_i(t) + x_j(t)} = 0 \quad &\text{ for all } i \in C_1, j \in C_2.
\end{align*}
Notably, bipartite consensus states can occur in systems where the adjacency matrix is allowed to take signed values, which indicate antagonistic relations between certain agents (corresponding to negative values $a_{ij}<0$) and cooperative relations with others (where $a_{ij}>0$). 
\vs
Combining these concepts, a \emph{periodic multiconsensus solution} exists when there is a partition into clusters of agents $\mathcal V = C_1\sqcup C_2 \sqcup\dots \sqcup C_k$ and a $T$-periodic solution to the closed loop dynamics such that for all $t\in \mathbb R,$ and for all $i,j\in C_p$ and all $p=1,\dots,k$, we have $x_i(t) = x_j(t)$\footnote{It is not necessary that different clusters oppose each other in this definition, but as we shall see, this will typically be the case when such periodic multiconsensus solutions arise from symmetry breaking bifurcation.}. A system is said to achieve periodic multiconsensus if there exists a periodic multiconsensus solution which is locally asymptotically stable.
\vs
These concepts have a natural and powerful interpretation in terms of cooperation and competition within populations of agents. If one takes the natural view of consensus as representing cooperation, bipartite consensus as representing antagonism, then periodic multiconsensus solutions where certain subpopulations have opposite signs can naturally represent competitive striving in a conservative or zero-sum scenario, where two synchronized subpopulations are ``chasing'' each other through the state space, each one alternatingly the winner and the loser. 
\vs
A common thread of study is to consider a parametrized MAS, usually with a known consensus equilibrium (e.g. at $x\equiv 0$) and study the parameter ranges under which this system achieves consensus. One linearizes the closed-loop dynamics around the known consensus equilibrium. Under the reasonable assumption that the spectrum of the linearized operator depends smoothly on the parameters, one can then obtain conditions on the parameters which guarantee that all eigenvalues have strictly negative real part, and thus the local asymptotic stability of the consensus equilibrium.
\vs
But what happens at the edge of these parameter ranges, where consensus breaks down? Cooperation might break down, factions could arise within newly antagonistic subpopulations of agents, and they could begin pursuing their respective goals at each other's expense, rather than through mutual cooperation. This is an intriguing possibility for modeling such breakdowns in human relations as well, both between individuals and among larger organizational structures. To study these phenomena, we must look beyond the asymptotic stability of stationary equilibria. 
\vs
Hopf bifurcation is a natural approach for this parallel question in dynamical systems, and is increasingly recognized as such in the study of multi-agent systems as well. However, the high degrees of symmetry present in homogeneous multi-agent systems with symmetric interaction topologies (describing a level initial playing field) is a significant obstacle for classical methods of Hopf bifurcation such as Lyapunov-Schmidt reduction to a center manifold. Equivariant Hopf bifurcation frameworks such as the twisted degree emerge as a natural tool to study such problems and obtain powerful results on periodic multiconsensus solutions, using the symmetric interactions as an asset rather than an obstacle. 
\vs
However, as noted previously, equivariant degree methods cannot tell us if periodic multiconsensus is actually \emph{achieved}, as this is a stability property. It can, however, tell us if certain configurations are possible within the system, and at which parameter values they will occur. The achieval of the actual periodic multiconsensus solution, as a matter of stability, can then be inferred via numerical simulation initialized at the place and parameter value indicated by the degree.

\subsection{Memory}
Now we will make the explicit link between the delay equation concepts studied in Section \ref{sec:intro:dist-delay} and multi-agent systems. Consider a protocol for agent $i$ written $u_i = r(x_i) + h_i(x)$, where $r(x_i)$ is a \emph{self-regulation} term which depends only on the state of agent $i$, and $h_i(x)$ is a network term which depends on the states $x_j$ for $j\in \mathcal N_i$. We say that a protocol possesses \emph{continuous memory} if either $r(x_i)$ or $h_i(x)$ depends on a weighted average of the agent's own past states over a time window. The simplest such term is a distributed delay operator of the form
\[
\int_0^\infty \phi(s) x_i(t-s)ds
\]
with delay kernel $\phi(s)$. We will consider the ``boxcar kernel'' situation where, for $\tau>0$, 
\[
\phi_\tau(s) = \begin{cases}
    1 &\text{ if }0\leq s<\tau\\
    0 &\text{ if }s>\tau,
\end{cases}
\]
which allows the continuous memory operator to be written
\[
\int_0^\tau x(t-s)ds.
\]
Clearly, the closed-loop dynamics of a MAS with a continuous memory protocol can be viewed as a distributed delay differential equation of retarded type. This can be used to model many types of agent interactions where agents remember their own past states but see only the current state of their neighbors, and we will consider both first-integrator and second-integrator MAS having this type of continuous memory.
\vs
We say that a protocol possesses \emph{trend memory} (or \emph{pseudoneutral memory} in a first-integrator context) if the protocol depends on the derivative of a continuous memory term, i.e. (taking the same delay kernel as above)
\[
\frac{d}{dt}\left[\int_o^\tau g(x(t-s))ds\right] = g(x(t)) - g(x(t-\tau)).
\]
where $g:\mathbb R^n \to \mathbb R^n$ is Lipschitz continuous. Notice that this results in a finite difference term of purely retarded type. An agent with this type of memory is comapring its current state (nonlinearly transformed by $g$) to its transformed state delayed by $\tau$. This is a trend indicator and measures how the moving continuous average has changed over time. In the limit of small $\tau$ it approximates $g'(x)\dot x$, providing a derivative-like feedback without the fundmanetal instability problems of a true derivative.
\vs
On the other hand, we say a protocol possesses \emph{momentum memory} (or \emph{neutral memory} in the first-integrator context) if the protocol depends on an actual delayed derivative term
\[
g(\dot x(t-\tau)),
\]
where $g$ again is Lipschitz continuous. The alternate terminologies of pseudoneutral memory and neutral memory are motivated by the fact that, for a single-integrator MAS, the closed loop dynamics of a MAS whose protocol possesses these types of memory becomes a pseudoneutral (cf. Section \ref{sec:prelim:pseudoneutral}) or a neutral equation respectively. 
\vs
In our analysis of MAS with trend memory which form pseudoneutral equations, we will use a ``neutral framework'' to analyze them. By this we mean that we will formulate them implicitly, with the pseudoneutral memory term inside the derivative, and we will analyze them with the Nussbaum-Sadovskii degree. It is easy to see from Section \ref{sec:prelim:n-s} that any admissible pair with a well-defined Leray-Schauder degree also possesses a well-defined Nussbaum-Sadovskii degree which coincides with it exactly, and so it is mathematically sound to take this approach. Formulating a pseudoneutral equation in this implicit form places constraints on the pseudoneutral memory term which are not strictly necessary from the point of view of taking the Leray-Schauder degree of a compact operator. In particular, we require the pseudoneutral memory term to be a contraction, which is not necessary for the solution operator to be compact if the pseudoneutral delay is expanded to the right-hand side as a finite difference.
\vs
We do so in spite of this for three primary reasons: It makes the connection between MAS with this type of memory and those having neutral memory more transparent; the additional requirements on the pseudoneutral memory term also prove necessary for the MAS to achieve consensus at the trivial equilibrium, which is obviously desirable from the point of view of studying the \emph{breakdown} of consensus; and because it allows us to formulate our analysis of the neutral memory MAS in an exactly parallel manner, which directly facilitates the comparison of the two models and their respective dynamics.

\chapter{Global Periodic Multiconsensus Bifurcation in Distributed Delay Multi-Agent Systems with Continuous Memory}\footnote{This chapter is based on joint work with Chaoquan Chen and Travis Hensley, originally published as \cite{crane2026}. The material has been reframed within the multi-agent systems context, extended with new stability results, and the application has been enhanced and expanded.}\label{chapter:distributed}
%  \author{Casey Crane$^*$} 
%  \address{Department of Mathematical Sciences, the University of Texas at Dallas,  Richardson, TX, 75080-3021, U.S.A. }
% \email{cmc102120@UTDallas.edu}

%  \author{ Travis Hensley} 
%  \address{Department of Mathematical Sciences, the University of Texas at
% Dallas,  Richardson, TX, 75080-3021, U.S.A. }
% \email{Tyler.Hensley@UTDallas.edu}

% %\author{Huafeng Xiao}
% %\address{School of Mathematics and Information Science, Guangzhou University Guangzhou, Guangdong 510006, China}
% %\email{huafeng@gzhu.edu.cn}

% %EndAName

% \subjclass{Primary:  34K13, 37J20, 37C81, 47H11 Secondary: 39A23 }
% \date{}
% \maketitle
% \footnote{\hskip.1cm$*$ Corresponding author.}

% \begin{abstract}
% In this paper, we study a system of coupled identical distributed delay differential equations. We transform the problem of the existence of critical points for Hopf bifurcation to that of zeros of an operator equation. Making use of the equivariant degree theory, we show the existence of a non-trivial unbounded symmetric branches of solutions emerging through Hopf bifurcation. Finally, we give an example to explain our results in the context of steady state oscillations in anti-windup PID controllers, and provide numerical simulations to illustrate the theoretical results.
% \end{abstract}

\section{Introduction}\label{sec:dist:introduction}
In this chapter, we consider a homogeneous multi-agent system with a $\Gamma_0$-symmetric interaction topology (where $\Gamma_0$ is a finite group) with state-dependent effective weights, where each agent protocol depends nonlinearly on a continuous memory term, that is, a nonlinear function of a continuous weighted average over a moving window of past states. Consider a network of $n$ agents, and denote the total state vector as $\bm x(t) = (x_1(t),\dots,x_n(t))^T\in \mathbb R^n$, where $x_i(t)$ is the state of agent $i$ at time $t$. Assume there is a fixed interaction topology between these agents, with interaction weights given by the nonlinear function $h:\mathbb R^n \to \mathbb R^n$, such that $h_i(\bm x(t))$ gives the effective weight for agent $i$ at time $t$ given the total state $\bm x(t)$. We consider single-integrator dynamics where each agent has the following protocol:
\begin{equation}\label{eq:dist:protocol}
u_i(x) := -ax_i - \alpha f\left(\int_0^1x_i(t-s)ds\right)-h_i(\bm x),
\end{equation}
where $\int_0^1x(t-s)ds$ is the continuous memory term with the delay length fixed at 1 for simplicity, $f:\mathbb R\to \mathbb R$ is a nonlinear function assumed to be differentiable at 0 with $f'(0) := b$, where we assume $b>0$, $a>0$ is a self-regulation parameter, and $h_i(\bm x)$ governs interactions with other agents. 
Then the closed-loop dynamics equation for agent $i$ can be written
\begin{equation}\label{eq:dist:basic1}
\dot x_i(t) = -ax_i - \alpha f\left(\int_0^1x_i(t-s)ds\right)-h_i(\bm x),
\end{equation}
and the closed-loop dynamics of the full system can be written
\begin{equation}\label{eq:dist:basic1.1}
\dot {\bm x}(t) = -a\bm x - \alpha \bm f\left(\int_0^1x(t-s)ds\right)-h(\bm x),
\end{equation}
where $\bm f:\mathbb R^n \to \mathbb R^n$ is defined as $\bm f(\bm x) = (f(x_1),f(x_2),\dots,f(x_n))$. We will also make the following assumptions:
\begin{enumerate}[label=($C_\arabic*$)]\setcounter{enumi}{0}\itemc
\item\label{dist:c1} the functions $\bm f$ and $h$ are odd, i.e. $\bm f(-\bm x)=-\bm f(\bm x)$, $h(- \bm x)=-h(\bm x)$,
\item\label{dist:c2} the functions $\bm f$ and $h$ are everywhere continuous, and are differentiable at 0, with $b:=f'(0) > 0$ and $Dh(0)=:C$, where $C$ is assumed to be a symmetric matrix,
\item\label{dist:c3} the function $h$ is $\Gamma_0$-equivariant.
\end{enumerate}
Note that the closed-loop dynamics equation \eqref{eq:dist:basic1.1} forms a symmetrically coupled distributed delay differential equation of retarded type.
\vs
Multi-agent systems of this type have many interesting applications. The continuous memory term has natural smoothing properties and is resilient against noise. If we view the state $x_i(t)$ as the distance from the current state of agent $i$ to some set point, then \eqref{eq:dist:protocol} also bears a strong resemblance to a protocol of the ubiquitous PI type, but modified with a finite time horizon, and can certainly be interpreted as such. We will return to this interpretation later in this chapter, and show how this protocol can also be viewed as another way of addressing the integral windup problem in PI controllers and PID controllers in general.
\vs
Taking a single agent on its own, it also has a natural economic interpretation as modeling an asset market where fundamentalist traders (represented by the continuous memory term) analyze an asset's market price relative to some underlying value, averaged over a time horizon, buying if it is undervalued (on average), and selling if it is overvalued. Under this interpretation, the $-ax_i$ term can be interpreted as a liquidity provider or market maker who instantly buys or sells under or overvalued assets. In this interpretation, the market itself may be viewed as an agent whose goal is price discovery. A multi-agent system would then correspond to multiple markets which exert a macroeconomic influence on each other's prices through the larger economy. We will study this interpretation in detail with the more intricate models studied in Chapters \ref{chapter:pseudo} and \ref{chapter:neutral}.
\vs
One interpretation of particular interest involves the control of drones or UAVs. Formation or swarming in UAVs or animals is one of the most classical applications of multi-agent systems, and we will show how \eqref{eq:dist:basic1} can be used to describe a formation control scheme which allows a symmetric formation of UAVs to spontaneously oscillate in complex symmetric patterns with zero mean displacement, for evasive or other purposes.
\vs
We will first state our main results here, in the language of multi-agent systems. We will then restate these theorems in a logically equivalent formulation in the language of dynamical systems, and prove them in this context. We will also relate those theorems back to these in the introduction. For the system \eqref{eq:dist:basic1.1} satisfying the assumptions \ref{dist:c1}---\ref{dist:c4}, we will prove the following main results:

\begin{theorem}\label{thm:dist:mas-asymp}
Denote by $\widehat \alpha$ the smallest $\alpha_{n,j}>0$ satisfying
\[
\alpha_{n,j} = \frac{-a_j \beta_{n,j}}{b\sin \beta_{n,j}}
\]
where $a_j := a + \mu_j$, $\mu_j$ is an eigenvalue of $C:=Dh(0)$, and $\beta_{n,j}$ is a solution of
\[
\beta_{n,j} = \frac{a_j(\cos \beta_{n,j} -1)}{\sin\beta_{n,j}}
\]
for some $j=0,1,2,\dots,r$. If $a_j>0$ for all $j=0,1,\dots,r$, then the trivial consensus at $x\equiv0$ is locally asymptotically stable for all $\alpha \in (0,\widehat\alpha)$.
\end{theorem}
\begin{theorem}\label{thm:dist:mas-local}
If $\beta_0>0$ satisfies
\[
\beta_0 = \frac{a_j(\cos \beta_{0} -1)}{\sin\beta_{0}}
\]
and $\alpha_0>0$ satisfies
\[
\alpha_0 = \frac{-a_j \beta_0}{b\sin \beta_{0}}
\]
for some $j$, and if 
\[
\alpha_0(a_j + a_j^2 + \beta_0^2 -\alpha_0b) \neq 0
\]
then there exists a connected branch of non-constant periodic multi-consensus solutions bifurcating from the trivial consensus branch at $\alpha=\alpha_0$, and this branch is global and unbounded.
\end{theorem}
\begin{theorem}\label{thm:dist:mas-sym}
If the conditions of Theorem \ref{thm:dist:mas-local} hold, and $(H)<\Phi_1^t(S^1\times \Gamma_0\times \mathbb Z_2;\mathscr E_{1,j})$ is a maximal twisted orbit type in $\mathscr E_{1,j}:=\{c\cos(t)+d\sin(t):c,d\in E(\mu_j)\}$, then there exists a connected branch of non-constant periodic solutions with symmetries at least $(H)$ bifurcating from the trivial solution at $\alpha=\alpha_0$, and this branch is global and unbounded.
\end{theorem}

The following four sections of this paper are structured as follows: In section \ref{dist:sec2}, we will study the closed-loop dynamic equations \eqref{eq:dist:basic1.1} as a system of distributed delay differential equations with symmetric coupling, compute the characteristic equation of the linearization of this system, show that consensus is achieved at $x \equiv 0$ (i.e. that $x \equiv 0$ is a locally asymptotically stable equilibrium), and find transversality conditions allowing us to determine the sign of the crossing number at every critical point. 
\vs
At the end of Section \ref{dist:sec2}, we will state a main theorem showing the bifurcation of the trivial consensus at $x \equiv 0$ into non-constant periodic multiconsensus solutions, global continuation of periodic multiconsensus branches, and an immediate equivariant corollary showing how information on the symmetries of these branches can also be obtained from local bifurcation invariants. 
\vs
We will then prove this theorem and corollary in Section \ref{dist:sec3}, by period normalizing the system and reformulating it as an operator equation between functional spaces. This will allow us to apply the tools of equivariant degree theory, namely the equivariant twisted Leray-Schauder degree. 
\vs
In Section \ref{dist:sec4}, we will perform a fixed-point reduction which strengthens the main result by allowing us to apply the main theorems even in the presence of degeneracies, and allowing us to characterize the symmetries of branches more globally.
\vs
Finally, in Section \ref{dist:sec5}, we will show how a MAS of this type can be used to control spontaneous tuneable oscillations with zero mean displacement in symmetric formations of UAVs. We also provide numerical simulations which confirm the theoretical results, showing the emergence of branches of periodic solutions with the predicted limit frequencies and symmetries at the predicted critical parameter values.
\vspace{32pt}
\section{Formulation of the problem}\label{dist:sec2}
We consider the closed-loop dynamics equation \eqref{eq:dist:basic1}. More precisely, let $\Gamma_0$ be a finite group and let $V:={\mathbb R}^n$. Assume that $\Gamma_0$ acts on $V$ by permuting components of vectors in $V$, i.e. $\Gamma_0 \leq S_n$ and for $\sigma \in \Gamma_0$,
\[
\sigma(x_1,x_2,\dots,x_n)^T = (x_{\sigma(1)},x_{\sigma(2)},\dots,x_{\sigma(n)})^T,\quad \bm x =(x_1,\dots,x_n)^T \in {\mathbb R}^n.
\]
Notice that, under \ref{dist:c1}---\ref{dist:c3}, the system \eqref{eq:dist:basic1.1} is $\Gamma_0 \times {\mathbb Z}_2$-symmetric (with the antipodal ${\mathbb Z}_2$ action on $V$). As $V$ is a $\Gamma_0 \times {\mathbb Z}_2$-representation, consider the $\Gamma_0 \times {\mathbb Z}_2$-isotypic decomposition
\begin{equation}\label{eq:dist:iso-decomp}
V = V_0 \oplus V_1 \oplus \dots \oplus V_r,
\end{equation}
where the $j$-th component $V_j$ is modeled on the $j$-th irreducible $\Gamma_0$-representation ${\mathcal V}_j$ with the antipodal ${\mathbb Z}_2$-action which makes it an irreducible $\Gamma_0 \times {\mathbb Z}_2$-representation denoted by ${\mathcal V}_j^-$.
\vs

In order to avoid unnecessary notational complexity in the analysis of the Hopf bifurcation problem for \eqref{eq:dist:basic1.1}, we make an additional simplifying assumption:

\begin{enumerate}[label=($C_\arabic*$)]\setcounter{enumi}{3}\itemc
\item\label{dist:c4} The spectrum $\sigma(C)$ of the matrix $C$ is given by $\sigma(C) = \{\mu_0,\mu_1,\dots,\mu_r\}$, and the eigenspaces $E(\mu_j)$ coincide with the isotypic components $V_j$ in \eqref{eq:dist:iso-decomp}.
\end{enumerate}

\vs
\subsection{Characteristic equation}
The characteristic operator for \eqref{eq:dist:basic1.1} is given by 
\begin{equation}\label{eq:dist:char1}
\triangle_\alpha(\lambda) := (\lambda + a)\text{Id} - \alpha b \tfrac{(e^{-\lambda}-1)}{\lambda}\text{Id} + C.
\end{equation}
Let $\Gamma := \Gamma_0 \times {\mathbb Z}_2$. Notice that this characteristic operator is $\Gamma$-equivariant, therefore one can introduce
% \begin{equation}\label{eq:dist:char2}
\[
    \triangle_{\alpha,j}(\lambda) := \triangle_\alpha(\lambda)_{\vert{\mathcal{V}_j^{\mathbb C}}}:\mathcal{V}_j^{\mathbb C} \to \mathcal{V}_j^{\mathbb C},
\]
    % \end{equation}
where $\mathcal{V}_j^{\mathbb C}$ denotes the complexification of the irreducible representation $\mathcal V_j$.
Observe that, under the assumptions \ref{dist:c1}--\ref{dist:c4}, one has
% \begin{equation}\label{eq:dist:char3}
\[
    \triangle_{\alpha,j}(\lambda) = \big((\lambda + a) - \alpha b\tfrac{(e^{-\lambda}-1)}{\lambda}+\mu_j\big)\text{Id},
\]
    % \end{equation}
i.e. 
% \begin{equation}\label{eq:dist:char4}
\[
\triangle_\alpha(\lambda) = \begin{bmatrix}
    \triangle_{\alpha,0}(\lambda) & &\makebox(0,0){\text{\huge0}}\\
    & \ddots & \\
    \text{\huge0}& & \triangle_{\alpha,r}(\lambda)
    \end{bmatrix},
\]
% \end{equation}
% \begin{equation}\label{eq:dist:char5}
\[
    {\det}_{\mathbb C}(\triangle_\alpha(\lambda)) = \prod_{j=0}^r {\det}_{\mathbb C}(\triangle_{\alpha,j}(\lambda))^{m_j},
\]
    % \end{equation}
where we define the {\it isotypic multiplicity} $m_j$ of the eigenvalue $\mu_j$ as
\begin{equation}\label{eq:dist:iso-mult}
m_j:=\frac{\text{dim}\;E(\mu_j)}{\text{dim}\;{\mathcal V}_j}
\end{equation}
and we have
\begin{equation}\label{eq:dist:charfunc}
P_{\alpha,j}(\lambda):={\det}_{\mathbb C}(\triangle_{\alpha,j}(\lambda)) = (\lambda + a) - \alpha b \tfrac{(e^{-\lambda}-1)}{\lambda} + \mu_j.
\end{equation}
Put $a_j := a + \mu_j$. Then $\lambda$ is a characteristic root for \eqref{eq:dist:basic1.1} if there exists a $j=0,1,\dots,r$ such that the equation \eqref{eq:dist:j-iso-char} is satisfied.
% \[
% (\lambda - a) + \alpha b \tfrac{(e^{-\lambda}-1)}{\lambda} + \mu_j = 0
% \]
% \[
% \lambda - (a + \mu_j) + \alpha b \tfrac{(e^{-\lambda}-1)}{\lambda} + \mu_j = 0
% \]
% \[
% \lambda^2 - (a + \mu_j)\lambda + \alpha b (e^{-\lambda}-1) + \mu_j = 0
% \]
\begin{equation}\label{eq:dist:j-iso-char}
\lambda^2 +  a_j \lambda - \alpha b (e^{-\lambda}-1) = 0.
\end{equation}
We will call $P_{\alpha,j}$ the \textit{characteristic function} for \eqref{eq:dist:basic1.1}, and \eqref{eq:dist:j-iso-char} the $j$-th {\it isotypic characteristic equation} for \eqref{eq:dist:basic1.1}.
In such a case, $\lambda$ satisfying \eqref{eq:dist:j-iso-char} is called the $j$-th isotypic characteristic root with corresponding isotypic multiplicity $m_j$.
\vs
% \subsection{Centers, critical points and crossing numbers}
% \begin{definition}\label{def:center} \rm
% Assume that the system \eqref{eq:dist:basic1.1} satisfies conditions \ref{dist:c1}--\ref{dist:c4}. A trivial solution $(\alpha_0,0)$ to \eqref{eq:dist:basic1.1} is called a \textit{center} if there exists $\beta_0 > 0$ (called its \textit{frequency}) such that $\det_{\mathbb C}\triangle_{\alpha_0}(i\beta_0)=0$. In addition, we will say that the center $(\alpha_0,0)$ with the frequency $\beta_0>0$ is \textit{isolated} if there exists $\epsilon > 0$ such that if $0 < \vert \alpha - \alpha_0 \vert + \vert \beta - \beta_0 \vert < \epsilon$ then for all $k=0,1,2,\dots$
% \[
% {\det}_{\mathbb C}\triangle_{\alpha}(i\beta) \not= 0.
% \]
% \end{definition}
% In what follows, for a center $(\alpha_0,0)$ with the limit frequency $\beta_0$, we will simply write $(\alpha_0,\beta_0,0)$, and we will refer to this as a \textit{critical point} for \eqref{eq:dist:basic1.1}.
% \vs
% The following theorem allows us to identify all the critical points for \eqref{eq:dist:basic1.1}. We denote this set of all critical points by $\Lambda$.
% \vs
Recalling the definitions of isolated centers and critical points from Section \ref{sec:prelim:two-parameter-bifurcation}, we have the following proposition regarding the critical set of \eqref{eq:dist:basic1.1}:
\begin{proposition}\label{prop:dist:crit1} 
The problem \eqref{eq:dist:basic1.1} admits an infinite number of critical points $(\alpha_0,\beta_0,0)$ which can be described as the solutions to the following system of equations
\begin{align}
\alpha_{n,j}&=\frac{- a_j \beta_{n,j}}{b\sin\beta_{n,j}}, \quad n=1,2,3,\dots, \; j=0,1,\dots,r,\label{eq:dist:a_n}\\
\beta_{n,j}&=\frac{ a_j (\cos\beta_{n,j}-1)}{\sin \beta_{n,j}},\label{eq:dist:b_n}
\end{align}
where $\beta_{n,j}$ denotes the sequence of solutions to the equation \eqref{eq:dist:basic1.1} (see Figure \ref{fig:dist:fig1}). In other words, the set
\begin{equation}\label{eq:dist:crit-set}
\Lambda := \{(\alpha_{n,j},\beta_{n,j},0) : n=1,2,\dots,\;j=0,1,\dots\,r,\; \alpha_{n,j}>0,\beta_{n,j}>0\}
\end{equation} is discrete.
\end{proposition}
\begin{proof}
Notice that for a characteristic equation given by \eqref{eq:dist:char1}
\[
\lambda=\mathfrak - a_j +\frac{\alpha b}{\lambda}(e^{-\lambda}-1)\quad \Leftrightarrow\quad \lambda^2+\lambda  a_j -\alpha b(e^{-\lambda}-1)=0.
\]
If $\lambda = i\beta$, then \eqref{eq:dist:char1} can be written as
\[
-\beta^2+i\beta a_j  - \alpha b(\cos\beta-i\sin\beta-1)=0,
\]
which leads to the following system of equations
\[
\begin{cases}
\beta^2+\alpha b(\cos\beta-1)=0\\
\beta  a_j +\alpha b\sin\beta=0
\end{cases}
\quad \Leftrightarrow\quad
\begin{cases}
\cos\beta=\frac{-\beta^2+\alpha b}
{\alpha b}\\
\sin\beta=\frac{-\beta a_j }{\alpha b},
\end{cases}
\]
Then, by substituting $\alpha=\frac{-\beta  a_j }{b\sin\beta}$ into the first equation, one gets
\[
\cos\beta=\frac{\beta^2+\frac{\beta  a_j }{\sin\beta}}{\frac{\beta  a_j }{\sin\beta}}\quad \Leftrightarrow\quad \frac{ a_j (\cos\beta-1)}{\sin \beta}=\beta.
\]

%or for the solutions $\alpha$ and $\lambda=i\beta$ to \eqref{eq:dist:char2}, i.e.
%\[
%-\beta^2+i\alpha \beta a- b(\cos\beta-i\sin\beta-1)=0,
%\]
%i.e. we need to solve the system of equations
%\begin{equation}\label{eq:dist:ch-2-2}
%\begin{cases}
%-\beta^2-b(\cos\beta-1)=0\\
%\alpha \beta a-b\sin\beta=0
%\end{cases}
%\quad \Leftrightarrow\quad
%\begin{cases}
%\cos\beta=\frac{\beta^2- b}
%{%b}\\
%\sin\beta=\frac{\alpha a \beta}{ b}.
%\end{cases}
%\end{equation}
%Let us consider first the system \eqref{eq:dist:ch-2-2}. For the system \eqref{eq:dist:ch-2-2}, in order to have solutions one needs 
%\[
%\left|\frac{\beta^2- b}b\right|\le 1\quad \Leftrightarrow\quad \beta^4-2b\beta^2\le 0,
%\]
%which implies $b\ge \frac 12$.
%Then, the equation $ \cos\beta=\frac{\beta^2- b}b$ has  a unique solution (remember $\beta>0$).

\begin{figure}[H]
	\centering
	\includegraphics[height=6cm]{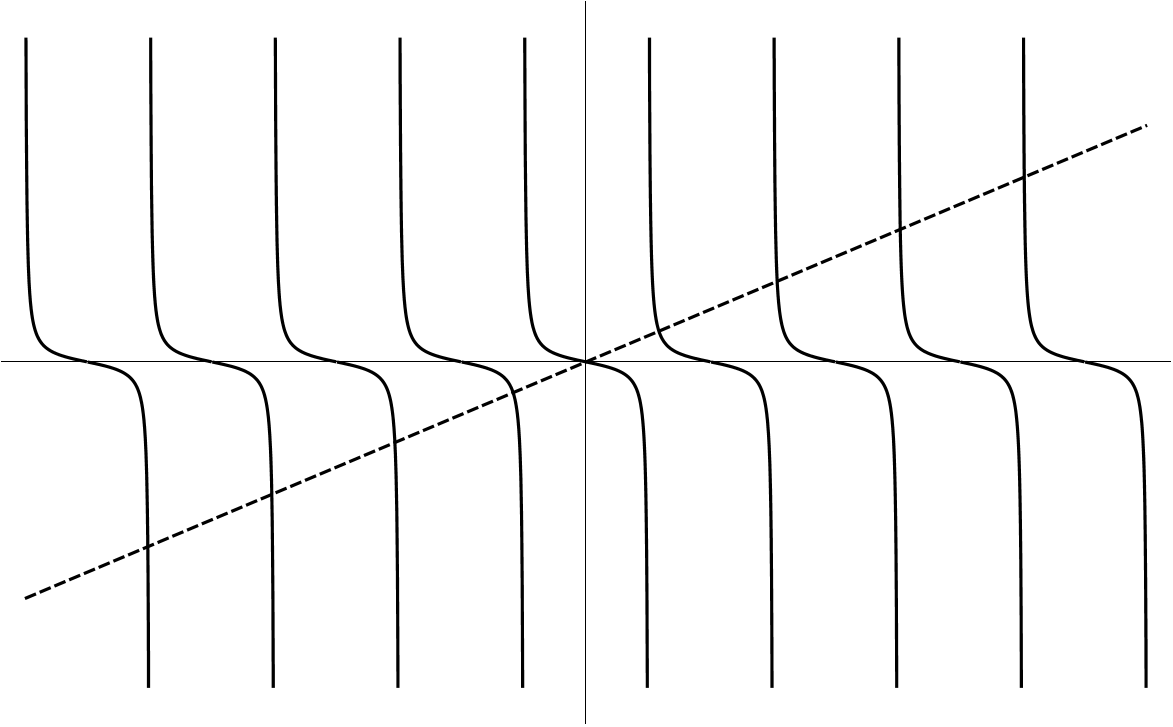} 
	\caption{Critical frequencies for \eqref{eq:dist:basic1}}\label{fig:dist:fig1}
	\rput(3.3,4.2){}
	\rput(-3.5,7.5){}
\end{figure}

By comparing the graphs of the functions $\vp(\beta)=\frac{ a_j  (\cos\beta-1)}{\sin \beta}$ and $\zeta(\beta)=\beta$ (see Figure \ref{fig:dist:fig1}) one can determine that there exists an infinite set of critical frequencies $\beta_{n,j}>0$ and corresponding values $\alpha_{n,j}:=\frac{- a_j \beta_{n,j}}{b\sin\beta_{n,j}}$. Notice that the values $\beta_n$ are regularly distributed, thus the set $\Lambda$ is discrete. 
\end{proof}
\vs
\begin{remark}
\rm Notice that the set $\Lambda$ is well defined. Indeed, since $b<0$ and we take only positive solutions to the equation $\zeta(\beta) = \beta$, one obtains that if $\beta_{n,j}>0$, then $\alpha_{n,j} > 0$. Therefore, the set $\Lambda$ is infinite, and $\alpha_{n,j}$, $\beta_{n,j} \to \infty$ as $n \to \infty$.
\end{remark}
By Proposition \ref{prop:dist:crit1},  each critical point  $(\alpha_0,\beta_0,0)\in \Lambda$  corresponds to the isolated center $(\alpha_o,0)$ with the critical frequency $\beta_o$. Thus for a sufficiently small $\delta>0$ and  $0<\max\{\vert \alpha-\alpha_0 \vert, \vert \beta-\beta_0 \vert \} \le \delta$, one has $(\alpha_0,k\beta_0,0)  \in \Lambda$ for all $k\in \bn$. We will call the set
\begin{equation}\label{eq:dist:isol-neigb}
U(\alpha_o,\beta_o):=\{(\alpha,\beta): \max\{\vert \alpha-\alpha_0 \vert, \vert \beta-\beta_0 \vert \} < \delta\}
\end{equation}
 an {\it isolating neighborhood} of $(\alpha_0,\beta_0,0)$.
\vs
Let $(\alpha_0,\beta_0,0) \in \Lambda$. Since the critical points in $\Lambda$ are isolated, there exist $\delta>0, \eta>0$ such that the set
\begin{equation}\label{eq:dist:crossnum-omega}
\Omega := \{u+i\beta \in {\mathbb C} : 0<u<\eta, \vert \beta - \beta_0 \vert < \delta\}
\end{equation}
satisfies 
\[
{\det}_{{\mathbb C}}\triangle_\alpha(u+i\beta) \not= 0,\] whenever
\begin{enumerate}
\item $0\leq\vert \alpha-\alpha_0 \vert \leq \delta$, $\vert \beta-\beta_0 \vert = \delta$ and $0\leq u \leq \eta$,
\item $0\leq\vert \alpha-\alpha_0 \vert \leq \delta$, $0\leq \vert \beta-\beta_0 \vert \leq \delta$ and $ u = \eta$,
\item $u=0$, $0\leq \vert \beta - \beta_0 \vert \leq \delta$, $0\leq \vert \alpha - \alpha_0 \vert \leq \delta$, except$(\alpha_0,\beta_0,0)$.
\end{enumerate}
In particular, the above conditions mean that the characteristic roots $\lambda_\alpha = u_\alpha + i\beta_\alpha$ for \eqref{eq:dist:basic1.1} can only enter or exit the region $\Omega$ when $\alpha=\alpha_0$ through the point $i\beta_0$.
The \textit{crossing number} $\mathfrak t(\alpha_0,\beta_0)$ can be conceptually understood as the net algebraic count of the eigenvalues $\lambda_\alpha$ exiting the region $\Omega$ (if $\mathfrak t(\alpha_0,\beta_0)>0$) or entering $\Omega$ (if $\mathfrak t(\alpha_0,\beta_0)<0$) when $\alpha$ is crossing $\alpha_0$.
\vs
Put $\alpha_\pm:=\alpha\pm\delta$ and denote by $\mathfrak t_\pm(\alpha_0,\beta_0)$ the number of characteristic values $\lambda_{\alpha_\pm}$ in $\Omega$ counted with their algebraic multiplicities. Then clearly $\mathfrak t(\alpha_0,\beta_0)$ = $\mathfrak t_-(\alpha_0,\beta_0) - \mathfrak t_+(\alpha_0,\beta_0)$.
\vs
Following this line of thinking, we introduce the concept of crossing numbers in a context suitable for equivariant analysis of Hopf bifurcation. More precisely, using the notation described above, we have the following definition:
\begin{definition}\rm
Assume $(\alpha_0,\beta_0,0)\in\Lambda$ and consider the set $\Omega$ given by \eqref{eq:dist:crossnum-omega}. Put $\alpha_\pm = \alpha_0 \pm \delta$. Then for $j=0,1,\dots,r$, we define the so-called $j$-th \textit{isotypic crossing number} $\mathfrak t^j(\alpha_0,\beta_0)$ to be
\[
\mathfrak t^j(\alpha_0,\beta_0):=\mathfrak t_-^j(\alpha_0,\beta_0)-\mathfrak t_+^j(\alpha_0,\beta_0),
\]
where $\mathfrak t_\pm^j(\alpha_0,\beta_0)$ denotes the number of $j$-th isotypic characteristic roots of \eqref{eq:dist:basic1.1} in $\Omega$ for $\alpha=\alpha_\pm$ counted with their $j$-th isotypic multiplicities.
\vs
For use in proving the global result in subsequent sections we also define, for $k\in\bbN$, the \textit{$k$-resonant $j$-th isotypic crossing numbers} 
\begin{equation}\label{eq:dist:k-res-crossnum}
\mathfrak t_k^j(\alpha_0,\beta_0) := \mathfrak t_j(\alpha_0,k\beta_0),\quad k\in\bbN.
\end{equation}
\end{definition}
\vs
\begin{remark}\rm
Let us point out that the $j$-th crossing number can also be represented by 
\[\mathfrak t^j(\alpha_0,\beta_0) = \text{deg}\left({\det}_{\mathbb C} \triangle^j_\alpha(i\beta),U(\alpha_0,\beta_0)\right),
\]
where $\text{deg}$ stands for the usual Brouwer degree.
\end{remark}
\vs
\subsection{Stability of the \texorpdfstring{$x\equiv 0$}{x ≡ 0} consensus equilibrium}
First, we observe that $x\equiv 0$ is clearly both an equilibrium and a consensus state (cf. Section \ref{sec:prelim:consensus}) for system \eqref{eq:dist:basic1.1}. Recall that the $x\equiv 0$ consensus (which we will also call the \emph{trivial consensus}) is said to be achieved in system \eqref{eq:dist:basic1.1} if this equilibrium is locally asymptotically stable.

\begin{proposition}
Given system \eqref{eq:dist:basic1.1} satisfying assumptions \ref{dist:c1}---\ref{dist:c3}, let $\widehat \alpha := \min_{\alpha \in \mathbb R}\{(\alpha,\beta,0)\in \:\Lambda\}$. If $a_j = a + \mu_j>0$ for all $j=0,1,\dots,r$, then for all $\alpha \in (0,\widehat \alpha)$, the trivial solution $x\equiv 0$ is locally asymptotically stable.
\end{proposition}
\begin{proof}
Since the roots of the characteristic equation \eqref{eq:dist:char1} are related to the spectrum of the linearization of \eqref{eq:dist:basic1.1} and depend continuously on $\alpha$, it is sufficient to show that for $\alpha=0$, all roots to \eqref{eq:dist:char1} have strictly negative real part. Since $\widehat \alpha$ is taken to be the minimal critical value of $\alpha$ across all $\Gamma_0$-isotypic components, this is sufficient to show the local asymptotic stability of the trivial equilibrium for all $\alpha \in (0,\widehat \alpha)$, and show that no steady-state bifurcation is possible for $\alpha \in(0,\widehat \alpha)$, since the set $\Lambda$ does not detect steady-state bifurcation points.
\vs
Setting $\alpha=0$ in \eqref{eq:dist:char1}, we obtain
\[
\lambda +  a_j = 0,
\]
and clearly no roots with positive real part exist if $a_j>0$ for all $j=0,1,\dots,r$. Now we consider \eqref{eq:dist:j-iso-char}, and show that it has no positive real roots for all $\alpha \in (0,\infty)$. We first take $\lambda = u+iv$ and separate the equation into real and imaginary parts, yielding
\begin{align*}
    u^2-v^2+ a_j u-\alpha b(e^{-u}\cos v-1)&=0\\
    2uv+ a_j  v+\alpha be^{-u}\sin v &= 0.
\end{align*}
Since we search for real roots, we set $v=0$ and obtain
\begin{align*}
    u^2+ a_j u-\alpha b(e^{-u}-1)&=0.
\end{align*}
Let $g_1(u) = u^2+a_ju$ and $g_2(u) = \alpha b (e^{-u}-1)$. Then a positive real root of \eqref{eq:dist:j-iso-char} exists if and only if the graphs of $g_1(u)$ and $g_2(u)$ intersect in the $u>0$ half plane of $\mathbb R^2$. However, $g_1(u)>0$ for all $u>0$, and $g_2(u)<0$ for all $u>0$, so this is impossible (and therefore steady-state bifurcation is impossible if $a_j>0$ for all $j=0,1,\dots,r$). This also suffices to prove Theorem \ref{thm:dist:mas-asymp}.
\end{proof}

\subsection{Transversality}
It is important to note that the local Krasnosel'skii theorem used to obtain local bifurcation in this dissertation does not actually require the kind of strict transversality condition that we prove in this section, in contrast to classical Hopf bifurcation theorems. All that is needed is a nonzero crossing number, which can itself be computed as a difference in the net count of eigenvalues in a small region around a limit frequency $i\beta_0$ in the positive half plane, counted with their isotypic multiplicities. However, from the point of view of obtaining general results on the signs of \emph{all} crossing numbers for use in the global continuation result, it is often actually more convenient to compute the stronger transversality condition, which immediately allows one to compute crossing numbers.
\vs
The equation \eqref{eq:dist:char1}, for $\lambda(\alpha):=u(\alpha)+iv(\alpha)$  can be written as the system
\begin{equation}\label{eq:dist:char-split}
\begin{cases}
u^2-v^2+ a_j u-\alpha b(e^{-u}\cos v-1)=0,\
    2uv+ a_j  v+\alpha be^{-u}\sin v = 0.
    \end{cases}
\end{equation}

Since we are interested in purely imaginary characteristic roots $\lambda = i \beta$, by substituting $u = 0$ and $v = \beta$ into the system \eqref{eq:dist:char-split}, we obtain
%  \begin{equation}\label{eq:dist:char3}
\[
 	\begin{cases}
		-\beta^2 - \alpha b(\cos\beta-1) = 0,\
		 a_j  \beta + \alpha b \sin \beta = 0,
	\end{cases}
\]
%  \end{equation}
which can be easily transformed to
\begin{equation}\label{eq:dist:char4}
\begin{cases}
	\cos \beta = \frac{-\beta^2}{\alpha b} + 1,\\
	\sin \beta = \frac{- a_j  \beta}{\alpha b}.
\end{cases}
\end{equation}
In order to determine the value of the crossing number associated with a purely imaginary characteristic root $\lambda_0 = i \beta_0$, we compute $\tfrac{d}{d \alpha} u(\alpha)$ by implicit differentiation. By differentiating the system \eqref{eq:dist:char-split} with respect to $\alpha$,
\begin{equation}\label{eq:dist:char5}
\begin{cases}
2u u' - 2v v' + a_j u' + \alpha b(e^{-u}\cos v \cdot u' +e^{-u}\sin v \cdot v') = b(e^{-u} \cos v -1),\\
2u v' + 2v u' + a_j v' + \alpha b(-e^{-u}\sin v \cdot u' +e^{-u}\cos v \cdot v') = -be^{-u} \sin v.
\end{cases}
\end{equation}
By substituting $\alpha = \alpha_0$, $u = 0$ and $v = \beta_0$ into the system \eqref{eq:dist:char5}, we obtain
% \begin{equation}\label{eq:dist:char6}
\[
\begin{cases}
u'( a_j +\alpha_0 b\cos \beta_0) - v'(2\beta_0-\alpha_0 b\sin \beta_0) = b(\cos \beta_0 - 1),\\
u'(2\beta_0-\alpha_0 b\sin \beta_0) + v'( a_j +\alpha_0 b\cos \beta_0) = -b\sin \beta_0.
\end{cases}
\]
% \end{equation}
Put $p= a_j  - \beta_0^2 + \alpha_0 b$, $q=2\beta_0 +  a_j \beta_0$. Then \eqref{eq:dist:char4} can be written
\begin{equation}\label{eq:dist:char7}
\begin{cases}
u'_0 p - v'_0 q = \frac{-\beta_0^2}{\alpha_0},\\
u'_0 q + v'_0 p = \frac{ a_j  \beta_0}{\alpha_0},
\end{cases}
\end{equation}
where $\frac{du}{d\alpha}(\alpha_0) =: u'_0$, and $\frac{dv}{d\alpha}(\alpha_0) =: v'_0$. Solving the linear system \eqref{eq:dist:char7} yields
\[
u'_0 = \frac{ a_j  \beta_0 q - \beta_0^2 p}{(p^2 + q^2) \alpha_0}.
\]
Therefore
\begin{equation}\label{eq:dist:du-dalpha}
% \begin{cases}
% \tfrac{du}{d\alpha} > 0 &\;\; \text{if }\;\alpha_0( a_j  +  a_j ^2 + \beta_0^2 -\alpha_0b) > 0\\ 
% \tfrac{du}{d\alpha} < 0 \quad &\;\;\text{otherwise} 
% \end{cases}.
\sign \tfrac{du}{d\alpha}(\alpha_0) = \sign(\alpha_0( a_j  +  a_j ^2 + \beta_0^2 -\alpha_0b))
\end{equation}
\vs

\begin{remark}
\rm It is possible that $\tfrac{du}{d\alpha}(\alpha_0)=0$, but since $\alpha_{n,j},\beta_{n,j} \to \infty$ as $n \to \infty$, such a situation can occur only for finitely many $n$, and moreover, for all but finitely many $(\alpha_{n,j},\beta_{n,j},0)\in \Lambda$, $\sign(\tfrac{du}{d\alpha}(\alpha_0)) = -1$. Since $\tfrac{du}{d\alpha}(\alpha_0)=0$ can occur for only finitely many $n$, we will consider this as a degeneracy. These degeneracies can be removed by fixed-point reduction which will be applied later in the proof of the global result. Nevertheless, this situation does not affect the local bifurcation result.
\end{remark}
The above result gives us conditions for the occurrence of local Hopf bifurcation at critical points, but cannot tell us anything about the global persistence of these bifurcations or how they relate to the inherent symmetries of the system. These bifurcating branches of solutions could be global and unbounded or they could be transient, merging back into the trivial branch and disappearing at some other critical points. The global unbounded branches are obviously of greater interest, but to show their existence we must introduce some further tools, namely the twisted equivariant degree. We will do so in order to prove the following main result, along with an associated corollary which provides information on the minimal symmetries of bifurcating solutions.
\begin{theorem}\label{thm:dist:main1}
Suppose that $f:{\mathbb R} \to {\mathbb R}$ is a continuous odd function such that $b:=f'(0)<0$, and assume that $\bm f:V \to V$ and $h:V\to V$ are continuous and satisfy \ref{dist:c1}---\ref{dist:c4}. Then for every critical point $(\alpha_0,\beta_0,0)\in \Lambda$ such that $\mathfrak t_j(\alpha_0,\beta_0)\neq0$ for some $j=0,1,\dots,r$, there exists a branch of non-constant periodic solutions to \eqref{eq:dist:basic1.1} bifurcating from $(\alpha_0,0)$ with the limit frequency $\beta_0$. Moreover, for sufficiently large $n$ and $(\alpha_{n,j},\beta_{n,j},0)\in \Lambda$ (given by \eqref{eq:dist:a_n}, \eqref{eq:dist:b_n}), the bifurcating branch is unbounded. 
\end{theorem}
\begin{corollary}\label{thm:dist:main2}
Under the same assumptions as Theorem \ref{thm:dist:main1}, if there exists $\lambda_0 \in \bm \Lambda$ (as defined in \eqref{eq:dist:crit-set-bold}) such that for some $j=0,1,\dots,r$, $k>0$, and $(H) \in \phi^t_1(G)$, we have $\mathfrak{t}_{k,j}(\lambda_0) \not= 0$ and $\text{coeff}^H(\text{deg}_{{\mathcal V}_{k,j}^-}) \not= 0$, then there exists a branch $\mathscr C$ of nontrivial solutions to \eqref{eq:dist:basic1.1} bifurcating from $(\lambda_0,0)$ with $(\lambda_0,0) \in \overline{\mathscr C^{(H)}}$.
\end{corollary}
This theorem and its corollary both follow from the equivariant Krasnosel'skii theorem as formulated in Section \ref{sec:prelim:krasnoselskii}. However, in order to apply this theorem to this problem, we must first reformulate it as an $S^1 \times \Gamma_0\times \mathbb Z_2$-equivariant compact perturbation of identity, which will be the focus of the following two sections. However, if we assume for the moment that the system \eqref{eq:dist:basic1.1} admits such a reformulation, then we can provide a partial proof of the above theorem and corollary, which will be fully completed following the fixed-point reduction in Section \ref{sec:dist:fpr}. 
\vs
First, observe that the sign of the derivative $\tfrac{du}{d\alpha}(\alpha_0)$ determines the sign of the crossing number $\mathfrak t_{k,j}(\alpha_{n,j},\beta_{n,j})$ as given by Lemma \ref{lem:prelim:sign-crossing-num}. This means that $\tfrac{du}{d\alpha}(\alpha_0)\neq0$ implies a nonzero crossing number, which, by the computational formula for the local bifurcation invariant given in Section \ref{sec:prelim:local-bif-twisted}, implies that if $\text{coeff}^H(\text{deg}_{{\mathcal V}_{k,j}^-}) \not= 0$, then $\text{coeff}^H(\omega_G(\alpha_{0},\beta_0)) \neq 0$, and the result follows. We also have that the sign of the crossing numbers is constant for all $n$ sufficiently large. This implies, by the equivariant Rabinowitz theorem, that the sum of local bifurcation invariants is nonzero, which gives the global unboundedness of the branch.
\vs
There are two remaining matters to complete the proof, which will also complete the proofs of Theorem \ref{thm:dist:mas-local} and Theorem \ref{thm:dist:mas-sym}. As mentioned earlier, we must actually show that this problem can be viewed as an $S^1 \times \Gamma_0 \times \mathbb Z_2$-equivariant compact perturbation of identity. Second, in order to show that the non-constant periodic branches are actually non-constant and periodic across their full extent in the parameter-function space, we require a fixed-point reduction to a space which contains only non-constant solutions. The first of these objectives will be obtained in Section \ref{dist:sec3} and \ref{dist:sec4}, and the second will be obtained in Section \ref{dist:sec5}.
\vs
\begin{remark}\normalfont
The above theorem and corollary contain no particular requirement on the multiplicity of eigenvalues in order to show the existence of periodic solutions. This is one of the great strengths of equivariant Hopf bifurcation results. While classical Hopf bifurcation results require a transversality condition and simple eigenvalues, the equivariant Hopf bifurcation theorem accounts for the higher multiplicity of eigenvalues which are forced due to symmetries. Furthermore, the crossing numbers used in the equivariant Hopf bifurcation theorem form a more general and robust condition for local bifurcation, and prove existence topologically with respect to the $G$-isotypic components. The requirement of simple eigenvalues is therefore replaced with the requirement that eigenvalues be $G$-isotypically simple, i.e. crossing on only one $G$-isotypic component, and even this requirement can be easily sidestepped using fixed-point reduction techniques. The condition \eqref{eq:dist:du-dalpha} is indeed necessary for transversality in classical Hopf bifurcation, but in this result is used indirectly to find the sign of crossing numbers, which is needed for the global result. 
\end{remark}
\section{Setting in functional spaces}\label{dist:sec3}
In order to apply the twisted equivariant degree method, the problem \eqref{eq:dist:basic1.1} must be reformulated in the appropriate functional spaces. More precisely, first we normalize the period in  \eqref{eq:dist:basic1.1}  by introducing the frequency $\beta>0$ as additional parameter, i.e.
we substitute  $u(t) := \bm x(\frac{p}{2\pi}t)$ and taking as a new parameter $\beta := \frac{2\pi}{p}$, we obtain the equation
\begin{equation}\label{eq:dist:basic2}
\dot u(t) = -\frac{a}{\beta}u(t) - \frac{\alpha}{\beta}\bm f \left(\int_0^1 u(t - \beta s) ds \right) + \frac{1}{\beta}h(u(t)).
\end{equation}
\vs
Then, we can reformulate equation \eqref{eq:dist:basic2} as a two-parameter bifurcation problem. First, we introduce the appropriate functional spaces and corresponding norms
\[
\mathscr{E} := C^1_{2\pi}({\mathbb R} ; {\mathbb R}^n), \quad \| \phi \| = \max\{\|\phi\|_\infty,\|\dot\phi\|_\infty\},\;\; \phi \in \mathscr{E},
\]
where $\| \cdot \|_\infty$ denotes the usual supremum norm. Then we introduce the following operators:
\begin{align*}
L&:\mathscr E\to C_{2\pi}({\mathbb R} ; {\mathbb R}^n) , \quad (Lu)(t):=\dot u (t),\\
N&:{\mathbb R}^2\times C_{2\pi}({\mathbb R} ; {\mathbb R}^n)\to C_{2\pi}({\mathbb R} ; {\mathbb R}^n), \quad N(\alpha,\beta,v)(t):=\frac{\alpha}{\beta}\bm f \left(\int_0^1 v(t - \beta s) ds \right) - \frac{1}{\beta}h(v(t)),
\end{align*}
where 
$u\in \mathscr E$ and $v\in  C_{2\pi}({\mathbb R} ; {\mathbb R}^n)$. 
Let $\bm j:\mathscr E\to C_{2\pi}({\mathbb R} ; {\mathbb R})$ denote the natural embedding of $\mathscr E$ into $C_{2\pi}({\mathbb R} ; {\mathbb R}^n)$. Observe that $N$ is a continuous map, the delay operator $\int_0^1v(t - \beta s)ds$ is a compact operator, and $\bm j$ is a compact operator. Then, \eqref{eq:dist:basic2} is equivalent to
\begin{equation}\label{eq:dist:basic3}
Lu=-\frac a \beta \bm j(u)-N(\alpha,\beta, \bm j (u)), \quad u\in \mathscr E,
\end{equation}
where $(\alpha,\beta)\in {\mathbb R}_+^2:=\{(u,v)\in {\mathbb R}^2: v>0\}$. 
Notice that $L+\bm j:\mathscr E\to C_{2\pi}({\mathbb R} ; {\mathbb R}^n)$ is an isomorphism, thus \eqref{eq:dist:basic3} is equivalent to the following operator equation:
\begin{equation}\label{eq:dist:basic4}
\mathscr F(\alpha,\beta,u):=u+(L+\bm j)^{-1} \left(N(\alpha, \beta, \bm j(u))+(\tfrac{a}{\beta} - 1)\bm j(u)\right)=0.
\end{equation}
Since $N(\alpha, \beta, \bm j(u))$ is the composition of a continuous operator and a compact operator, it is compact, and so $N(\alpha, \beta, \bm j(u))+(\tfrac{a}{\beta} - 1)\bm j(u)$ is a sum of a compact operators and is also compact. Therefore $\mathscr F:{\mathbb R}^2_+ \times \mathscr E \to \mathscr E$ is a completely continuous field (i.e. a compact perturbation of identity). We also note that the equation \eqref{eq:dist:basic1.1} has no nonzero periodic solutions for $\alpha = 0$. 
\vs
Put $G := \Gamma_0 \times {\mathbb Z}_2 \times S^1$ and notice that the space $\mathscr E$ is a natural $G$-representation, where the three factors of $G$ act on $u \in \mathscr E$ by permuting indices, changing sign, and shifting phase, respectively. That is, an element of $G$ acts on $u(t)$ as follows:
\[
(\gamma,\pm 1,e^{i\theta})u(t) := \pm \gamma u(t+\theta).
\]
One can also easily see that under the assumptions $(C_1)-(C_3)$ and following the above period normalization, $\mathscr F$ is $G$-equivariant.
\vs
Notice that the map $\mathscr F$ is differentiable with respect to $u$ at all points $(\alpha,\beta,0) \in {\mathbb R}^2_+ \times \mathscr E$, and
\begin{equation}\label{eq:dist:lin1}
D_u\mathscr F(\alpha,\beta,0)u(t):=u(t)+(L+\bm j)^{-1} \left(\frac{\alpha b}{\beta} \int_0^1 u(t - \beta s) ds - \tfrac{1}{\beta}Cu(t) +(\tfrac{a}{\beta} - 1)u(t)\right).
\end{equation}
\vs
Put $\mathscr A (\alpha,\beta) := D_u \mathscr F(\alpha,\beta,0)$. Since $\mathscr A (\alpha,\beta)$ is $S^1$-equivariant, it preserves Fourier modes in the space $\mathscr E$. That means for the subspaces 
\[
\mathscr{E}_k := \{\cos(kt)c+\sin(kt)d: c,
\,d\in V\}, \quad k=0,1,2,3,\dots
\]
one has $\mathscr A(\alpha,\beta) \mathscr{E}_k\subset \mathscr{E}_k$. Put 
\[
\mathscr A_k(\alpha,\beta):=\mathscr A(\alpha,\beta)|_{\mathscr{E}_k}.
\]

We now need to describe the $G$-isotypic decomposition of $\mathscr E$. To do so, we first describe the $S^1$-isotypic decomposition of $\mathscr E$, then further decompose it according to the $\Gamma_0 \times {\mathbb Z}_2$-isotypic decomposition described earlier in \eqref{eq:dist:iso-decomp}. Let $\varphi_k:S^1 \to S^1$ denote the so-called $k$-folding homomorphism, given by
\begin{equation}\label{eq:dist:k-folding}
\varphi_k(e^{i\theta}) := e^{ik\theta}.
\end{equation}
Then each irreducible $S^1$-representation $\cW_k$ can be described as the image of $\mathscr E$ (taken as a natural $S^1$-representation) under $\varphi_k$ for all $k \in \bbN$, and we will refer to this as the $k$-folded $S^1$-action. Thus we can view each subspace $\mathscr E_k$ as an $S^1$-isotypic component of $\mathscr E$ modeled on the irreducible $S^1$-representation $\cW_k$.
Now take the $\Gamma_0 \times {\mathbb Z}_2$-isotypic decomposition of $V$ given in \eqref{eq:dist:iso-decomp}, defined for $j=0,1,\dots,r$, and put
\begin{align*}
\mathscr E_{k,j} &:= \{\cos(kt)c + \sin(kt)d:c,d\in V_j\} &k>0,\\
\mathscr E_{0,j} &:= V_j &k=0.
\end{align*}
Then we have that $\mathscr E_{k,j}$ is the $G$-isotypic component of $\mathscr E$ modeled on the irreducible ${\mathcal V}^-_{k,j}$ representation, where ${\mathcal V}^-_{k,j}$ (as a $G$-representation) is ${\mathcal V}_j \otimes \cW_k$ with the antipodal ${\mathbb Z}_2$-action, and the $k$-folded $S^1$-action. Notice that ${\mathcal V}^-_{0,j} = V_j = {\mathcal V}^-_j$ (with the antipodal ${\mathbb Z}_2$-action), and observe that for $k > 0$ one has the following $G$-equivalence:
\[
\mathscr E_{k,j} \equiv V^{\mathbb C}_{k,j},
\]
where $V^{\mathbb C}_{k,j}$ denotes the complexification of $V_j$ with the $k$-folded $S^1$-action.
Consequently, we have the following $G$-isotypic decomposition of $\mathscr E$:
\vs
\begin{equation}\label{eq:dist:g-iso-decomp}
\mathscr E = \overline{\mathop{\bigoplus_{k = 0}^\infty}\mathop{\bigoplus_{j=0}^r} \mathscr E_{k,j}}.
\end{equation}
\vs
Put
\[
\mathscr A_{k,j}(\alpha,\beta) := \mathscr A(\alpha,\beta)\vert_{\mathscr E_{k,j}}:\mathscr E_{k,j} \to \mathscr E_{k,j}
\]
and notice that 
\begin{align}
\mathscr A_{k,j}(\alpha,\beta)&=\frac{1}{\beta(1+ik)}\triangle _\alpha(i k\beta) \vert_{V^{\mathbb C}_{k,j}}\quad k>0,\label{eq:dist:A_k}\\
\mathscr A_{0,j}(\alpha,\beta)&=\left(\frac \alpha\beta b+\frac{ a_j }{\beta}\right)\vert_{V_0}\label{eq:dist:A_0}.
\end{align}
\vs
\section{System \texorpdfstring{\eqref{eq:dist:basic1.1}}{(4.1.3)} as a two-parameter bifurcation problem}\label{dist:sec4}
The method of the twisted equivariant degree and its application to symmetric Hopf bifurcation problems of various types is described in detail in  \cite{BalanovEtAl2025}. For completeness, let us cite the main results relevant to this problem.
\vs
System \eqref{eq:dist:basic1.1} is equivalent to the following $G$-equivariant two-parameter bifurcation problem
\begin{equation}\label{eq:dist:operator-eq}
\mathscr F(\alpha,\beta,v)=0, \quad (\alpha,\beta,v)\in {\mathbb R}^2_+\times \mathscr E.
\end{equation} 
where $\mathscr F:{\mathbb R}^2_+ \times \mathscr E \to \mathscr E$ is given by \eqref{eq:dist:basic4}, with $\mathscr A(\alpha,\beta) := D_u\mathscr F(\alpha,\beta,0)$. One can easily verify that conditions $(C_1)-(C_4)$ lead to the following properties of the map $\mathscr F$:
\vs
\begin{enumerate}[label=(${\mathscr B}_\arabic*$)]
\item\label{bif0} $\mathscr F$ is a $G$-equivariant compact perturbation of identity;
\item\label{bif1} $\mathscr F(\alpha,\beta,0)=0$ for every $(\alpha,\beta)\in {\mathbb R}^2_+$; 
\item\label{bif2} for  all $(\alpha,\beta)\in {\mathbb R}^2_+$,  the derivative $\mathscr A(\alpha,
\beta):=D_v\mathscr F(\alpha,\beta,0):\mathscr E\to \mathscr E$ exists,
depends continuously on  $(\alpha,\beta)$ and  
\begin{equation}\label{eq:dist:A-twisted}
\lim_{(\alpha',\beta',v)\to (\alpha',\beta',0)} \frac{\|\mathscr F(\alpha',\beta',v)-\mathscr A(\alpha',\beta')v\|}{\|v\|}=0.
\end{equation}
\end{enumerate}
We recall the definitions of bifurcation points, branching points, and branches from Section \ref{sec:prelim:krasnoselskii}, along with the branch of trivial solutions $M$ and nontrivial solutions $\mathscr S$. Then we have the following necessary condition for $(\alpha,\beta,0)\in M$ to be a bifurcation point.
% \vs
% Elements of the set
% \begin{equation}\label{eq:dist:trivialM-twisted}
% M:=\{(\alpha,\beta,0): (\alpha,\beta,0)\in {\mathbb R}^2_+\times \mathscr E\}
% \end{equation}
% are called {\it trivial solutions} to \eqref{eq:dist:operator-eq} and   
% $\mathscr S$  denotes the set of all  {\it nontrivial solutions} to  \eqref{eq:dist:operator-eq}, i.e.
% \begin{equation}\label{eq:dist:nontr-twisted}
% \mathscr S:=\{(\alpha,\beta,v)\in {\mathbb R}^2_+\times \mathscr E: \mathscr F(\alpha,\beta,v)=0 \;\text{ and }\; v\not=0\}.
% \end{equation}
% \vs
% \begin{definition}
% \rm $(i)$: A point $(\alpha,\beta,0) \in M$ is called a {\it bifurcation point} for \eqref{eq:dist:basic4} if in every neighborhood $\Omega \subset {\mathbb R}^2_+ \times \mathscr E$ of $(\alpha,\beta,0)$, there is a nontrivial solution $(\alpha,\beta,u)\in\Omega$.\\
% $(ii)$: A nonempty set $\mathscr C \subset \mathscr S$ is called a branch of nontrivial solutions to \eqref{eq:dist:basic4} if for some connected component $\mathscr D$ of $\overline{\mathscr S}$, one has $\mathscr C = \mathscr S \cap \mathscr D$. In addition, if $(\alpha,\beta,0) \in \overline{\mathscr C}$, we call $(\alpha,\beta,0)$ a {\it branching point}. Notice that any branching point is also a bifurcation point.
% \end{definition}
% \vs
% The following proposition provides a necessary condition for $(\alpha,\beta,0) \in M$ to be a bifurcation point.
\begin{proposition}
Under the above assumptions, if $(\alpha_0,\beta_0,0) \in M$ is a bifurcation point, then
\begin{equation}\label{eq:dist:A-crit}
\mathscr A(\alpha_0,\beta_0):\mathscr E \to \mathscr E
\end{equation}
is not an isomorphism.
\end{proposition}
\begin{proof}
Suppose that $(\alpha_0,\beta_0,0) \in M$ is a bifurcation point but $\mathscr A(\alpha_0,\beta_0)$ is an isomorphism. Then there exists $c > 0$ such that 
\[
\|\mathscr A(\alpha_0,\beta_0)v\| \geq c\|v\|, \quad \forall v \in \mathscr E.
\]
By the continuity of $\mathscr A(\alpha,\beta)$, there also exists a neighborhood $U \subset \mathbb R_+^2$ and constant $c' >0$ such that
\[
||\mathscr A(\alpha,\beta)v||, \geq c'||v|| \quad \forall v \in \mathscr E.
\]
By assumption \ref{bif2}, for $\varepsilon = \frac{c'}{2}$, there exists $\delta >0$ such that $\|(\alpha,\beta) - (\alpha_0,\beta_0)\|<\delta$ and $\|v\| < \delta$, then
\[
\|\mathscr F(\alpha,\beta,v) - \mathscr A(\alpha,\beta)v\| \leq \frac{c'}{2}\|v\|.
\]
We then put $U':=U \cap \{(\alpha,\beta)\in \mathbb R_+^2 : \|(\alpha,\beta)-(\alpha_0,\beta_0)\|<\delta\}$, take $V := \{v\in \mathscr E\setminus \{0\} : \|v\|< \delta\}$, and put $\Omega := U' \times V$. Then for all $(\alpha,\beta,v)\in \Omega$, by the triangle inequality we have
\[
\|\mathscr F(\alpha,\beta,v)\| \geq \|\mathscr A(\alpha,\beta)v\| - \|\mathscr F(\alpha,\beta,v) - \mathscr A(\alpha,\beta)v\| \geq c'\|v\| - \frac{c'}{2}\|v\| \geq \frac{c'}{2}\|v\| > 0.
\]
Therefore there exists a neighborhood of $(\alpha_0,\beta_0,0) \in M$ in which \eqref{eq:dist:operator-eq} has no nontrivial solutions, and thus we have a contradiction.
\end{proof}
\vs
This leads us to the following definition for the critical set of \eqref{eq:dist:operator-eq}:
\[
\Lambda := \{(\alpha,\beta,0) \in M:\mathscr A(\alpha,\beta):\mathscr E \to \mathscr E \text{ is not an isomorphism}\}.
\]
\vs
By \eqref{eq:dist:A_k}, $(\alpha,\beta,0) \in \Lambda$ if and only if there exists some $k \in \bbN$ such that  ${\det}_{{\mathbb C}}\triangle_\alpha(ik \beta) = 0$, which implies that $(\alpha,k\beta,0)$ is a critical point (in the sense of an isolated center, as defined in Section \ref{sec:prelim:two-parameter-bifurcation}). Consequently, for the system \eqref{eq:dist:basic1.1}, the critical set $\Lambda$ is given by \eqref{eq:dist:crit-set}.
\vs
The set $\Lambda$ provides the candidate $(\alpha,\beta,0)$ for the possible bifurcation point of \eqref{eq:dist:operator-eq}. As we are not interested in studying the bifurcation of steady-state solutions, we exclude from $\Lambda$ the points $(\alpha_0,\beta_0,0)$ such that $\det\mathscr A_0(\alpha_0,\beta_0)=0$.
\vs
Consequently, using \eqref{eq:dist:A_k} and \eqref{eq:dist:A_0}, we define 
\begin{equation}\label{eq:dist:crit-set-bold}
\bm{\Lambda} := \left\{(\alpha_0,\beta_0,0) \in \Lambda: \exists_{k\in\bn}\;\; {\det}_{\mathbb C} \triangle_{\alpha_0}(ik\beta_0)=0 \text{ and } \alpha_0b+a_j \not= 0, j=0,1,\dots,r\right\}.
\end{equation}
Since ${\det}_{\mathbb C} \triangle_{\alpha_0}(ik\beta_0)=0$ implies that $(\alpha_0,k\beta_0,0)\in \Lambda$, it follows that $\bm \Lambda \subset \Lambda$.
\vs
Let us point out that by Proposition \ref{prop:dist:crit1}, any $(\alpha_0,\beta_0,0) \in \bm \Lambda$ is an isolated point. Then, exactly following the construction in Section \ref{sec:prelim:local-bif-twisted}, for each $(\alpha_0,\beta_0,0)\in \bm \Lambda$, we can construct an isolated $G$-invariant neighborhood $\Omega(\alpha_0,\beta_0,0)$, construct a continuous $G$-equivariant auxiliary function $\eta$ and define the complemented map $\mathscr F_\eta$, and obtain the local bifurcation invariant $\omega_G(\alpha_0,\beta_0) = \gdeg(\mathscr F_\eta,\Omega(\alpha_0,\beta_0,0))$.

\vs
We also have the following computational formula for $\omega_G(\lambda_0)$, formulated earlier in Section \ref{sec:prelim:local-bif-twisted}, and proven in a general setting in \cite{BalanovEtAl2025}.
\begin{theorem}\label{thm:dist:global-bif}
Suppose that $f:{\mathbb R} \to {\mathbb R}$ is a continuous odd function such that $b:=f'(0)<0$ and assume that $\mathbf f:V \to V$ is given by \eqref{eq:dist:basic1.1}, and $h:V\to V$ is continuous and $\Gamma_0$-equivariant and satisfies \ref{dist:c1}---\ref{dist:c4}. If $(\alpha_0,\beta_0,0) \in \bm \Lambda$ is an isolated critical point, then we have:
% \[
% \omega_G(\lambda_0) = \Gamma\text{\rm-deg}(\mathscr A_0, \lambda_0, B(V))\cdot \sum_{k,j}\mathfrak t_k^j(\lambda_0)\text{\rm deg}_{\mathcal V_{k,j}^-}
% \]
\[
\omega_G(\lambda_0) =\sum_{k,j}\mathfrak t_{k,j}(\alpha_0,\beta_0)\deg_{\mathcal V_{k,j}^-}
\]
where, if $(\alpha_0,k \beta_0,0)\in\Lambda$ (see \eqref{eq:dist:crit-set}), then $\mathfrak t_{k,j}(\lambda_0) = \pm m_j$ (see \eqref{eq:dist:iso-mult}), and the sign of $m_j$ is described by \eqref{eq:dist:du-dalpha} as given by Lemma \ref{lem:prelim:sign-crossing-num}.
\end{theorem}
This is the way we relate the transversality condition proved earlier to the crossing numbers, and therefore to the nontriviality of certain orbit types in the local bifurcation invariant, which are obtained from the twisted basic degrees. There is only one final point to be addressed to complete the proof of the main results, and this is to ensure that global branches of non-constant periodic solutions \emph{remain} non-constant and periodic globally, i.e. they do not bifurcate to a steady-state solution or otherwise vanish. To do this, we will apply the Krasnosel'skii and Rabinowitz theorems in a reduced space whose solutions are also solutions of \eqref{eq:dist:basic4}. 
\subsection{Fixed-point reduction}\label{sec:dist:fpr}
Now we will describe the aforementioned fixed-point reduction in detail, and show how it not only provides the necessary conditions to complete the proof of the main result, but also how it provides information about the symmetries of unbounded solutions described in the symmetric corollary to the main result. We will use this fixed-point reduction method to find solutions whose isotypic crossing numbers can be computed explicitly and shown to have the same sign. In conjunction with the above explicit formula for computation of the local bifurcation invariant, this establishes the main result.

For a given $\kappa \in \bbN$, put $\bm K := \{\bm e\}\times {\mathbb Z}_{2\kappa}^d \leq G$, where $\bm e \in \Gamma_0$ is the neutral element, $\xi := e^{\tfrac{i\pi}{\kappa}}$, and 
\[
{\mathbb Z}_{2\kappa}^d := \{(1,1),(-1,\xi),(1,\xi^2),\dots,(-1,\xi^{2\kappa-1})\}.
\]
Then we reduce the problem \eqref{eq:dist:operator-eq} to the fixed-point space $\mathscr E^{\bm K}$, i.e. instead of \eqref{eq:dist:operator-eq}, we study the equation
\begin{equation}\label{eq:dist:fp-op-eq}
    \mathscr F^{\bm K}(\alpha,\beta,u)=0,\quad (\alpha,\beta,u)\in {\mathbb R}^2_+ \times \mathscr E^{\bm K}.
\end{equation}
The isotypic components $\mathscr E_{k,j}^{\bm K}$ of the fixed-point space $\mathscr E^{\bm K}$ can easily be described as
\[
\mathscr {E}_{k,j}^{\bm K} := \begin{cases}
\mathscr E_{k,j} \quad &k \text{ is an odd multiple of $\kappa$}\\
0 \quad &\text{otherwise}
\end{cases}
\]
Therefore
\[
\mathscr E^{\bm K} = \overline{\bigoplus_{l=1}^{\infty}\bigoplus_{j=0}^r\mathscr E_{(2l-1)\kappa,j}}.
\]
\vs
Since ${\bm K}$ is normal in $G$, and $G_0:=G/{\bm K} = \Gamma_0 \times S^1$, the problem \eqref{eq:dist:fp-op-eq} is a two-parameter $\Gamma_0 \times S^1$-equivariant bifurcation problem.

One can easily notice that any solution to \eqref{eq:dist:fp-op-eq} is also a solution to \eqref{eq:dist:operator-eq} (and consequently to \eqref{eq:dist:basic1.1}), thus by establishing a global bifurcation result for \eqref{eq:dist:fp-op-eq} we also obtain a global bifurcation result for \eqref{eq:dist:basic1.1}.

By taking $k$ sufficiently large, one can guarantee that none of the issues listed above affect the problem \eqref{eq:dist:fp-op-eq}. Moreover, the set of critical points for $\mathscr F^{\bm K}$ can be described as follows:
\[
{\bm \Lambda}^{\bm K} := \{(\alpha_0,\beta_0,0):\exists_{l \in \bbN} \;{\det}_{\mathbb C} \triangle_{\alpha_0}(i(2l-1)\kappa\beta_0)=0\} 
\]
Clearly, the set ${\bm \Lambda}^{\bm K}$ is infinite, and every point $(\alpha_0,\beta_0,0)\in {\bm \Lambda}^{\bm K}$ is isolated. Furthermore, we have that all non-zero crossing numbers satisfy $\mathfrak t_j(\alpha_0,\beta_0) = -m_j$.
\vs
Applying the formula for the local bifurcation invariant provided in Theorem (\ref{thm:dist:global-bif}), we then have
\[
\sum_{i=1}^N \omega_{G_0}(\alpha_i,\beta_i) = \sum_{i=1}^N\sum_{j=0}^r\sum_{l=1}^{\infty}\mathfrak t_{(2l-1)\kappa,j}(\alpha_i,\beta_i) \text{deg}_{{\mathcal V}_{(2l-1)\kappa,j}} 
\]
(here $\text{deg}_{{\mathcal V}_{(2l-1)\kappa,j}}$ denotes the basic $G_0$-degree corresponding to the irreducible $G_0$-representation ${\mathcal V}_{(2l-1)\kappa,j}$ with trivial ${\mathbb Z}_2$-action), and for every $i=1,\dots,N$
\[
\mathfrak t^j_{(2l-1)\kappa}(\alpha_i,\beta_i)=\begin{cases}
-m_j \quad&\text{for some }l,j\\ 
0 \quad&\text{otherwise}
\end{cases}.
\]
We also point out that for fixed $j$, each orbit type present in the basic degree $\text{deg}_{{\mathcal V}^-_{k,j}}$ uniquely corresponds to an orbit type in $\text{deg}_{{\mathcal V}^-_{1,j}}$, related by the $k$-folding homomorphism \eqref{eq:dist:k-folding}. In particular, if we have
\[
\text{deg}_{{\mathcal V}^-_{k,j}} = c_1(H_1) + \dots + c_s(H_s),
\]
then
\[
\text{deg}_{{\mathcal V}^-_{1,j}} = c_1(\tilde{\varphi_k}(H_1)) + \dots + c_s(\tilde{\varphi_k}(H_s)),
\]
where $\tilde{\varphi_k}:G \to G$ is the natural extension of $\varphi_k: S^1 \to S^1$ such that for $(\gamma,\kappa,\xi) \in G = \Gamma_0 \times {\mathbb Z}_2 \times S^1$, $\tilde{\varphi_k}((\gamma,\kappa,\xi))=(\gamma,\kappa,\varphi_k(\xi))$. An important consequence of this relationship is that orbit types are disjoint across different ${\mathcal V}^-_{k,j}$ (at least for fixed $j$), differing in their folding but not in their coefficient or sign. So it is indeed sufficient to show that $\mathfrak t_k^j$ is of constant sign to establish the above result about the local bifurcation invariant.

Therefore, it is not possible for any branch bifurcating from $(\alpha_0,\beta_0,0)$ to be bounded, nor to collapse to a steady-state solution. This completes the proof of Theorem \ref{thm:dist:main1} and Corollary \ref{thm:dist:main2}, and therefore also completes the proofs of Theorem \ref{thm:dist:mas-local} and Theorem \ref{thm:dist:mas-sym}.

\vs
We will now show a plausible and concrete application for multi-agent systems of the type of \eqref{eq:dist:basic1.1}, demonstrate how Theorem \ref{thm:dist:mas-local} and Theorem \ref{thm:dist:mas-sym} can be applied to this problem to show the emergence of periodic multiconsensus states following consensus-breaking bifurcation of the $\bm x\equiv 0$ consensus solution, and validate and enhance these theoretical results through numerical simulations showing multiple stable periodic multiconsensus solutions.
\section{Application: tunable spontaneous oscillations in symmetric formations of UAVs control}\label{dist:sec5}
Consider an octahedral formation of UAVs, where each UAV is situated at a vertex of the octahedron. Suppose that the center of the octahedron is an invisible set point, or a ``mothership'' which is escorted by this surrounding formation. There are many reasons why one might desire periodic movement or oscillation in a drone formation, either when the formation is in motion or at rest. For example, if UAVs are scanning or sensing the environment, symmetric periodic motions could give them better sensor coverage. For military applications, such maneuvers might aid in evasion. However, it is also very important that the formation retains cohesion throughout such periodic motion. In particular, we desire that
\begin{itemize}
    \item There should be zero mean displacement from the steady-state formation, 
    \item Oscillations should be bounded in amplitude,
    \item The frequency and amplitude of oscillations should be tunable in some way,
    \item The formation can return to a stable steady-state formation, and these regimes can be controlled with a single parameter.
\end{itemize}
In this section, we will show a protocol which achieves this, under certain assumptions. 
For simplicity and clarity, we assume that the low-level attitude and velocity controllers in each UAV are sufficiently fast that the translational dynamics can be approximated by single-integrator dynamics across three orthogonal $x,y,z$-axes. We will call the coordinate system with $(0,0,0)$ as the center of the octahedron the ``formation-local'' coordinates. Then we define the set point for UAV $i$ as $p_i$ in $\mathbb R^3$, where
\begin{align*}
    p_1&:=(\delta,0,0)\\
    p_2&:=(0,0,\delta)\\
    p_3&:=(-\delta,0,0)\\
    p_4&:=(0,0,-\delta)\\
    p_5&:=(0,\delta,0)\\
    p_6&:=(0,-\delta,0),
\end{align*}
and $\delta>0$ stands for the distance from the formation origin to the vertices of the octahedron, and so determines the ``size'' of the formation. For UAV $i$, we define the ``UAV-local'' coordinates centered around its respective set point $p_i$, such that $x_i(t),y_i(t),z_i(t)$ define the displacement of UAV $i$ from the set point $p_i$ along each of the three UAV-local coordinate axes. We also set the UAV-local axes parallel to the formation-local axes. In other words, the ``forward'' direction for each drone is parallel to the ``forward'' direction of the formation as whole, and likewise with other axes. If $(\tilde x_i(t),\tilde y_i(t),\tilde z_i(t)) \in \mathbb R^3$ represents the position of UAV $i$ in formation-local coordinates, then $r_i(t):= (x_i(t),y_i(t),z_i(t)) := (\widetilde x_i(t),\widetilde y_i(t),\widetilde z_i(t))-p_i$ represents its position in UAV-local coordinates. Then we define 
\begin{align*}
\bm x(t) &:= (x_1(t),\dots,x_6(t))^T \in \mathbb R^6,\\
\bm y(t) &:= (y_1(t),\dots,y_6(t))^T \in \mathbb R^6,\\ 
\bm z(t) &:= (z_1(t),\dots,z_6(t))^T \in \mathbb R^6,\\
r_i(t) &:=(x_i(t),y_i(t),z_i(t))^T \in \mathbb R^3\\
\bm r(t)&:=(\bm x(t),\bm y(t), \bm z(t))^T \in \mathbb R^{18}.
\end{align*}
Then $\bm x \equiv 0, \bm y \equiv 0, \bm z \equiv 0$ is a trivial consensus state, and corresponds to all UAVs maintaining the octahedral formation. For each UAV $i$ on each local coordinate axis, we define the UAV control protocol $u_i(r_i):=(u_i^x(\bm x)^T, u_i^y(\bm y)^T,u_i^z(\bm z)^T)^T$, where
\begin{align*}
    u_i^x(x_i) := -a_xx_i - \alpha f_x\left(\int_0^1x_i(t-s)ds\right)-h_i^x(\bm r)\\
    u_i^y(y_i) := -a_yy_i - \alpha f_y\left(\int_0^1y_i(t-s)ds\right)-h_i^y(\bm r)\\
    u_i^z(z_i) := -a_zz_i - \alpha f_z\left(\int_0^1z_i(t-s)ds\right)-h_i^z(\bm r).
\end{align*}
Then we can put
\[
 u_i(\bm r) := -ar_i - \alpha f\left(\int_0^1r(t-s)ds\right)-h_i(\bm r),
\]
where $f:\mathbb R\to \mathbb R^3$ is defined $f(r_i) = (f_x(x_i)^T,f_y(y_i)^T,f_z(z_i)^T)^T$, $h_i:\mathbb R^6 \to \mathbb R^3$ is defined $h_i(\bm r) :=(h_i^x(x_i)^T,h_i^y(y_i)^T,h_i^z(z_i)^T)^T$, and $\bm a = (a_x,a_y,a_z)^T$. For simplicity, we will assume that the axes are decoupled and have identical dynamics, i.e.  i.e. $a := a_x = a_y = a_z$, $f_x = f_y = f_z$, and the axial coupling functions $h^x, h^y,h^z : \mathbb R^6 \to \mathbb R^6$ are identical and depend only on $\bm x, \bm y,$ and $\bm z$, respectively. Then we can write $h = h^x = h^y = h^z$, and the axis-independent protocol
\[
u(\bm x) = -a\bm x - \alpha \bm f\left(\int_0^1\bm x(t-s)ds\right)-h(\bm x)
\]
which gives
\begin{align*}
u(\bm x) &:= -a\bm x - \alpha \bm f\left(\int_0^1\bm x(t-s)ds\right)-h(\bm x),\\
u(\bm y) &:= -a\bm y - \alpha \bm f\left(\int_0^1\bm y(t-s)ds\right)-h(\bm x),\\
u(\bm z) &:= -a\bm z - \alpha \bm f\left(\int_0^1\bm z(t-s)ds\right)-h(\bm x).
\end{align*}
This yields the following three decoupled closed-loop dynamics equations
\begin{equation}\label{eq:dist:UAV-dynamics}
\begin{aligned}
\dot {\bm x}(t) &:= -a\bm x(t) - \alpha \bm f\left(\int_0^1\bm x(t-s)ds\right)-h(\bm x(t)),\\
\dot {\bm y}(t) &:= -a\bm y(t) - \alpha \bm f\left(\int_0^1\bm y(t-s)ds\right)-h(\bm y(t)),\\
\dot {\bm z}(t) &:= -a\bm z(t) - \alpha \bm f\left(\int_0^1\bm z(t-s)ds\right)-h(\bm z(t)).
\end{aligned}
\end{equation}
Note that if one of these equations satisfies \ref{dist:c1}---\ref{dist:c4}, then they all satisfy these assumptions, so without loss of generality, we can study $\dot {\bm x} = u(\bm x)$ and apply this analysis equally to each axis. Note that the symmetry group of rigid motions of an octahedron is given by $\mathbb O\cong S_4 \leq S_6$. We view $\mathbb O$ as a subgroup of $S_6$ which acts on vertices by permuting their indices. Then the entire closed loop dynamics equation is $S_4$ equivariant on each axis.
\begin{remark}
One could also consider protocols of the same type controlling the pitch, yaw, and roll of each UAV. In this case, the assumptions that there is a sufficiently fast closed-loop controller which allows us to abstractly manipulate the UAV's cartesian coordinates is unnecessary, but the coupling arrangement and symmetries become much more intricate. Likewise, we could also allow for cross-coupling between the axes in this example. This could result in richer symmetric formation oscillations, but again complicates the isotypic decomposition and group theory. 
\end{remark}
Each protocol for a single axis of a single UAV is similar to a coupled PI controller, but with integration taken over a finite time horizon. This addresses potential integral windup issues which occur in PI and PID controllers, and also allows for stable periodic motion. In this case, we can think of the $-ax_i$ term as pushing the displacement towards 0, the memory term $f(\int_0^1x_i(t-s)ds)$ as creating stable periodic oscillations as $x_i(t)$ crosses 0, and the coupling term as controlling the phase relationship between the oscillations of different UAVs based on their relative position. In a practical application, the coupling term could be based on each UAV sensing its neighbors and their relative displacement through either active or passive sensors or communication, or could be supplied to them externally by a coordinator or controller UAV which tracks the entire formation. 

\label{fig:dist:formation}\rm
\begin{figure}[h]
	\begin{center}
		\scalebox{.6}
		{

			\begin{tikzpicture}
				\node at (2.5,0){$x$}; ;
				\filldraw [red,solid](-1,-1)--(-4,0)--(0,-5);
				\filldraw [yellow,solid](4,0)--(0,-5)--(-1,-1);
				\filldraw [green,solid](-4,0)--(-1,-1)--(0,5);
				\filldraw [lightblue,solid](-1,-1)--(4,0)--(0,5);
				%red diagonal%

				%green diagonal

				%lightblue diagonal

				%yellow diagonal
                \draw (-1.5,0)--(1.5,0);
                \draw (-0.5,-0.5)--(0.3,0.3);
                \draw (0,-2)--(0,2);
                \node at (2,0){\huge $\tilde x_+$};
                \node at (0.6,0.6){\huge$\tilde y_+$};
                \node at (0.15,2.4){\huge$\tilde z_+$};

                \draw [blue] (3,0)--(5,0);
                \draw [blue] (3.5,-0.5)--(4.5,0.5);
                \draw [blue] (4,-0.5)--(4,0.5);
                \node at (2,0){\huge $\tilde x_+$};
                \node at (0.6,0.6){\huge$\tilde y_+$};
                \node at (0.15,2.4){\huge$\tilde z_+$};
                
				\draw (0,5)--(-4,0)--(0,-5)--(4,0)--(0,5);
				\draw (-4,0)--(-1,-1)--(4,0);
				\draw [dashed](4,0)--(1,1)--(-4,0);
				\draw (0,5)--(-1.,-1)--(0,-5);
				\draw [dashed](0,5)--(1.,1)--(0,-5);
				\fill (-4,0) circle (3pt) (-1,-1) circle (3pt) (0,5) circle (3pt) (4,0) circle (3pt) (0,-5) circle (3pt) (1,1)circle (3pt);                   
				%  \rput(-1.9,-2){\Large\bf $1$}
				% \rput(1,2){\Large\bf $3$}
				%\rput(-2,2){\Large\bf $2$}
				\node at (-.5,-1.3){\huge $p_4$};
				\node at (0,-5.5){\huge $p_6$};
				\node at (4.8,0.3){\huge $p_1$};
				\node at (0,5.4){\huge $p_5$};
				\node at (-4.5,0){\huge $p_3$};
				\node at (1.5,1.4){\huge $p_2$};
				% \rput(1,-2){\Large\bf $4$}

			\end{tikzpicture}
		}
		\caption{$\widetilde x_+,\widetilde y_+,$ and $\widetilde z_+$ represent the positive axis directions of the formation-local coordinate axes $\widetilde x,\widetilde y,$ and $\widetilde z$. The parallel UAV-local coordinate axes for UAV 1 at $p_1$ are shown in blue.}
		\label{fig:dist:oct}
	\end{center}
\end{figure}
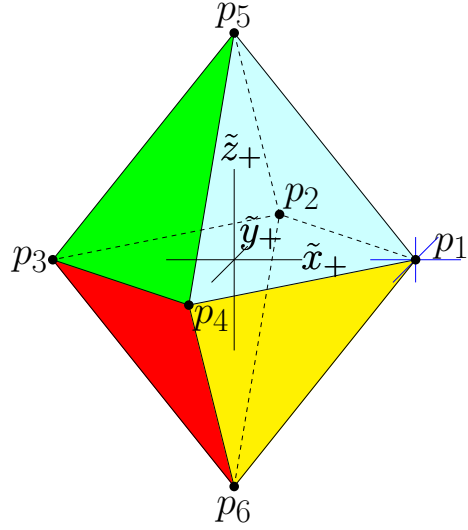

We put $C:= Dh(0)$ and write
\[
C:=Dh(0)=  \left[
    \begin{array}
    [c]{cccccc}
   0 &  d & 0 & d & d &  d \\
     d  &     0 &  d  &0  & d  &  d  \\
    0&  d  &    0 &   d &   d  &   d  \\
     d  &0 &  d &     0 &   d  &  d  \\
     d  &  d  &  d  &   d &   0&0\\
  d  &  d  & d  &  d  & 0 &    0 
    \end{array}
    \right]
\]
Note that this is the undirected adjacency matrix of the octahedron shown in Figure \ref{fig:dist:formation}. This means that each UAV will only be influenced by the positions of its ``neighbors'' connected by edges of the octahedron. $C$ admits the eigenspaces
\begin{align*}
E(\mu_0)&:= \text{span}\{ (1,1,1,1,1,1)^T\}, \quad \mu_0:=4d,\\
E(\mu_2)&:= \text{span}\{ (-1,0,-1,0,1,1)^T,(-1,1,-1,1,0,0)^T\}, \quad \mu_2:=-2d,\\
E(\mu_3)&:= \text{span}\{ (-1,0,1,0,0,0)^T,(0,0,0,0,-1,1)^T,(0,-1,0,1,0,0)^T\},
\quad \mu_3=0.
\end{align*}
\vs
Although the octahedron has 6 vertices, its group of rigid motions which transform vertices while preserving orientation are isomorphic to the group $S_4$. Therefore when we refer to the {\it octahedral symmetry group} $S_4$, we understand it to mean an embedding of $S_4$ in $S_6$, where $S_6$ acts on vectors in $V$ by permuting their components. This embedding can be chosen arbitrarily without affecting the isotypic decomposition, so we choose the following correspondences between representatives of conjugacy classes in $S_4$ and $S_6$.
\begin{gather} 
%\label{eq:dist:octahedral_permutation}
(1,2)\;\;\; \leftrightarrow \;\; (1,4)(2,3)(5,6),\quad (1,2,3)\;\;\leftrightarrow \;\; (1,4,6)(3,2,5), \nonumber \\ (1,2)(3,4)\;\; \leftrightarrow \;\; (1,4)(2,3),\quad (1,2,3,4)\;\;\leftrightarrow \;\; (1,2,3,4), \nonumber
\end{gather}
then the vertices of the octahedron (see Figure \ref{fig:dist:formation}) are permuted in the following way
\begin{align} \label{def:S4-action}
\forall_{\sigma\in S_4\le S_6}\;\;\sigma(u_{1},u_{2},u_{3},u_{4},u_{5},u_6)^T=
(u_{\sigma(1)},u_{\sigma(2)},u_{\sigma(3)},u_{\sigma(4)},u_{\sigma(5)},u_{\sigma(6)})^T,\quad  u\in V.
\end{align}
In this way, $V$ is an orthogonal $S_4$-representation with respect to the action \eqref{def:S4-action}. The character table of irreducible $S_4$-representations of $V$ can be easily computed.
\begin{table}[h]
\centering
\begin{tabular}{|c|ccccc|}
\hline
con. classes &$(1)$ & $(1,2)$&$(1,2)(3,4)$ & $(1,2,3)$&$(1,2,3,4)$\\ \hline $\chi_0$ &$1$ & $1$ & $1$ & $1$ &$1$\\ $\chi_1$ &$1$ & $-1$ & $1$ & $1$ &$-1$\\ $\chi_2$ &$2$ & $0$ & $2$ & $-1$ &$0$\\ $\chi_3$ &$3$ & $-1$ & $-1$ & $0$ &$1$\\ $\chi_4$ &$3$ & $1$ & $-1$ & $0$ &$-1$\\ \hline 
$\chi_V$ &$6$ & $0$ & $2$ & $0$ &$2$\\ \hline
\end{tabular}
\caption{Character Table for $S_4$.}
\end{table}
\vs
Let $\Gamma = S_4 \times {\mathbb Z}_2$. Inspection of the above character table immediately yields the $S_4 \times {\mathbb Z}_2$-isotypic decomposition of $V$
\[
V = V_0 \oplus V_2 \oplus V_3,
\]
Observe that $E(\mu_0)=V_0$, $E(\mu_2)=V_2$, $E(\mu_3)=V_3$, and that the $\Gamma$-isotypic multiplicities of $V_i$ for all $i=0,2,3$ is 1. Our system can be reformulated as a $G$-symmetric two-parameter bifurcation problem \eqref{eq:dist:basic4} with $G=S^1 \times S_4 \times {\mathbb Z}_2$, and the theoretical results obtained in section 4 can be applied. In our case
\begin{align*}
    a_0 &= a+4d,\\
    a_2 &= a-2d,\\
    a_3 &= a,
\end{align*}
with the critical set $\bm \Lambda$ given by \eqref{eq:dist:crit-set-bold}. Consequently, our main global bifurcation result in Theorem \ref{thm:dist:main1} holds for this system.
\vs
The most interesting aspect of this example is that we are actually able to analyze the spatio-temporal symmetric properties of the bifurcating solutions. To simplify the procedure, we use the GAP software system which allows us to conveniently compute all basic twisted degrees for the group $G = S^1 \times S_4 \times {\mathbb Z}_2$. These computations yield the following basic degrees for our system
\begin{align*}
    \text{deg}_{{\mathcal V}_{k,0}^-} = &-( {\mathbb Z}_{2k} {}^{{\mathbb Z}_k}\times {}^{S_4} S_4^p)\\
    \text{deg}_{{\mathcal V}_{k,2}^-} = &( {\mathbb Z}_{6k} {}^{{\mathbb Z}_k}\times {}^{V_4} A_4^p) + ({\mathbb Z}_{2k} {}^{{\mathbb Z}_k}\times {}^{D_4^{\hat{d}}} D_4^p) + ( {\mathbb Z}_{2k} {}^{{\mathbb Z}_k}\times {}^{D_4} D_4^p) - ( {\mathbb Z}_{2k} {}^{{\mathbb Z}_k}\times {}^{V_4} V_4^p)\\
    \text{deg}_{{\mathcal V}_{k,3}^-} =&({\mathbb Z}_{2k} {}^{{\mathbb Z}_k}\times {}^{D_4^z} D_4^p) + ({\mathbb Z}_{2k} {}^{{\mathbb Z}_k}\times {}^{D_3^z} D_3^p) + ({\mathbb Z}_{2k} {}^{{\mathbb Z}_k}\times {}^{D_2^d} D_4^z)+({\mathbb Z}_{4k} {}^{{\mathbb Z}_k}\times {}^{{\mathbb Z}_2^-} {\mathbb Z}_4^p)+\\ &({\mathbb Z}_{6k} {}^{{\mathbb Z}_k}\times {\mathbb Z}_3^p) - ({\mathbb Z}_{2k} {}^{{\mathbb Z}_k}\times {}^{{\mathbb Z}_2^-} D_2^p) - ({\mathbb Z}_{2k} {}^{{\mathbb Z}_k}\times {}^{D_1^z} D_1^p)
\end{align*}

Furthermore, one can identify 7 types of symmetries corresponding to the following orbit types of maximal kind $(H)$ in $A_1^t(G;\mathscr{E}_{k,j})$:

\begin{align*}
    \mathfrak M_{k,0} = \{&( {\mathbb Z}_{2k} {}^{{\mathbb Z}_k}\times {}^{S_4} S_4^p)\}\\
    \mathfrak M_{k,2} = \{&( {\mathbb Z}_{6k} {}^{{\mathbb Z}_k}\times {}^{V_4} A_4^p), ({\mathbb Z}_{2k} {}^{{\mathbb Z}_k}\times {}^{D_4^{\hat{d}}} D_4^p), ( {\mathbb Z}_{2k} {}^{{\mathbb Z}_k}\times {}^{D_4} D_4^p)\}\\
    \mathfrak M_{k,3} = \{&({\mathbb Z}_{2k} {}^{{\mathbb Z}_k}\times {}^{D_4^z} D_4^p) , ({\mathbb Z}_{2k} {}^{{\mathbb Z}_k}\times {}^{D_3^z} D_3^p) , ({\mathbb Z}_{2k} {}^{{\mathbb Z}_k}\times {}^{D_2^d} D_4^z),\\&({\mathbb Z}_{4k} {}^{{\mathbb Z}_k}\times {}^{{\mathbb Z}_2^-} {\mathbb Z}_4^p), ({\mathbb Z}_{6k} {}^{{\mathbb Z}_k}\times {\mathbb Z}_3^p)\}
\end{align*}

\begin{figure}[H]
    \includegraphics[height=9cm]{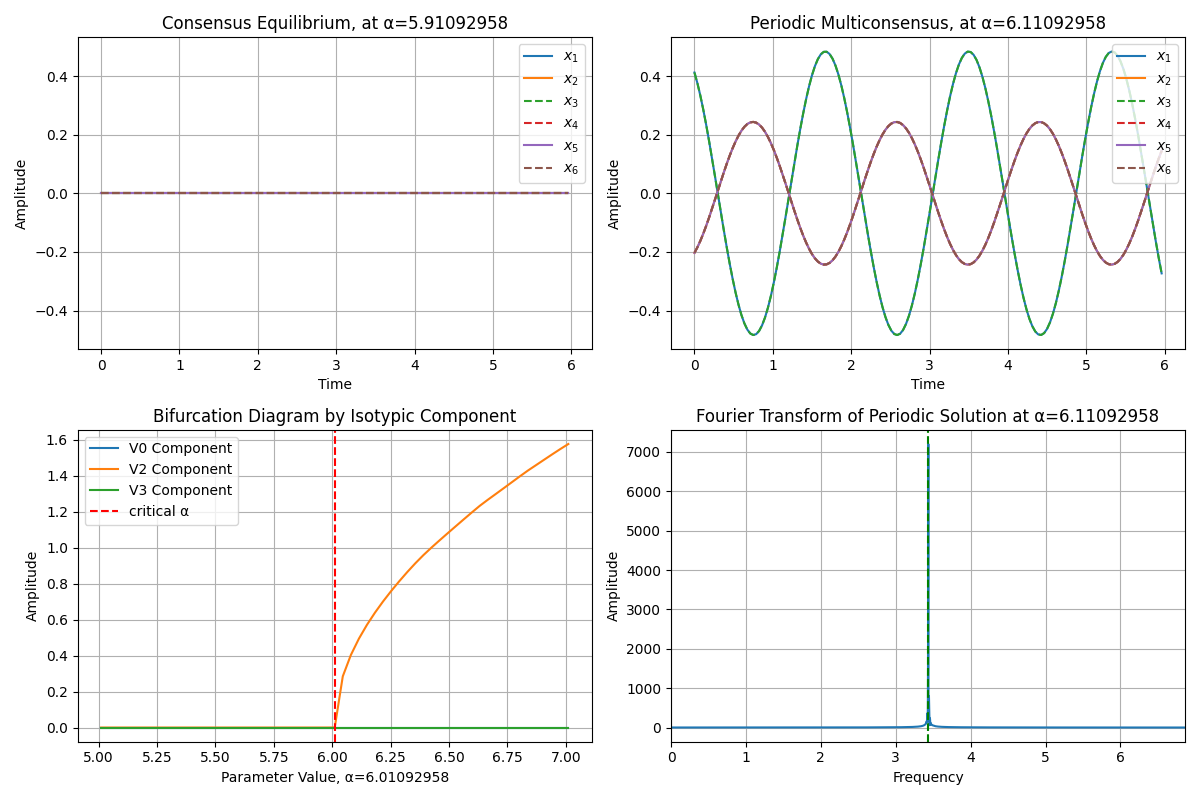} 
	\caption{Representative of the standing wave solution type, initialized with a small perturbation in the direction $(-1,-1,-1,-1,2,2)^T$.}\label{fig:dist:standing}
\end{figure}

\begin{figure}[H]
    \includegraphics[height=9cm]{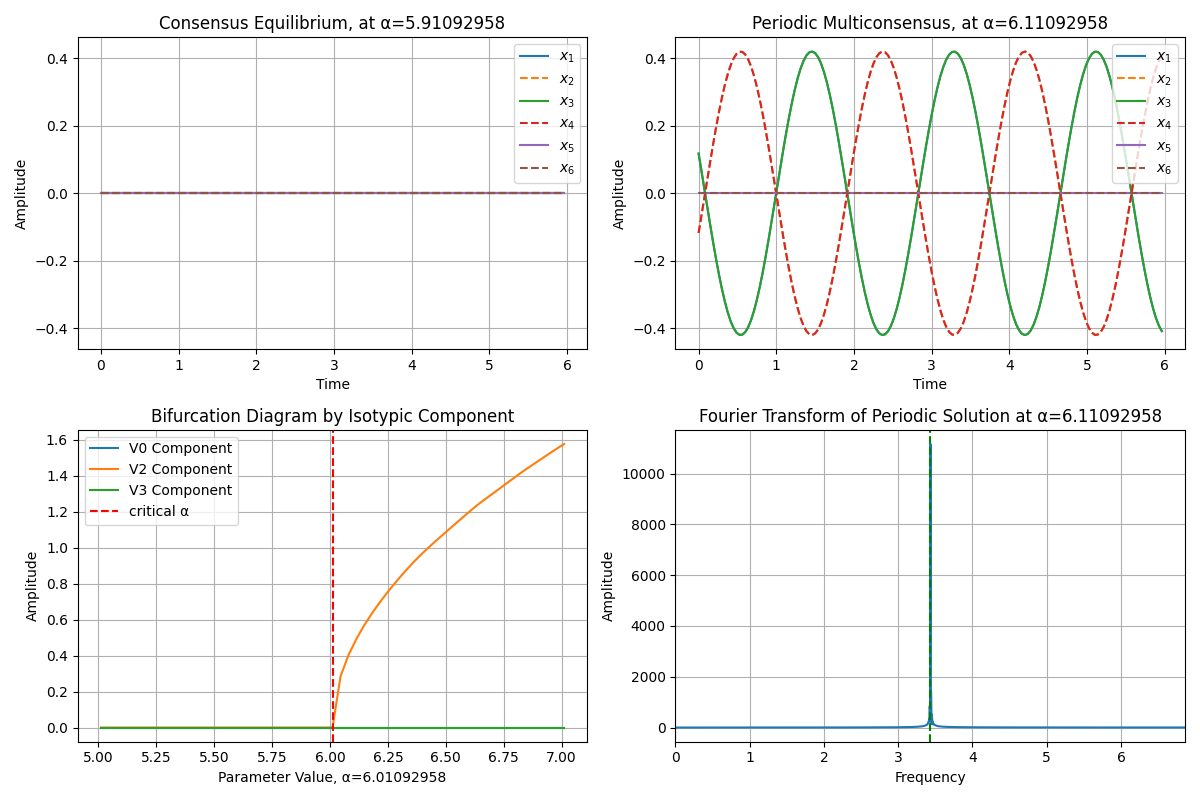} 
	\caption{Representative of the twisted wave solution type, initialized with a small perturbation in the direction $(-1,0,-1,0,1,1)^T$.}\label{fig:dist:twisted}
\end{figure}

\begin{figure}[H]
    \includegraphics[height=9cm]{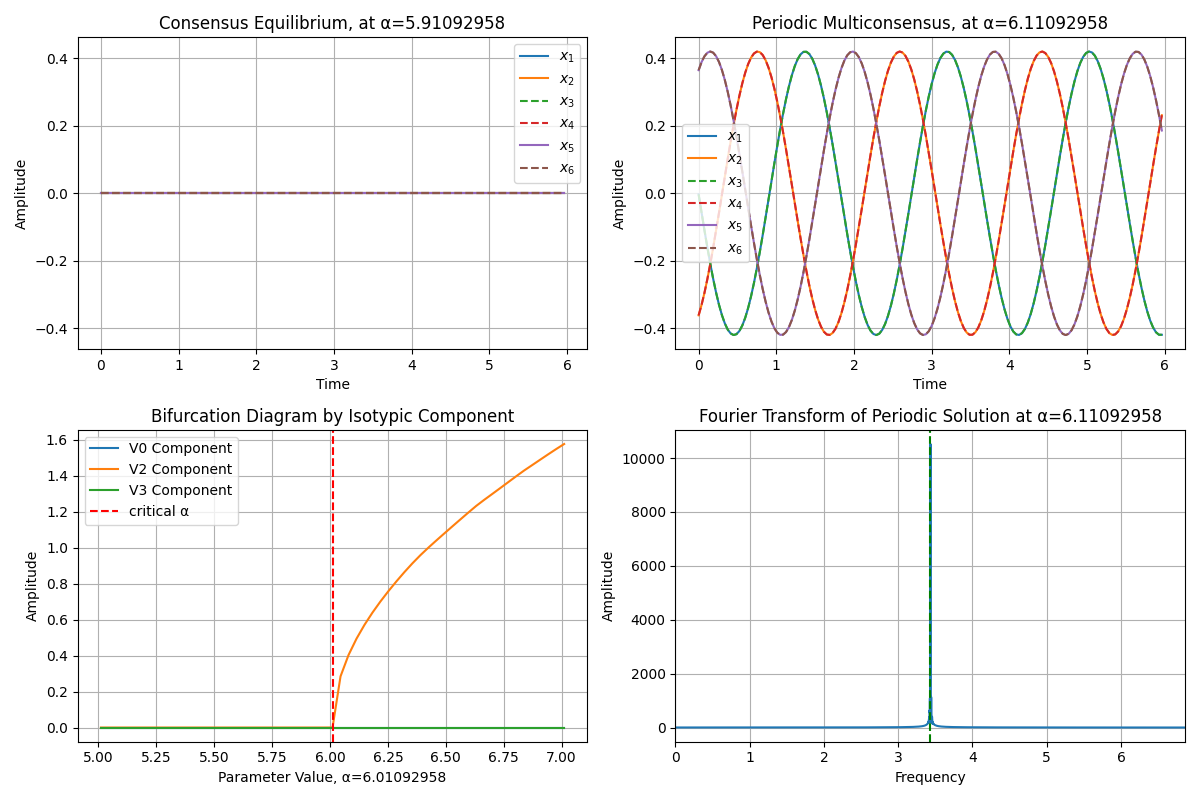} 
	\caption{Representative of the standing wave solution type, initialized with a small perturbation in the direction $(-0.9,0.6,-0.9,0.6,0.3,0.3)^T$.}\label{fig:dist:traveling}
\end{figure}

By Theorem \ref{thm:dist:global-bif}, this implies that for each of the three decoupled dynamics equations for each spatial axis given in \eqref{eq:dist:UAV-dynamics}, if $\alpha$ passes through $\alpha_{1,j}$ with $j=0,2,3$, then for each of the above orbit types in $(H)\in \mathfrak M_{1,j}$, there exists a bifurcating branch of non-constant periodic solutions with symmetries at least $(H)$. The above orbit types are given in the amalgamated notation, which is explained in detail in Appendix \ref{app:amalgamated}.
\vs
To show these solutions numerically, to show what types of spontaneous pattern formation is possible in synchronized formation oscillations, and to show how Theorem \ref{thm:dist:main1} and Corollary \ref{thm:dist:main2} can be applied concretely, we will now choose some specific parameter values. Put $b = 1$, $a = 1$, $d = 1/4$, and choose $f(x_i) = \tanh(bx_i)$, and $h(\bm x):= \tanh\left(\sum_{j=1}^8c_{ij}x_j\right),$ where $c_{ij}$ is an entry of $C$. Then we have
\begin{align*}
    a_0 &= a+4d = 2\\
    a_2 &= a-2d = \frac{1}{2}\\
    a_3 &= a = 1
\end{align*}
Since $a_2<a_3<a_0$, this implies that Hopf bifurcation will occur first on the $V_2$ isotypic component, and there exist solutions with isotropy types (in the first Fourier mode)
\[
\mathfrak M_{1,2} = \{( {\mathbb Z}_{6} {}^{{\mathbb Z}_1}\times {}^{V_4} A_4^p), ({\mathbb Z}_{2} {}^{{\mathbb Z}_1}\times {}^{D_4^{\hat{d}}} D_4^p), ( {\mathbb Z}_{2} {}^{{\mathbb Z}_1}\times {}^{D_4} D_4^p)\}
\]
Indeed, solving the coincidence equations on each isotypic component, we obtain:
\begin{align*}
    \alpha_{1,0} &\approx10.231716, &\quad \beta_{1,0} &\approx 4.0575156, \\
    \\
    \alpha_{1,2} &\approx6.0109295, &\quad \beta_{1,2} &\approx 3.4310143, \\
    \\
    \alpha_{1,3} &\approx7.2461785, &\quad \beta_{1,3} &\approx 3.6731944.
\end{align*}
Note that the $V_1$ subspace is two-dimensional and contains three orbit types of maximal kind. After bifurcation occurs, a periodic multiconsensus solution corresponding to each of these emerges, but the degree theory is insufficient to tell us which is stable. In order to ``select'' the branch corresponding to a particular orbit type, we must give the zero solution a ``kick'' in the direction of the subspace fixed by the orbit type of interest. This system can be seen to exhibit multistability, with each of the maximal orbit types in $\mathfrak M_{1,j}$ forming a stable periodic multiconsensus solution. The clustering behavior of the agents differs between these maximal orbit types. Across three axes, this implies \textbf{27 distinct (up to conjugacy) oscillating formations} satisfying our desired requirements. 
\vs
The different orbit maximal orbit types can be interpreted as three different phase-correlations between each agent's periodic motion. They can be broadly classed as ``standing wave,'' ``twisted wave,'' and ``travelling wave'' symmetries. $( {\mathbb Z}_{2} {}^{{\mathbb Z}_1}\times {}^{D_4} D_4^p)$ is the standing wave symmetry. It divides the agents into two clusters, $K_1,K_2$, which, for one representative of the conjugacy class, can be written $K_1:=\{x_1,x_2,x_3,x_4\},$ and $K_2:=\{x_5,x_6\}.$ All $x_i$ in $K_j$ are locked in phase, but $K_1$ and $K_2$ are in antiphase. Moreover, the oscillations in $K_2$ have twice the amplitude of those in $K_1$. This can be understood by the fact that $(-1,-1,-1,-1,2,2)^T\in E(\mu_2)$ is fixed by $( {\mathbb Z}_{2} {}^{{\mathbb Z}_1}\times {}^{D_4} D_4^p)$, and a kick in this direction will send agents into this stable configuration after bifurcation. A representative of this stable branch can be seen in Figure \ref{fig:dist:standing}.
\vs
$({\mathbb Z}_{2} {}^{{\mathbb Z}_1}\times{}^{D_4^{\hat{d}}} D_4^p)$ is the twisted wave symmetry. It divides agents into clusters $K_1,K_2,K_3$ which, for one representative of the conjugacy class, can be written $K_1:=\{x_1,x_3\}$, $K_2:=\{x_2,x_4\}$, $K_3:=\{x_5,x_6\}$. To summarize their phase relationship, let $K_i(t)$ denote the amplitude of cluster $i$ at time $t$. Then if $K_1(t)=1$, then $K_2(t)=-1$, and $K_3(t)=0$. The vectors $(-1,0,-1,0,1,1)^T$ and $(-1,1,-1,1,0,0)^T \in E(\mu_2)$ are each fixed by different conjugate isotropy groups in this orbit type and yield different conjugate branches. A representative of this stable branch can be seen in Figure \ref{fig:dist:twisted}.
\vs
Finally, the $( {\mathbb Z}_{6} {}^{{\mathbb Z}_1}\times {}^{V_4} A_4^p)$ orbit type corresponds to a ``travelling wave'' solution, where agents form three clusters $K_1,K_2,K_3$ which, for one representative, can be written $K_1:=\{x_1,x_3\}, K_2:=\{x_2,x_4\}, K_3:=\{x_5,x_6\}$. Each cluster is phase locked, and differs from the other two clusters in phase by $2\pi/3$. This has a two-dimensional fixed-point subspace, and this regime can be intered by a perturbation along any generic vector in $V_2$ not belonging to one of the other two fixed-point subspaces. A representative of this stable branch can be seen in Figure \ref{fig:dist:traveling}.
\vs
Since the perturbation vectors leading to the different branches differ significantly, it is easy to select for the desired branch. One can imagine different situations where different symmetric branches along different spatial axes are more desirable for different purposes. It is also possible to induce a phase shift between different axes in order to generate more sophisticated three dimensional spatio-temporal patterns. Another advantage of this protocol is that the formation of these oscillations (or quieting them back to the stable formation consensus) is controlled by a single parameter across all UAVs, and does not require any centralized coordination to orchestrate these complex movements. The parameter $\alpha$ not only controls the bifurcation, but also its amplitude. This can be used to guarantee that the maximum amplitude of oscillations is sufficient to avoid any collisions of UAVs.  
\vs
Finally, we note that the frequency of the oscillations is tunable by changing the delay length in the distributed delay $\int_0^1x_i(t-s)ds$. By setting $\int_0^\tau x_i(t-s)ds$ and solving the associated transcendental equations, one can see that as $\tau$ increases, critical frequencies decrease. This effect will be explored in more detail in Chapters \ref{chapter:pseudo} and \ref{chapter:neutral}, where multiple interacting delay lengths are considered.


\vspace{24pt}

\chapter{Global Periodic Multiconsensus Bifurcation in Pseudoneutral Multi-Agent Systems with Trend Memory}\label{chapter:pseudo}
 \author{Casey Crane}

\section{Introduction}\label{sec:pseudo:introduction}
Consider a MAS of $n$ homogeneous agents whose protocol depends on a \emph{trend memory} term $\frac{d}{dt}\left[\int_0^{\tau_1}g(x(t-s))ds\right] = g(x(t)) - g(x(t-s))$. This describes an agent which keeps track of the instantaneous rate of change of a moving continuous average over past states. We will consider a MAS whose agent protocols incorporate this type of trend memory, as well as continuous memory with possibly different time horizons. Leaving the trend memory term in implicit form in the left-hand side, the closed-loop dynamics of such a MAS for agent $i$ can be written
\[
\frac{d}{dt}\left[x_i - \int_0^{\tau_1} g(x_i(t-s))ds\right] = -ax_i - \alpha f\left( \int_0^{\tau_2} x_i(t-s)ds \right) + h_i(x), \;\quad x_i(t) \in \mathbb R.
\]
The components of the protocol consist of
\begin{itemize}
    \item A linear self-regulation term $-ax_i$, with strength $a>0$.
    \item A nonlinear interaction term $h_i(x)$ representing agent interaction in a fixed interaction topology with symmetry group $\Gamma_0$ and with state-dependent effective weights.
    \item A nonlinearly transformed continuous memory term $\int _0^{\tau_2}x_i(t-s)ds$ with gain strength $\alpha>0$, where we assume $f$ is differentiable at $0$ and denote $b:= f'(0)$.
    \item A trend memory term, written in implicit form as $\frac{d}{dt}\left[\int_0^{\tau_1}g(x_i(t-s))ds\right]$, where $g$ is a nonlinear function differentiable at $0$ and we denote $\gamma := g'(0)$. 
\end{itemize}
This forms a delay equation of \emph{pseudoneutral} type, as defined in Section \ref{sec:prelim:pseudoneutral}. Although this equation is of purely retarded type, we will analyze it in this chapter using the framework designed for neutral equations. We choose to do so for three primary reasons:
\begin{enumerate}
    \item To show the connection between trend memory and momentum memory in MAS dynamics in parallel with the connection between neutral and pseudoneutral equations.
    \item To develop the framework needed for neutral equations in more detail, which simplifies the development of the neutral model with momentum memory in the next chapter.
    \item Because the conditions needed for neutral equations to be analyzed using the equivariant Nussbaum-Sadovskii degree (i.e. to be of Hale class) coincide with certain fundamental stability conditions for pseudoneutral equations.
\end{enumerate}
This last point, which we will explore in detail, is of particular importance. Since we are interested in studying the breakdown of consensus states into periodic multiconsensus, we must first establish that there is a nonempty parameter range where the reference consensus at $x\equiv 0$ is achieved. It turns out that the conditions on the operator formulation of the left-hand side which allow the system to be formulated as a condensing perturbation of identity are also necessary to guarantee that there are parameter values (of trend memory parameter $\gamma$, delay terms $\tau_1, \tau_2, $ continuous memory gain $a$, and agent interaction strength $\mu_j$) such that the $x\equiv 0$ solution is locally asymptotically stable.
\vs
Assuming $n$ homogeneous agents with identical protocols (up to the symmetry of the interaction term), we put $\bm x(t):= (x_1(t),\dots,x_n(t))^T\in \mathbb R^n$, and 
if we put $\bm f,\bm g, \bm h:\mathbb R^n \to \mathbb R^n$ as $\bm f(\bm x) := (f(x_1),\dots,f(x_n))^T$, $\bm g(\bm x) := (g(x_1),\dots,g(x_n))^T$, and $\bm h(\bm x) := (h_1(x),\dots,h_n(x))^T \in \mathbb R^n$, 
then the closed-loop dynamics of the full system can be written
\begin{equation}\label{eq:pseudo:basic-sys}
\frac{d}{dt}\left[\bm x - \int_0^{\tau_1} \bm g(\bm x(t-s))ds\right] = -a \bm x - \alpha \bm f\left( \int_0^{\tau_2} \bm x(t-s)ds \right) - \bm h(\bm x).
\end{equation}
\vs
To account for the symmetries of the agent interaction topology, let $\Gamma_0\leq S_n$ be a finite symmetry group which acts by permuting the indices of $\bm x$. Let $V:= \mathbb R^n$, and put $\Gamma := \Gamma_0 \times \mathbb Z_2$, where $\mathbb Z_2$ acts antipodally, i.e. the action of $(\sigma,\pm1)$ on $\bm x \in V$ is given by
\[
    (\sigma,\pm1)\bm x :=  (\sigma,\pm1)(x_1,\dots,x_n)^T=\pm(x_{\sigma(1)},\dots,x_{\sigma(n)})^T.
\]
We will make the following assumptions about $\bm f, \bm g,$ and $\bm h$:
\begin{enumerate}[label=($C_\arabic*$)]\setcounter{enumi}{0}\itemc
\item\label{c1} $\bm f$, $\bm g$, and $\bm h$ are $\Gamma$-equivariant (and therefore are odd functions).
\item\label{c2} $\bm f$, $\bm g$, and $\bm h$ are continuous functions differentiable at 0, and $b := f'(0), \gamma := g'(0), C:= Dh(0)$ satisfy:
    \begin{enumerate}
        \item $b\neq0$.
        \item $0 < \gamma < 1$.
        \item C is a $\Gamma$-equivariant symmetric matrix.
    \end{enumerate}
\item\label{c3} $\bm g$ is $\kappa$-Lipschitzian with $\kappa<1$, i.e. 
\[
\exists_{\kappa \in [0,1)}\quad \forall_{\varphi,\psi \in C(\bbR^n;\bbR)} \quad ||\bm g(\varphi) - \bm g(\psi)|| \leq \kappa||\varphi - \psi||_\infty
\]
\end{enumerate}
The assumption \ref{c3} and the requirement that $0<\gamma<1$ are not strictly necessary if the trend memory is written explicitly as a finite difference, but in practical terms, such systems where the trend memory operator is not contractive have serious stability problems, even if their operator equation formulations are still compact perturbations of identity. 
\vs
We will refer to $x\equiv 0$ as the \emph{reference consensus} or the \emph{consensus equilibrium}, and we will provide conditions for its local asymptotic stability and thus the achieval of consensus in system \eqref{eq:pseudo:basic-sys}. We will then show conditions for the breakdown of consensus into periodic multiconsensus regimes using the framework of symmetric global Hopf bifurcation. We will also show conditions for these multiconsensus branches to be unbounded, and classify them according to their spatio-temporal symmetries, corresponding to different configurations of subpopulations of agents cooperating and competing with one another. We will also provide an illustrative example of a MAS simulating a coupled asset market with multiple trading strategies corresponding to the two types of memory present in this system. The question of which multiconsensus states are actually achieved (i.e. which are at least locally asymptotically stable) cannot be answered with equivariant degree techniques alone, and so in this example we will show how this can be addressed using numerical methods. 
\vs
Since the closed-loop dynamics equation \eqref{eq:pseudo:basic-sys} takes the form of a delay differential equation, we will take a ``bilingual'' approach to our results, stating them here in terms of symmetric periodic multiconsensus solutions emerging from a consensus-breaking bifurcation, and then restating them and proving them throughout the chapter in their logically equivalent form as results on global symmetric Hopf bifurcation of non-constant periodic solutions in families of delay differential equations of the form \eqref{eq:pseudo:basic-sys}, regardless of their intepretation.
Specifically, for a system of the form \eqref{eq:pseudo:basic-sys} satisfying \ref{c1} -- \ref{c3}, we will prove the following theorems:
\begin{theorem}\label{thm:pseudo:mas-asymp}
Denote by $\alpha_0$ the smallest $\alpha>0$ satisfying \eqref{eq:pseudo:coincidence}. If $0<\gamma\tau_1\leq 1, b>0$, and $a_j>0$ for all $j$, then $\bm x\equiv0$ is a locally asymptotically stable consensus equilibrium for all $\alpha\in(0,\alpha_0)$.  
\end{theorem}
\begin{theorem}\label{thm:pseudo:mas-local}
If $\beta_0>0$ satisfies
\[
\beta_0 = \gamma\sin(\beta_0\tau_1) + \frac{(\gamma(\cos(\beta_0\tau_1)-1)+a_j)(\cos(\beta_0\tau_2)-1)}{\sin(\beta_0\tau_2)}
\]
and $\alpha_0>0$ satisfies
\[
\alpha_0 = \frac{-\gamma\beta_0(\cos(\beta_0\tau_1)-1)-\beta_0 a_j}{b\sin(\beta_0\tau_2)}
\]
for some $j$, and if 
\[
{\rho_j}(\alpha_0,\beta_0):= p(\alpha_0,\beta_0)b(\cos(\beta_0\tau_2)-1)+ q(\alpha_0,\beta_0)b\sin(\beta_0\tau_2) \neq 0,
\]
where
\begin{align*}
p(\alpha,\beta) &:= \gamma(\cos(\beta\tau_1)-1) + a_j - \gamma\beta\tau_1\sin(\beta\tau_1) + \alpha b \tau_2 \cos(\beta\tau_2)\\
q(\alpha,\beta) &:= \gamma\sin(\beta\tau_1)+\gamma\beta\tau_1\cos(\beta\tau_1)+\alpha b \tau_2\sin(\beta\tau_2)-2\beta
\end{align*}
is satisfied at $\alpha_0,\beta_0$, then there exists a connected global branch of non-constant periodic multi-consensus solutions bifurcating from the trivial consensus branch at $\alpha=\alpha_0$.
\vs
\end{theorem}
\begin{theorem}\label{thm:pseudo:mas-global}
If the conditions of Theorem \ref{thm:pseudo:mas-local} hold and additionally $|\gamma(\tau_1-\tau_2)|<|\tau_2(\gamma-a_j)-2|$, then the bifurcating branch is global and unbounded.
\end{theorem}
\begin{theorem}\label{thm:pseudo:mas-sym}
If the conditions of Theorem \ref{thm:pseudo:mas-local} hold, and $(H)<\Phi_1^t(S^1\times \Gamma_0\times \mathbb Z_2;\mathscr E_{1,j})$ is a maximal twisted orbit type in $\mathscr E_{1,j}:=\{c\cos(t)+d\sin(t):c,d\in E(\mu_j)\}$, then there exists a connected branch of non-constant periodic solutions with symmetries at least $(H)$ bifurcating from the trivial solution at $\alpha=\alpha_0$.
Moreover, if the conditions of Theorem \ref{thm:pseudo:mas-global} are satisfied, then this branch is unbounded.
\end{theorem}
We consider the gain coefficient $\alpha$ of the continuous memory term as our bifurcation parameter. We do this so as to facilitate an easier, more direct, and therefore more interesting and useful comparison with the related systems studied in Chapter \ref{chapter:distributed} and Chapter \ref{chapter:neutral}. There are obviously many potentially interesting choices of bifurcation parameter for this MAS. Other canonical choices for MAS with time-lags of various types are the delay terms themselves. The overall methodology for any bifurcation parameter is the same, and one would only need to change the proofs of the signs of crossing numbers (i.e. transversality conditions), the stability of the reference consensus, and the derivation of the transcendental equations governing the critical set. This is not to say that such reformulations would be trivial, but none of the core theoretical assumptions would need to be changed in any way.  
\vs
This chapter will be structured as follows: In Section \ref{sec:pseudo:introduction:overview}, we will give a general treatment of how neutral functional differential equations (informally, neutral equations or NFDEs) are treated using equivariant degree methods. In particular, this concerns conditions for an implicitly formulated neutral or pseudoneutral equation to be a condensing perturbation of identity and therefore admissible for the Nussbaum-Sadovskii degree. For neutral equations, this coincides with the equation being of Hale class, and for pseudoneutral equations, it has spectral implications for the linearized pseudoneutral delay operator.
\vs
Following this general treatment, we will then linearize system \eqref{eq:pseudo:basic-sys} in Section \ref{sec:pseudo:crit}, compute its characteristic equation, derive a transversality condition used later to compute crossing numbers to show global bifurcation, and prove conditions for local asymptotic stability of the reference consensus. In Section \ref{sec:pseudo:funcspace}, we will reformulate this specific system as a fixed point equation in line with the general method outlined in Section \ref{sec:pseudo:introduction:overview} as a fixed point equation and show that assumptions \ref{c1}---\ref{c3} are sufficient to guarantee that it is a condensing perturbation of identity. In Section \ref{sec:pseudo:two-parameter}, we will formulate the Hopf bifurcation problem for system \eqref{eq:pseudo:basic-sys} as a two-parameter bifurcation problem and show how the relevant topological bifurcation invariants can be computed using spectral information derived from the characteristic equation. In Section \ref{sec:pseudo:fixed-point-reduction}, we will complete the proof of the main result by showing how degeneracies can be resolved by restricting the problem to an appropriate fixed point space which guarantees that global branches of non-constant periodic solutions remain non-constant and periodic across their full extent.
\vs
This will complete the first part of this chapter and conclude the abstract results. In Section \ref{sec:pseudo:example}, we will show an interesting application of a MAS of type \eqref{eq:pseudo:basic-sys}, modeling coupled asset markets where traders use a combination of price fundamental strategies and momentum trading strategies, and show the bifurcation of the consensus state to symmetric periodic multiconsensus. All of the above theorems will be applied, and we will provide numerical simulations which not only support the theoretical predictions, but also demonstrate how numerical methods can work in concert with the wealth of information provided by equivariant degree methods.
We will finish with some concluding remarks in Section \ref{sec:pseudo:conclusion} bridging this chapter and the following one, which will use much of the machinery developed in this chapter in parallel to study a system of truly neutral type.
\vs

\subsection{Preliminaries}\label{sec:pseudo:introduction:overview}
The degree we will use to analyze system \eqref{eq:pseudo:basic-sys} is the equivariant twisted Nussbaum-Sadovskii degree (ETNS degree). To use this degree, the system of interest must first be period normalized, reformulated as an operator equation between functional spaces, and then linearized. To streamline further exposition and simplify Chapter \ref{chapter:neutral}, we will first show this process for a generic neutral or pseudoneutral system where the neutral or pseudoneutral delay operator is written implicitly inside the derivative. Following this, we will also give some essential definitions regarding the representation theory of the relevant groups on this functional space and its $G$-isotypic decomposition, and the twisted orbit types used to classify spatio-temporal symmetries of solutions.

\vs

\subsubsection{NFDE preliminaries and reformulation as a condensing perturbation of identity}
Take the system 
\begin{equation}\label{eq:pseudo:fde2}
\frac d{dt}\big[x(t)- \bm k (x_t)\big]=\bm f(\alpha,x_t), \quad x(t)\in \br^n,
\end{equation}
and, for a fixed number $r>0$, consider the Banach space $\bm C_r:=C([-r,0];\br^n)$ equipped with the standard sup-norm, assume that $\bm f(\alpha,x):\bbR \times \bm C_r \to \br^n$ is a continuous map and $\bm k:\bm C_r\to \br^n$ is continuous, and take $\alpha$ as the bifurcation parameter. 
% Such a system is a functional differential equation of neutral type, or an NFDE for short.
\vs
We further assume that $\bm f$ is differentiable at $0$ with respect to $x$ (and $D_x \bm f(\alpha,0)$ is continuous with respect to $\alpha$), that $\bm k$ is continuously differentiable and $\kappa$-Lipschitzian with $\kappa <1$, and that $\norm{D_x\bm k}<1$ for all $x \in \bm C_r$. Finally, we assume without loss of generality that $x\equiv0$ is a trivial solution to \eqref{eq:pseudo:fde2}, i.e. $\bm k(0) = 0$ and $\bm f(\alpha,0)=0$.
\vs
We wish to reformulate \eqref{eq:pseudo:fde2} as a fixed point equation on a functional space. We first perform period normalization, using the substitution $x(t) = u(\beta t)$ where $\beta:=\frac{2\pi}{p}$ and $p>0$ is a known period of $x(t)$, to obtain
\begin{equation}\label{eq:pseudo:fde-norm}
    \frac d{dt}\big[u(t) - \bm k(u_{t,\beta})\big] = \frac{1}{\beta}\bm f(\alpha,u_{t,\beta}), \;\text{ and } u_{t,\beta}(s) := u(t+\beta s)
\end{equation}
Put $\mathscr{E} := C^1_{2\pi}(\bbR;\bbR^n)$ and $||\phi|| =\max\{||\phi||_\infty,||\dot{\phi}||_\infty\}$,
where $||\phi||_\infty$ denotes the usual supremum norm on $\phi$. Let $\bm j: \mathscr{E} \to C_{2\pi}(\bbR;\bbR^n)$ denote the natural embedding of $\mathscr{E}$ into $C_{2\pi}(\bbR;\bbR^n)$. Let $\bbR_+^2 := \bbR \times \bbR_+$, where $\bbR_+ := \beta>0$, and put 
\[
\mathscr{D}:= \{(\alpha,\beta,u)\in \bbR_+^2 \times C^1_{2\pi}(\bbR;\bbR^n): u+N_{\bm{k}}(u) \in C^1_{2\pi}(\bbR;\bbR^n)\}
\]
where $N_{\bm{k}}:\bbR^2_+ \times C_{2\pi}(\bbR;\bbR^n) \to C_{2\pi}(\bbR;\bbR^n)$, and $N_{\bm{k}}(\alpha,\beta,u)(t) := \bm k(u_{t,\beta})$. 

Notice that $C_{2\pi}^1(\bbR;\bbR^n) \mathop{\hookrightarrow}\limits^{\bm j}C_{2\pi}(\bbR;\bbR^n)$ is a compact embedding, and put 
\[
N_{\bm f}(\alpha,\beta,u)(t) := \frac{1}{\beta} \bm f(\alpha,u_{t,\beta}).
\]
Define $L:\mathscr D \subset \mathscr E \to C_{2\pi}(\bbR;\bbR^n)$ as $(Lu)(t) := \dot u(t)$. Then the equation \eqref{eq:pseudo:fde-norm} can be written as
\begin{equation}\label{eq:pseudo:fde-func1}
L(u-N_{\bm k}(\alpha,\beta,\bm j(u))) = N_{\bm f}(\alpha,\beta,u),
\end{equation}
where $N_{\bm f}(\alpha,\beta,u) \in \mathscr D$. From \eqref{eq:pseudo:fde-func1}, we obtain
\begin{align*}
L(u-N_{\bm k}(\alpha,\beta,u)))+\bm j(u) - N_{\bm k}(\alpha,\beta,\bm j(u))&= N_{\bm f}(\alpha,\beta,u) + \bm j(u) - N_{\bm k}(\alpha,\beta,\bm j(u))\\
(L + \bm j)(u-N_{\bm k}(\alpha,\beta,u))&=N_{\bm f}(\alpha,\beta,u) + \bm j(u) - N_{\bm k}(\alpha,\beta,\bm j(u)),
\end{align*}
from which we obtain the fixed point operator equation
\begin{equation}\label{eq:pseudo:fde-fp-op-eq}
\mathscr F(\alpha,\beta,u) = u - N_{\bm k}(\alpha,\beta,u) - (L + \bm j)^{-1}\Big(N_{\bm f}(\alpha,\beta,u) + j(u) - N_{\bm k}(\alpha,\beta,\bm j(u))\Big)
\end{equation}
Clearly, $\mathscr F$ is a Darbo operator, i.e. it is a sum of a completely continuous and a contractive function, and therefore $\mathscr F$ is a condensing perturbation of identity (cf. Section \ref{sec:prelim:n-s}). Clearly, $(\alpha,\beta,u)\in\bbR_+^2\times \mathscr E$ satisfying $\mathscr F(\alpha,\beta,u) = 0$ are solutions to \eqref{eq:pseudo:fde2}. This allows us to view \eqref{eq:pseudo:fde2} as a two-parameter bifurcation problem which can be analyzed using the ETNS degree.
\vs

\subsubsection{Symmetry groups, irreducible representations, and $G$-isotypic decomposition of $\mathscr E$}
Now we will introduce some standard preliminaries regarding the structure of irreducible $G$-representations, $G$-isotypic components, and how $G$ acts on the functional space $\mathscr E$. Let $\Gamma_0$ be a finite group representing the coupling symmetries of $\eqref{eq:pseudo:fde2}$, and put $\Gamma:=\Gamma_0\times\bbZ_2$. Suppose $V:=\bbR^n$ is an orthogonal $\Gamma$-representation. If $\bm f,\bm k$ are $\Gamma$-equivariant, then $\mathscr F$ is $G$-equivariant, where $G := S^1 \times \Gamma$, and $S^1$ acts on the period-normalized equation by phase shifting on Fourier modes.
\vs
V has a $\Gamma$-isotypic decomposition given by
\[
V = V_0 \oplus V_1 \oplus \dots \oplus V_r
\]
where $V_i$ is the (real) $\Gamma$-isotypic component modeled on the irreducible $\Gamma$-representation $\cV_i^-$, which is the irreducible $\Gamma_0$-representation $\cV_i$ with the antipodal $\bbZ_2$-action. Clearly, if $\cV_i$ is an irreducible $\Gamma_0$-representation of \emph{real type}, then
\[
\cV_i^{\mathbb C} := \mathbb C \otimes_\bbR \cV_i
\]
is an irreducible complex $\Gamma_0$-representation, and is called the \emph{complexification} of $\cV_i$. On the other hand, if $\cV_i$ is of \emph{complex type}, then it is also an irreducible complex $\Gamma_0$-representation. Therefore
\[
V^\mathbb C = V_0^{\mathbb C} \oplus V_1^{\mathbb C} \oplus \dots \oplus V_r^{\mathbb C},
\]
where, for $\cV_i$ of real type
\[
U_i := \cV_i^{\mathbb C}
\]
is a complex $\Gamma_0$-isotypic component modeled on $\cV_i^{\mathbb C}$, and if $\cV_i$ is of complex type, then
\[
V_i^{\mathbb C} = U_i \oplus U_i'
\]
where $U_i$ is a complex $\Gamma_0$-isotypic component modeled on $\cV_i$ and $U_i'$ is a complex $\Gamma_0$-isotypic component modeled on $\overline{\cV_i}$. For simplicity, in this paper we will assume all the irreducible $\Gamma_0$-representations $\cV_i$ are of real type. 
\vs
$S^1$ has a two-dimensional real irreducible representation $\mathcal W_k$ for each $k \in \mathbb N$, where 
\[
\mathcal W_k := \{p\cos(kt) + q\sin(kt): p,q\in \mathbb R\}
\]
and $\theta \in S^1$ acts on $u \in \mathcal W_k$ as $\theta u(t) = u(t+k\theta)$. Then all irreducible $G$-representation are given by
\[
\mathcal W_{k,j} := \mathcal W_k \otimes \mathcal V_j^-
\]
where $(\theta,\sigma,\pm1) \in G$ acts on $u = (u_1(t),\dots,u_n(t))^T \in \mathcal W_{k,j}$ as
\[
(\theta,\sigma,\pm1)(u_1(t),\dots,u_n(t))^T = (u_{\sigma(1)}(t+k\theta),\dots,u_{\sigma(n)}(t+k\theta))^T,
\]
where $\sigma \in \Gamma_0$ is viewed as a permutation of the indices of vectors in $\mathbb R^n$. Since every $u \in \mathscr E$ admits a Fourier decomposition, the $G$-isotypic component $W_{k,j}$ (modeled on the irreducible representation $\mathcal W_{k,j})$ is equivalent as a G-space to
\[
\mathscr E_{k,j} := \{c\cos(kt)+d\sin(kt):c,d \in V_j\}
\]
and so the $G$-isotypic decomposition of $\mathscr E$ is given by
\[
\mathscr E := \overline{\bigoplus_{k=0}^\infty \bigoplus_{j=0}^r \mathscr E_{k,j}}
\]
\vs
\subsubsection{Orbit types}\label{sec:pseudo:prelim-orbittypes}
We denote by $\Phi(G)$ the set of all conjugacy classes of subgroups of $G$. A conjugacy class of subgroups $(H)$ is called an \textit{orbit type} if there is some $u \in \mathscr E\setminus\{0\}$ such that $hu = u$ for all $h \in H$. We denote by $\Phi(G;\mathscr E)$ the set of all orbit types of $G$ on $\mathscr E$. These sets admit a natural partial order given by subconjugacy: $(H) \leq (K)$ if and only if there exists some $H' \in (H)$ and $K'\in (K)$ such that $H' \leq K'$. 
\vs
We say that $(H)$ is a \textit{maximal orbit type} if $(H)\leq (K)$ implies $(K)=(H)$ or $(K)=(G)$. On the other hand, a vector $u \in \mathscr E$ has the \textit{isotropy group} $G_u := \{g\in G:gu=u\}$, and we call $(G_u)$ the \textit{isotropy type} of $u$. So we may talk about the orbit types of $G$ on $\mathscr E$, or about the isotropy types of some particular $u\in \mathscr E$, and these are dual notions.
\vs
The ETNS degree is a form of twisted degree, which can detect solutions according to their spatio-temporal symmetries. In particular, it can detect non-constant periodic solutions whose isotropy $H \leq G$ satisfies $\dim W_G(H) = 1$, where $W_G(H) = N_G(H)/H$ is the Weyl group of $H$ in $G$ and $N_G(H)$ is the normalizer of $H$ in $G$.
\vs
These are so-called \textit{twisted orbit types}, and take the form
\[
({\mathbb Z_{km}}^{\mathbb Z_k} \times {}^{H_0}H) := \{(z,h)\in \mathbb Z_{km} \times H: \theta(h)=z^m\}
\]
where $\theta:\Gamma \to S^1$ is a homomorphism with $\ker \theta = H_0$. We denote the set of twisted orbit types by $\Phi_1^t(G;\mathscr E)$. 
\vs

We will denote the set of maximal twisted orbit types in $\Phi_1^t(G;\mathscr E_{k,j})$ as $\mathfrak M_{k,j}$, and the maximal twisted orbit types in $\Phi_1^t(G;\mathscr E_{1,j})$ simply as $\mathfrak M_j$. For more details on the amalgamated notation used above to describe twisted subgroups, see Appendix \ref{app:amalgamated}. For a detailed construction of $\Phi_1^t(G)$ and the twisted degree, see \cite{BalanovEtAl2025}.
\vs
We also note that maximal twisted orbit types in $G$ may be very easily computed in GAP for any finite $\Gamma_0$ and any representation. Concrete examples of such computations will be illustrated in Section \ref{sec:pseudo:example}. 

\vs
Viewing system \eqref{eq:pseudo:basic-sys} as an FDE, under assumptions \ref{c1} -- \ref{c3}, it clearly meets the conditions outlined in Section \ref{sec:pseudo:introduction:overview}, and its period normalization and functional space reformulation will be along similar lines. To apply the equivariant Krasnosel'skii and Rabinowitz theorems, we must find a computable formula for the local bifurcation invariants, and show that those invariants are well-defined in the first place. 
\vs
This means showing that critical values of $\alpha$ are isolated, and studying the spectrum of our linearized system. We will do this in two stages. First, we will linearize the system \eqref{eq:pseudo:basic-sys} and study its characteristic operator equation. Then we will perform the same type of reformulation as in Section \ref{sec:pseudo:introduction:overview}, linearize the period-normalized operator equation, and describe its spectrum in terms of the characteristic operator equation we obtain here. This will provide us with all necessary information to compute local bifurcation invariants and thus ultimately to prove the main results.

\section{Linearization, spectral analysis, and critical points}\label{sec:pseudo:crit}

\vs
Linearizing \eqref{eq:pseudo:basic-sys} around $x=0$ we obtain
\begin{equation}\label{eq:pseudo:basic-sys-lin}
\frac{d}{dt} \left[\bm x - \gamma \int_0^{\tau_1} \bm x(t-s)ds\right] + a \bm x + \alpha b \int_0^{\tau_2} \bm x(t-s)ds + C \bm x = 0
\end{equation}
\vs
The characteristic operator equation for \eqref{eq:pseudo:basic-sys-lin} is given by
\begin{equation}\label{eq:pseudo:char-op}
\triangle_{\alpha}(\lambda):=\left(\lambda + \gamma(e^{-\lambda\tau_1}-1) + a - \alpha b \frac{e^{-\lambda\tau_2} -1}{\lambda}\right)\text{Id} + C
\end{equation}
\begin{definition}[Centers and limit frequencies]\normalfont
The $\bm x \equiv 0$ solution to \eqref{eq:pseudo:basic-sys-lin} at $\alpha=\alpha_0$ is called a \emph{center} if there exists a corresponding $\beta_0>0$ (called the \textit{limit frequency}) such that $\det_{\mathbb C}\triangle_{\alpha_0}(i\beta_0) =0$. If, in addition, there exists $\varepsilon>0$ such that $0<|\alpha-\alpha_0|+|\beta-\beta_0|<\varepsilon$ implies
\[
{\det}_{\mathbb C}\triangle_{\alpha}(ik\beta) \not=0 
\]
for all $k\in\mathbb N$, then it is called an \emph{isolated center}.
\end{definition}

\vs
Since the system \eqref{eq:pseudo:basic-sys} is $\Gamma_0$-equivariant, its linearization \eqref{eq:pseudo:basic-sys-lin} is also $\Gamma_0$-equivariant and therefore so is $\triangle_\alpha(\lambda)$. As such, it can be decomposed along the $\Gamma_0$-isotypic components of $V$ by restricting it to the complexified irreducible $\Gamma_0$-representations $\cV_j^{\mathbb C}$ on which the $\Gamma_0$-isotypic components are modeled. These will be called $\Gamma_0$-isotypic characteristic equations, denoted
\begin{equation}\label{eq:pseudo:char-op-iso-restrict}
\triangle_{\alpha,j}(\lambda):= \triangle_\alpha(\lambda)_{|\cV^{\mathbb C}_j}:\cV^{\mathbb C}_j\to \cV^{\mathbb C}_j
\end{equation}
Put $a_j := a + \mu_j$, where $\mu_j$ is the $j$-th eigenvalue of $C$ whose eigenspace $E(\mu_j)$ corresponds to the $j$-th $\Gamma_0$-isotypic component. Then under the assumptions \ref{c1} -- \ref{c3}, \eqref{eq:pseudo:char-op-iso-restrict} can be written
\begin{equation}\label{eq:pseudo:char-op-iso}
\triangle_{\alpha,j}(\lambda):=\left(\lambda + \gamma(e^{-\lambda\tau_1}-1) + a_j - \alpha b \frac{e^{-\lambda\tau_2} -1}{\lambda}\right)\text{Id},
\end{equation} 
and \eqref{eq:pseudo:char-op} and its determinant can be expressed in terms of the $\Gamma_0$-isotypic characteristic equations as a block matrix:
\[
\triangle_\alpha(\lambda) = \begin{bmatrix}
    \triangle_{\alpha,0}(\lambda) & &\makebox(0,0){\text{\huge0}}\\
    & \ddots & \\
    \text{\huge0}& & \triangle_{\alpha,r}(\lambda)
\end{bmatrix},
\]
\begin{equation}\label{eq:pseudo:char-block-det}
{\det}_{\mathbb C}(\triangle_\alpha(\lambda)) = \prod_{j=0}^r {\det}_{\mathbb C}(\triangle_{\alpha,j}(\lambda))^{m_j}
\end{equation}
Here, $m_j$ denotes the \textit{isotypic multiplicity} of the eigenvalue $\mu_j$, defined
\begin{equation}\label{eq:pseudo:def-iso-mult}
m_j := \frac{\text{dim }E(\mu_j)}{\text{dim }\cV_j}
\end{equation}
\subsection{The critical set}
We define the \textit{characteristic quasipolynomial} corresponding to $\triangle_{\alpha,j}(\lambda)$, written
\begin{equation}\label{eq:pseudo:char-root}
P_j(\alpha,\lambda) := \lambda^2 + \gamma(e^{-\lambda\tau_1}-1)\lambda + a_j\lambda -\alpha b(e^{-\lambda\tau_2}-1).
\end{equation}  
From \eqref{eq:pseudo:char-block-det} we can clearly see that $\lambda \in {\mathbb C}\setminus\{0\}$ is a characteristic root of \eqref{eq:pseudo:char-op} if and only if, for some $j=0,\dots,r$ and some $\alpha$, we have $P_j(\alpha,\lambda)=0$. Furthermore, the linearized system \eqref{eq:pseudo:basic-sys-lin} admits a center if and only if, for some $j$, there exists $\alpha_0$ and $\beta_0$ such that $P_j(\alpha_0,i\beta_0) = 0$.
We will now show that the system \eqref{eq:pseudo:basic-sys-lin} admits an infinite number of isolated centers.
Setting $\lambda = i\beta$ and substituting into \eqref{eq:pseudo:char-root}, we obtain 
\[
-\beta^2 + i\beta\gamma(\cos(\beta\tau_1) - i\sin(\beta\tau_1) -1) + i\beta a_j - \alpha b (\cos(\beta\tau_2) - i\sin(\beta\tau_2) - 1) = 0
\]
which yields the system
\begin{align*}
-\beta^2 + \gamma\beta\sin(\beta\tau_1) - \alpha b (\cos(\beta\tau_2) - 1) &= 0\label{eq:pseudo:system-sincos}\\
\gamma\beta(\cos(\beta\tau_1)-1) + \beta a_j + \alpha b \sin(\beta\tau_2) &= 0.\nonumber
\end{align*}
Solving the second equation for $\alpha$, we get
\begin{equation}\label{eq:pseudo:alpha-relation}
\alpha = \frac{-\gamma\beta(\cos(\beta\tau_1)-1)-\beta a_j}{b\sin(\beta\tau_2)}.
\end{equation}
Substituting this into the first equation and dividing out $\beta$, we obtain the relation
\begin{equation}\label{eq:pseudo:beta-relation}
\beta = \gamma\sin(\beta\tau_1) + \frac{(\gamma(\cos(\beta\tau_1)-1)+a_j)(\cos(\beta\tau_2)-1)}{\sin(\beta\tau_2)}.
\end{equation}

\begin{figure}[tbp]
    \centering
    \begin{minipage}{0.45\textwidth}
        \centering
        \includegraphics[width=\linewidth]{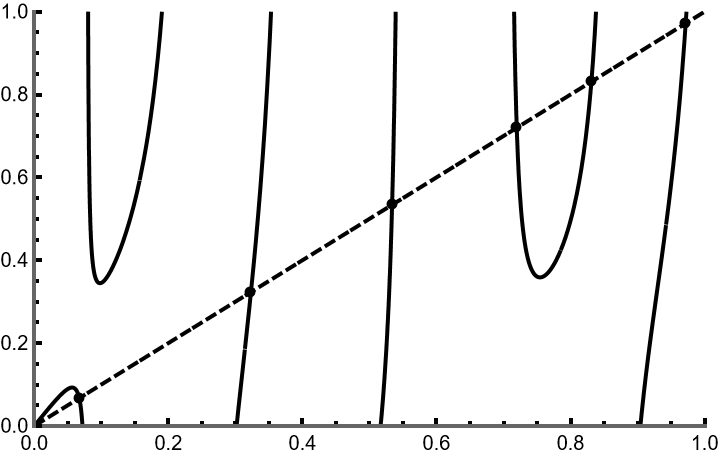}
        \subcaptiontext{a}{Limit frequencies for $0<\beta<1$}
    \end{minipage}\hfill
    \begin{minipage}{0.45\textwidth}
        \centering
        \includegraphics[width=\linewidth]{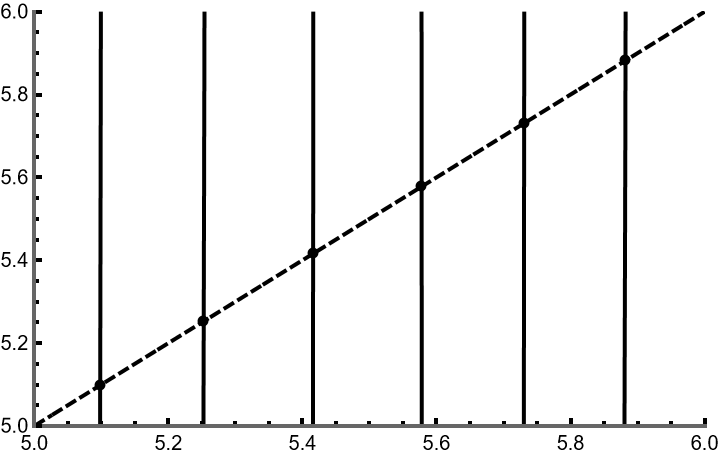}
        \subcaptiontext{b}{Limit frequencies for $5<\beta<6$}
    \end{minipage}
    \caption{Example plots of limit frequencies of \eqref{eq:pseudo:basic-sys} showing behavior at relatively lower and higher values of $\beta$. Here $\gamma = 0.6$, $\tau_1 = 9$, $\tau_2 = 40$, $a_j=0.17$, and $b=0.4$. As $\beta$ grows larger, $\beta_{n+1,j} - \beta_{n,j}$ approaches $\frac{2\pi}{\tau_2}$.}
    \label{fig:pseudo:beta-relation}
\end{figure}
A solution $\beta_0$ to the transcendental coincidence problem \eqref{eq:pseudo:beta-relation} corresponds to a characteristic root of \eqref{eq:pseudo:char-op-iso}. Such a root, through \eqref{eq:pseudo:alpha-relation}, corresponds to an isolated center $\alpha_0$ of \eqref{eq:pseudo:basic-sys} at $x=0$ with limit frequency $\beta_0$. To further analyze this coincidence problem and describe the set of all such centers and limit frequencies, put
\begin{align*}
\phi_j(\beta) &:= \gamma\sin(\beta\tau_1) + \frac{(\gamma(\cos(\beta\tau_1)-1)+a_j)(\cos(\beta\tau_2)-1)}{\sin(\beta\tau_2)}\\
\zeta(\beta) &:= \beta
\end{align*}
and note that the intersection of the graphs of $\phi_j$ and $\zeta$ in $\mathbb R_+^2$ forms a discrete and monotonically increasing sequence of values, as shown in Figure \ref{fig:pseudo:beta-relation}. This can be seen from the fact that $\phi_j$ is periodic and has a bounded numerator. As $\beta$ grows large, intersections can only occur once per period, near the points where $\sin(\beta\tau_2)$=0, which clearly form a discrete set. Therefore we obtain a sequence $\alpha_{n,j}$ of isolated centers of \eqref{eq:pseudo:basic-sys}, and corresponding characteristic roots $\beta_{n,j}$ satisfying
\begin{equation}\label{eq:pseudo:coincidence}
\begin{aligned}
\alpha_{n,j} &= -\beta_{n,j}\frac{\gamma(\cos(\beta_{n,j}\tau_1) -1) + a_j}{b\sin(\beta_{n,j}\tau_2)}\\
\beta_{n,j} &= \gamma\sin(\beta_{n,j}\tau_1) + \frac{(\gamma(\cos(\beta_{n,j}\tau_1)-1)+a_j)(\cos(\beta_{n,j}\tau_2)-1)}{\sin(\beta_{n,j}\tau_2)}
\end{aligned}
\end{equation}
A pair of these values together with the trivial solution $(\alpha_{n,j},\beta_{n,j},0)$ is called a \emph{critical point} of our system \eqref{eq:pseudo:basic-sys}, and the set 
\begin{equation}\label{eq:pseudo:crit-set}
\Lambda := \{(\alpha_{n,j},\beta_{n,j},0) : n\in{\mathbb N},\; j=0,1,\dots,r, \;\alpha_{n,j},\beta_{n,j}>0\}    
\end{equation}
is called the \emph{critical set} for \eqref{eq:pseudo:basic-sys}. 
\vs

\subsection{Transversality}
Taking $\lambda(\alpha) = u(\alpha) + iv(\alpha)$ and substituting into \eqref{eq:pseudo:char-root}, we obtain the system
\begin{align*}
u^2 - v^2 + \gamma (u (e^{-u\tau_1}\cos(v\tau_1)-1) + ve^{-u\tau_1}\sin(v\tau_1)) + u a_j - \alpha b (e^{-u\tau_2}\cos(v\tau_2)-1) &= 0\\
2uv + \gamma (v(e^{-u\tau_1}\cos(v \tau_1) - 1) - u e^{-u\tau_1}\sin(v\tau_1)) + v a_j + \alpha b e^{-u\tau_2}\sin(v\tau_2) &= 0
\end{align*}
\vs
\noindent
Differentiating this with respect to $\alpha$ gives
\begin{align*}
&\begin{cases}
\begin{aligned}
&2uu' - 2vv' - \gamma u'\tau_1e^{-u\tau_1}(u\cos(v\tau_1) + v\sin(v\tau_1)) -\gamma u'\\
+ &\gamma e^{-u\tau_1}(u'\cos(v\tau_1)-uv'\sin(v\tau_1)\tau_1 + v'\sin(v\tau_1)+vv'\cos(v\tau_1)\tau_1)\\
+ &u'a_j + \alpha b \tau_2 e^{-u\tau_2} (u'\cos(v\tau_2) + v'\sin(v\tau_2)) = b(e^{-u\tau_2}\cos(v\tau_2)-1)
\end{aligned}
\end{cases}\\
\\
&\begin{cases}
\begin{aligned}
&2vu' + 2uv' - \gamma u' \tau_1 e^{-u\tau_1}(v\cos(v\tau_1) - u\sin(v\tau_1)) -\gamma v'\\
+ &\gamma e^{-u\tau_1}(v'\cos(v\tau_1) - vv'\sin(v\tau_1)\tau_1 - u'\sin(v\tau_1) - uv'\cos(v\tau_1)\tau_1)\\
+ &v'a_j + \alpha b \tau_2 e^{-u\tau_2}(v'\cos(v\tau_2) - u'\sin(v\tau_2)) = -be^{-u\tau_2}\sin(v\tau_2)\\
\end{aligned}
\end{cases}
\end{align*}
where $u' := \frac{d}{d\alpha}u(\alpha)$ and $v' := \frac{d}{d\alpha}v(\alpha)$. Since we are interested in the behavior of purely imaginary roots, we set $u = 0$ and $v = \beta$ to obtain
\begin{subequations}\label{eq:pseudo:u'v'}
\begin{align}
&\begin{cases}
\begin{aligned}
&u'(\gamma(\cos(\beta\tau_1)-1)+a_j-\gamma\beta\tau_1\sin(\beta\tau_1)+\alpha b \tau_2 \cos(\beta\tau_2))\\+ &v'(\gamma\sin(\beta\tau_1)+\gamma\beta\tau_1\cos(\beta\tau_1)+\alpha b \tau_2\sin(\beta\tau_2)-2\beta) = b(\cos(\beta\tau_2)-1)
\end{aligned}
\end{cases}\\ \nonumber
\\ 
&\begin{cases}
\begin{aligned}
&u'(-\gamma\sin(\beta\tau_1)-\gamma\beta\tau_1\cos(\beta\tau_1)-\alpha b \tau_2 \sin(\beta\tau_2) + 2\beta)
\\+&v'(\gamma(\cos(\beta\tau_1)-1)+a_j-\gamma\beta\tau_1\sin(\beta\tau_1)+\alpha b \tau_2 \cos(\beta\tau_2)) = -b\sin(\beta\tau_2)
\end{aligned}
\end{cases}
\end{align}
\end{subequations}
Put
\begin{align*}
p &:= \gamma(\cos(\beta\tau_1)-1) + a_j - \gamma\beta\tau_1\sin(\beta\tau_1) + \alpha b \tau_2 \cos(\beta\tau_2)\\
q &:= \gamma\sin(\beta\tau_1)+\gamma\beta\tau_1\cos(\beta\tau_1)+\alpha b \tau_2\sin(\beta\tau_2)-2\beta
\end{align*}
Then \eqref{eq:pseudo:u'v'} can be written as the linear system
\begin{align*}
u'p + v'q &= b(\cos(\beta\tau_2) -1)\\
-u'q + v'p &= b\sin(\beta\tau_2)
\end{align*}
which can be solved to obtain
\begin{equation}\label{eq:pseudo:pq}
\tfrac{\partial}{\partial \alpha}P_j(\alpha,i\beta)=u' = \frac{pb(\cos (\beta\tau_2) -1) + qb \sin (\beta\tau_2)}{p^2+q^2}
\end{equation}
Put
\[
{\rho_j} := pb(\cos (\beta\tau_2) -1) + qb \sin (\beta\tau_2)
\]
\vs
This is sufficient to state and prove the following local bifurcation theorem:
\begin{theorem}\label{thm:pseudo:dde-local}
Let $(\alpha_0,\beta_0)\in \Lambda$. If $\rho_j(\alpha_0,\beta_0)\neq 0$ for some $j=0,1,\dots,r$, then there is a continuous branch of non-constant periodic solutions bifurcating from the trivial branch at $\alpha = \alpha_0$. 
\end{theorem}
\begin{proof}
If $\rho_j(\alpha_0,\beta_0)\neq 0$ for some $j$, then by Lemma \ref{lem:prelim:sign-crossing-num}, the $\Gamma_0$-isotypic crossing number $\mathfrak t_j(\alpha_0,\beta_0)\neq 0$, which implies that the local bifurcation invariant $\omega_G(\alpha_0,\beta_0)$ (applied to the two-parameter period-normalized problem formualated in the next sections) is non-trivial, which guarantees the local bifurcation of a continuum of non-constant periodic solutions by the equivariant Krasnosel'skii theorem. This also suffices to prove Theorem \ref{thm:pseudo:mas-local}.
\end{proof}

Clearly, the sign of ${\rho_j}$ determines the sign of $\tfrac{\partial}{\partial \alpha}P_j(\alpha,i\beta)$. For the purposes of proving the unboundedness of branches of solutions, we are mostly interested in determining the sign of $\tfrac{\partial}{\partial \alpha}P_j(\alpha,i\beta)$ for very large critical values of $\alpha$. 
 
This leads to the following proposition:
\begin{proposition}\label{prop:pseudo:crossing-num}
Let $P_j(\alpha,\lambda)$ be defined as in \eqref{eq:pseudo:char-root} and put $\lambda(\alpha) = u(\alpha)+iv(\alpha)$. If 
\[
|\gamma(\tau_1 - \tau_2)| < |\tau_2(\gamma - a_j)-2|
\]
then
\[
\sum_{n=1}^{\infty}\sign\big(\tfrac{\partial}{\partial \alpha}P_j(\alpha_{n,j},i\beta_{n,j})\big) \not=0
\]
\end{proposition}
\begin{proof}
The transcendental coincidence problem \eqref{eq:pseudo:beta-relation} implies that $\beta_{n,j}$ is monotone increasing. Since the numerator of \eqref{eq:pseudo:beta-relation} is bounded, large values of $\beta_{n,j}$ occur when $\sin(\beta_{n,j}\tau_2)$ is very small, i.e. when $\beta_{n,j} = \frac{2\pi l}{\tau_2}+\epsilon$ for some $l\in{\mathbb N}$ and $0<|\epsilon|\ll1$. This allows us to use the small angle approximations to analyze the proportional relationships between the quantities which make up ${\rho_j}$. Indeed, for $\beta_{n,j}$ sufficiently large, we have
\begin{align*}
&\sin(\beta_{n,j}\tau_2) \approx \epsilon\tau_2 &\quad &\cos(\beta_{n,j}\tau_2) \approx 1-\frac{1}{2}(\epsilon\tau_2)^2\\
&\sin(\beta_{n,j}\tau_1) = \sin(2\pi l\tfrac{\tau_1}{\tau_2}+\epsilon\tau_1) &\quad &\cos(\beta_{n,j}\tau_1) = \cos(2\pi l\tfrac{\tau_1}{\tau_2}+\epsilon\tau_1)
\end{align*}
Since $|\frac{1}{2}(\epsilon\tau_2)^2| \ll |\epsilon\tau_2|$ for sufficiently small $\epsilon$, we have
\[
\begin{aligned}
\sign({\rho_j}) &= \sign(qb\epsilon\tau_2)\\
&=\sign\big(b\epsilon\tau_2(\gamma\sin(2\pi l\tfrac{\tau_1}{\tau_2}+\epsilon\tau_1) + \gamma\beta\tau_1\cos(2\pi l\tfrac{\tau_1}{\tau_2}+\epsilon\tau_1) + \alpha b\epsilon\tau_2^2 - 2\beta)\big).
\end{aligned}
\]
Substituting \eqref{eq:pseudo:alpha-relation} into the second equation:
\begin{align*}
\sign(b\epsilon\tau_2(\gamma\sin(2\pi l\tfrac{\tau_1}{\tau_2}+\epsilon\tau_1) + \gamma\beta(\tau_1-\tau_2)\cos(2\pi l\tfrac{\tau_1}{\tau_2}+\epsilon\tau_1)+\beta\tau_2(\gamma-a_j) - 2\beta)).
\end{align*}
Since $\beta$ is assumed to be very large (and in particular $\beta \gg \gamma$), we need only consider terms with a factor of $\beta$. Therefore, we see that the expression
\[
\gamma(\tau_1-\tau_2)\cos(2\pi l \tfrac{\tau_1}{\tau_2} + \epsilon\tau_1) +\tau_2(\gamma-a_j) - 2
\]
determines the sign of $\tfrac{\partial}{\partial \alpha}P_j(\alpha_{n,j},i\beta_{n,j})$ for all $\alpha_{n,j}$ sufficiently large. Therefore by the boundedness of $\cos(2\pi l\tfrac{\tau_1}{\tau_2}+\epsilon\tau_1)$, if 
\[
|\gamma(\tau_1 - \tau_2)| < |\tau_2(\gamma - a_j) -2|
\]
then $\tfrac{\partial}{\partial \alpha}P_j(\alpha_{n,j},i\beta_{n,j})$ must have a constant sign for all critical $\alpha_{n,j}$ sufficiently large.
\end{proof}
\vs
We also wish to find conditions on the system and on the bifurcation parameters which guarantee local asymptotic stability of the trivial solution (i.e. consensus). This requires us to analyze \eqref{eq:pseudo:char-op} for a given value of $\alpha$, and show that all roots have negative real part. Due to the infinitude of roots of \eqref{eq:pseudo:char-op}, this can be difficult for arbitrary values of $\alpha$. 
\vs
However, if we can find the first bifurcation point $\alpha_0$ where the trivial solution changes stability, and if we can find the asymptotic stability conditions when $\alpha=0$, then since the spectrum of \eqref{eq:pseudo:basic-sys-lin} depends continuously on $\alpha$, this stability will hold for all $\alpha\in(0,\alpha_0)$. Since \eqref{eq:pseudo:alpha-relation} has infinitely many discrete roots, we can take the first positive root as the first bifurcation point, provided that there is no smaller value of $\alpha$ which admits a zero eigenvalue. This leads to the following proposition:
\vs

\begin{proposition}\label{prop:pseudo:asymptotic-stability}
Consider the system \eqref{eq:pseudo:basic-sys} satisfying \ref{c1}-\ref{c3}. Denote by $\alpha_0$ the smallest positive value of $\alpha$ satisfying \eqref{eq:pseudo:alpha-relation} for any $j$. If $0<\gamma\tau_1 \leq 1$, $b>0$, and $a_j>0$ for all $j$, then the trivial solution is locally asymptotically stable for all $\alpha \in (0,\alpha_0)$, and there is no steady-state bifurcation for any $\alpha>0$.
\end{proposition}

\begin{proof}
Take equation \eqref{eq:pseudo:char-op-iso} and set $\alpha = 0$. This yields
\[
\lambda + \gamma(e^{-\lambda \tau_1}-1) + a_j = 0.
\]
Substituting $\lambda = u+iv$ and separating real and imaginary parts obtains the system
\begin{align*}
    &u + \gamma(e^{-u\tau_1}\cos(v\tau_1)-1) +a_j = 0\\
    &v-\gamma e^{-u\tau_1}\sin(v\tau_1) = 0
\end{align*}
Obviously any solution to this system must satisfy the second equation, which is equivalent to the coincidence problem
\[
v=\gamma e^{-u\tau_1}\sin(v\tau_1).
\]
Clearly $v=0$ satisfies this equation. By the mean value theorem, this coincidence problem only admits solutions for $v>0$ if $1<\gamma \tau_1 e^{-u\tau_1}$, which implies $e^{u\tau_1}<\gamma\tau_1$, and so $u< \ln(\gamma\tau_1)/\tau_1$. Therefore if $\ln(\gamma\tau_1) \leq 0$, $u < 0$, and so $0\leq\gamma\tau_1 \leq 1$ suffices to guarantee that the coincidence problem has no nonzero solutions. 
\vs
When $v=0$, the first equation becomes
\[
u +\gamma e^{-u\tau_1}-\gamma + a_j = 0.
\]
Let $r(u) := u +\gamma e^{-u\tau_1}-\gamma$. Then $r(0)=0$. Since $\tfrac{d^2r}{du^2} = \tau_1^2 \gamma e^{-u\tau_1}\geq0$ for all $u$,  $r(u)$ is a convex function. We also have $\tfrac{dr}{du}>0$ at $u=0$ if $\gamma\tau_1\leq1$. This means that $r(u)$ takes a single minimum value which must be in the $u<0$ half-plane. So if $r(u)$ has another nonzero root, then this root must be negative. Therefore, if $a_j >0$, then the upwards translation of the graph of $r(u)$ by $a_j$ must cause the root at $u=0$ to move into the negative half plane. Therefore, real roots of $r(u)+a_j=0$ (if they exist) must all be negative.
\vs
Finally, we must address the possibility of steady-state bifurcation for $\alpha \in (0,\alpha_0)$. This corresponds to a purely real eigenvalue crossing zero and becoming positive. A value of $\alpha$ where this happens would not appear as a solution to \eqref{eq:pseudo:alpha-relation}, as this relation assumes that $\lambda = i\beta$ with $\beta>0$. 
\vs
Taking equation \eqref{eq:pseudo:char-op-iso}, substituting $\lambda = u$, and moving the term proportional to $\alpha$ to the right hand side, we obtain the coincidence problem
\[
r(u) = \frac{\alpha b(e^{-u\tau_2} -1)}{u}.
\]
Therefore, the existence of a real positive eigenvalue for some $\alpha>0$ corresponds to a real positive solution to this coincidence problem. By the continuity of eigenvalues as functions of $\alpha$, if we can show that no such solutions exist, then we have demonstrated that real eigenvalues cannot cross zero and that steady state bifurcation is therefore impossible (thus also justifying our consideration of the first $\alpha_0$ where Hopf bifurcation occurs as the boundary of the parameter range where the trivial solution is locally asymptotically stable).
\vs
Since $r(u)$ is convex and both roots (if they exist) are negative, this means that $r(u)>0$ for all $u > 0$. On the other hand, if $u>0$ then $e^{-u\tau_2}-1<0$. This implies that if $b>0$, then $\alpha b(e^{-u\tau_2}-1)/u < 0$ for all $u>0$. Therefore, the graphs of $r(u)$ and $\alpha b(e^{-u\tau_2}-1)/u$ cannot intersect on the $u\geq0$ half-plane, and so steady-state bifurcation is impossible for all $\alpha>0$.
\end{proof}
\vs
Note that this suffices to prove Theorem \ref{thm:pseudo:mas-asymp}.

\begin{remark}\normalfont
The requirements that $0 \leq\gamma\tau_1\leq 1$, $a_j>0$, and $b>0$ have very natural interpretations in our system, viewed both as a general NFDE dynamical system and as a MAS. The term $\gamma\tau_1$ can be viewed as a type of upper bound on the feedback introduced by the neutral delay term near the trivial solution, viewed proportionally to $x(t)$. If $\gamma\tau_1>1$, this means that the neutral delay term is not dissipative and can cause unbounded growth of the derivative, or high frequency instabilities. This also corresponds to the $\kappa$-Lipschitz condition imposed on $G$ through condition \ref{c3}. 
\vs
In a MAS context, this suggests that an agent's memory of its own past derivative must not have a higher weight than the contribution of the rest of its state function (including its own memories of past states). The requirement that $a_j>0$ corresponds to an attractive coupling between oscillators, or a cooperative agent interaction in the MAS context. This reflects the intuition that consensus should reflect a fundamentally cooperative relationship between agents, which may destabilized by the weight of continuous state memory. Finally, $b>0$ simply ensures that the continuous memory term should promote a return to consensus, otherwise the consensus solution would be fundamentally unstable. 
    
\end{remark}

\vs
\section{Functional space reformulation}\label{sec:pseudo:funcspace}

In this section, we reformulate \eqref{eq:pseudo:basic-sys} as a fixed point problem in functional spaces, and show that it is a condensing perturbation of identity. We linearize this operator and describe the spectrum of this linearization in terms of the spectrum of \eqref{eq:pseudo:basic-sys-lin} by leveraging the relationship between $G$-isotypic components of $\mathscr E$ and $\Gamma_0$-isotypic components of $V$. (Note that in this section and for the remainder of the paper, we will use $u(t)$ to refer to a period-normalized function, not to be confused with its usage as the real part of an eigenvalue in the previous section.) 
\vs
As with \eqref{eq:pseudo:fde-norm}, we set $u(t):=\bm x(\tfrac{p}{2\pi}t)$ and $\beta = \tfrac{2\pi}{p}$ and substitute into system \eqref{eq:pseudo:basic-sys} to obtain the period-normalized system
\begin{equation}\label{eq:pseudo:basic-sys-norm}
\frac{d}{dt}\left[u- \int_0^{\tau_1} \bm g(u(t-\beta s))ds\right] = -\frac{a}{\beta}u - \frac{\alpha}{\beta}\bm f\left(\int_0^{\tau_2} u(t-\beta s)ds\right) - \frac{1}{\beta}\bm h(u) 
\end{equation}

We define the space $\mathscr E$ and the operators $L, \bm j, N_{\bm f}$ and $N_{\bm k}$ as in Section \ref{sec:pseudo:introduction:overview} and put
\begin{align*}
N_{\bm k}(\alpha,\beta,u) &= \int_0^{\tau_1}\bm g(u(t-\beta s))ds\\
N_{\bm f}(\alpha,\beta,u) &= -\frac{\alpha}{\beta}u - \frac{\alpha}{\beta}\bm f\left(\int_0^{\tau_2} u(t-\beta s)ds\right) - \frac{1}{\beta}\bm h(u) 
\end{align*}
\vs
yielding
\begin{equation}\label{eq:pseudo:op-sys}
\mathscr F(\alpha,\beta,u) = u - N_{\bm k}(\alpha,\beta,u) - (L + \bm j)^{-1}\Big(N_{\bm f}(\alpha,\beta,u) + j(u) - N_{\bm k}(\alpha,\beta,\bm j(u))\Big)
\end{equation}
The assumptions \ref{c1}--\ref{c3} guarantee that $\mathscr F$ is $G$-equivariant. Since $\bm f$ is continuous and the delay operator $\int_0^{\tau_2}u(t-\beta s)ds$ is compact, their composition is compact, and so the Nemytskii operator $N_{\bm f}$ is compact. Because $\bm g$ is $\kappa$-Lipschitzian with $\kappa<1$, $N_{\bm k}$ is condensing. Since the sum of a condensing operator and a compact operator is condensing, and because $L+j$ is an isomorphism, $\mathscr{F}:\bbR^2_+ \times \mathscr{E} \to \mathscr{E}$ is a $G$-equivariant condensing perturbation of identity, and thus can be analyzed using the ETNS degree.

\vs
Following this reformulation, the problem of finding periodic solutions to \eqref{eq:pseudo:basic-sys} is now equivalent to finding $(\alpha,\beta,u) \in \bbR_+^2 \times \mathscr{E}$ such that $\mathscr{F}(\alpha,\beta,u) = 0$. $\mathscr F$ is also differentiable at every $(\alpha,\beta,0) \in \bbR_+^2 \times \mathscr E$. Put $\mathscr A(\alpha,\beta) := D_u \mathscr F(\alpha,\beta,0)$. To simplify notation, we also define a distributed delay operator:
\begin{align*}
K_{\tau}u(t) &:= \int_0^{\tau}u(t-\beta s)ds
\end{align*}
Then we have
\[
\begin{aligned}
\mathscr A(\alpha,\beta)u &= u - \gamma K_{\tau_1}u - (L+\bm j)^{-1}\left(-\frac{\alpha b}{\beta}K_{\tau_2}u - \frac{1}{\beta}Cu - \Big(\frac{a}{\beta} - 1\Big)u- \gamma K_{\tau_1}u\right)  \\
&= u - \gamma K_{\tau_1}u + (L+\bm j)^{-1}\left(\frac{\alpha b}{\beta}K_{\tau_2}u + \frac{1}{\beta}Cu + \Big(\frac{a}{\beta} - 1\Big)u+ \gamma K_{\tau_1}u\right).
\end{aligned}
\]
The $G$-equivariance of $\mathscr F$ implies that $\mathscr A$ is also $G$-equivariant. Accordingly, its eigenspaces correspond to $G$-isotypic components, and we define the restricted maps $\mathscr A_{k,j}(\alpha,\beta) := \mathscr A_{|\mathscr E_{k,j}}(\alpha,\beta):\mathscr E_{k,j} \to \mathscr E_{k,j}$. 
\vs
Note that the characteristic operator equation of the period normalized system \eqref{eq:pseudo:basic-sys-norm} is given by
\[
\frac{1}{\beta}\triangle_{\alpha}(\lambda)
\]
and the characteristic operator equation for $(L+\bm j)$ is given by
\[
(L+\bm j)e^{\lambda t} = (\lambda + 1)e^{\lambda t}
\]
Combining this with the fact that $u(t) = \bm x(\frac{t}{\beta})$ and taking $\lambda = ik\beta$, we have
\begin{align}
\mathscr A_{k,j}(\alpha,\beta) &= \frac{1}{ik\beta + \beta}\triangle_{\alpha,j}(ik\beta)_{|\mathscr E_{k,j}}\quad &k>0\label{eq:pseudo:char-Ak}\\
\mathscr A_{0,j}(\alpha,\beta) &= \left(\frac{\alpha b \tau_2 +a_j }{\beta}\right)_{|\mathscr E_{0,j}}\quad &k=0\label{eq:pseudo:char-A0}
\end{align}
\section{Two-parameter bifurcation}\label{sec:pseudo:two-parameter}
Here we will define the required $G$-homotopy invariant quantities (namely the crossing number and the local bifurcation invariant) and show how they are computed. Using the equivariant versions of the existence theorem of M.A. Krasnosel'skii and the unbounded continuation argument of the Rabinowitz alternative, we will prove Hopf bifurcation of global branches and conditions for their unboundedness in the absence of degeneracies, i.e. in the case that eigenvalue crossings are $G$-isotypically simple. 
\vs
We can view the system \eqref{eq:pseudo:basic-sys-norm} as a two-parameter bifurcation problem given by
\[
\mathscr F(\alpha,\beta,0) = 0, \quad (\alpha,\beta,u) \in \bbR^2_+ \times \mathscr E
\]
Indeed, one can easily see that from the functional space reformulation shown above that the following facts hold:
\begin{enumerate}
    \item $\mathscr F$ is a condensing $G$-equivariant field.
    \item For all $(\alpha,\beta) \in \bbR^2_+$, $\mathscr F(\alpha,\beta,0)=0$. 
    \item For all $(\alpha,\beta) \in \bbR^2_+$, the derivative $\mathscr A(\alpha,\beta) := D_u \mathscr F(\alpha,\beta,u)$ exists, depends continuously on $(\alpha,\beta)$, and for any $(\alpha_0,\beta_0,u_0) \in \bbR^2_+ \times \mathscr E$ we have
    \[
    \lim_{(\alpha_0,\beta_0,u_0) \to (\alpha_0,\beta_0,0)}\frac{||\mathscr F(\alpha_0,\beta_0,u_0) - \mathscr A(\alpha_0,\beta_0)u_0||}{||u_0||} = 0
    \]
\end{enumerate}
We now refer to the definitions of branches, branching points, and bifurcation points given in Section \ref{sec:prelim:krasnoselskii}. From these definitions, we obtain a necessary condition for $(\alpha_0,\beta_0,0)\in \mathscr M$ to be a bifurcation point. Namely, if $(\alpha_0,\beta_0,0)$ is a branching point, then $\mathscr A(\alpha_0,\beta_0):\mathscr E \to \mathscr E$ is not an isomorphism. Equations \eqref{eq:pseudo:char-Ak} and \eqref{eq:pseudo:char-A0} imply that this requires either $\det_{\mathbb C} \triangle_{\alpha_0}(ik\beta_0) = 0$ for some $k\in {\mathbb N}$, or $\tfrac{\alpha_0 b \tau_2 + a_j}{\beta_0} = 0$. This provides the connection between the critical set for the system \eqref{eq:pseudo:basic-sys}, given by \eqref{eq:pseudo:crit-set}, and the critical set of the period normalized system \eqref{eq:pseudo:basic-sys-norm}. 
\vs
However, if $\tfrac{\alpha_0 b \tau_2 + a_j}{\beta_0} = 0$, this corresponds to a steady-state bifurcation. Since we are interested in showing the bifurcation of non-constant periodic solutions, we view this as a degeneracy. Excluding such points leads us to the non-constant critical set for \eqref{eq:pseudo:basic-sys-norm}, denoted $\widetilde\Lambda$: 

\[
\widetilde{\Lambda} := \left\{(\alpha,\beta,0) \in \mathscr M:\;\; \exists_{k \in {\mathbb N}} \;\;{\det}_{\mathbb C} \triangle_\alpha(ik\beta) = 0 \text{ and } \frac{\alpha b \tau_2 + a_j}{\beta} \not= 0\right\}.
\]
Following the construction in Section \ref{sec:prelim:local-bif-twisted}, we can construct isolated $G$-invariant neighborhoods around each isolated point, create complemented maps, and obtain the computational formula for the local bifurcation invariant used by the equivariant Krasnosel'skii theorem:
\[
\omega_G(\alpha_0,\beta_0) = \text{\rm $\Gamma$-deg}(\mathscr A_{0}(\alpha_0,\beta_0),B(V))\cdot \sum_{k=1}^\infty\sum_{j=0}^r\mathfrak t_{k,j}(\alpha_0,\beta_0)\text{\rm deg}_{\cV_{k,j}^-}.
\]
Since $\widetilde \Lambda$ filters out points where steady-state bifurcation occurs, we can write this in a simpler way as
\[
\omega_G(\alpha_0,\beta_0) =  \sum_{k=1}^\infty\sum_{j=0}^r\mathfrak t_{k,j}(\alpha_0,\beta_0)\text{\rm deg}_{\cV_{k,j}^-}.
\]
Recall the equivariant Krasnosel'skii theorem guarantees local bifurcation if this local bifurcation invariant has some twisted orbit type with a non-trivial coefficient. Since these coefficients are given in terms of crossing numbers, and crossing numbers are given in terms of the transversality indicator function $\rho_j(\alpha,\beta)$, this is enough to establish local bifurcation.
\vs
For global bifurcation, the equivariant Rabinowitz alternative, as formulated in Section \ref{sec:prelim:two-parameter-bifurcation}, is also sufficient to guarantee the existence of branches which are unbounded in $\mathbb R^2_+\times \mathscr E$, taken with the conditions we have already obtained on the sum of $\rho_j(\alpha_{n,j},\beta_{n,j})$. However, this says nothing about their symmetries across this full global extent. In other words, the Rabinowitz theorem on its own is not sufficient to guarantee that branches of solutions which may be non-constant and periodic at the point of bifurcation will not collapse to constant solutions or undergo some secondary bifurcation which destroys their interesting symmetries at some point in their global extent (or maybe even very soon after bifurcation). To do this, we must restrict the problem to a fixed-point subspace which cannot contain constant solutions. Applying the Rabinowitz theorem in this space, we obtain unbounded global branches of non-constant periodic solutions which must remain non-constant and periodic across their full extent (although they could experience other nonlocal symmetry-breaking or period-doubling bifurcations), and these are guaranteed to be solutions to the original problem as well.

\section {Fixed point reduction}\label{sec:pseudo:fixed-point-reduction}

\vs

Fix $\kappa \in {\mathbb N}$ and consider the subgroup of $G$ generated by 
\[
\bm K := \left<\left(\bm e,-1,e^{\tfrac{i\pi}{\kappa}}\right)\right> \leq \Gamma_0 \times \bbZ_2 \times S^1
\]
We then examine \eqref{eq:pseudo:op-sys} in the fixed point space $\mathscr E^{\bm K}$, i.e. we are interested in solutions to the system
\begin{equation}\label{eq:pseudo:func-sys-fixedpoint}
\mathscr F^{\bm K}(\alpha,\beta,u) = 0,\quad (\alpha,\beta,u)\in \bbR^2_+ \times \mathscr E^{\bm K}    
\end{equation}
We can see that any solution to \eqref{eq:pseudo:func-sys-fixedpoint} must be an odd function and is also a solution to \eqref{eq:pseudo:op-sys}, and the isotypic decomposition of $\mathscr E^{\bm K}$ is given by
\[
\mathscr E^{\bm K}_{k,j} := \begin{cases}
    \mathscr E_{k,j} \quad &\text{ if $k=(2l-1)\kappa$ for some $l\in{\mathbb N}$}\\
    0 \quad &\text{ otherwise}
\end{cases} 
\]
which immediately yields
\[
\mathscr E^{\bm K} = \overline{\bigoplus_{l=1}^\infty \bigoplus_{j=0}^r \mathscr E_{(2l-1)\kappa,j}}
\]
Clearly $\bm K$ is normal in $G$, and since $G_0 := G/\bm K = S^1 \times \Gamma_0$, the above system is a two-parameter $S^1 \times \Gamma_0$-equivariant bifurcation problem, and the set of critical points for \eqref{eq:pseudo:func-sys-fixedpoint} can be described as follows:
\[
\widetilde{\Lambda}^{\bm K}:=\{(\alpha_0,\beta_0,0) \in \bbR^2_+ \times \mathscr E: \exists l \in {\mathbb N}\;\;{\det}_{\mathbb C}(\triangle_{\alpha_0}(i(2l - 1)\kappa\beta_0))=0\}
\]
We also define
\[
\widetilde{\Lambda}_k^{\bm K}:=\{(\alpha_0,\beta_0,0) \in \bbR^2_+ \times \mathscr E: \exists j\in {\mathbb N} \;\;{\det}_{\mathbb C}(\triangle_{\alpha_0,j}(i(2k - 1)\kappa\beta_0))=0\}
\]
and so
\[
\widetilde\Lambda^{\bm K} = \bigcup_{k=1}^\infty \widetilde\Lambda_k^{\bm K}
\]
\vs
This set is clearly infinite and all critical points in it are isolated, and the signs of crossing numbers are as given in Proposition \ref{lem:prelim:sign-crossing-num}.
\vs
Thus we obtain the following formula for the local bifurcation invariant in the $\bm K$-fixed point reduction:
\[
\sum_{i=1}^N \omega_G(\alpha_i,\beta_i) = \sum_{i=1}^N\sum_{j=0}^r \sum_{l=1}^\infty \mathfrak t_{(2l-1)\kappa,j}(\alpha_i,\beta_i) \text{deg}_{\cV_{(2l-1)\kappa,j}}
\]

Now we can present the main global existence theorem. It is derived from the Rabinowitz alternative and proved for completely continuous maps in \cite{BalanovEtAl2025}, and here is trivially extended to condensing maps.
\begin{theorem}
Let $G=\Gamma_0\times\bbZ_2\times S^1$ and let $\mathscr F:\bbR_+^2\times \mathscr E \to \mathscr E$ be a $G$-equivariant condensing map satisfying assumptions \ref{c1} -- \ref{c3}. If $\widetilde\Lambda_k^{\bm K}$ is discrete and finite for some $k=1,2,\dots$, and for some orbit type $(H)\in\Phi^t_1(G;\cV^-_{(2k-1)\kappa,j})$ we have
\[
\sum_{(\alpha_0,\beta_0,0)\in\widetilde\Lambda_k^{\bm K}} \text{coeff}^H(\omega_G(\alpha_0,\beta_0)) \not=0
\]
where $\text{coeff}^H$ stands for the coefficient of $(H)$, then there exists an unbounded global branch $\mathscr C'\subset \mathscr S^H$ with $\mathscr C'\cap\widetilde\Lambda_k^{\bm K} \not=\emptyset$.
\end{theorem}
The proof of this theorem follows exactly from the proof in \cite{BalanovEtAl2025}. It is also worth mentioning that, if $(\alpha_0,\beta_0,0)\in \widetilde\Lambda_k^{\bm K}$, then $(\alpha_0,(2k-1)\kappa\beta_0,0)\in \widetilde\Lambda\subset \Lambda$. In other words, if a maximal orbit type is detected on a higher mode, then there is a corresponding orbit type on the first mode. Therefore, the fixed point reduction does not impede our detection of solutions, and we are justified in only considering the first mode.
This theorem, which permits us to use the local and global bifurcation results even in the presence of the aforementioned degeneracies and guarantees the detected periodic solutions are non-constant, taken together with Prop \ref{prop:pseudo:crossing-num} and Lemma \ref{lem:prelim:sign-crossing-num}, conclude the proof of Theorem \ref{thm:pseudo:mas-local}, Theorem \ref{thm:pseudo:mas-global}, and Theorem \ref{thm:pseudo:mas-sym}.

\section{Application: coupled asset markets with fundamentalist and chartist traders}\label{sec:pseudo:example}
Heterogeneous agent models (HAMs) have been very popular in modeling market dynamics. They consider the effects of different populations of traders with different trading strategies. This framework was first established by Beja and Goldman \cite{beja1980price}, and has been extended in numerous ways into continuous time delay systems. Indeed, there are many models in the HAM framework whose closed-loop dynamic equations closely resemble the equations studied in this dissertation (cf. \cite{chiarella1992dynamics,dibeh2007dynamics,he2010dynamics,matsumoto2016heterogeneous,dobrescu2016asset,guerrini2014heterogeneous}). There has been a great deal of sustained interest in these types of models, because they are good models of a very fundamental economic question: On one hand, the competing forces of supply and demand in the free market should result in price discovery, where an asset has a single stable value within a given market (not including outside forces or pressures). On the other hand, we regularly see certain boom-bust cycles in real market data, and oscillatory price movements at many different time scales.
\vs
We must also make a certain distinction of terminology. In HAM literature, the ``agents'' in question are the different traders or other forces within the market, and they are heterogeneous because they use different trading strategies. However, throughout this dissertation, we have focused on multi-agent systems with homogeneous agents with identical protocols (up to symmetry). In this context, we consider the market itself to be the agent, and its ``goal'' is price discovery. Its protocol consists of components representing the market maker, fundamentalist traders, and chartist traders, and a coupling term which captures the fact that markets in related sectors can exert influences on each other. Therefore, we can think of this system as a kind of MAS of HAMs. We will uses the term ``agent'' to refer to the markets themselves, to stay consistent with previous chapters, and we will refer to the different populations of traders within each market simply as ``traders'', but it should be understood that this example could also be formulated purely within a HAM context with some simple differences in terms and labels.
\vs
\begin{figure}[tbp]
\centering
\begin{minipage}{0.4\textwidth}
% \resizebox{\linewidth}{!}{%    
\centering
\begin{tikzpicture}[
    scale=2.5,
    x={(1cm,0cm)},          % x direction → right
    y={(0cm,1cm)},          % y direction → up
    z={(-0.5cm,-0.5cm)},    % z direction → back/left
    vertex/.style={circle, draw=black, fill=black, inner sep=1.5pt},
    edge/.style={solid, thick, black},
    facediag/.style={dashed, thick, black},
    spacediag/.style={dotted, thick, black},
    every label/.style={font=\tiny},  % optional: style for all labels
]

% Define cube vertices in desired order
\coordinate (x1) at (0,1,1);   % top left back
\coordinate (x2) at (1,1,1);   % top right back
\coordinate (x3) at (0,1,0);   % top left front
\coordinate (x4) at (1,1,0);   % top right front
\coordinate (x5) at (0,0,1);   % bottom left back
\coordinate (x6) at (1,0,1);   % bottom right back
\coordinate (x7) at (0,0,0);   % bottom left front
\coordinate (x8) at (1,0,0);   % bottom right front

% Draw vertices (just the dots)
\foreach \i in {1,...,8}
    \node[vertex] at (x\i) {};

% Add labels separately
% Top face: labels above
\node[above=3pt] at (x1) {$x_3$};
\node[above=3pt] at (x2) {$x_4$};
\node[above=3pt] at (x3) {$x_1$};
\node[above=3pt] at (x4) {$x_2$};
% Bottom face: labels below
\node[below=3pt] at (x5) {$x_7$};
\node[below=3pt] at (x6) {$x_8$};
\node[below=3pt] at (x7) {$x_5$};
\node[below=3pt] at (x8) {$x_6$};

% Draw edges (solid)
\foreach \i/\j in {1/2, 1/3, 1/5, 2/4, 2/6, 3/4, 3/7, 4/8, 5/6, 5/7, 6/8, 7/8}
    \draw[edge] (x\i) -- (x\j);

% Draw face diagonals (dashed)
\foreach \i/\j in {1/4, 2/3, 5/8, 6/7, 3/8, 4/7, 1/6, 2/5, 1/7, 3/5, 2/8, 4/6}
    \draw[facediag] (x\i) -- (x\j);

% Draw space diagonals (dotted)
\foreach \i/\j in {1/8, 2/7, 3/6, 4/5}
    \draw[spacediag] (x\i) -- (x\j);

\end{tikzpicture}
% }
\end{minipage}%
\begin{minipage}{0.45\textwidth}
\centering
\small
\resizebox{\linewidth}{!}{%
    \(
    C = \begin{pmatrix}
        0 & c_1 & c_1 & c_2 & c_1 & c_2 & c_2 & c_3 \\
        c_1 & 0 & c_2 & c_1 & c_2 & c_1 & c_3 & c_2 \\
        c_1 & c_2 & 0 & c_1 & c_2 & c_3 & c_1 & c_2 \\
        c_2 & c_1 & c_1 & 0 & c_3 & c_2 & c_2 & c_1 \\
        c_1 & c_2 & c_2 & c_3 & 0 & c_1 & c_1 & c_2 \\
        c_2 & c_1 & c_3 & c_2 & c_1 & 0 & c_2 & c_1 \\
        c_2 & c_3 & c_1 & c_2 & c_1 & c_2 & 0 & c_1 \\
        c_3 & c_2 & c_2 & c_1 & c_2 & c_1 & c_1 & 0
    \end{pmatrix}
    \)%
    }
\end{minipage}
\caption{Cube-coupled system with coupling matrix $C$. Solid, dashed, and dotted lines indicate interactions with coupling strength given by $c_1,c_2$, and $c_3$ respectively.}\label{fig:pseudo:cube}
\end{figure}
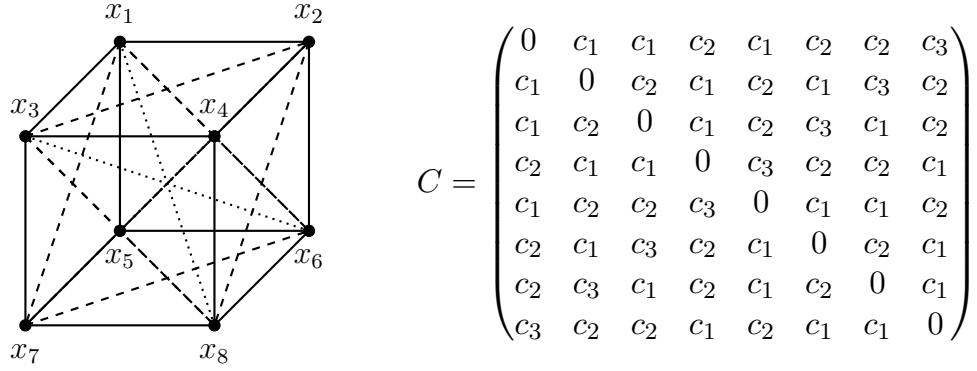
We consider a network of 8 coupled homogeneous identical HAMs. Put $\bm x(t) \in \mathbb R^8$, and denote by $x_i(t)$ the difference of the price of the $i$-th asset from its fundamental value. Then if $x_i(t)>0$, it is overvalued, and if $x_i(t)<0$, it is undervalued. We consider each market as being populated by two types of trader: fundamentalist traders, who react to a weighted average of asset price deviations over a finite time horizon $\tau_1$, and chartists, who react to a weighted average of the rate of change of asset prices over some finite time horizon $\tau_2$. Then the state dynamics of the $i$-th asset market are given by: 
\[
\frac{d}{dt}\left[x_i - \int_0^{\tau_1}\tanh(\gamma x_i(t-s))ds\right]=-ax_i - \alpha \tanh\left(b\int_0^{\tau_2}x_i(t-s)ds\right) - h(x) \]
where the $-ax_i$ term represents a market maker or liquidity provider which reacts instantly to mispricing and forces the asset's price towards its fundamental value, $\alpha$ represents the strength or aggressiveness of fundamentalist traders, and $h(x) := \tanh\left(\sum_{j=1}^8c_{ij}x_j\right)$ represents saturable coupling between markets (i.e. cross-market arbitrage), and the saturation function $\tanh$ is chosen to represent saturable liquidity.

\vs
In this case, the distributed delay terms may also be thought of as smoothing out the noise or high-frequency price fluctuations often seen in real markets. This forms a multi-agent system where the goal of each agent is price discovery. The existence of consensus at $x\equiv0$ has a natural interpretation as representing an efficient, arbitrage-free market, where all 8 assets are priced at their true value. Multiconsensus represents situations where different clusters of assets are undergoing boom and bust cycles, which could also be viewed as the formation and collapse of pricing bubbles. 
\vs
Different assets or markets can be identified with certain characteristics, such as their region (e.g. North America vs. Europe), their sector (e.g. technology vs. energy), or their size (e.g. high market capitalization vs. low market capitalization). Markets with more of these characteristics in common tend to be more strongly correlated. We will consider three different properties which take binary values, $(p_1,p_2,p_3)\in \{0,1\}^3$. If two markets share two properties and differ in one, we will give them a higher amount of correlation. If they share one property and differ in two, they get a medium amount of correlation, and if they differ in every property, they get the lowest amount of correlation. Arranging the $2\times 2\times 2=8$ possibilities as the vertices of a cube, we can use an edge to represent sharing two properties, a face diagonal to represent sharing one property, and a space diagonal to represent sharing no properties (cf. Figure \ref{fig:pseudo:cube}).
\vs
To choose parameters, we will suppose that chartists look at trends over the past month, $\tau_1 = 20$, and fundamentalists over the past quarter, $\tau_2=60$. We will set $a = 0.5$, $b =0.2$, $\gamma = 0.04$, $c_1 = 0.15$, $c_2 = 0.05$, and $c_3 = 0.01$, and $L_{\text{sat}} = 2$. This leads to the following linearized system:
\begin{equation}\label{eq:pseudo:example-sys-lin}
\frac{d}{dt} \left[x - 0.2 \int_0^{20} x(t-s)ds\right] + 0.5x + 0.2\alpha \int_0^{60} x(t-s)ds + C = 0
\end{equation}
where
\[
C = \begin{bmatrix}
0 & 0.15 & 0.15 & 0.05 & 0.15 & 0.05 & 0.05 & 0.01 \\
0.15 & 0 & 0.05 & 0.15 & 0.05 & 0.15 & 0.01 & 0.05 \\
0.15 & 0.05 & 0 & 0.15 & 0.05 & 0.01 & 0.15 & 0.05 \\
0.05 & 0.15 & 0.15 & 0 & 0.01 & 0.05 & 0.05 & 0.15 \\
0.15 & 0.05 & 0.05 & 0.01 & 0 & 0.15 & 0.15 & 0.05 \\
0.05 & 0.15 & 0.01 & 0.05 & 0.15 & 0 & 0.05 & 0.15 \\
0.05 & 0.01 & 0.15 & 0.05 & 0.15 & 0.05 & 0 & 0.15 \\
0.01 & 0.05 & 0.05 & 0.15 & 0.05 & 0.15 & 0.15 & 0
\end{bmatrix},
\]
We note that system $\eqref{eq:pseudo:example-sys-lin}$ satisfies assumptions \ref{c1} -- \ref{c3}. Before taking its isotypic decomposition to apply the main result, we will first make some brief remarks on the group theory of the symmetry group of rigid motions of a cube/octahedron, denoted $\mathbb O$.
\vs
First we note that since $C$ can be written as a weighted sum of adjacency matrices of three undirected cubic graphs corresponding to cubic edges, face diagonals, and space diagonals, respectively. Therefore, $C$ is clearly $\mathbb O$-equivariant. To better describe the action of $\mathbb O$ on our 8 markets, we will consider $\mathbb O \cong S_4 \leq S_8$, the full symmetry group of the 8 vertices. $\mathbb O$ has order 24 and can be generated by a rotation of order 4 about an axis connecting the centers of two opposite faces, a rotation of order 3 about an axis connecting two space diagonal opposite vertices, and a rotation of order 2 about an axis between the midpoints of two opposite edges. We note the following correspondence between elements of $\mathbb O$ (written as permutations in $S_4$) and elements in $S_8$:
\begin{table}[tbp]
\centering
\small  % reduce font size for the whole table
\begin{tabular}{lll}
Rotation type & $S_4$ representative & $S_8$ representative\\
\hline
Face rotation   & (1234)       & (1243)(5687) \\
Vertex rotation & (142)        & (283)(167)   \\
Edge rotation   & (34)         & (18)(27)(34)(56) \\
\end{tabular}
\caption{Generators of $S_4$ with corresponding generators in $S_8$.}
\end{table}
Let $\mathbb O$ act by permuting vertices in $V:=\mathbb R^8$. Then $V$ is an isometric $\mathbb O$-representation with the following character table:
\begin{table}[tbp]
\centering
\small
\begin{tabular}{|c|ccccc|}
\hline
con. classes & $(1)$ & $(12)$ & $(12)(34)$ & $(123)$ & $(1234)$ \\ \hline
$\chi_0$ & $1$ & $1$ & $1$ & $1$ & $1$ \\
$\chi_1$ & $1$ & $-1$ & $1$ & $1$ & $-1$ \\
$\chi_2$ & $2$ & $0$ & $2$ & $-1$ & $0$ \\
$\chi_3$ & $3$ & $-1$ & $-1$ & $0$ & $1$ \\
$\chi_4$ & $3$ & $1$ & $-1$ & $0$ & $-1$ \\ \hline
$\chi_V$ & $8$ & $0$ & $0$ & $2$ & $0$ \\ \hline
\end{tabular}
\caption{Irreducible characters of $S_4$ with character of $V$.}
\label{table:pseudo:char}
\end{table}
Let $\Gamma = \mathbb O \times \bbZ_2$. Then $V$ is also an isometric $\Gamma$-representation, and from the above character table we immediately obtain the isotypic decomposition
\[
V = V_0 \oplus V_1 \oplus V_3 \oplus V_4
\]
where each $V_i$ has isotypic multiplicity 1 and is modeled on the $\cV_i^-$ irreducible representation, i.e., the $\cV_i$ irreducible $\mathbb O$-representation equipped with the antipodal $\bbZ_2$-action. Let $G = S^1 \times \mathbb O \times \bbZ_2$. Then our system \eqref{eq:pseudo:example-sys-lin} can be reformulated as a $G$-symmetric two-parameter bifurcation problem as described in Section \ref{sec:pseudo:funcspace} and our main results can be applied. We compute the eigenspaces of $C$ by looking for $S_4$-invariant subspaces corresponding to the irreducible representations in the isotypic decomposition:
\begin{align*}
    E(\mu_0) = \text{span}(&(1,1,1,1,1,1,1,1)^T) \quad &\mu_0 &= 3c_1 +3c_2 +c_3\\
    \\
    E(\mu_1) = \text{span}(&(1,-1,-1,1,-1,1,1,-1)^T) \quad &\mu_1 &= -3c_1 + 3c_2 - c_3\\
    \\
    E(\mu_3) = \text{span}(&(1,-1,-1,1,1,-1,-1,1)^T, \quad &\mu_3 &= -c_1 - c_2 + c_3\\
    &(1,1,-1,-1,-1,-1,1,1)^T,\\
    &(1,-1,1,-1,-1,1,-1,1)^T)\\
    \\
    E(\mu_4) = \text{span}(&(1,1,1,1,-1,-1,-1,-1)^T,\quad &\mu_4 &= c_1 - c_2 - c_3\\&(1,1,-1,-1,1,1,-1,-1)^T,\\
    &(1,-1,1,-1,1,-1,1,-1)^T)\\
\end{align*}
By computing the traces of conjugacy classes of permutations in $\mathbb O$ on these eigenspaces and comparing them to Table \ref{table:pseudo:char}, one immediately finds $E(\mu_0) = V_0$, $E(\mu_1) = V_1$, $E(\mu_3) = V_3$, and $E(\mu_4) = V_4$. Therefore we have
\begin{align*}
    a_0 &= a + 3c_1 + 3c_2+c_3 &= 0.5 +0.45 + 0.15 + 0.01 &= 1.11\\
    a_1 &= a - 3c_1+3c_2-c_3 &= 0.5 -0.45 + 0.15 - 0.01 &= 0.19\\
    a_3 &= a - c_1-c_2+c_3 &= 0.5 -0.15 - 0.05 + 0.01 &= 0.31\\
    a_4 &= a + c_1-c_2-c_3 &= 0.5 +0.15 - 0.05 - 0.01 &= 0.59
\end{align*}

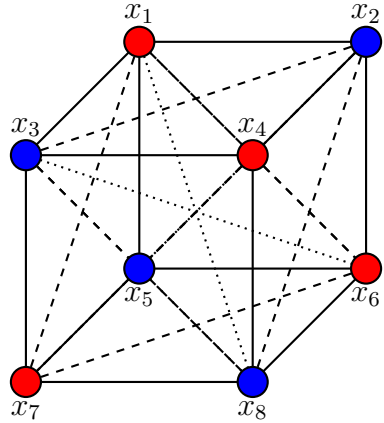
\begin{figure}
\centering
\begin{tikzpicture}[
    scale=3,
    x={(1cm,0cm)},          % x direction → right
    y={(0cm,1cm)},          % y direction → up
    z={(-0.5cm,-0.5cm)},    % z direction → back/left
    vertex/.style={circle, draw=black, fill=black, inner sep=1.5pt},
    edge/.style={solid, thick, black},
    facediag/.style={dashed, thick, black},
    spacediag/.style={dotted, thick, black},
    every label/.style={font=\tiny},  % optional: style for all labels
]

% Define cube vertices in desired order
\coordinate (x1) at (0,1,1);   % top left back
\coordinate (x2) at (1,1,1);   % top right back
\coordinate (x3) at (0,1,0);   % top left front
\coordinate (x4) at (1,1,0);   % top right front
\coordinate (x5) at (0,0,1);   % bottom left back
\coordinate (x6) at (1,0,1);   % bottom right back
\coordinate (x7) at (0,0,0);   % bottom left front
\coordinate (x8) at (1,0,0);   % bottom right front

% % Draw vertices (just the dots)
% \foreach \i in {1,...,8}
%     \node[vertex] at (x\i) {};

% Add labels separately
% Top face: labels above
\node[above=3pt] at (x1) {$x_3$};
\node[above=3pt] at (x2) {$x_4$};
\node[above=3pt] at (x3) {$x_1$};
\node[above=3pt] at (x4) {$x_2$};
% Bottom face: labels below
\node[below=3pt] at (x5) {$x_7$};
\node[below=3pt] at (x6) {$x_8$};
\node[below=3pt] at (x7) {$x_5$};
\node[below=3pt] at (x8) {$x_6$};

% Draw edges (solid)
\foreach \i/\j in {1/2, 1/3, 1/5, 2/4, 2/6, 3/4, 3/7, 4/8, 5/6, 5/7, 6/8, 7/8}
    \draw[edge] (x\i) -- (x\j);

% Draw face diagonals (dashed)
\foreach \i/\j in {1/4, 2/3, 5/8, 6/7, 3/8, 4/7, 1/6, 2/5, 1/7, 3/5, 2/8, 4/6}
    \draw[facediag] (x\i) -- (x\j);

% Draw space diagonals (dotted)
\foreach \i/\j in {1/8, 2/7, 3/6, 4/5}
    \draw[spacediag] (x\i) -- (x\j);

\node[circle, thick, draw=black, fill=blue, inner sep=4pt] at (x1) {};
\node[circle, thick, draw=black, fill=blue, inner sep=4pt] at (x4) {};
\node[circle, thick, draw=black, fill=blue, inner sep=4pt] at (x6) {};
\node[circle, thick, draw=black, fill=blue, inner sep=4pt] at (x7) {};
\node[circle, thick, draw=black, fill=red, inner sep=4pt] at (x2) {};
\node[circle, thick, draw=black, fill=red, inner sep=4pt] at (x3) {};
\node[circle, thick, draw=black, fill=red, inner sep=4pt] at (x5) {};
\node[circle, thick, draw=black, fill=red, inner sep=4pt] at (x8) {};

\end{tikzpicture}
\caption{Clustering in periodic multiconsensus states on the $V_1$ component, where the trivial consensus first loses stability. Red circles correspond to agents in cluster $K_1 := \{x_1,x_4,x_6,x_7\}$, and blue circles correspond to agents in cluster $K_2 := \{x_2,x_3,x_5,x_8\}$. Agents in $K_1$ and $K_2$ are in antiphase. This corresponds to alternating overvaluation and undervaluation of assets across the respective markets.}\label{fig:pseudo:clusters}
\end{figure}

Now we can compute the critical set on each isotypic component. While it is possible to obtain an unlimited number of critical limit frequencies and parameter values on each isotypic component, we will show only the first critical point on each component.
\vs
By numerically solving equations \eqref{eq:pseudo:coincidence} with high machine precision, we obtain the following approximate values for $\alpha_{n,j},\beta_{n,j}$, and crossing number $\mathfrak t(\alpha_{n,j},\beta_{n,j})$:
\begin{align*}
V_0:\quad&\alpha_{1,0} = 4.24124372, & \beta_{1,0} = 0.10259439, &\quad\quad\mathfrak t(\alpha_{1,0},\beta_{1,0}) = -1\\
V_1:\quad&\alpha_{1,1} = 0.09529711, & \beta_{1,1} = 0.09239073, &\quad\quad\mathfrak t(\alpha_{1,1},\beta_{1,1}) = -1\\
V_3:\quad&\alpha_{1,3} = 0.27955407, & \beta_{1,3} = 0.09705997, &\quad\quad\mathfrak t(\alpha_{1,3},\beta_{1,3}) = -1\\
V_4:\quad&\alpha_{1,4} = 1.12332838, & \beta_{1,4} = 0.10070033, &\quad\quad\mathfrak t(\alpha_{1,4},\beta_{1,4}) = -1\\
\end{align*}
Here we can see that the first bifurcation point occurs on the $V_1$ isotypic component. Inspecting the corresponding one-dimensional eigenspace of $C$, we can see that this corresponds to a situation where edge-adjacent markets are in anti-phase, and face-diagonal markets are in phase. In other words, when trivial consensus breaks down to periodic multiconsensus on this isotypic component, despite a positive correlation between the prices of these assets, one cluster of assets will be experiencing a pricing bubble while the other experiences a crash, and these boom-bust cycles will alternate due to the competing influences of fundamentalist and chartist traders. For one representative of this conjugacy class, we can denote the cluster $K_1 := \{x_1,x_4,x_6,x_7\}$, and $K_2 := \{x_2,x_3,x_5,x_8\}$. Then for all $x_i,x_j$ in $K_l$, $x_i(t) = x_j(t)$. If we denote by $\phi(K_i)(t)$ the phase of all agents in cluster $K_l$ at time $t$, then the phase relationship between the clusters is $\phi (K_1)(t) = -\phi(K_2)(t)$. The spatial configuration of clusters can be seen in Figure \ref{fig:pseudo:clusters}.

\begin{figure}[tbp]
    \centering
    \begin{minipage}{0.45\textwidth}
        \centering
        \includegraphics[width=\linewidth]{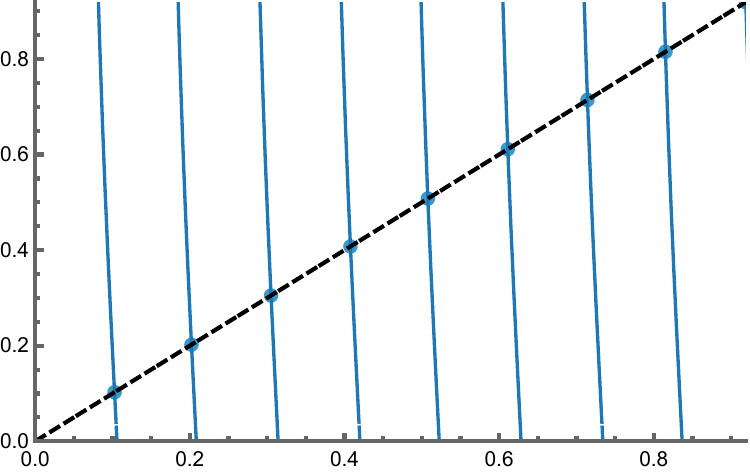}
        \subcaptiontext{a}{Limit frequencies $\beta_{n,0}$ on  $V_0$}
    \end{minipage}\hfill
    \begin{minipage}{0.45\textwidth}
        \centering
        \includegraphics[width=\linewidth]{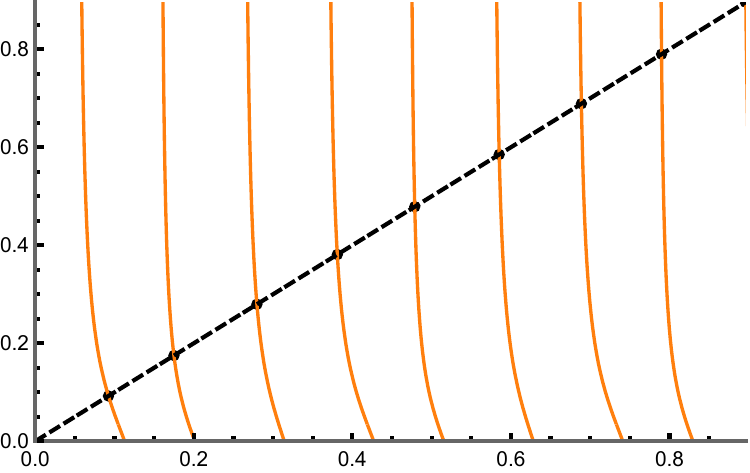}
        \subcaptiontext{b}{Limit frequencies $\beta_{n,1}$ on  $V_1$}
    \end{minipage}
    
    \vspace{0.5cm}
    
    \begin{minipage}{0.45\textwidth}
        \centering
        \includegraphics[width=\linewidth]{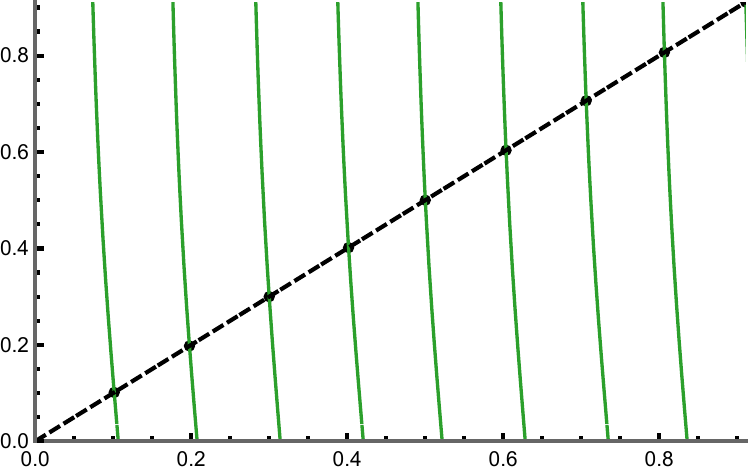}
        \subcaptiontext{c}{Limit frequencies $\beta_{n,3}$ on  $V_3$}
    \end{minipage}\hfill
    \begin{minipage}{0.45\textwidth}
        \centering
        \includegraphics[width=\linewidth]{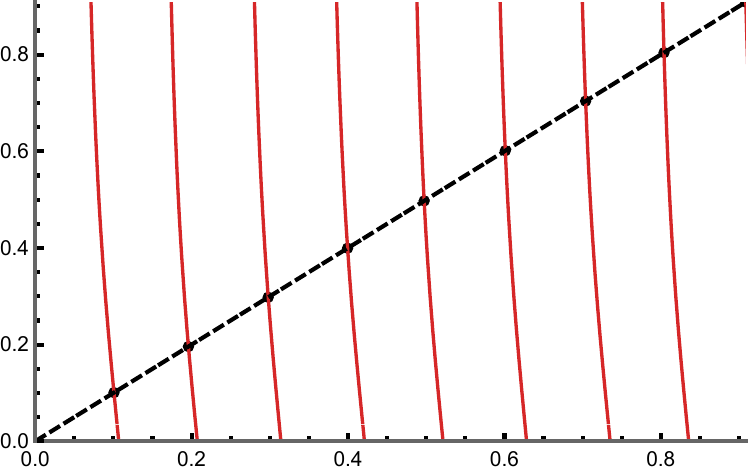}
        \subcaptiontext{d}{Limit frequencies $\beta_{n,4}$ on  $V_4$}
    \end{minipage}
    \caption{Coincidence plots of \eqref{eq:pseudo:coincidence} across isotypic components. Each crossing corresponds to a limit frequency on the corresponding isotypic component. Even for small values of $\beta$, one can see how the spacing between successive limit frequencies becomes increasingly regular as $\beta$ increases.}
    % \label{fig:pseudo:beta-coinc}
\end{figure}
% Second figure: alpha coincidence
\begin{figure}[tbp]
    \centering
    \begin{minipage}{0.45\textwidth}
        \centering
        \includegraphics[width=\linewidth]{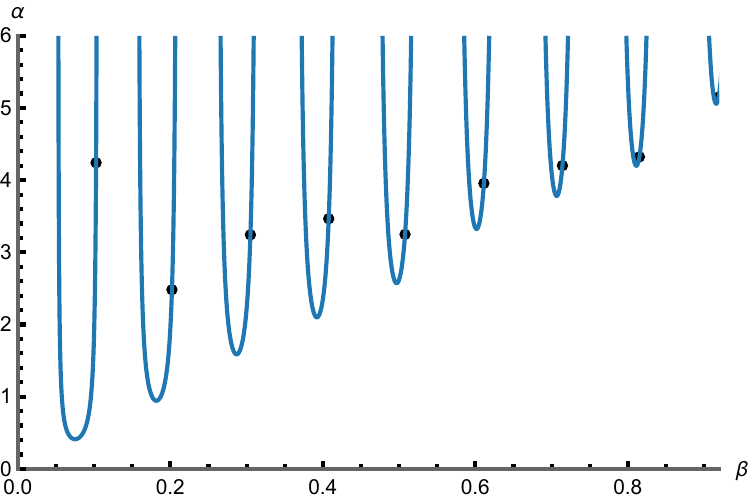}
        \subcaptiontext{a}{Critical points $\alpha_{n,0}$ on  $V_0$}
    \end{minipage}\hfill
    \begin{minipage}{0.45\textwidth}
        \centering
        \includegraphics[width=\linewidth]{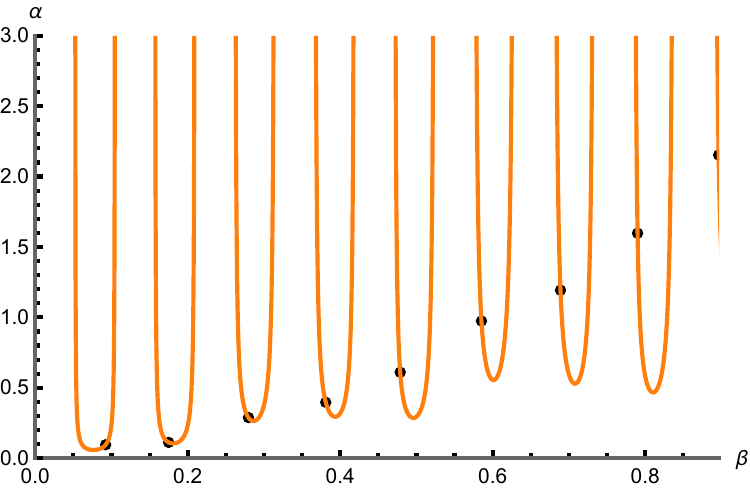}
        \subcaptiontext{b}{Critical points $\alpha_{n,1}$ on  $V_1$}
    \end{minipage}
    
    \vspace{0.5cm}
    
    \begin{minipage}{0.45\textwidth}
        \centering
        \includegraphics[width=\linewidth]{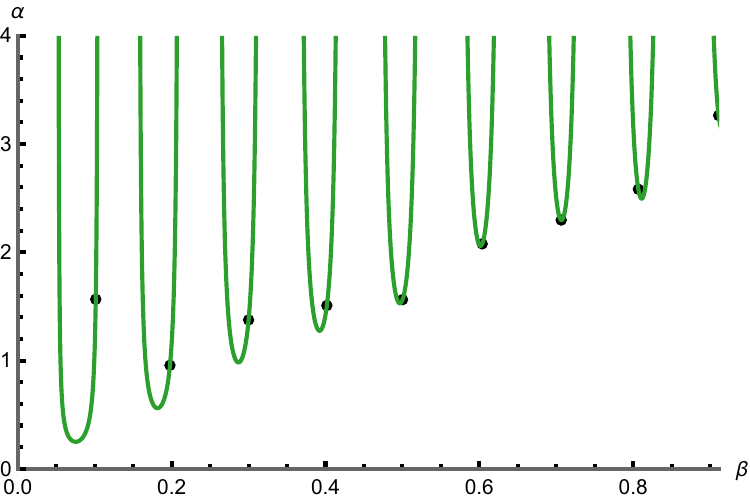}
        \subcaptiontext{c}{Critical points $\alpha_{n,3}$ on  $V_3$}
    \end{minipage}\hfill
    \begin{minipage}{0.45\textwidth}
        \centering
        \includegraphics[width=\linewidth]{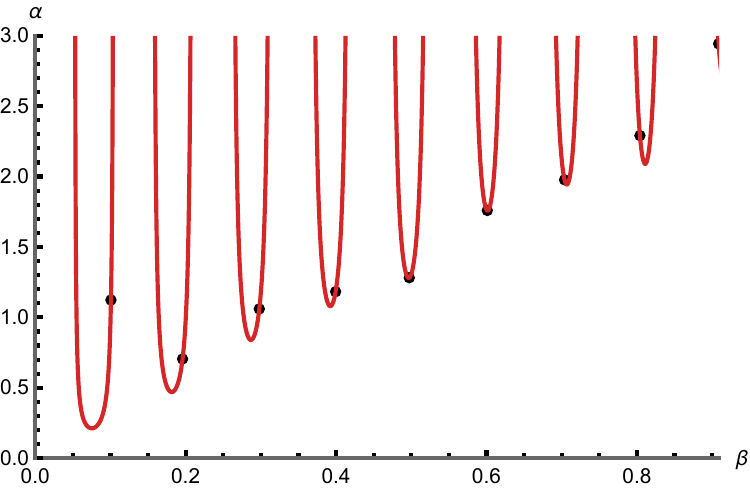}
        \subcaptiontext{d}{Critical points $\alpha_{n,4}$ on  $V_4$}
    \end{minipage}
    \caption{Critical points $\alpha_{n,j}$ corresponding to limit frequencies $\beta_{n,j}$ on each isotypic component, with $0<\alpha_{n,j}<1$. Note that although limit frequencies are always monotonically increasing, critical values of $\alpha$ may not be monotonic, especially for small values of $\alpha$.}
    % \label{fig:pseudo:alpha-coinc}
\end{figure}

\begin{figure}[tbp]
    \centering
    \includegraphics[width=\linewidth]{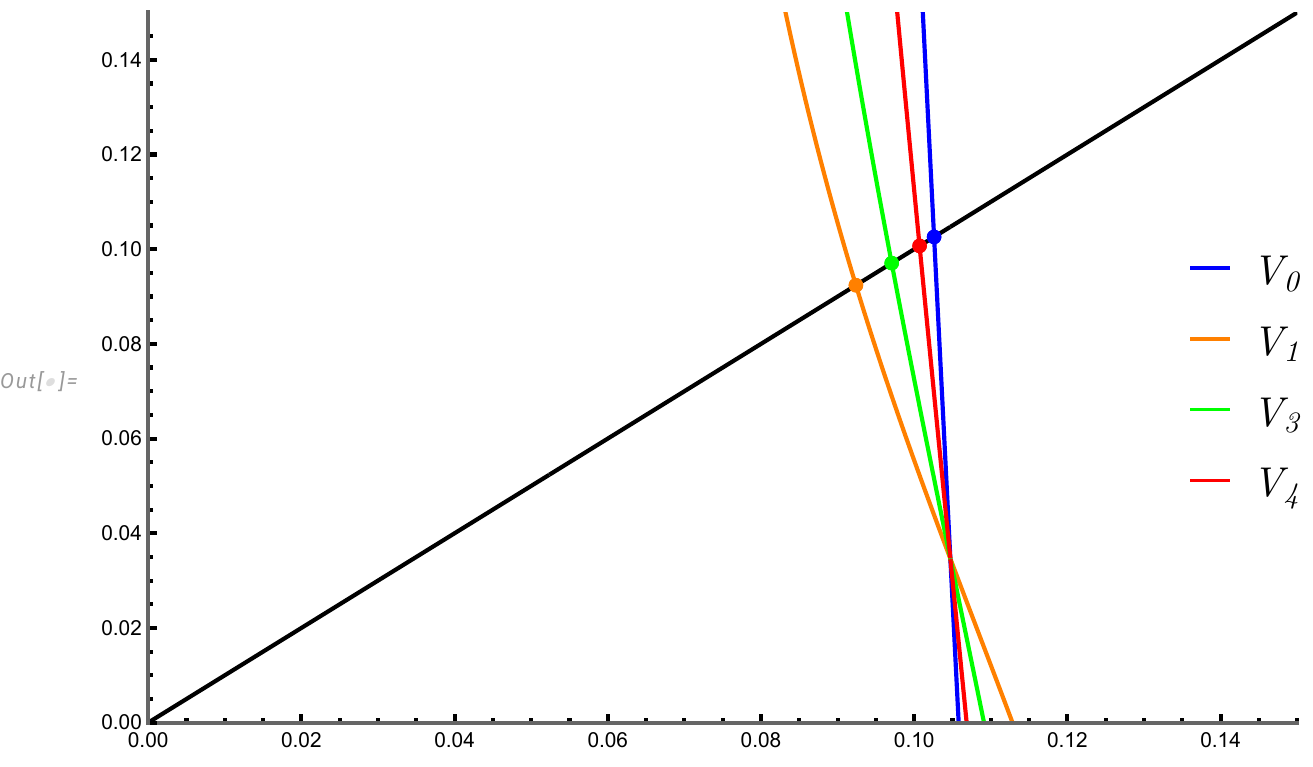}
    \caption{Overlay plot of the coincidence equations for $\beta$ near the first point of intersection across the four $S_4\times \mathbb Z_2$-isotypic components.}
    % \label{fig:pseudo:beta-coinc}
\end{figure}

\vs
For each $V_j$, the corresponding irreducible $G$-representations are given by $\mathcal V_{m,j}^-$. Using GAP with the Equideg package, we can compute the twisted basic degrees $\deg_{\mathcal V_{m,j}^-}$, yielding
\begin{align*}
    \text{deg}_{\cV_{m,0}^-} = &-( \bbZ_{2m} {}^{\bbZ_m}\times {}^{S_4} S_4^p)\\
    \text{deg}_{\cV_{m,1}^-} = &-( \bbZ_{2m} {}^{\bbZ_m}\times {}^{S_4^-} S_4^p)\\
    \text{deg}_{\cV_{m,3}^-} = &(\bbZ_{2m} {}^{\bbZ_m}\times {}^{D_4^z} D_4^p) + (\bbZ_{2m} {}^{\bbZ_m}\times {}^{D_3^z} D_3^p) + (\bbZ_{2m} {}^{\bbZ_m}\times {}^{D_2^d} D_2^p)+(\bbZ_{4m} {}^{\bbZ_m}\times {}^{\bbZ_2^-} \bbZ_4^p)+\\ 
    &(\bbZ_{6m} {}^{\bbZ_m}\times \bbZ_3^p) - (\bbZ_{2m} {}^{\bbZ_m}\times {}^{\bbZ_2^-} D_2^p) - (\bbZ_{2m} {}^{\bbZ_m}\times {}^{D_1^z} D_1^p)\\
    \text{deg}_{\cV_{m,4}^-} =&(\bbZ_{2m} {}^{\bbZ_m}\times {}^{D_4^d} D_4^p) + (\bbZ_{2m} {}^{\bbZ_m}\times {}^{D_3} D_3^p) + (\bbZ_{2m} {}^{\bbZ_m}\times {}^{D_2^d} D_2^p)+(\bbZ_{4m} {}^{\bbZ_m}\times {}^{\bbZ_2^-} \bbZ_4^p)+\\ 
    &(\bbZ_{6m} {}^{\bbZ_m}\times \bbZ_3^p) - (\bbZ_{2m} {}^{\bbZ_m}\times {}^{\bbZ_2^-} D_2^p) - (\bbZ_{2m} {}^{\bbZ_m}\times {}^{D_1} D_1^p)
\end{align*}
GAP can also be used to identify which of these orbit types are maximal, which yields:
\begin{align*}
    \mathfrak M_0 = \{&( \bbZ_{2m} {}^{\bbZ_m}\times {}^{S_4} S_4^p)\}\\
    \mathfrak M_1 = \{&( \bbZ_{2m} {}^{\bbZ_m}\times {}^{S_4^-} S_4^p)\}\\
    \mathfrak M_3 = \{&(\bbZ_{2m} {}^{\bbZ_m}\times {}^{D_4^z} D_4^p),(\bbZ_{2m} {}^{\bbZ_m}\times {}^{D_3^z} D_3^p),(\bbZ_{2m} {}^{\bbZ_m}\times {}^{D_2^d} D_2^p),(\bbZ_{4m} {}^{\bbZ_m}\times {}^{\bbZ_2^-} \bbZ_4^p),\\ 
    &(\bbZ_{6m} {}^{\bbZ_m}\times \bbZ_3^p)\}\\
    \mathfrak M_4 =\{&(\bbZ_{2m} {}^{\bbZ_m}\times {}^{D_4^d} D_4^p),(\bbZ_{2m} {}^{\bbZ_m}\times {}^{D_3} D_3^p),(\bbZ_{2m} {}^{\bbZ_m}\times {}^{D_2^d} D_2^p),(\bbZ_{4m} {}^{\bbZ_m}\times {}^{\bbZ_2^-} \bbZ_4^p),\\ 
    &(\bbZ_{6m} {}^{\bbZ_m}\times \bbZ_3^p)\}
\end{align*}
Hence, by Theorems \ref{thm:pseudo:mas-local}, \ref{thm:pseudo:mas-global}, and \ref{thm:pseudo:mas-sym}, if a bifurcation point occurs on the component $V_j$, then there is a distinct global branch of non-constant periodic solutions having symmetries $(H)$ for each $(H)\in \mathfrak M_j$. Since, in this example, the first bifurcation point occurs on the $V_1$ component, there must be a branch of non-constant periodic solutions having symmetries at least $( \bbZ_{2m} {}^{\bbZ_m}\times {}^{S_4^-} S_4^p)$. To determine if this solution represents multiconsensus, we will numerically simulate the system with a perturbation initialized on the $V_1$ subspace for values of $\alpha$ near $\alpha_{1,1} = 0.09529711$. 
\begin{figure}[h]\label{fig:pseudo:bif-diagram}
    \centering
    \includegraphics[width=\linewidth]{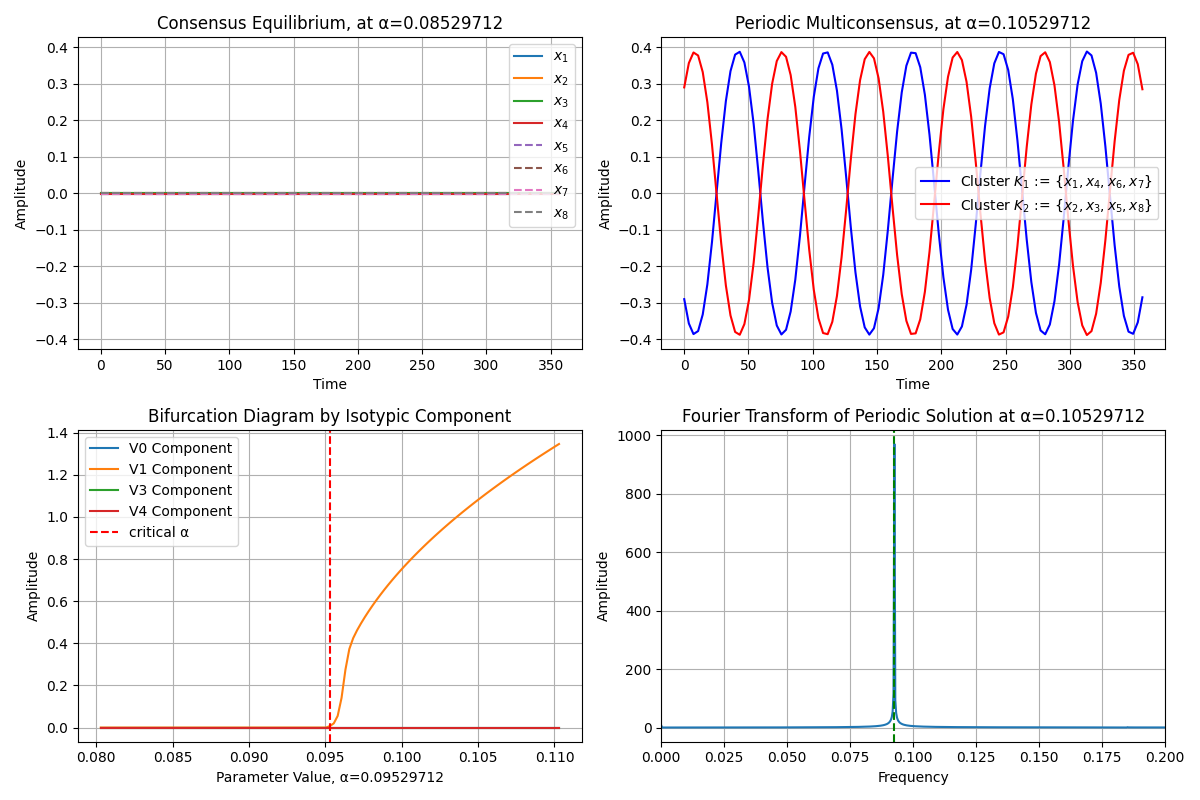}
    \caption{Time series of solution (after discarding initial transient) initialized by a perturbation near $x\equiv0$ on either side of $\alpha_{1,1}$. The dashed red line in the bottom left figure represents the value of $\alpha_{1,1}$, where bifurcation is predicted to occur, and the dashed green line in the bottom right plot represents the value of the limit frequency $\beta_{1,1}$ compared with the fast Fourier transform of the time series data.}
\end{figure}
The simulation results illustrate the local stability of the $x\equiv 0$ consensus for $\alpha<\alpha_{1,1}$, as predicted by Theorem \ref{thm:pseudo:mas-asymp}. For $\alpha > \alpha_{1,1}$, small perturbations are attracted to a locally stable non-constant periodic solution. Taking a fast Fourier transform of the time series data of the periodic solution, we also see that it has a limit frequency corresponding to $\beta_{1,1}$. By projecting the time series data onto the isotypic components, and taking its maximal value as we sweep across $\alpha$ around the bifurcation point, we see that the amplitude of the periodic solution grows as $\alpha$ grows. Since random perturbations of $x\equiv 0$ are all attracted to this same periodic solution for $\alpha>\alpha_{1,1}$ near the bifurcation point, we infer numerically that this solution is locally asymptotically stable and hence represents periodic multiconsensus. This corresponds to a cycle of price bubbles having a period of approximately 10.8 days, during which one cluster of assets is overvalued and the opposite cluster is correspondingly undervalued. 

\section{Conclusions and discussion}\label{sec:pseudo:conclusion}
In this chapter, we showed a consensus-breaking bifurcation to periodic multiconsensus in multi-agent systems with continuous memory and trend memory. We then showed an economic application as a network of heterogeneous agent models, where each asset market can be thought of as an agent whose protocol represents two populations of traders following different memory-dependent trading strategies, and whose goal is price discovery.
\vs
Since the closed-loop dynamics equations for this system formed a delay equation with distributed delays and a pseudoneutral delay, this can also be viewed as a symmetric global Hopf bifurcation result for these types of systems. In Section \ref{sec:pseudo:introduction:overview}, we established the general framework for taking a neutral or pseudoneutral delay equation, with the neutral or pseudoneutral delay written in implicit form, and formulating it as a condensing perturbation of identity which can then be analyzed using the twisted Nussbaum-Sadovskii degree. As mentioned before, this framework is not strictly necessary for the pseudoneutral equation studied in this chapter. If the pesudoneutral term is differentiated, then by the fundamental theorem of calculus, it can be written as a finite difference of retarded delays. This yields a delay equation of purely retarded type, which forms a compact perturbation of identity and can be treated using the Leray-Schauder degree. 
\vs
However, we wish to stress that using this framework for this problem makes absolutely no difference in terms of the detected orbit types or corresponding solutions. The only difference is that the framework used here imposes the condition that $\gamma\tau_1<1$, as this is necessary for the operator $\bm k_g$ to be a contraction and therefore for $\mathscr F$ to be a condensing perturbation of identity. If we wrote the pseudoneutral delay as a finite difference of retarded delays, we could relax this condition. However, the condition that $\gamma \tau_1 < 1$ is also closely related to the local asymptotic stability of the $x\equiv0$ consensus equilibrium. Indeed, Proposition \ref{prop:pseudo:asymptotic-stability} implies that for $\gamma\tau_1 \geq 1$, there may be eigenvalues with positive real part for all $\alpha$, and therefore the $x\equiv0$ equilibrium is never stable. Not only does this mean that the $x\equiv 0$ consensus is never achieved in this case, it also means that studying Hopf bifurcation from this equilibrium is not very interesting from a practical point of view. 
\vs
It is useful to see the direct relationship between this condition and the analogous condition for a neutral system. If we replace the pseudoneutral delay with a neutral delay, then the analogous condition imposed on the linear system by requiring the operator $\bm k_g$ (as given in Section \ref{sec:pseudo:introduction:overview}) to be a contraction is that $\gamma < 1$. In this case, which will be explored in the next chapter, violation of this condition leads to an infinite chain of eigenvalues with positive real part, preventing the operator $\mathscr F$ from being a condensing perturbation of identity, and making the entire system ill-posed as an initial value problem. This not only illustrates another way in which pseudoneutral systems can be viewed as a kind of family of approximations of neutral systems, but also shows some of the practical reasons for analyzing both types of system using the same overarching framework.

\chapter{Global Periodic Multiconsensus Bifurcation in Neutral Multi-Agent Systems with Momentum Memory}\label{chapter:neutral}
 \author{Casey Crane}

\section{Introduction}\label{sec:neutral:introduction}
Consider a multi-agent system of $n$ homogeneous agents with single-integrator dynamics, whose open loop dynamics are given by
\[
\dot{\bm x} =  \bm u_i,\quad \bm x(t) = (x_1,\dots,x_n)^T\in\mathbb R^n.
\]
We define an identical protocol for each agent
\[
u_i:= -ax_i - \alpha f\left( \int_0^{\tau_2} x_i(t-s)ds \right) + \frac{d}{dt}\left[g(x_i(t-\tau_1))\right] - h_i(\bm x),
\]
which consists of the following terms:
\begin{itemize}
\item A self-regulation term $-ax_i(t)$, representing the agent's tendency to return to a reference point at $x\equiv0$.
\item An interaction rule given by some nonlinear coupling function $h_i:\mathbb R^n \to \mathbb R$ which defines a $\Gamma_0$-symmetric fixed interaction topology with state-dependent effective weights.
\item A \emph{continuous memory term}, consisting of a nonlinear function of a continuous weighted average of the agent's own state history, $-\alpha f(\int_0^{\tau_2}x_i(t-s)ds)$, where $\alpha$ is taken as a bifurcation parameter.
\item A \emph{momentum memory term}, consisting of a derivative of the nonlinearly transformed state of the agent at a fixed point in the past, written in derivative form as $\frac{d}{dt}\left[g(x_i(t-\tau_1))\right]$. 
\end{itemize}
Then the closed loop dynamics of a single agent can be written
\[
\frac{d}{dt}\left[x_i - g(x_i(t-\tau_1))\right] = -ax_i - \alpha f\left( \int_0^{\tau_2} x_i(t-s)ds \right) - h_i(x), \;\quad x_i(t) \in \mathbb R.
\]
The presence of the trend memory term makes this a neutral functional differential equation. If we group the identical components of the protocol together and put $\bm f,\bm g, \bm h:\mathbb R^n \to \mathbb R^n$ as $\bm f(\bm x) := (f(x_1),\dots,f(x_n))^T$, $\bm g(\bm x) := (g(x_1),\dots,g(x_n))^T$, and $\bm h(\bm x) := (h_1(x),\dots,h_n(x))^T \in \mathbb R^n$, respectively, then the closed-loop dynamics of the full system can be written in vector form as
\begin{equation}\label{eq:neutral:basic-sys}
\frac{d}{dt}\left[\bm x - \bm g(\bm x(t-\tau_1))\right] = -a \bm x - \alpha \bm f\left( \int_0^{\tau_2} \bm x(t-s)ds \right) - \bm h(\bm x).
\end{equation}
\vs
To account for the symmetries of the agent interactions, let $\Gamma_0\leq S_n$ be a finite symmetry group which acts by permuting the indices of $\bm x$. Let $V:= \mathbb R^n$, and put $\Gamma := \Gamma_0 \times \mathbb Z_2$, where $\mathbb Z_2$ acts antipodally, i.e. the action of $(\sigma,\pm1)$ on $\bm x \in V$ is given by
\[
    (\sigma,\pm1)\bm x :=  (\sigma,\pm1)(x_1,\dots,x_n)^T=\pm(x_{\sigma(1)},\dots,x_{\sigma(n)})^T.
\]
We will make identical assumptions about $\bm f, \bm g,$ and $\bm h$ as in Chapter \ref{chapter:pseudo}, reproduced here for convenience:
\begin{enumerate}[label=($C_\arabic*$)]\setcounter{enumi}{0}\itemc
\item\label{neutral:c1} $\bm f$, $\bm g$, and $\bm h$ are $\Gamma$-equivariant (and therefore are odd functions).
\item\label{neutral:c2} $\bm f$, $\bm g$, and $\bm h$ are continuous functions differentiable at 0, and $b := f'(0), \gamma := g'(0), C:= Dh(0)$ satisfy:
    \begin{enumerate}
        \item $b\neq0$.
        \item $0 < \gamma < 1$.
        \item C is a $\Gamma$-equivariant symmetric matrix.
    \end{enumerate}
\item\label{neutral:c3} $\bm g$ is $\kappa$-Lipschitzian with $\kappa<1$, i.e. 
\[
\exists_{\kappa \in [0,1)}\quad \forall_{\varphi,\psi \in C(\bbR^n;\bbR)} \quad ||\bm g(\varphi) - \bm g(\psi)|| \leq \kappa||\varphi - \psi||_\infty
\]
\end{enumerate}
\vs
We again refer to the equilibrium $\bm x\equiv0$ as the \textit{reference consensus} (or \textit{consensus equilibrium}), and we will prove the following theorems establishing conditions for local and global symmetric multiconsensus solution bifurcating from the reference consensus, and their classification according to their spatio-temporal symmetries:
\begin{theorem}\label{thm:neutral:mas-asymp}
Let $\widehat \alpha_j := \min_{n\in \mathbb N}\{\alpha\in\mathbb R:(\alpha_{n,j},\beta_{n,j},0)\in \Lambda\}$, and $\widehat \alpha := \min_{j=0,1,\dots,r} \{\widehat a_j\}$. Then if $a_j > 0$ for all $j=0,1,\dots,r$, $0\leq\lambda<1$, and $b>0$, then the reference consensus equilibrium at $x\equiv 0$ is locally asymptotically stable for all $\alpha \in (0,\widehat \alpha)$.
\end{theorem}
\begin{theorem}\label{thm:neutral:mas-local}
If $\beta_0>0$ satisfies
\[
\beta_0 = \frac{(\beta_0\gamma\sin(\beta_0\tau_1)-a_j)(\cos(\beta_0\tau_2)-1)}{(1-\gamma\cos(\beta_0\tau_1))\sin(\beta_0\tau_2)}
\]
and $\alpha_0>0$ satisfies
\[
\alpha_0 = \frac{\beta_0^2\gamma\sin(\beta_0 \tau_1) - \beta_0 a_j}{b\sin(\beta_0 \tau_2)}
\]
for some $j=0,1,\dots,r$, and if $\alpha_0,\beta_0$ satisfy 
\[
\rho_j(\alpha_0,\beta_0):=p(\alpha_0,\beta_0)b(\cos(\beta_0\tau_2)-1)+q(\alpha_0,\beta_0)b\sin(\beta_0\tau_2) \neq 0
\]
where
\begin{align*}
p(\alpha,\beta) &:= 2\beta\gamma\sin(\beta\tau_1) - \beta^2\gamma\tau_1 \cos(\beta\tau_1) + a_j - \alpha b \tau_2 \cos(\beta\tau_2),\\
q(\alpha,\beta) &:= 2\beta\gamma\cos(\beta\tau_1)-\beta^2\gamma\tau_1\sin(\beta\tau_1)-2\beta+\alpha b\tau_2\sin(\beta\tau_2),
\end{align*}
then there exists a connected global branch of non-constant periodic multiconsensus solutions bifurcating from the reference consensus at $\alpha=\alpha_0$.
\vs
\end{theorem}
\begin{theorem}\label{thm:neutral:mas-global}
Fix $j=0,1,\dots,r$, and denote by $\beta_{n,j}$ the sequence of solutions of 
\[
\beta_{n,j} = \frac{(\beta_{n,j}\gamma\sin(\beta_{n,j}\tau_1)-a_j)(\cos(\beta_{n,j}\tau_2)-1)}{(1-\gamma\cos(\beta_{n,j}\tau_1))\sin(\beta_{n,j}\tau_2)}
\]
and $\alpha_{n,j}$ the corresponding values
\[
\alpha_{n,j} = \frac{\beta_{n,j}^2\gamma\sin(\beta_{n,j} \tau_1) - \beta_{n,j} a_j}{b\sin(\beta_{n,j} \tau_2)}.
\]
If $\tau_1/\tau_2 \in \mathbb Q$, and if $(\widetilde \beta_{l,j})_{l=1}^P$ is a $P$-periodic sequence of roots of 
    \[
    1 = \frac{(\gamma\sin(\beta\tau_1)-\varepsilon)(1-\cos(\beta\tau_2))}{(1-\gamma\cos(\beta\tau_1))\sin(\beta\tau_2)},
    \]
    where $\varepsilon>0$ is taken sufficiently small, such that $(\widetilde \beta_{l,j})_{l=1}^P$ satisfies
    \[
    \sum_{l=1}^P\sign \Upsilon_j(\widetilde \beta_{l,j}) \neq 0
    \]
    where 
    \[
    \Upsilon_j(\beta):= \gamma\tau_1\cos(\beta\tau_1) + \tau_2\cos(\beta\tau_2) -\gamma\tau_1\cos(\beta(\tau_1-\tau_2))-\gamma\tau_2\cos(\beta(\tau_1+\tau_2)),
    \]
    or if there exists some $K>0$ such that $\sign \rho_j(\alpha_{K+i,j},\beta_{K+1,j}) = \sign \Upsilon_j(\widetilde \beta_{l,j})$ for all $l\in \mathbb N$, where indices of $\widetilde \beta_{l,j}$ are taken mod $P$, and
    \[
    \sum_{i=1}^K\sign \rho_j(\alpha_{i,j},\beta_{i,j}) \neq 0,
    \]
    then there exists a global unbounded branch of non-constant periodic multiconsensus solutions bifurcating from the reference consensus at $x\equiv 0$.
\end{theorem}
\begin{theorem}\label{thm:neutral:mas-sym}
If the conditions of Theorem \ref{thm:neutral:mas-local} hold, and $(H)\in \Phi_1^t(S^1\times \Gamma_0\times \mathbb Z_2;\mathscr E_{k,j})$ is a maximal twisted orbit type which fixes some $u(t) \in \mathscr E_{1,j}:=\{c\cos(t)+d\sin(t):c,d\in E(\mu_j)\}$, then there exists a connected component of non-constant periodic solutions with symmetries at least $(H)$ bifurcating from the trivial solution at $\alpha=\alpha_0$.
Moreover, if the conditions of Theorem \ref{thm:neutral:mas-global} are satisfied, then this branch is unbounded.
\end{theorem}
\vs
As before, we will take a ``bilingual'' approach to proving these theorems. In particular, we will state and prove logically equivalent formulations of these theorems in terms of local and global Hopf bifurcation of non-constant periodic solutions from the trivial equilibrium for the system of neutral functional differential equations which governs the closed-loop dynamics of this MAS. Accordingly, these results are valid for any such system of NFDEs, irrespective of its interpretation.
\vs
This chapter will build on the framework established in Chapter \ref{chapter:pseudo}. As discussed in greater detail in that chapter, the closed loop dynamics of the MAS studied in that chapter formed a delay equation which, in spite of its outward appearance when the trend memory term was written inside the derivative, was actually of purely retarded type. For this reason we coined the designation ``pseudoneutral'' for that type of equation (cf. Sec.  \ref{sec:prelim:pseudoneutral}). This equation, however, is truly of neutral type. We will show that very similar requirements prevail for its stability as those which were proven for the pseudoneutral case. Although the trend memory term in that protocol exhibited ``derivative-like'' behavior, especially for small values of $\tau_1$, the fact that it was a finite difference prevented it from fully exhibiting the most problematic behaviors associated with neutral equations as a class. That is not so in this case.
\vs
However, we still maintain the requirement that the trend memory term (which we will also occasionally refer to as the neutral memory or neutral delay term when the NFDE context is to be emphasized) is contractive, as this is necessary for equation \eqref{eq:neutral:basic-sys} to be of Hale class, to have an essential spectrum contained in the negative half-plane, and for it to be possible for only finitely many eigenvalues to have non-negative real part. While violation of this condition merely made the zero conensus unstable for all parameter values in the pseudoneutral equation (but did not impede use of twisted equivariant degree techniques \emph{per se}), here this condition is absolutely essential in order to guarantee that the operator equation reformulation is a condensing perturbation of identity. 
\vs
Using the framework developed in Chapter \ref{chapter:pseudo} allows us to dispense with many of the preliminaries and move directly to the linearization of the equation around the $x\equiv 0$ equilibrium and analysis of its spectral properties.

\section{Linearization, spectral analysis, and critical points}\label{sec:neutral:crit}

\vs
Linearizing \eqref{eq:neutral:basic-sys} around $x=0$ gives
\begin{equation}\label{eq:neutral:basic-sys-lin}
\frac{d}{dt} \left[\bm x - \gamma \bm x(t-\tau_1)\right] + a \bm x + \alpha b \int_0^{\tau_2} \bm x(t-s)ds + C \bm x = 0,
\end{equation}
\vs
with corresponding characteristic operator equation
\begin{equation}\label{eq:neutral:char-op}
\triangle_{\alpha}(\lambda):=\left(\lambda(1-\gamma e^{-\lambda\tau_1}) + a - \alpha b \frac{e^{-\lambda\tau_2} -1}{\lambda}\right)\text{Id} + C
\end{equation}
Taking $\mu_j$ as an eigenvalue of the $\Gamma_0$-equivariant matrix $C$ and restricting \eqref{eq:neutral:char-op} to complexified irreducible $\Gamma_0$-representations, we obtain the $\Gamma_0$-isotypic characteristic equations
% \begin{definition}[Centers and limit frequencies]\normalfont
% The $\bm x \equiv 0$ solution to \eqref{eq:neutral:basic-sys-lin} at $\alpha=\alpha_0$ is called a \emph{center} if there exists a corresponding $\beta_0>0$ (called the \textit{limit frequency}) such that $\det_{\mathbb C}\triangle_{\alpha_0}(i\beta_0) =0$. If, in addition, there exists $\varepsilon>0$ such that $0<|\alpha-\alpha_0|+|\beta-\beta_0|<\varepsilon$ implies
% \[
% {\det}_{\mathbb C}\triangle_{\alpha}(ik\beta) \not=0 
% \]
% for all $k\in\mathbb N$, then it is called an \emph{isolated center}.
% \end{definition}

\begin{equation}\label{eq:neutral:char-op-iso-restrict}
\triangle_{\alpha,j}(\lambda):= \triangle_\alpha(\lambda)_{|\cV^{C}_j}:\cV^{C}_j\to \cV^{C}_j
\end{equation}
\begin{figure}[tbp]
    \centering
    \begin{minipage}{0.45\textwidth}
        \centering
        \includegraphics[width=\linewidth]{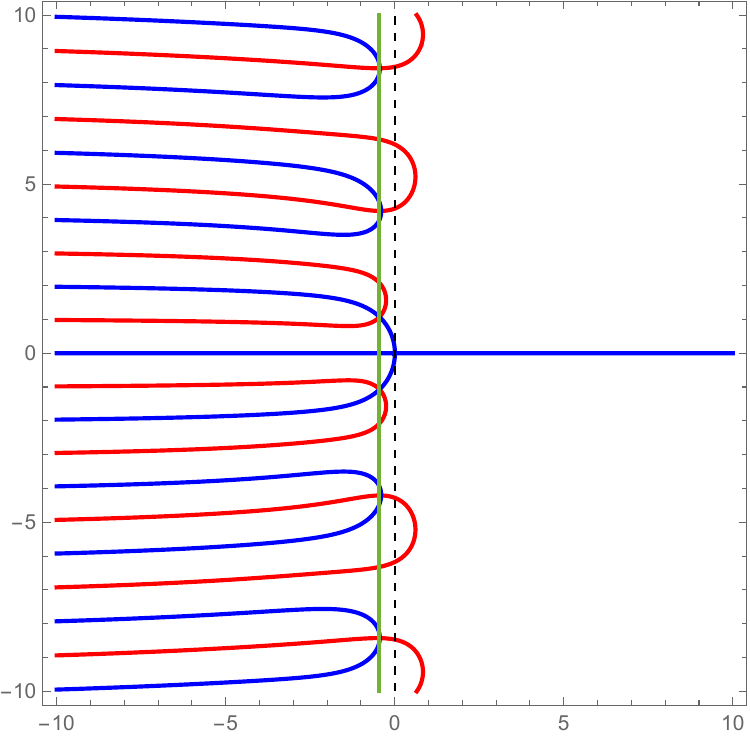}
        \subcaptiontext{a}{Characteristic roots for \eqref{eq:neutral:basic-sys-lin}}
    \end{minipage}\hfill
    \begin{minipage}{0.45\textwidth}
        \centering
        \includegraphics[width=\linewidth]{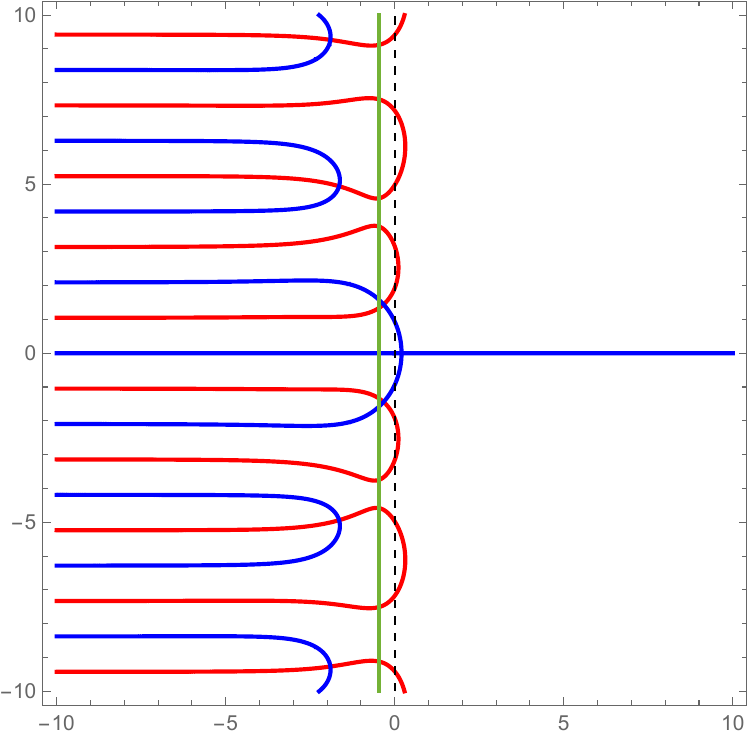}
        \subcaptiontext{b}{Characteristic roots for \eqref{eq:pseudo:basic-sys-lin}}
    \end{minipage}
    \caption{Plots of roots of the real and imaginary parts of the characteristic equations for the neutral and pseudoneutral systems, taken with the same parameter values $\gamma,a_j = 0.5,\alpha,b=1,\tau_1 = 1.5,\tau_2=1.01$. Red contours are $\operatorname{{Re}} \triangle_{\alpha,j}(\lambda)=0$ level curves, blue contours are $\operatorname{Im} \triangle_{\alpha,j}(\lambda)=0$ level curves, the dashed vertical line represents the imaginary axis, and the green vertical line represents the accumumlating line of the essential spectrum at $\operatorname{Im} \lambda = \tfrac{\ln \gamma}{\tau_1}$. Notice how eigenvalues (where the red and blue curves intersect) accumulate on this line for the neutral equation but fall away from it as $\abs \lambda$ increases for the pseudoneutral equation. If $\gamma>1$, this line lies in the positive half-plane.}
    \label{fig:neutral:char-root-comparison}
\end{figure}

We again put $a_j := a + \mu_j$, where $\mu_j$ is the $j$-th eigenvalue of $C$ whose eigenspace $E(\mu_j)$ corresponds to the $j$-th $\Gamma_0$-isotypic component. Then under the assumptions \ref{neutral:c1} -- \ref{neutral:c3}, \eqref{eq:neutral:char-op-iso-restrict} can be written
\begin{equation}\label{eq:neutral:char-op-iso}
\triangle_{\alpha,j}(\lambda):=\left(\lambda(1-\gamma e^{-\lambda\tau_1}) + a_j - \alpha b \frac{e^{-\lambda\tau_2} -1}{\lambda}\right)\text{Id},
\end{equation} 
and \eqref{eq:neutral:char-op} can be written as the block matrix:
\[
\triangle_\alpha(\lambda) = \begin{bmatrix}
    \triangle_{\alpha,0}(\lambda) & &\makebox(0,0){\text{\huge0}}\\
    & \ddots & \\
    \text{\huge0}& & \triangle_{\alpha,r}(\lambda)
\end{bmatrix},
\]
\begin{equation}\label{eq:neutral:char-block-det}
{\det}_{\mathbb C}(\triangle_\alpha(\lambda)) = \prod_{j=0}^r {\det}_{\mathbb C}(\triangle_{\alpha,j}(\lambda))^{m_j}
\end{equation}
where $m_j$ denotes the $\Gamma_0$-isotypic multiplicity of the eigenvalue $\mu_j$, defined
\begin{equation}\label{eq:neutral:def-iso-mult}
m_j := \frac{\text{dim }E(\mu_j)}{\text{dim }\cV_j}
\end{equation}
\subsection{The critical set}
We define the \textit{characteristic quasipolynomial} corresponding to $\triangle_{\alpha,j}(\lambda)$, written
\begin{equation}\label{eq:neutral:char-root}
P_j(\alpha,\lambda) := \lambda^2(1-\gamma e^{-\lambda\tau_1}) + a_j\lambda -\alpha b(e^{-\lambda\tau_2}-1).
\end{equation}  
Clearly, $\lambda \in {\mathbb C}\setminus\{0\}$ is a characteristic root of \eqref{eq:neutral:char-op} if and only if, for some $j=0,\dots,r$ and some $\alpha$, we have $P_j(\alpha,\lambda)=0$. Furthermore, the linearized system \eqref{eq:neutral:basic-sys-lin} admits a center if and only if, for some $j$, there exists $\alpha_0$ and $\beta_0$ such that $P_j(\alpha_0,i\beta_0) = 0$.
\vs
We will now show that the system \eqref{eq:neutral:basic-sys-lin} admits an infinite number of isolated centers.
Setting $\lambda = i\beta$ and substituting into \eqref{eq:neutral:char-root}, we obtain 
% \[
% -\beta^2 + i\beta\gamma(\cos(\beta\tau_1) - i\sin(\beta\tau_1) -1) + i\beta a_j - \alpha b (\cos(\beta\tau_2) - i\sin(\beta\tau_2) - 1) = 0
% \]
\[
-\beta^2(1-\gamma(\cos(\beta\tau_1)-i\sin(\beta\tau_1)) + i\beta a_j - \alpha b(\cos(\beta\tau_2)-i\sin(\beta\tau_2)-1).
\]
Splitting into real and imaginary parts yields
\begin{align*}
-\beta^2(1-\gamma\cos(\beta\tau_1))+ \alpha b (1-\cos(\beta\tau_2)) &= 0\label{eq:neutral:system-sincos}\\
-\beta^2\gamma\sin(\beta\tau_1) + \beta a_j + \alpha b \sin(\beta\tau_2) &= 0.\nonumber
\end{align*}
Solving the second equation for $\alpha$, we get
\begin{equation}\label{eq:neutral:alpha-relation}
% \alpha = \frac{-\gamma\beta(\cos(\beta\tau_1)-1)-\beta a_j}{b\sin(\beta\tau_2)}.
\alpha = \frac{\beta^2\gamma\sin(\beta \tau_1) - \beta a_j}{b\sin(\beta \tau_2)}
\end{equation}
Substitution into the first equation and solving for $\beta$ yields
\begin{equation}\label{eq:neutral:beta-relation}
% \beta = \gamma\sin(\beta\tau_1) + \frac{(\gamma(\cos(\beta\tau_1)-1)+a_j)(\cos(\beta\tau_2)-1)}{\sin(\beta\tau_2)}.
\beta = \frac{(\beta\gamma\sin(\beta\tau_1)-a_j)(\cos(\beta\tau_2)-1)}{(1-\gamma\cos(\beta\tau_1))\sin(\beta\tau_2)}
\end{equation}

\begin{figure}[tbp]
    \centering
    \begin{minipage}{0.45\textwidth}
        \centering
        \includegraphics[width=\linewidth]{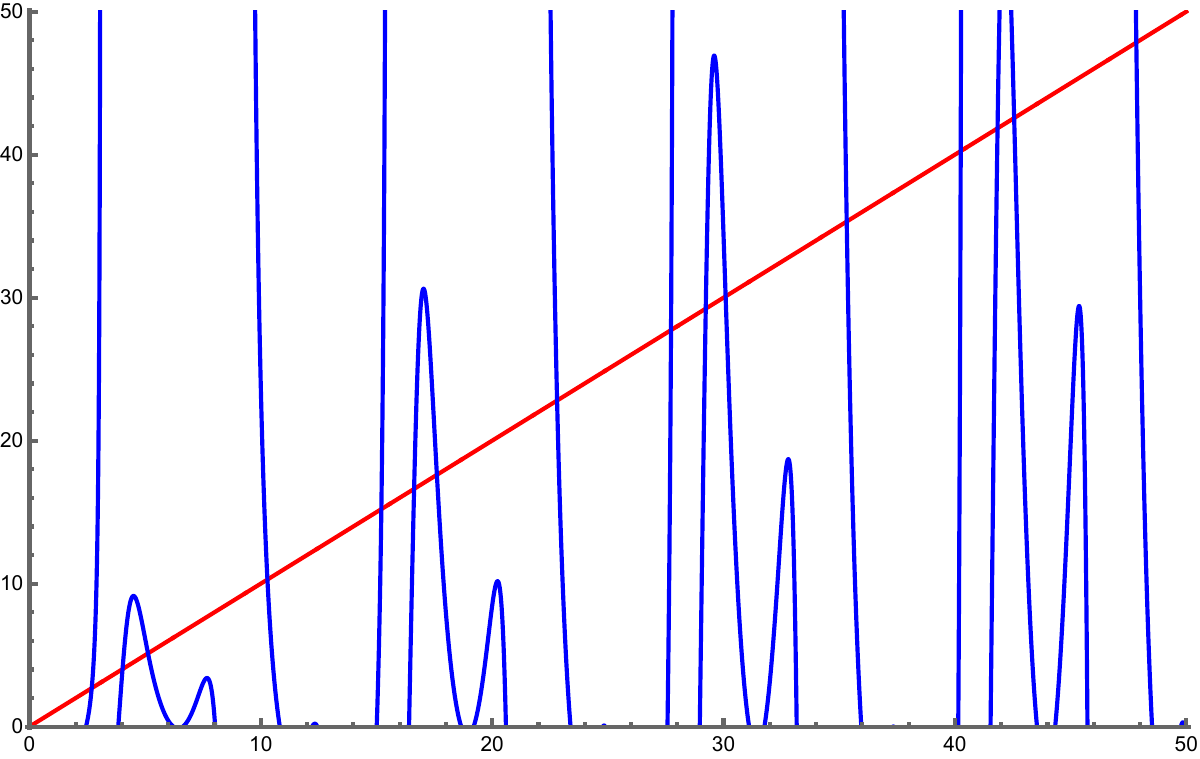}
        \subcaptiontext{a}{Coincidence plot for \eqref{eq:neutral:beta-relation} for the neutral system}
    \end{minipage}\hfill
    \begin{minipage}{0.45\textwidth}
        \centering
        \includegraphics[width=\linewidth]{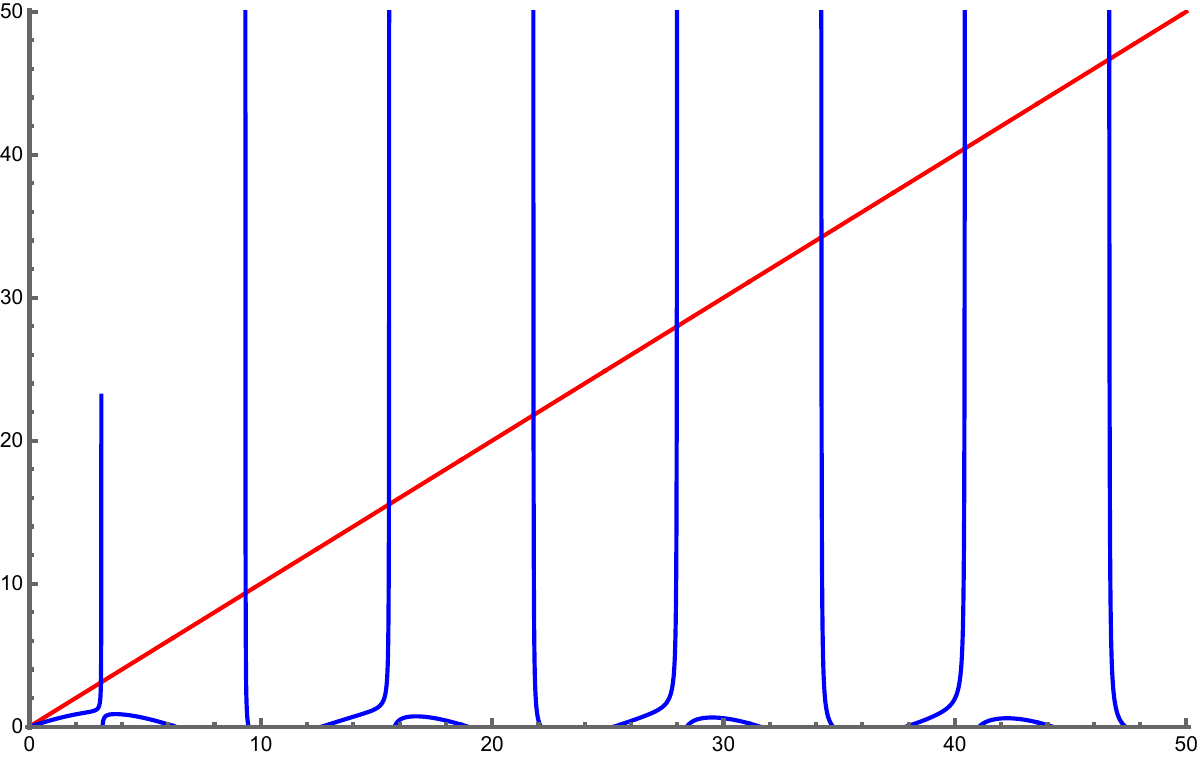}
        \subcaptiontext{b}{Coincidence plot for \eqref{eq:pseudo:beta-relation} for the pseudoneutral system}
    \end{minipage}
    \caption{Comparison of the transcendental coincidence equations for $\beta$ for the neutral and pseudoneutral systems. Both are plotted with the same parameter values of $\gamma,a_j = 0.5,\alpha,b=1,\tau_1 = 1.5,\tau_2=1.01$. The red and blue curves in each plot represent the left-hand and right-hand sides of the coincidence equations, respectively. Notice how the blue curve for the neutral equation scales with $\beta$, causing intersections at similar points in each period. This reflects the propagation of instabilities inherent in neutral equations.}
    \label{fig:neutral:coinc-comparison}
\end{figure}

A solution $\beta_0$ to \eqref{eq:neutral:beta-relation} corresponds to a purely imaginary characteristic root of \eqref{eq:neutral:char-op-iso}. Plugging this root into \eqref{eq:neutral:alpha-relation} yields a corresponding $\alpha_0$, and $(\alpha_0,\beta_0)$ is an isolated center of \eqref{eq:neutral:basic-sys-lin} at $x=0$ with limit frequency $\beta_0$. That this set is discrete and isolated can be seen from the plots of the coincidence problem \eqref{eq:neutral:beta-relation} on the left side of Figure \ref{fig:neutral:coinc-comparison}, and from the periodicity of the right-hand side. In particular, in order to have non-isolated intersections of these two graphs, we would require that their derivatives (with respect to $\beta$) agree on a set of non-zero measure. This implies a set of non-zero measure in $\beta$ on which the derivative
\[
\frac{d}{d\beta}\left[\frac{(\beta\gamma\sin(\beta\tau_1)-a_j)(1-\cos(\beta\tau_2))}{(1-\gamma\cos(\beta\tau_1))\sin(\beta\tau_2)}\right] = 1,
\]
but this is a meromorphic function, and its derivative with respect to $\beta$ is also meromorphic, and consists of sums and products of scaled and shifted trigonometric functions, so solutions to \eqref{eq:neutral:beta-relation} must be isolated. Since one can form a sequence strictly increasing in $\beta$ with corresponding values of $\alpha$ sharing an index for each $j=0,1,\dots,r$, we are justified in labeling the sequence of solutions to \eqref{eq:neutral:beta-relation} as $\beta_{n,j}$ and corresponding values of $\alpha$ obtained from \eqref{eq:neutral:alpha-relation} as $\alpha_{n,j}$ (although this sequence may not be monotonic in $\alpha$).
\vs
A triple $(\alpha_{n,j},\beta_{n,j},0)$ such that $x\equiv 0$ is an isolated center of \eqref{eq:neutral:basic-sys-lin} at $\alpha=\alpha_0$ with limit frequency $\beta_0$ is called a \emph{critical point} of system \eqref{eq:neutral:basic-sys}, and the set 
\begin{equation}\label{eq:neutral:crit-set}
\Lambda := \{(\alpha_{n,j},\beta_{n,j},0) : n\in{\mathbb N},\; j=0,1,\dots,r, \;\alpha_{n,j},\beta_{n,j}>0\}    
\end{equation}
is called the \emph{critical set} for \eqref{eq:neutral:basic-sys}. 
\vs

\subsection{Stability of the trivial equilibrium}
The monotonicity of the critical set in $\beta$ means that, although it is not monotonic in $\alpha$, we can nonetheless choose the smallest value of $\alpha_{n,j}$ for each $j=0,1,\dots,r$, which we denote $\widehat \alpha_j := \min_{n\in \mathbb N}\{\alpha\in\mathbb R:(\alpha_{n,j},\beta_{n,j},0)\in \Lambda\}$. Then the continuous dependence of characteristic roots on $\alpha$ implies that if the $x\equiv 0$ reference consensus can be shown to be locally asymptotically stable for $\alpha =0$, then this stability holds up until either $\alpha = \widehat \alpha_j$ or until a steady-state bifurcation occurs for some $\alpha_* < \widehat \alpha_j$. This leads to the following theorem on the stability of the reference consensus:
\begin{theorem}\label{thm:neutral:trivial-consensus-stability}
Let $\widehat \alpha_j := \min_{n\in \mathbb N}\{\alpha\in\mathbb R:(\alpha_{n,j},\beta_{n,j},0)\in \Lambda\}$, and $\widehat \alpha := \min_{j=0,1,\dots,r} \{\widehat a_j\}$. Then if $a_j > 0$ for all $j=0,1,\dots,r$, $0\leq\lambda<1$, and $b>0$, then the trivial equilibrium at $x\equiv 0$ is stable for all $\alpha \in (0,\widehat \alpha)$.
\end{theorem}
\begin{proof}
For $\alpha=0$, the negative power of $\lambda$ disappears in the $\Gamma_0$-isotypic charactersistic equation, and so characteristic roots must satisfy
\begin{equation}\label{eq:neutral:stability-alpha-zero}
\lambda(1-\gamma e^{-\lambda \tau_1}) = -a_j.
\end{equation}
First, we must show that there are no purely imaginary roots when $\alpha=0$ (i.e. $x\equiv0$ is not itself a center when $\alpha=0$). Put $\lambda = i\beta$. Then we obtain
\[
\beta(1-\gamma(\cos(\beta\tau_1) - i\sin(\beta\tau_1))) =-a_j.
\]
Which gives the system
\begin{align*}
\beta(1-\gamma\cos(\beta\tau_1)) &= -a_j,\\
\beta\gamma\sin(\beta\tau_1) &=0.
\end{align*}
This implies that $\beta\tau_1 = \pi k$ for some $k\in \mathbb N$. Suppose $\beta\tau_1 = 2\pi k$. Then the first equation gives
\[
2\pi k(1-\gamma) = -a_j,
\]
which cannot be satisfied if $a_j>0$ and $\gamma<1$, for any $k\in \mathbb N$. On the other hand, if $\beta = \pi(2k-1)$, then we have
\begin{align*}
\pi(2k-1)(1+\gamma) &= -a_j,\\
\end{align*}
which also has no solutions for $a_j>0$, $\gamma<1$, and $k\in \mathbb N$. Note that the same argument suggests that purely imaginary roots for $\alpha=0$ \emph{are} possible for $a_j < 0$, but only for highly non-generic values of $\gamma$ and $a_j$, as might be expected.

Now we turn to the question of complex roots with positive real part: if $\re \lambda >0$, then this is only possible if $1<\re (\gamma e^{-\lambda \tau_1})$, which implies $\ln\gamma>\re \lambda \tau_1$ and so $\re(\lambda) < \tfrac{\ln\gamma}{\tau_1}$. Therefore, if $\tfrac{\ln \gamma}{\tau_1}<0$, we have a contradiction, and all roots must have negative real part. Since $\tau_1>0$, this immediately yields the requirement $\gamma<1$.
\vs
Now we will consider the situation of steady-state bifurcation, which requires us to look at the full $\Gamma_0$-isotypic characteristic equation for all $\alpha$. Substituting $\lambda = r$, we obtain
\[
r(1-\gamma e^{-r\tau_1}) -\alpha b \frac{e^{-r\tau_2}-1}{r}=-a_j.
\]
If $r>0$, then $e^{-r\tau_1}<1$ and $e^{-r\tau_2} < 1$. Therefore, if $a_j>0$ and $\gamma<1$, then once again, the left-hand side is strictly positive, and so no solutions exist. Notice that this implies that $a_j>0$ and $\gamma<1$ are actually sufficient to guarantee that steady-state bifurcation \emph{never} occurs for any $\alpha >0$.
\end{proof}
This also suffices to prove Theorem \ref{thm:neutral:mas-asymp}.

\subsection{Transversality}
We will now attain some conditions analogous to standard Hopf bifurcation transversality conditions, which will be used later to aid in computing the crossing numbers. Here it is more convenient to use the characteristic quasipolynomial. Taking $\lambda(\alpha) = u(\alpha) + iv(\alpha)$ and substituting into \eqref{eq:neutral:char-root}, we obtain the system
\begin{align*}
(u^2 - v^2)(1-\gamma e^{-u\tau_1}\cos(v\tau_1)) - 2\gamma u v \sin(v\tau_1) + a_ju- \alpha be^{-u\tau_2}(\cos(v\tau_2) -1 ) &= 0\\
(u^2-v^2)\gamma e^{-u\tau_1}\sin(v\tau_1) + 2uv(1-\gamma \cos(v\tau_1)) + a_jv +\alpha b e^{-u\tau_2}\sin(v\tau_2)&= 0
\end{align*}
\vs
\noindent
Differentiating with respect to $\alpha$, setting $u=0$ and $v = \beta$, and factoring out $u' := \tfrac{du}{d\alpha}$ and $v' := \tfrac{dv}{d\alpha}$, we obtain
\begin{subequations}\label{eq:neutral:u'v'}
\begin{align}
&\begin{cases}
\begin{aligned}
&u'(2\beta\gamma\sin(\beta\tau_1) - \beta^2\gamma\tau_1 \cos(\beta\tau_1) + a_j - \alpha b \tau_2 \cos(\beta\tau_2))\\+ &v'(2\beta\gamma\cos(\beta\tau_1)-\beta^2\gamma\tau_1\sin(\beta\tau_1)-2\beta+\alpha b\tau_2\sin(\beta\tau_2)) = b(\cos(\beta\tau_2)-1)
\end{aligned}
\end{cases}\\ \nonumber
\\ 
&\begin{cases}
\begin{aligned}
&u'(-2\beta\gamma\cos(\beta\tau_1) + \beta^2\gamma\tau_1\sin(\beta\tau_1)+2\beta -\alpha b \tau_2 \sin(\beta \tau_2))
\\+&v'(2\beta\gamma\sin(\beta\tau_1) - \beta^2\gamma\tau_1 \cos(\beta\tau_1) + a_j - \alpha b \tau_2 \cos(\beta\tau_2)) = -b\sin(\beta\tau_2)
\end{aligned}
\end{cases}
\end{align}
\end{subequations}
Put
\begin{align*}
p &:= 2\beta\gamma\sin(\beta\tau_1) - \beta^2\gamma\tau_1 \cos(\beta\tau_1) + a_j - \alpha b \tau_2 \cos(\beta\tau_2)\\
q &:= 2\beta\gamma\cos(\beta\tau_1)-\beta^2\gamma\tau_1\sin(\beta\tau_1)-2\beta+\alpha b\tau_2\sin(\beta\tau_2)
\end{align*}
Then \eqref{eq:neutral:u'v'} can be written as the linear system
\begin{align*}
u'p + v'q &= b(\cos(\beta\tau_2) -1)\\
-u'q + v'p &= -b\sin(\beta\tau_2)
\end{align*}
which can be solved to obtain
\begin{equation}\label{eq:neutral:pq}
\tfrac{\partial}{\partial \alpha}P_j(\alpha,i\beta)=u' = \frac{pb(\cos (\beta\tau_2) -1) + qb \sin (\beta\tau_2)}{p^2+q^2}
\end{equation}
Put
\begin{equation}\label{eq:neutral:rho}
\rho_j(\alpha,\beta) := pb(\cos (\beta\tau_2) -1) + qb \sin (\beta\tau_2)
\end{equation}
This is sufficient to state and prove the local bifurcation theorem in the NFDE context (and hence to prove Theorem \ref{thm:neutral:mas-local}). 
\begin{theorem}\label{thm:neutral:nfde-local}
Let $(\alpha_0,\beta_0)\in \Lambda$. If $\rho_j(\alpha_0,\beta_0)\neq 0$, then there exists a connected branch of non-constant periodic solutions bifurcating from the trivial branch at $\alpha=\alpha_0$.
\end{theorem}
\begin{proof}
If $\rho_j(\alpha_0,\beta_0)\neq 0$, then by Lemma \ref{lem:prelim:sign-crossing-num}, this implies the corresponding $\Gamma_0$-isotypic crossing number $\mathfrak t_j(\alpha_0,\beta_0)\neq 0$, which in turn implies that the local bifurcation invariant $\omega_G(\alpha_0,\beta_0)\neq 0$ (cf. Section \ref{sec:prelim:local-bif-twisted}).  
\end{proof}

As was shown in Section \ref{sec:prelim:krasnoselskii}, the condition that $\rho_j(\alpha_0,\beta_0)\neq 0$ is sufficient to guarantee that the associated crossing number is non-zero, implying that the local bifurcation invariant is non-trivial, and thus local bifurcation is obtained through the equivariant Krasnosel'skii theorem. We
\vs
Then the sign of $\rho_j$ evaluated at a critical point $(\alpha_{n,j},\beta_{n,j},0)$ is sufficient to determine the $k$-resonant $\Gamma$-isotypic crossing number at that point. Depending on the nature of the application, there are three primary ways in which $\rho_j$ might be used:
\begin{enumerate}
    \item We are only interested in the first Hopf bifurcation where the trivial consensus equilibrium loses stability, denoted $(\widehat\alpha, \widehat \beta)\in \Lambda$ (cf. Thm. \ref{thm:neutral:trivial-consensus-stability}), and $\rho_j(\widehat \alpha, \widehat \beta)$ can be evaluated numerically.
    \item We are interested in showing that some bifurcating branch is global and unbounded over some finite parameter range $\alpha\in(\alpha_1,\alpha_2)$ which is physically meaningful. Since $\Lambda$ is discrete, $\Lambda\cup(\alpha_1,\alpha_2)\times\mathbb R_+ \times \mathscr E$ is finite, and so $\sign\rho_j(\alpha,\beta)$ can be numerically computed and summed across this parameter range.
    \item We are interested in showing that bifurcating branches are unbounded over an arbitrary parameter range, or one which is not known \emph{a priori}, or we want guarantees that \emph{all} branches are unbounded, etc.
\end{enumerate}
In the third case, we require more analytical results on the sign of $\rho_j$ at critical points. Unlike the case for the pseudoneutral equation, here it is much more difficult to obtain any kind of conditions for $\rho_j$ to have constant sign for generic choices of $\tau_1$ and $\tau_2$. Note that although $\rho_j$ is quasiperiodic, it is actually only evaluated at the critical values $(\alpha,\beta)\in\Lambda$ which satisfy \eqref{eq:neutral:beta-relation} and \eqref{eq:neutral:alpha-relation}. We can provide strong criteria for the special case when $\tau_1 = \tau_2$. We will also transform $\rho_j$ into an expression which is more useful both for numerical analysis and which allows us to obtain conditions for $\beta$ sufficiently large.

\begin{proposition}
    If $\tau_1=\tau_2, a_j>0$, and $b>0$, then $\rho_j(\alpha_{n,j},\beta_{n,j})<0$ for all $(\alpha_{n,j},\beta_{n,j},0)\in \Lambda$.
\end{proposition}
\begin{proof}
If $\tau :=\tau_1=\tau_2$, then all relevant equations have only one period. Moreover, \eqref{eq:neutral:beta-relation} can be satisfied only near the beginning of this period, assuming $a_j > 0$. To see this, notice that for $\beta\tau \approx \pi(2k-1) +\varepsilon$, for $\varepsilon>0$ sufficiently small and $k\in \mathbb N$, \eqref{eq:neutral:beta-relation} satisfies
\[
\frac{(\beta\gamma\sin(\beta\tau)-a_j)(1-\cos(\beta\tau))}{(1-\gamma\cos(\beta\tau_1))(\sin(\beta\tau_2)}\approx \frac{(-\beta\epsilon - a_j)(2-\tfrac{1}{2}\varepsilon^2)}{(1+\gamma-\gamma\tfrac{1}{2}\varepsilon^2)(-\epsilon)} \approx \frac{2a_j}{1+\gamma}.
\]
Substituting $\beta\tau = \pi(2 k-1)+\varepsilon$ into $\rho_j$ and again taking the small angle approximation, we obtain
\begin{align*}
&(2\beta\gamma\epsilon - \beta\gamma(\pi(2k-1)+\varepsilon)(-1+\tfrac{1}{2}\varepsilon^2) + a_j + \alpha b \tau(1-\tfrac{1}{2}\varepsilon^2))(-2+\tfrac{1}{2}\varepsilon^2)\\
+&(2\beta\gamma(-1+\tfrac{1}{2}\varepsilon^2) - \beta\gamma(\pi(2k-1)+\varepsilon)(-\varepsilon)-2\beta -\alpha b\tau \varepsilon)(-\varepsilon).
\end{align*}
Discarding all terms of order $\varepsilon$ or smaller, this yields
\[
-2(\beta\gamma\pi(2k-1) +a_j +\alpha b \tau), 
\]
which is negative for all $a_j>0$ and $b>0$. 
\end{proof}
\begin{proposition}
    Fix $j=0,1,\dots,r$. If $\tau_1/\tau_2 \in \mathbb Q$, and if $(\widetilde \beta_{l,j})_{l=1}^P$ is a $P$-periodic sequence of roots of 
    \[
    1 = \frac{(\gamma\sin(\beta\tau_1)-\varepsilon)(1-\cos(\beta\tau_2))}{(1-\gamma\cos(\beta\tau_1))\sin(\beta\tau_2)},
    \]
    where $\varepsilon>0$ is taken sufficiently small, such that $(\widetilde \beta_{l,j})_{l=1}^P$ satisfies
    \[
    \sum_{l=1}^P\sign \Upsilon_j(\widetilde \beta_{l,j}) \neq 0
    \]
    where 
    \[
    \Upsilon_j(\beta):= \gamma\tau_1\cos(\beta\tau_1) + \tau_2\cos(\beta\tau_2) -\gamma\tau_1\cos(\beta(\tau_1-\tau_2))-\gamma\tau_2\cos(\beta(\tau_1+\tau_2)),
    \]
    or if there exists some $K>0$ such that $\sign \rho_j (\alpha_{K+i,j},\beta_{K+i,j}) = \sign \Upsilon_j(\widetilde \beta_{l,j})$ for all $l\in \mathbb N$, where indices of $\widetilde \beta_{l,j}$ are taken mod $P$, and
    \[
    \sum_{i=1}^K\sign \rho_j(\alpha_{i,j},\beta_{i,j}) \neq 0,
    \]
    then
    \[
    \sum_{(\alpha_{n,j},\beta_{n,j},0)\in \Lambda}\sign \rho_j(\alpha_{n,j},\beta_{n,j}) \neq 0.
    \]
\end{proposition}
\begin{proof}
Since we are only interested in the value of $\rho_j$ at critical points, we can substitute \eqref{eq:neutral:alpha-relation} and \eqref{eq:neutral:beta-relation} to obtain a more useful expression solely in terms of $\beta$, which agrees with $\rho_j$ when evaluated at critical points. First, we substitute $\alpha \mapsto \tfrac{\beta^2\gamma\sin(\beta\tau_1)-\beta a_j}{b\sin(\beta\tau_2)}$ in the expressions for $p$ and $q$, obtaining
\begin{align*}
    \widetilde p &:= 2\beta\gamma\sin(\beta\tau_1) - \beta^2\gamma\tau_1 \cos(\beta\tau_1) + a_j - \tau_2 \cos(\beta\tau_2)\frac{(\beta^2\gamma\sin(\beta\tau_1) - \beta a_j)}{\sin(\beta\tau_2)}\\
\widetilde q &:= -2\beta(1-\gamma\cos(\beta\tau_1))-\beta^2\gamma\tau_1\sin(\beta\tau_1)+ \tau_2(\beta^2\gamma\sin(\beta\tau_1) - \beta a_j)
\end{align*}
Then we have the expression
\begin{equation}\label{eq:neutral:widetilde-pq}
\widetilde pb(\cos(\beta\tau_2)-1) + \widetilde qb(\sin(\beta\tau_2)),
\end{equation}
which is equal to $\rho_j(\alpha_{n,j},\beta_{n,j})$ for any critical $(\alpha_{n,j},\beta_{n,j})\in \Lambda$. This eliminates dependence on $\alpha$ but now has an undesirable asymptote in the $\widetilde p$ term. We can rearrange \eqref{eq:neutral:beta-relation} to obtain
\[
\beta(1-\gamma\cos(\beta\tau_1)) = \frac{(\beta\gamma\sin(\beta\tau_1)-a_j)(1-\cos(\beta\tau_2))}{\sin(\beta\tau_2)},
\]
which can be substituted into \eqref{eq:neutral:widetilde-pq} to obtain (after dividing out $b$)
\begin{align*}
\widetilde \rho_j(\beta) :=\;\; &(2\beta\gamma\sin(\beta\tau_1) - \beta^2\gamma\tau_1 \cos(\beta\tau_1) + a_j)(\cos(\beta\tau_2)-1)\\
+&(-2\beta(1-\gamma\cos(\beta\tau_1)) - \beta^2\gamma\tau_1 \sin(\beta\tau_1) + \tau_2(\beta^2\gamma\sin(\beta\tau_1) - \beta a_j))\sin(\beta\tau_2)\\
+&\beta^2\tau_2 \cos(\beta\tau_2)(1-\gamma\cos(\beta\tau_1)).
\end{align*}
This expression, though unwieldy, depends only on $\beta$ and possesses no vertical asymptotes. Importantly, this means that the terms proportional to $\beta^2$ dominate for large values of $\beta$. Gathering these terms and performing trigonometric reduction, we obtain
\begin{equation}\label{eq:neutral:upsilon}
\Upsilon_j(\beta):
= \gamma\tau_1\cos(\beta\tau_1) + \tau_2\cos(\beta\tau_2) -\gamma\tau_1\cos(\beta(\tau_1-\tau_2))-\gamma\tau_2\cos(\beta(\tau_1+\tau_2))
\end{equation}
Unlike $\widetilde \rho_j$, this is a periodic function of $\beta$ (provided $\tau_1/\tau_2 \in \mathbb Q$). On the other hand, note that \eqref{eq:neutral:beta-relation} can be written
\[
1 = \frac{(\gamma\sin(\beta\tau_1)-\tfrac{a_j}{\beta})(1-\cos(\beta\tau_2))}{(1-\gamma\cos(\beta\tau_1))\sin(\beta\tau_2)}.
\]
Put $\varepsilon\approx o(\tfrac{a_j}{\beta})$. Then as $\varepsilon\to0$, solutions to the periodic equation
\begin{equation}\label{eq:neutral:largebetalimit}
1 = \frac{(\gamma\sin(\beta\tau_1)-\varepsilon)(1-\cos(\beta\tau_2))}{(1-\gamma\cos(\beta\tau_1))\sin(\beta\tau_2)}
\end{equation}
converge to solutions of \eqref{eq:neutral:beta-relation} taken mod $T$, where $T$ is the period of \eqref{eq:neutral:largebetalimit} (which is finite if $\tau_!/\tau_2\in\mathbb Q$). These two facts taken together imply that if $\tau_1/\tau_2 \in \mathbb Q$, then eventually, for $\beta$ large enough, its crossing numbers will repeat periodically. Therefore, if a representative periodic sequence of solutions $(\widetilde \beta_{l,j})_{l=1}^P$ to \eqref{eq:neutral:largebetalimit} can be found such that 
\[
\sum_{l=1}^P \sign\Upsilon_j(\widetilde \beta_{l,j}) \neq 0, 
\]
then there exists a global branch of solutions unbounded for all $\alpha>0$. On the other hand, if 
\[
\sum_{l=1}^P \sign\Upsilon_j(\widetilde \beta_{l,j}) = 0, 
\]
but there exist $\beta_{0,j},\beta_{1,j},\dots,\beta_{K,j}$ such that for all $n>K$, $\sign\rho_j(\alpha_{n,j},\beta_{n,j})$ forms a $P$-periodic sequence, i.e. for all $n>K$ there exists an $l'=1,\dots,P$ with $\sign \rho_j(\alpha_{n,j},\beta_{n,j}) = \sign \Upsilon_j(\beta_{n,j} \text{ (mod T)}) = \sign \Upsilon_j(\widetilde \beta_{l',j})$, and if
\[
\sum_{n=1}^K \sign \rho_j(\alpha_{n,j},\beta_{n,j})\neq 0,
\]
then there also exists a global branch of solutions unbounded for all $\alpha>0$. Note the above statement can be put more intuitively as saying that if the periodic sequence sums to zero, but for the first few critical values of $\beta$ before the periodic behavior begins, these values do not sum to zero, then the full sum taken over all critical points must also have a nonzero sum.
\vs
If $\tfrac{\tau_1}{\tau_2}=\tfrac{c}{d}$, where $\tfrac{c}{d}$ is assumed to be in reduced form with $\gcd(c,d)=1$, then both \eqref{eq:neutral:largebetalimit} and \eqref{eq:neutral:upsilon} are $T$-periodic, where
\[
T = \frac{2\pi c}{\tau_1} = \frac{2\pi d}{\tau_2}.
\]
Therefore, the general procedure for obtaining unboundedness over all $\alpha$ is to numerically solve \eqref{eq:neutral:largebetalimit} for $\beta\in(0,T)$ to obtain the sequence $(\widehat \beta_{l,j})_{l=1}^P$, and then compute
\[
\sum_{l=1}^P \sign \Upsilon_j(\widetilde \beta_{l,j}).
\]
In the degenerate case where this sum equals zero, one can then compute the first few critical points and substitute them into $\rho_j$ until the periodic pattern begins to emerge, and check this sum (numerical tests show that typically this periodic behavior emerges very quickly). 
\vs
However, there is a final important caveat: The equation \eqref{eq:neutral:largebetalimit} must be solved (numerically) with $\varepsilon>0$. In the limit where $\varepsilon =0$, certain degeneracies appear when taking roots of \eqref{eq:neutral:largebetalimit} and applying them to $\Upsilon_j$ which are not present when taking actual critical points from \eqref{eq:neutral:beta-relation} and applying them to $\rho_j$, due to perturbation by lower order terms. In particular, numerical experiments show that solutions to \eqref{eq:neutral:largebetalimit} with $\varepsilon = 0$ produce constant values in their evaluation under $\Upsilon_j$. It is sufficient to take $\varepsilon \ll a_j$ to ensure that the behavior at high values of $\beta$ is being faithfully captured. 
\end{proof}
\begin{remark}
Since the dependence of solutions to \eqref{eq:neutral:beta-relation} on $a_j$ flattens asymptotically as $\beta \to \infty$, for $\varepsilon\ll a_j$ sufficiently small (and hence for $\beta$ sufficiently large), all solutions to \eqref{eq:neutral:beta-relation}, when taken mod $T$ as described above, will converge to the same periodic sequence of solutions of \eqref{eq:neutral:largebetalimit} for every $\Gamma$-isotypic component. This means that the boundedness or unboundedness of branches is generally the same across $\Gamma$-isotypic components---except in the case where all but finitely many initial crossing numbers cancel periodically---and $a_j$ primarily determines on which $\Gamma$-isotypic component bifurcation first occurs (and hence what types of symmetries are guaranteed for the first stability-exchanging bifurcation). 
\end{remark}
\vs
Now we will reformulate the system \eqref{eq:neutral:basic-sys} in a functional setting where the $S^1$-action is well-defined, which allows us to view Hopf bifurcation of periodic multiconsensus solutions as a two-parameter bifurcation problem which can be addressed using the equivariant twisted Nussbaum-Sadovskii degree.
\vs
\section{Functional space reformulation}\label{sec:neutral:funcspace}

As before, we first period-normalize system \eqref{eq:neutral:basic-sys}. We do this by first assuming that solutions exist with some arbitrary period $p>0$. We convert this to a frequency $\beta := \tfrac{2\pi}{p}$, and rescale time to $\widetilde t := \beta t$, so a function which is $p$-periodic in $t$ is $2\pi$-periodic in $\widetilde t$. We then set $u(t) := \bm x(\widetilde t) = \bm x(\beta t)$. Note that we are freeing the notation $u(t)$ to exclusively refer to the period normalized function of time, disregarding its other use in the previous section. This yields the period-normalized system
\begin{equation}\label{eq:neutral:basic-sys-norm}
\frac{d}{dt}\left[u- \bm g(u(t-\beta \tau_1))\right] = -\frac{a}{\beta}u - \frac{\alpha}{\beta}\bm f\left(\int_0^{\tau_2} u(t-\beta s)ds\right) - \frac{1}{\beta}\bm h(u) 
\end{equation}

We define the space $\mathscr E$, the operators $L, \bm j, N_{\bm f}$, and $N_{\bm k}$ as in Section \ref{sec:pseudo:introduction:overview} and, in parallel with the functional-analytic framework of Chapter \ref{chapter:pseudo}, Section \ref{sec:pseudo:funcspace}, we put
\begin{align*}
N_{\bm k}(\alpha,\beta,u) &= \bm g(u(t-\beta \tau_1))\\
N_{\bm f}(\alpha,\beta,u) &= -\frac{a}{\beta}u - \frac{\alpha}{\beta}\bm f\left(\int_0^{\tau_2} u(t-\beta s)ds\right) - \frac{1}{\beta}\bm h(u) 
\end{align*}
\vs
which yields the operator equation
\begin{equation}\label{eq:neutral:op-sys}
\mathscr F(\alpha,\beta,u) = u - N_{\bm k}(\alpha,\beta,u) - (L + \bm j)^{-1}\Big(N_{\bm f}(\alpha,\beta,u) + j(u) - N_{\bm k}(\alpha,\beta,\bm j(u))\Big)
\end{equation}
Recall $\Gamma := \Gamma_0\times \mathbb Z_2$, and put $G:=S^1\times \Gamma$. Then by \ref{neutral:c1}--\ref{neutral:c3}, $\mathscr F$ is $G$-equivariant, where $S^1$ acts naturally on $\mathscr E$ by shifting time. Since $\bm f$ is continuous and $\int_0^{\tau_2}u(t-\beta s)ds$ is a compact operator, their composition is compact, and so the Nemytskii operator $N_{\bm f}$ is compact. Because $\bm g$ is $\kappa$-Lipschitzian with $\kappa<1$, $N_{\bm k}$ is condensing. Note that the sum of a condensing operator and a compact operator is condensing, and so because $L+j$ is an isomorphism, $\mathscr{F}:\bbR^2_+ \times \mathscr{E} \to \mathscr{E}$ is a $G$-equivariant condensing perturbation of identity.

\vs
A $\tfrac{2\pi}{\beta_0}$-periodic function $x(t)$ is a solution of \eqref{eq:neutral:basic-sys} for $\alpha=\alpha_0$ if and only if $\mathscr F(\alpha_0,\beta_0,u) = 0$, where $u(t) = x(\beta_0 t)$. Note that by \ref{neutral:c1}---\ref{neutral:c3}, $\mathscr F$ is also Fr\'echet differentiable at $u=0$, and we put $\mathscr A(\alpha,\beta) := D_u \mathscr F(\alpha,\beta,0)$. Analogously to Chapter \ref{chapter:pseudo}, we define the continuous memory operator and the momentum memory operator (which acts as a translation operator)
\begin{align*}
K_{\tau_2}u(t) &:= \int_0^{\tau_2}u(t-\beta s)ds,\\
\widehat K_{\tau_1}u(t) &:= u(t-\beta\tau_1).
\end{align*}
Note that both $K$ and $\widehat K$ are bounded linear operators, and 
\begin{align*}
D_u N_{\bm k}(\alpha,\beta,0) &= \gamma \widehat K,\\
D_uN_{\bm f}(\alpha,\beta,0) &= -\frac{a}{\beta}\id -\frac{\alpha b}{\beta}K_{\tau_2}- \frac{1}{\beta}C.
\end{align*}
Then we have
\[
\begin{aligned}
\mathscr A(\alpha,\beta)u &= u - \gamma \widehat K_{\tau_1}u - (L+\bm j)^{-1}\left(-\frac{\alpha b}{\beta}K_{\tau_2}u - \frac{1}{\beta}Cu - \Big(\frac{a}{\beta} - 1\Big)u- \gamma \widehat K_{\tau_1}u\right)  \\
&= u - \gamma \widehat K_{\tau_1}u + (L+\bm j)^{-1}\left(\frac{\alpha b}{\beta}K_{\tau_2}u + \frac{1}{\beta}Cu + \Big(\frac{a}{\beta} - 1\Big)u+ \gamma \widehat  K_{\tau_1}u\right).
\end{aligned}
\]
Since \eqref{eq:neutral:op-sys} is $G$-equivariant, $\mathscr A$ is a $G$-equivariant linear operator. Its eigenspaces correspond to $G$-isotypic components, which are invariant under $\mathscr A$, and we put $\mathscr A_{k,j}(\alpha,\beta) := \mathscr A_{|\mathscr E_{k,j}}(\alpha,\beta):\mathscr E_{k,j} \to \mathscr E_{k,j}$. 
\vs
We now need to relate $\mathscr A$ to the characteristic equation of \eqref{eq:neutral:basic-sys-lin} and thereby bridge the gap between imaginary roots of $\triangle_\alpha(i\beta)=0$ (which correspond, in a classical Hopf bifurcation context, to bifurcating periodic solutions with limit frequency $\beta$), and the period-normalization frequency parameter $\beta$ as it appears here. First, we note that linearizing the system \eqref{eq:neutral:basic-sys-norm} around $u=0$ yields the linear period normalized system
\[
\frac{d}{dt}\left[u- \gamma(u(t-\beta \tau_1))\right] +\frac{a}{\beta}u + \frac{\alpha b}{\beta}\int_0^{\tau_2} u(t-\beta s) - \frac{1}{\beta}\bm C =0,
\]
whose characteristic equation is clearly $\widetilde \triangle_{\alpha}(\lambda) = \tfrac{1}{\beta}\triangle_\alpha(\lambda)$. Pick $k>0$ and let $u(t) =e^{ikt}v\in \mathscr E_{k,j}$. Then
\begin{align*}
    \widehat K_{\tau_1}u(t) &= e^{ik\beta\tau_1}e^{ikt}v\\
    K_{\tau_2}u(t) &= \frac{1-e^{-ik\beta\tau_2}}{ik\beta}e^{ikt}v\\
    (L+\bm j) u &= (ik+1)e^{ikt}v \Rightarrow (L+\bm j)^{-1}e^{ikt}v = \frac{1}{ik+1}e^{ikt}v.
\end{align*}
Then we obtain
\[
\mathscr A(\alpha,\beta)(e^{ikt}v) = \frac{(ik+1)(1-\gamma e^{-ik\beta\tau_1}) + \tfrac{\alpha b}{\beta}\tfrac{1-e^{ik\beta\tau_2}}{ik\beta}+\tfrac{1}{\beta}\mu_j + \frac{a}{\beta}-1 + \gamma e^{-ik\beta\tau_1}}{ik+1}e^{ikt}v.
\]
Multiplying the numerator and denominator by $\beta$ and canceling terms yields
\[
\mathscr A(\alpha,\beta)(e^{ikt}v) = \frac{ik\beta(1-\gamma e^{-ik\beta\tau_1}) + \alpha b\tfrac{1-e^{ik\beta\tau_2}}{ik\beta}+\mu_j + a}{\beta + ik\beta}e^{ikt}v.
\]
But note that the numerator of this expression is exactly $\triangle_{\alpha,j}(ik\beta)$, and so we have
\[
\mathscr A(\alpha,\beta)e^{ikt}v = \frac{1}{\beta+ik\beta}\triangle_\alpha(ik\beta)e^{ikt}v.
\]
On the other hand, for $k=0$, take $u(t) = v \in \mathscr E_{0,j}$. Then
\begin{align*}
    \widehat K_{\tau_1}v &= v\\
    K_{\tau_2}v &= v\\
    (L+\bm j) v &=  (L+\bm j)^{-1}v = v,
\end{align*}
and the same computation gives
\[
\mathscr A(\alpha,\beta)v = (1-\gamma)v + \left(\frac{\alpha b}{\beta}\tau_2 + \frac{1}{\beta}C + \frac{a}{\beta} - 1 + \gamma\right)v = \frac{1}{\beta}(\alpha b \tau_2 +\mu_j +a)v,
\]
and so
\[
\mathscr A(\alpha,\beta)v = \frac{1}{\beta}\triangle_{\alpha,j}(0)v  = \left(\frac{\alpha b\tau_2 + a_j}{\beta}\right)v.
\]
This gives us the following forms for $\mathscr A_{k,j}(\alpha,\beta)$:
\begin{align}
\mathscr A_{k,j}(\alpha,\beta) &= \frac{1}{ik\beta + \beta}\triangle_{\alpha,j}(ik\beta)\id_{\mathscr E_{k,j}}\quad &k>0\label{eq:neutral:char-Ak},\\
\mathscr A_{0,j}(\alpha,\beta) &= \left(\frac{\alpha b \tau_2 +a_j }{\beta}\right)\id_{\mathscr E_{0,j}}\quad &k=0.\label{eq:neutral:char-A0}
\end{align}

\section{Two-parameter bifurcation}\label{sec:neutral:two-parameter}
Notice that a purely imaginary root of the original characteristic equation $\triangle_\alpha(i\beta)=0$ corresponds to a zero eigenvalue of $\mathscr A$, and thus a parameter value where $\mathscr A(\alpha,\beta)$ fails to be an isomorphism. This is the natural setting which allows us to view Hopf bifurcation as a two-parameter bifurcation problem. We can therefore use the Krasnosel'skii/Rabinowitz paradigm established in previous chapters to show the existence of bifurcating non-constant periodic solutions. Here, we will show this in the absence of certain degeneracies. These degenerate cases are then addressed using fixed point reduction techniques.
\vs
The functional space reformulation above clearly gives us the following facts, in parallel with Chapter \ref{chapter:pseudo}:
\begin{enumerate}
    \item $\mathscr F$ is a $G$-equivariant condensing perturbation of identity.
    \item $\mathscr F(\alpha,\beta,0)=0$ for all $\alpha,\beta \in \mathbb R^2_+$ (i.e. the trivial branch exists). 
    \item $\mathscr F$ is Fr\'echet differentiable at $u=0$ for all $(\alpha,\beta)\in\mathbb R^2_+$, and $\mathscr A(\alpha,\beta):= D_u\mathscr F(\alpha,\beta,0)$ depends continuously on $\alpha$ and $\beta$ (i.e. the linearization $\mathscr A$ exists and is well-defined everywhere on the trivial branch, and depends continuously on parameters).
\end{enumerate}
\vs
Under these conditions, we can apply the equivariant Krasnosel'skii and Rabinowitz theorems to show local bifurcation of branches of solutions and their global continuation, using the spectral data and transversality conditions already obtained. We use the same definitions and notation as before for the trivial and nontrivial sets, i.e. 
\begin{align*}
\mathcal M &:= \{(\alpha,\beta,0) \in \bbR^2_+ \times \mathscr E:\quad \mathscr F(\alpha,\beta,0)=0\}\\
\mathscr S &:= \{(\alpha,\beta,u)\in\bbR^2_+ \times \mathscr E:\quad\mathscr F(\alpha,\beta,0)=0,\quad u\not\equiv0\}
\end{align*}
Bearing in mind the definitions of bifurcation points, branching points, and branches of solutions, we obtain another definition of the critical set: A point $(\alpha_0,\beta_0,0) \in \mathbb R^2_+\times \mathscr E$ is called a critical point of \eqref{eq:neutral:op-sys} if $\mathscr A(\alpha_0,\beta_0)$ fails to be an isomorphism. This implies that for some $k\geq0$, $j=0,1,\dots,r$, $\mathscr A_{k,j}$ fails to be an isomorphism, which in turn implies that $\det_\mathbb C\triangle_{\alpha_0,j}(ik\beta)=0$ for some $k>0$, or that $\alpha b \tau_2 + a_j = 0$. We therefore put
\[
\Lambda_0 := \{(\alpha_0,\beta_0,0)\in \mathbb R^2_+\times \mathscr E: \mathscr A(\alpha_0,\beta_0)\text{ is not an isomorphism}\}.
\]
Clearly $\Lambda \subseteq \Lambda_0$, and moreover, for any $(\alpha_0,\beta_0,0)\in \Lambda$ and any $k\in \mathbb N$, we have $(\alpha_0,\beta_0/k,0)\in \Lambda_0$. However, we are not interested in steady-state bifurcation, and such points where $\mathscr A_{0,j}$ is not an isomorphism are degeneracies, from our point of view. So we filter these out to obtain a refined critical set
\[
\widetilde \Lambda := \{(\alpha,\beta,0) \in  \mathbb R^2_+\times \mathscr E:\;\; \exists_{k \in {\mathbb N}} \;\;{\det}_{\mathbb C} \triangle_\alpha(ik\beta) = 0 \text{  \;and\;  } \alpha b \tau_2 + a_j \not= 0\}
\]
Since $\widetilde{\Lambda}$ is discrete, for some sufficiently small $\varepsilon,\delta > 0$, one can find an \textit{isolated $G$-invariant neighborhood} $\Omega(\alpha_0,\beta_0,0)$ given by
\[
\Omega(\alpha_0,\beta_0,0) := \{(\alpha,\beta,u) \in \bbR^2_+ \times \mathscr E:\;\; |(\alpha,\beta)-(\alpha_0,\beta_0)| < \delta,\;\; ||u|| < \varepsilon\}
\]
such that
\[
\Omega(\alpha_0,\beta_0,0) \cap \overline{\mathscr S} = \emptyset.
\]
\vs
As before, we construct a continuous $G$-equivariant auxiliary function $\eta:\bbR^2_+ \times \mathscr E \to \bbR$ satisfying
\[
\begin{cases} \eta(\alpha,\beta,0)<0 &\text{ if } \; |(\alpha,\beta)-(\alpha_0,\beta_0)|  = \delta,\\
\eta(\alpha,\beta,v)>0 &\text{ if } \; |(\alpha,\beta)-(\alpha_0,\beta_0)|\le \delta \text{ and }\; \|v\|=\varepsilon,\\
\end{cases}
\]
\vs
and define $\mathscr F_\eta : \overline{\Omega(\alpha_0,\beta_0,0)} \to \mathscr \bbR \times \mathscr E$ by 
\[
\mathscr F_\eta(\alpha,\beta,u) = (\eta(\alpha,\beta,u),\mathscr F(\alpha,\beta,u)).
\]
Then $(\mathscr F_\eta,\Omega(\alpha_0,\beta_0,0))$ is an admissible $G$-pair, and we can define the local bifurcation invariant $\omega_G(\alpha_0,\beta_0)\in A_1^t(G)$ by
\[
\omega_G(\alpha_0,\beta_0) := G\text{-deg}(\mathscr F_\eta, \Omega(\alpha_0,\beta_0,0)),
\]
where $G$-deg is the twisted $G$-equivariant Nussbaum-Sadovskii degree. 
\subsection{Crossing numbers}
As before, the crossing number $\mathfrak t_{k,j}(\alpha_0,\beta_0)$ is defined as the net algebraic count of eigenvalues, counted with their $G$-isotypic multiplicity, which pass from the positive half-plane to the negative half-plane as $\alpha$ passes through $\alpha_0$. It is easy to see from this definition that, at a $G$-isotypically simple critical point $(\alpha_0,\beta_0,0)\in \widetilde \Lambda$ (i.e. where eigenvalues cross on only one $G$-isotypic component $\mathscr E_{k,j}$), we have
\[
\mathfrak t_{k,j}(\alpha_0,\beta_0) = -\rho_j(\alpha_0,k\beta_0)m_j,
\]
where $\rho_j(\alpha,\beta)$ is given by \eqref{eq:neutral:rho} and $m_j = \tfrac{\dim \mathscr E_{k,j}}{\dim \mathcal W_{k,j}}$ is the $G$-isotypic multiplicity. One should take care to note that the above formula holds only if $(\alpha_0,k\beta_0,0)\in \Lambda$. On the other hand, if a critical point is not $G$-isotypically simple, which can occur if $(\alpha_0,\beta_0,0)\in \lambda$ and $(\alpha_0,k\beta_0,0)\in \Lambda$ for some $k>1$, or if $\det_\mathbb C \triangle_{\alpha_0,j_1}(ik\beta_0) = 0 = \det_\mathbb C \triangle_{\alpha_0,j_2}(ik\beta_0)$ for some $j_1\neq j_2$, then put $\widehat \Lambda(\alpha_0) = \{(\beta,k,j)\in \mathbb R_+\times\mathbb N\times\{0,1,\dots,r\}: {\det}_{\mathbb C}\triangle_{\alpha_0,j}(ik\beta) = 0 \}$. Then the crossing number can be computed through the sum 
\[
\mathfrak t_{k,j}(\alpha_0,\beta_0) = \sum_{\beta_0,k,j}-\rho_j(\alpha_0,k\beta_0)m_j.
\]
However, we generally use fixed point reduction to deal with these degeneracies.
\vs
We have the following formula for computing the local bifurcation invariant:
\[
\omega_G(\alpha_0,\beta_0) = \prod_{j=0}^r\Gamma\text{-deg}(\mathscr A_{0,j}({\alpha_0,\beta_0),B(V))}\cdot\sum_{k=1}^\infty\sum_{j=0}^r \mathfrak t_{k,j}(\alpha_0,\beta_0)\deg_{\mathcal W_{k,j}},
\]
where $\Gamma$-deg is the $\Gamma$-equivariant Nussbaum-Sadovskii degree, which detects solutions on the zero mode, and $\deg_{\mathcal W_{k,j}}$ is the twisted basic degree on the irreducible $G$-representation $\mathcal W_{k,j}$ on which the $G$-isotypic component $\mathscr E_{k,j}$ is modeled. We can then apply the equivariant Krasnosel'skii theorem and Rabinowitz alternative to show bifurcation and global unbounded continuation of branches of periodic solutions to \eqref{eq:neutral:basic-sys}, in the absence of degeneracies. To address these degenerate cases, and to guarantee that the detected periodic solutions are non-constant, we require fixed point reduction.

% \begin{remark}
%     Compared to classical Hopf bifurcation theorems, the degeneracy conditions for equivariant Hopf bifurcation using this method are very relaxed. We do not require eigenvalues to be simple, nor even $\Gamma$-isotypically simple, and it is still possible to show existence, global continuation, and spatio-temporal symmetric classification of branches of periodic solutions even when eigenvalues cross on several $\Gamma$-isotypic components simultaneously. 
% \end{remark}

\section {Fixed point reduction}\label{sec:neutral:fixed-point-reduction}
It's important to note that many situations which would be considered degenerate in classical Hopf bifurcation theories are not problematic for two-parameter twisted degree methods such as this. For example, eigenvalues failing to be simple, or even failing to be $G$-isotypically simple, do not pose any problems for these methods. It is entirely possible that for a critical point $(\alpha_0,\beta_0,0)\in \Lambda$, we have  ${\det}_{\mathbb C}\triangle_{\alpha_0,j}(ik_1\beta_0)=0$ and ${\det}_{\mathbb C}\triangle_{\alpha_0,j}(ik_2\beta_0)=0$ with $j_1\neq j_2$ or $k_1 \neq k_2$. As shown above, the $G$-isotypic crossing number $\mathfrak t_{k,j}(\alpha_0,\beta_0)$ can be computed even in this case, and therefore so can the local bifurcation invariant, and the equivariant Krasnosel'skii theorem gives the local bifurcation of solutions with symmetries corresponding to each twisted orbit type which appears in $\omega_G(\alpha_0,\beta_0)$. Rather, the issues and potential degeneracies pertain to the global picture.
\vs
In particular, it is possible that a periodic solution which bifurcates locally collapses to a constant solution due to a secondary bifurcation, or intersects with a branch of steady-state solutions nonlocally. In general, even if one knows from the Krasnosel'skii theorem that certain symmetric solutions bifurcate at a given critical point, and knows from the Rabinowitz alternative that a branch continues globally, it is not possible to exclude these situations without further assumptions. In particular, if we choose a subgroup and an associated fixed point space where constant solutions cannot occur, and restrict the problem to this space, then clearly solutions to this restricted problem are also solutions to \eqref{eq:neutral:op-sys}.
\vs
This section proceeds in parallel to the development of the fixed point reduction in Chapter \ref{chapter:pseudo}. Fix $\widehat kappa \in {\mathbb N}$\footnote{For convenience, in this section we are freeing the notation $\kappa$ from its earlier use as the Lipschitz constant of the map $\bm g$} and define $\bm K\leq G$ as
\[
\bm K := \left<\left(e^{\tfrac{i\pi}{\kappa}},-1,e_{\Gamma_0}\right)\right> \leq S^1 \times \Gamma_0 \times \bbZ_2,
\]
where $e_{\Gamma_0}$ is the neutral element in $\Gamma_0$. We then restrict \eqref{eq:neutral:op-sys} to the fixed point space $\mathscr E^{\bm K}$ to obtain
\begin{equation}\label{eq:neutral:func-sys-fixedpoint}
\mathscr F^{\bm K}(\alpha,\beta,u) := \mathscr F(\alpha,\beta,u)\vert_{\mathscr E^{\bm K}}=0,\quad (\alpha,\beta,u)\in \bbR^2_+ \times \mathscr E^{\bm K}    
\end{equation}
Clearly any $(\alpha,\beta,u)$ satisfying $\mathscr F^{\bm K}(\alpha,\beta,u)=0$ also satisfies $\mathscr F(\alpha,\beta,u)=0$, and must be a non-constant periodic function. Indeed, we can take the $G$-isotypic decomposition of $\mathscr E^{\bm K}$ in terms of the $G$-isotypic components of $\mathscr E$, and see that it consists of only odd multiples of $\mathscr E_{\kappa, j}$, i.e. we have
\[
\mathscr E^{\bm K}_{k,j} := \begin{cases}
    \mathscr E_{k,j} \quad &\text{ if $k=(2l-1)\kappa$ for some $l\in{\mathbb N}$}\\
    0 \quad &\text{ otherwise}
\end{cases} 
\]
and
\[
\mathscr E^{\bm K} = \overline{\bigoplus_{l=1}^\infty \bigoplus_{j=0}^r \mathscr E_{(2l-1)\kappa,j}}
\]
Since $\bm K$ is a normal subgroup, and since $G_0 := G/\bm K \cong S^1 \times \Gamma_0$, $\mathscr F^{\bm K}(\alpha,\beta,u)=0$ is itself a two-parameter $S^1 \times \Gamma_0$-equivariant bifurcation problem, and the critical set for \eqref{eq:neutral:func-sys-fixedpoint} can be related to the critical set for \eqref{eq:neutral:op-sys} as:
\[
\tilde{\Lambda}^{\bm K}:=\{(\alpha_0,\beta_0,0) \in \bbR^2_+ \times \mathscr E: \exists l \in {\mathbb N}\;\;{\det}_{\mathbb C}(\triangle_{\alpha_0}(i(2l - 1)\kappa\beta_0))=0\}.
\]
Note that the condition which excluded steady-state bifurcation in the original definition of $\widetilde \Lambda$ is unnecessary here, as $\mathscr E^{\bm K}$ already excludes the zero mode. We also define
\[
\tilde{\Lambda}_l^{\bm K}:=\{(\alpha_0,\beta_0,0) \in \bbR^2_+ \times \mathscr E: \exists j\in {\mathbb N} \;\;{\det}_{\mathbb C}(\triangle_{\alpha_0,j}(i(2l - 1)\kappa\beta_0))=0\}
\]
and so
\[
\tilde\Lambda^{\bm K} = \bigcup_{l=1}^\infty \tilde\Lambda_l^{\bm K}
\]
\vs
That this set is infinite and isolated follows from the same properties of $\widetilde \Lambda$. The computation of the local bifurcation invariant for the fixed point reduced problem also takes the simpler form (free of any products) given by
\[
\sum_{i=1}^N \omega_G(\alpha_i,\beta_i) = \sum_{i=1}^N\sum_{j=0}^r \sum_{l=1}^\infty \mathfrak t_{(2l-1)\kappa,j}(\alpha_i,\beta_i) \text{deg}_{\mathcal W_{(2l-1)\kappa,j}}.
\]
Now we can provide a more comprehensive global bifurcation theorem which guarantees that solutions are non-periodic, which concludes the main abstract results.
\begin{theorem}\label{thm:neutral:nfde-global1}
Let $G=\Gamma_0\times\bbZ_2\times S^1$ and let $\mathscr F:\bbR_+^2\times \mathscr E \to \mathscr E$ be a $G$-equivariant condensing map satisfying assumptions \ref{neutral:c1} -- \ref{neutral:c3}. If $\tilde\Lambda_k^{\bm K}$ is discrete and finite for some $k=1,2,\dots$, and for some orbit type $(H)\in\Phi^t_1(G;\cV^-_{(2k-1)\kappa,j})$ we have
\[
\sum_{(\alpha_0,\beta_0,0)\in\tilde\Lambda_k^{\bm K}} \text{coeff}^H(\omega_G(\alpha_0,\beta_0)) \not=0
\]
where $\text{coeff}^H$ stands for the coefficient of $(H)$, then there exists an unbounded global branch $\mathscr C'\subset \mathscr S^H$ of non-constant periodic solutions with symmetries at least $(H)$, and with $\mathscr C'\cap\tilde\Lambda_k^{\bm K} \not=\emptyset$.
\end{theorem}
\begin{proof}
As stated above, since solutions to \eqref{eq:neutral:func-sys-fixedpoint} are solutions to \eqref{eq:neutral:op-sys}, and are non-constant and periodic, the equivariant Krasnosel'skii theorem and equivariant Rabinowitz alternative suffice to prove this result. Here, the equivariant twisted degree used to resolves to the Nussbaum-Sadovskii degree, but this makes absolutely no difference in the application of either of these theorems, as the Nussbaum-Sadovskii degree of a condensing perturbation of identity ultimately reduces to the Leray-Schauder degree, which applies equally to the equivariant version of both. The main result presented in the introduction therefore simply combines this theorem with the already obtained results on the signs of the crossing numbers in terms of \eqref{eq:neutral:rho} or \eqref{eq:neutral:upsilon} in the large $\beta$ limit, and so this suffices to prove Theorem \ref{thm:neutral:mas-global} and Theorem \ref{thm:neutral:mas-sym}.
\end{proof}

\section{Application: coupled asset markets with fundamentalists and momentum traders}\label{sec:neutral:application}
\begin{figure}[tbp]
    \centering
        \includegraphics[width=\linewidth]{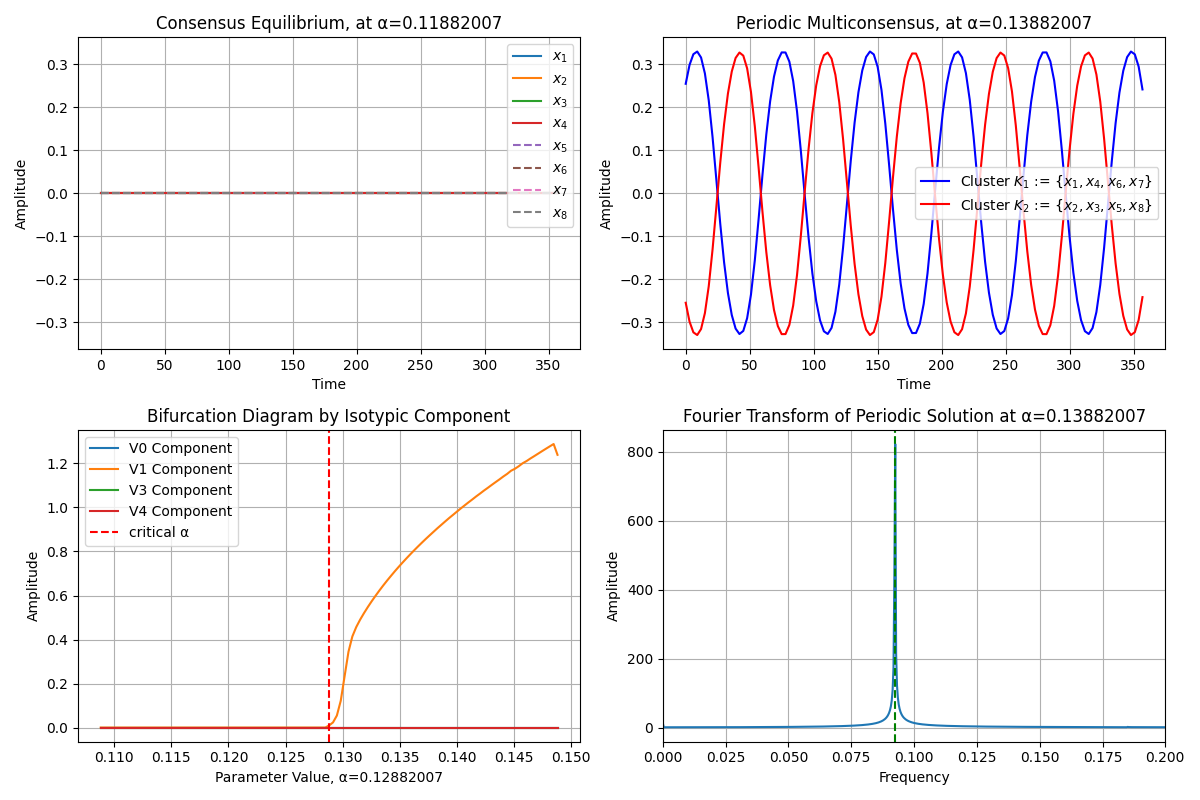}

    \caption{The top two figures show the evolution of the solutions (discarding initial transients) after a long simulation time ($t>50000$), for values of $\alpha$ slightly less than and slightly greater than $\alpha_{1,1}$. The bottom right plot shows the frequency spectrum of the steady-state solution, computed via fast Fourier transform, and normalized to $\beta$ by multiplying frequencies by $2\pi$. The bottom left plot shows the maximum amplitude of the solution projected onto all 4 $\Gamma_0$-isotypic components, computed by sweeping $\alpha$ along a mesh of 60 points near $\alpha_{1,1}$.}
    \label{fig:neutral:iso-comp-1-bif}
\end{figure}

Here we will follow the same type of application as formulated in Chapter \ref{chapter:pseudo}. However, the neutral delay has a different interpretation in this context. While the trend memory of the pseudoneutral delay offered an obvious interpretation as chartist traders who trade based on analysis of chart trends, the momentum memory of the neutral delay corresponds to momentum traders, who trade based on instantaneous derivatives at some point in the past. Clearly, the latter can be viewed as a type of limiting case of the former, as the averaging period over which trends are considered approaches 0 and the pseudoneutral term approaches a neutral term in its corresponding limit. We present this in order to allow a direct comparison of these two types of models taken with the same parameters, which illustrates many of the challenges in both the analysis and numerical simulation of neutral equations compared to their pseudoneutral counterparts of purely retarded type. 
\vs
% In this case, it also makes sense to choose $\tau_1$ to be substantially smaller than the choice in Chapter \ref{chapter:pseudo}, as momentum traders would consider much more recent data. We consider a model very similar to the one formulated in Chapter \ref{chapter:pseudo}, to facilitate a more direct comparison. 
% We again define the saturable piecewise-linear function
% \[
% L:= \begin{cases}
%     L_{\text{sat}}&\quad x\geq L_{\text{sat}}\\
%     x&\quad -1<x<1\\
%     -L_{\text{sat}}&\quad x \leq -L_{\text{sat}},
% \end{cases}
% \]
% used as a computationally convenient saturation function, representing finite liquidity in the market. 
Let $x_i(t)$ represent the value of asset $i$ relative to its fundamental, where $x_i>0$ represents overvaluation and $x_i<0$ represents undervaluation. For $\bm x = (x_1,x_2,\dots,x_8)^T\in \mathbb R^8$, we consider the following closed-form dynamics for 8 coupled homogeneous markets (each of which may be viewed as a HAM):
\[
\frac{d}{dt}\left[\bm x(t) - \tanh(\gamma \bm x(t-\tau_1))\right]=-a \bm x(t) - \alpha \tanh\left(b\int_0^{\tau_2}\bm x(t-s)ds\right) - h(\bm x(t)) \]
where $-a \bm x$ represents a market maker which buys or sells to push the asset's price towards its fundamental value, $\alpha$ represents the aggressiveness of fundamentalist traders given by $\tanh\left(b\int_0^{\tau_2}\bm x(t-s)ds\right)$, and $h(\bm x) := \tanh\left(\sum_{j=1}^8c_{ij}x_j\right)$ (where $c_{ij}$ is an entry of $C$) represents influence between markets which share certain properties. An edge connection on the cubic graph represents assets which share two out of three properties, a face diagonal represents assets which share one out of three properties, and a space diagonal represents assets which differ in every property. We use $\tanh$ to model saturative effects of, for example, capital or liquidity limitations. The momentum memory term $\frac{d}{dt}[\tanh(\gamma\bm x(t-\tau_1))]$ represents algorithmic traders, or high-frequency traders, who use price information from the market to make trades as quickly as possible based on computed or approximated price derivatives.
\vs
We will take mostly similar parameters as in Chapter \ref{chapter:pseudo}, but since the momentum memory term now represents algorithmic or high-frequency traders who are reacting to market information as soon as it is available, after a short delay in acquiring that information, processing it, and placing their trades, we set $\tau_1$ to a much lower value of half a day. Therefore, we set $a = 0.5, b=0.2, \gamma = 0.2, \tau_1 = 0.5, \tau_2=60, c_1 =0.15, c_2=0.05,c_3=0.01$, where the coupling matrix $C:=Dh(0)$ is given by
\[
C = \begin{bmatrix}
0 & 0.15 & 0.15 & 0.05 & 0.15 & 0.05 & 0.05 & 0.01 \\
0.15 & 0 & 0.05 & 0.15 & 0.05 & 0.15 & 0.01 & 0.05 \\
0.15 & 0.05 & 0 & 0.15 & 0.05 & 0.01 & 0.15 & 0.05 \\
0.05 & 0.15 & 0.15 & 0 & 0.01 & 0.05 & 0.05 & 0.15 \\
0.15 & 0.05 & 0.05 & 0.01 & 0 & 0.15 & 0.15 & 0.05 \\
0.05 & 0.15 & 0.01 & 0.05 & 0.15 & 0 & 0.05 & 0.15 \\
0.05 & 0.01 & 0.15 & 0.05 & 0.15 & 0.05 & 0 & 0.15 \\
0.01 & 0.05 & 0.05 & 0.15 & 0.05 & 0.15 & 0.15 & 0
\end{bmatrix}.
\]
Linearizing the system around the $\bm x\equiv 0$ consensus with the above parameters, we obtain the linear system
\begin{equation}\label{eq:neutral:example-sys-lin}
\frac{d}{dt} \left[x - 0.2 x(t-\tfrac{1}{2})\right] + 0.5x + 0.2\alpha \int_0^{60} x(t-s)ds + C = 0.
\end{equation}
Since we have the same interaction topology as in Chapter \ref{chapter:pseudo}, we set $\Gamma_0 := \mathbb O \cong S_4 <S_8$, where $\mathbb O$ is the symemtry group of rigid motions of the cube, which is isomorphic to the symmetric group $S_4$, viewed as a subgroup of $S_8$ which acts by permuting the $8$ vertices of the cube. Then $V:=\mathbb R^8$ is a natural isometric orthogonal $\Gamma_0$-representation, and we have the same $\Gamma_0$-isotypic decomposition of $V$, given by
\[
V = V_0 \oplus V_1 \oplus V_3 \oplus V_4,
\]
where each isotypic component $V_j$ has isotypic multiplicity 1 and is modeled on the real absolutely irreducible $\mathbb O$-representation $\mathcal V_j$ as given by Table \ref{table:pseudo:char}. 
Numerically solving \eqref{eq:neutral:beta-relation} and \eqref{eq:neutral:alpha-relation} and searching for the smallest value of $\alpha$ gives the pair $\alpha_{1,1} = 0.1288200, \beta_{1,1}=0.0923093$ on the $V_1$ component. Substitution into $\rho_1$ implies that $t_{1,1}(\alpha_{1,1},\beta_{1,1})=-1$. A negative crossing number, indicating an eigenvalue whose real part crosses from negative to positive at $\alpha = \alpha_{1,1}$, is exactly what is expected from a bifurcation from a stable consensus solution. If we perturb the zero equilibrium with a small value for $\alpha>\alpha_{1,1}$ and let the system evolve, it should attract towards a non-constant periodic solution on the $V_1$ subspace. Numerical simulations, shown in Figure \ref{fig:neutral:iso-comp-1-bif}, verify that a periodic non-constant solution does emerge on the predicted isotypic component, at the predicted critical value of $\alpha$, and with the predicted limit frequency\footnote{Since $\beta=\tfrac{2\pi}{p}$ and the frequency of a $p$-periodic solution is $\omega = \tfrac{1}{p}$, the true frequency predicted by $\beta$ is actually $\tfrac{\beta}{2\pi}$. To make the results easier to see graphically, we instead multiply the computed frequencies from the fast Fourier transform by $2\pi$, which is equivalent to dividing $\beta$ by $2\pi$.}.

\section{Hopf-Hopf bifurcation and consensus breakdown to chaotic multiconsensus}
\begin{figure}[tbp]
    \centering
    \begin{minipage}{0.45\textwidth}
        \centering
        \includegraphics[width=\linewidth]{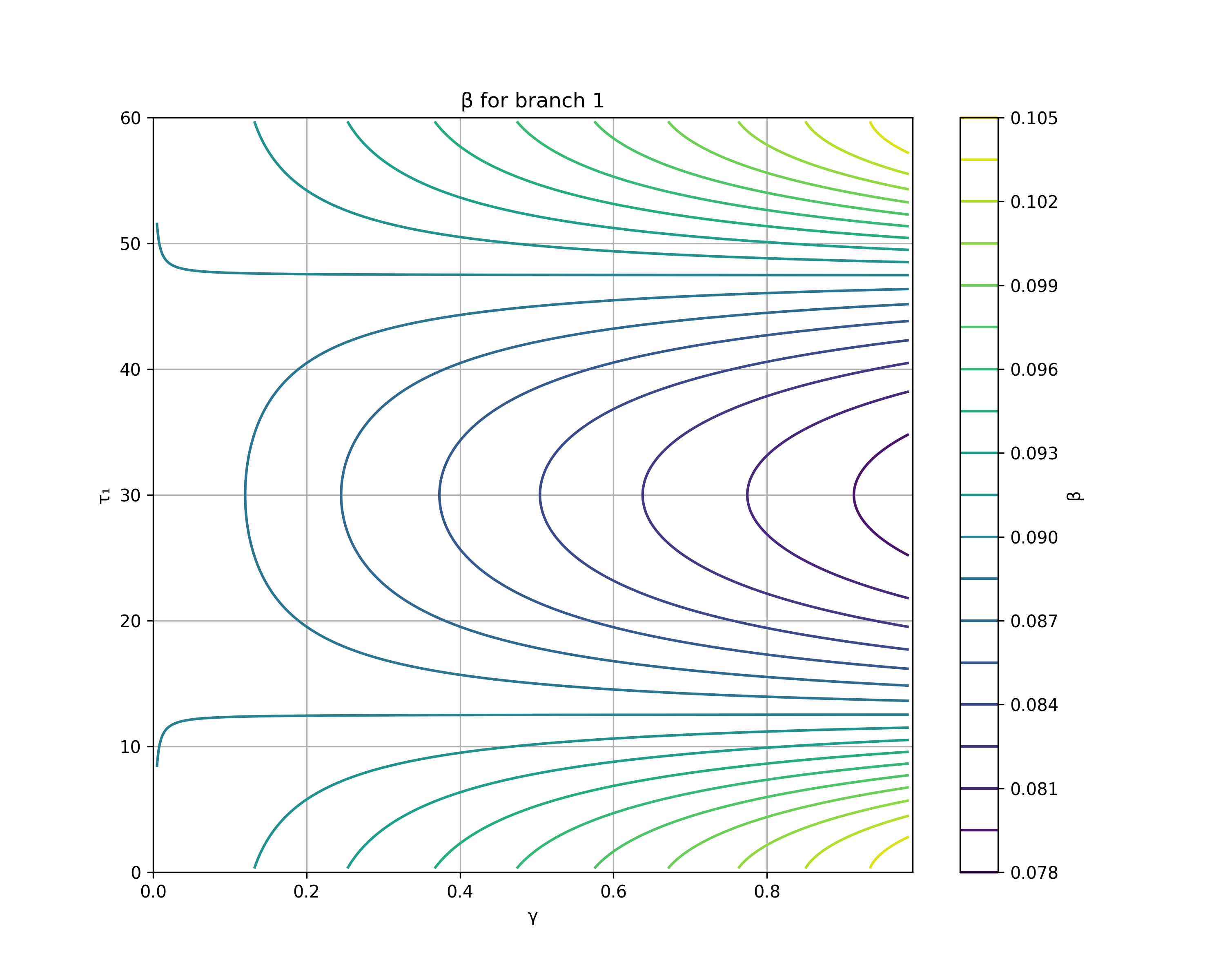}
        \subcaptiontext{a}{\footnotesize Contour plot showing the first critical value of $\beta$ on the $V_1$ isotypic component, $\beta_{1,1}(\gamma,\tau_1)$, plotted as a function of $\gamma$ and $\tau_1$.}
    \end{minipage}\hfill
    \begin{minipage}{0.45\textwidth}
        \centering
        \includegraphics[width=\linewidth]{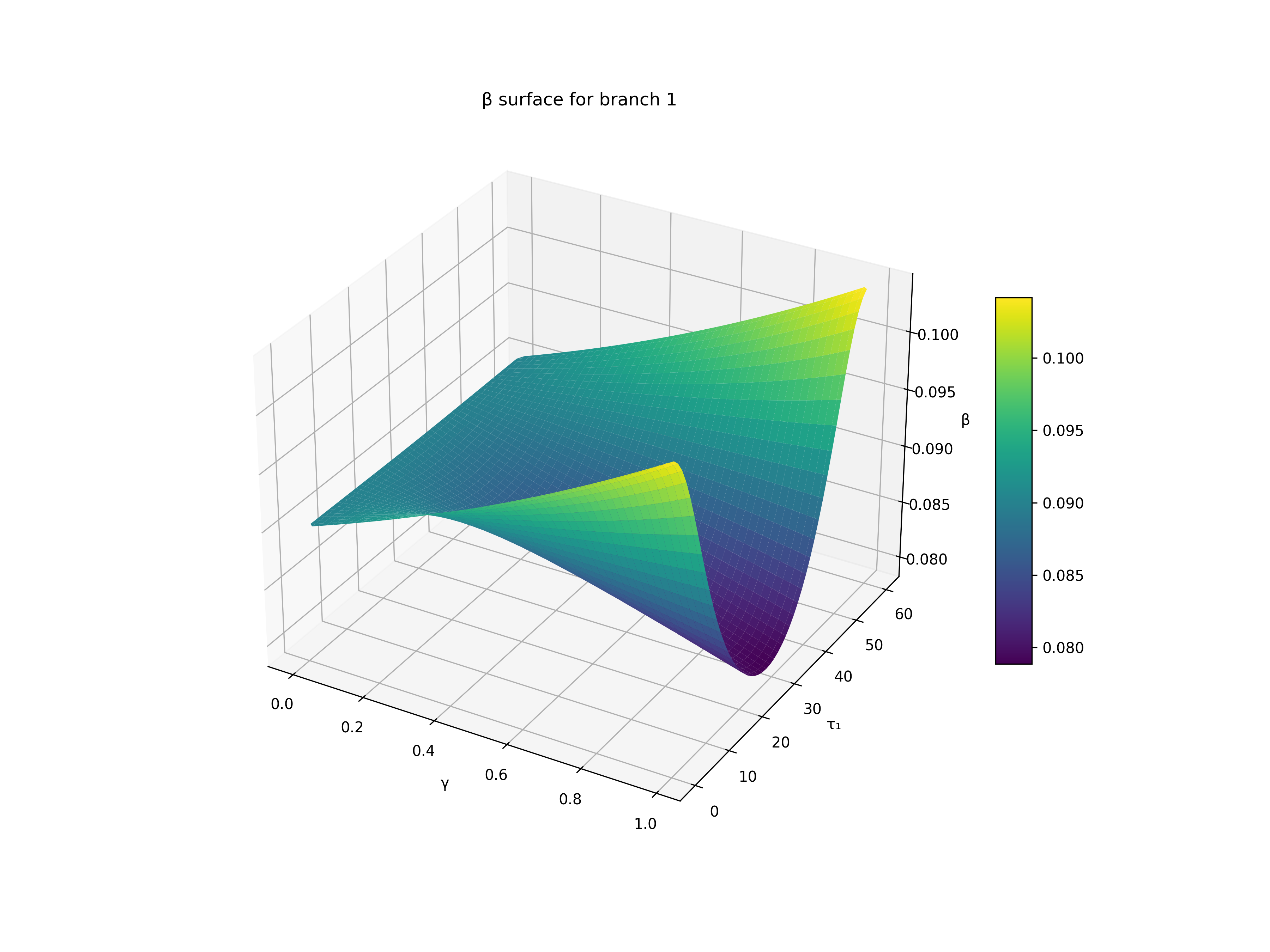}
        \subcaptiontext{b}{\footnotesize 3D surface plot showing the $\beta_{1,1}$ limit frequency surface as a function of $\gamma$ and $\tau_1$}
    \end{minipage}
    \caption{Parameter study of the dependence of the first $V_1$ limit frequency $\beta_{1,1}$ on the parameters $\gamma$ and $\tau_1$, for $0\leq \gamma <1$, and $0\leq \tau_1 \leq 60$, with $a=0.5, b=0.2,c_1=0.15,c_2=0.10,c_3=0.5$.}
    \label{fig:neutral:beta-gamma-tau1}
\end{figure}

\begin{figure}[tbp]
    \centering
    \begin{minipage}{0.45\textwidth}
        \centering
        \includegraphics[width=\linewidth]{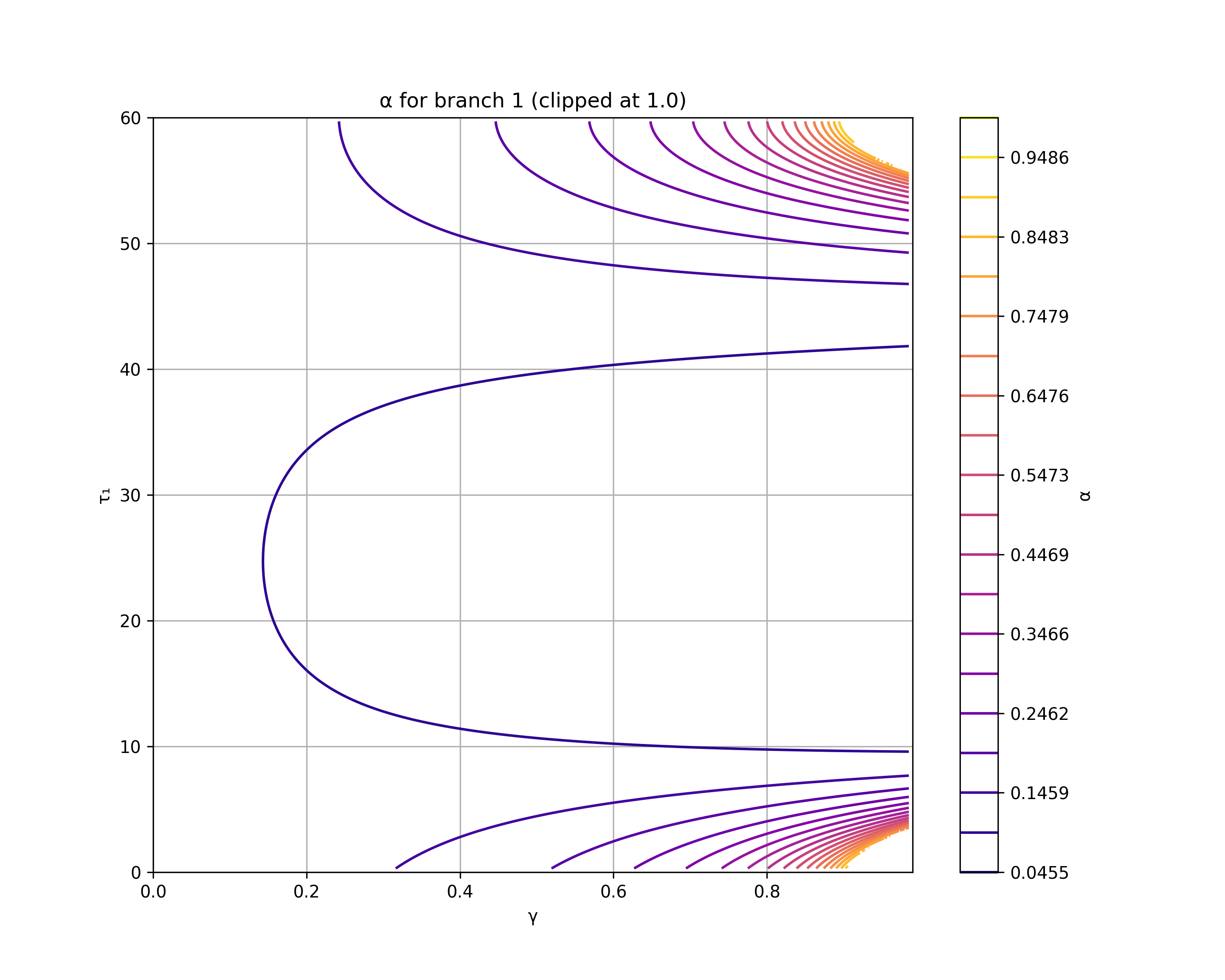}
        \subcaptiontext{a}{\footnotesize Contour plot showing the first critical value of $\alpha$ on the $V_1$ isotypic component, $\alpha_{1,1}(\gamma,\tau_1)$, plotted as a function of $\gamma$ and $\tau_1$. We consider only $\alpha<1$ to give more regularly spaced contours, but $\alpha_{1,1}$ can take larger values near the corners.}
    \end{minipage}\hfill
    \begin{minipage}{0.45\textwidth}
        \centering
        \includegraphics[width=\linewidth]{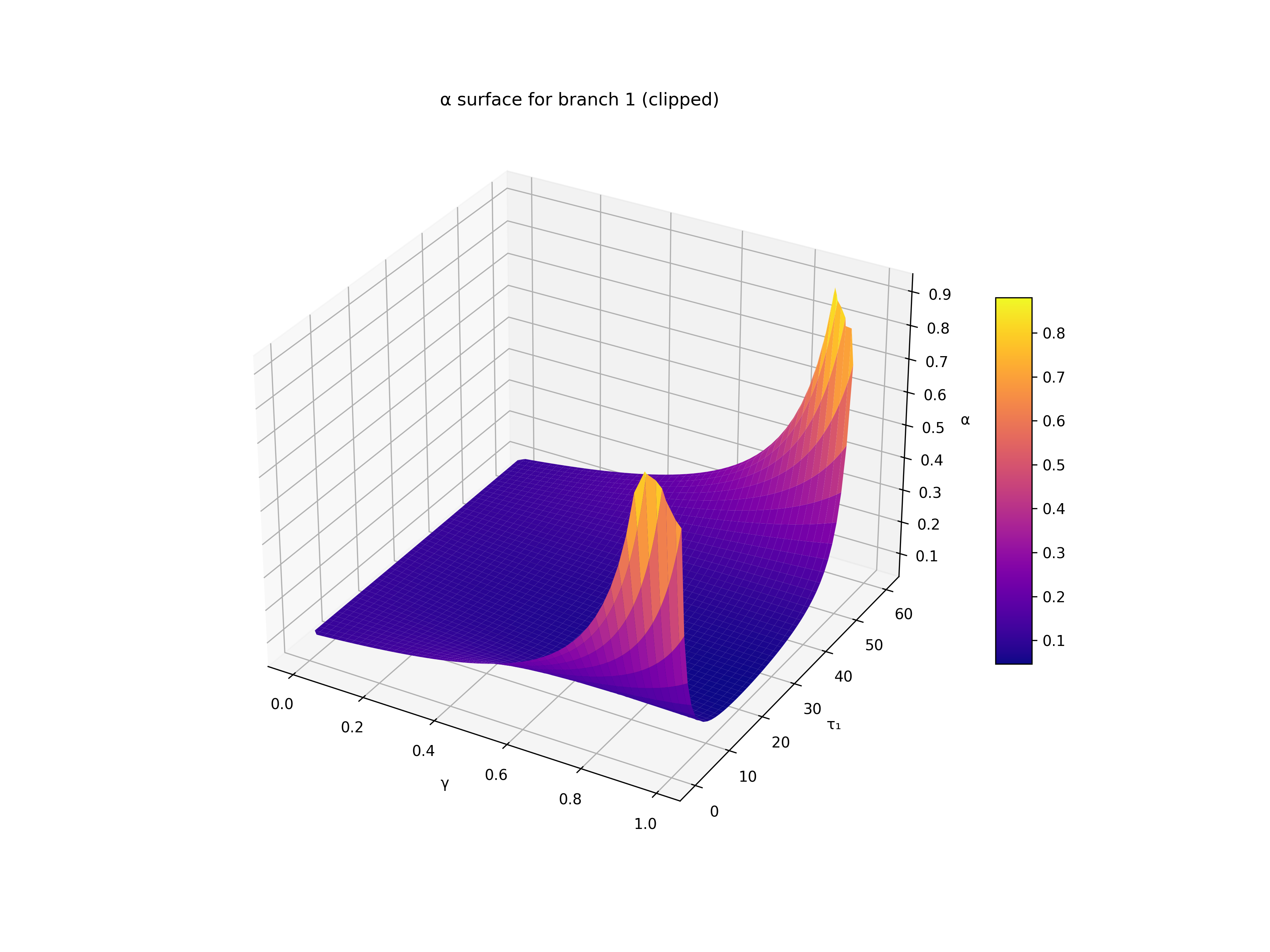}
        \subcaptiontext{b}{\footnotesize 3D surface plot showing the $\alpha_{1,1}$ Hopf bifurcation surface as a function of $\gamma$ and $\tau_1$}
    \end{minipage}
    \caption{Parameter study of the dependence of the first $V_1$ limit frequency $\alpha_{1,1}$ on the parameters $\gamma$ and $\tau_1$, for $0\leq \gamma <1$, and $0\leq \tau_1 \leq 60$, with $a=0.5, b=0.2,c_1=0.15,c_2=0.10,c_3=0.5$. Critical points correspond to taking $\beta$ from Figure \ref{fig:neutral:beta-gamma-tau1} and alpha from subplot \textbf{(a)} for the same pair of $(\gamma,\tau_1)$ coordinates.}
    \label{fig:neutral:alpha-gamma-tau1}
\end{figure}

The first bifurcation point is typically of greatest interest, since it is at this point that the trivial consensus first loses stability, and this periodic orbit is often more attractive to trajectories even after subsequent Hopf bifurcations from the trivial consensus equilibrium have occurred. As we have seen, the values $a_j$ taken from the parameter $a$ and the eigenvalues $\mu_j$ of the coupling matrix $C$ determine on which isotypic component Hopf bifurcation will first occur. For $c_3 < c_2 <c_1$, this will occur on the $V_1$ isotypic component. Even without performing a full separate bifurcation analysis using other parameters, we can explore the stability region of the trivial consensus in the wider parameter space and the effects of other parameters on the frequency of the initial multiconsensus Hopf bifurcation numerically. 
\vs
By numerically solving \eqref{eq:neutral:beta-relation} and \eqref{eq:neutral:alpha-relation} for the smallest limit frequency $\beta_{1,1}$ on the $V_1$ isotypic component, and seeing how the location of the associated critical point $(\alpha_{1,1},\beta_{1,1},0)\in \Lambda $ changes as a function of $\gamma$ and $\tau_1$ within the ranges $0\leq \gamma <1$ and $0 \leq \tau_1 \leq 60$, we can visualize the multiconsensus bifurcation surface, as seen in subplot \textbf{(b)} of Figure \ref{fig:neutral:alpha-gamma-tau1}. The corresponding surface of limit frequencies $\beta_{1,1}$ can be seen in Figure \ref{fig:neutral:beta-gamma-tau1}. This shows that for the momentum memory delay $\tau_1$ has a relatively small effect on both the limit frequency and the critical value of $\alpha$ if $\gamma<0.1$. When $\gamma > 0.5$, we begin to see a pronounced effect where the value of $\alpha_{1,1}$ corresponding to $\beta_{1,1}$ increases dramatically. 

However, note that these figures track the value of $\alpha_{1,1}$ associated with $\beta_{1,1}$ over the range $0\leq\gamma<1, 0\leq\tau_1\leq\tau_2$. As mentioned before, the dependence of critical $\alpha$ on $\beta$ is not necessarily monotonic. This means that for certain values of $\gamma$ and $\tau_1$, higher limit frequencies $\beta_{n,1}$ might correspond to lower critical parameter values, i.e. such that $\alpha_{n,1}<\alpha_{1,1}$. This fact also implies the existence of certain critical curves on the $(\gamma,\tau_1)$ plane where, for example, $\alpha_{1,1} = \alpha_{2,1}$. For these parameter values, a single critical value of $\alpha$ can correspond to multiple critical frequencies. Such values of $\alpha$ correspond to \emph{double Hopf points} where codimension-2 bifurcations can occur. The dynamics of the bifurcation at these double Hopf points also depends on the relationship between the two frequencies $\beta_{1,1}(\gamma,\tau_1)$ and $\beta_{1,2}(\gamma,\tau_1)$, particularly on whether or not these frequencies are commensurate.
\vs
\begin{figure}[tbp]
    \centering
        \centering
        \includegraphics[width=\linewidth]{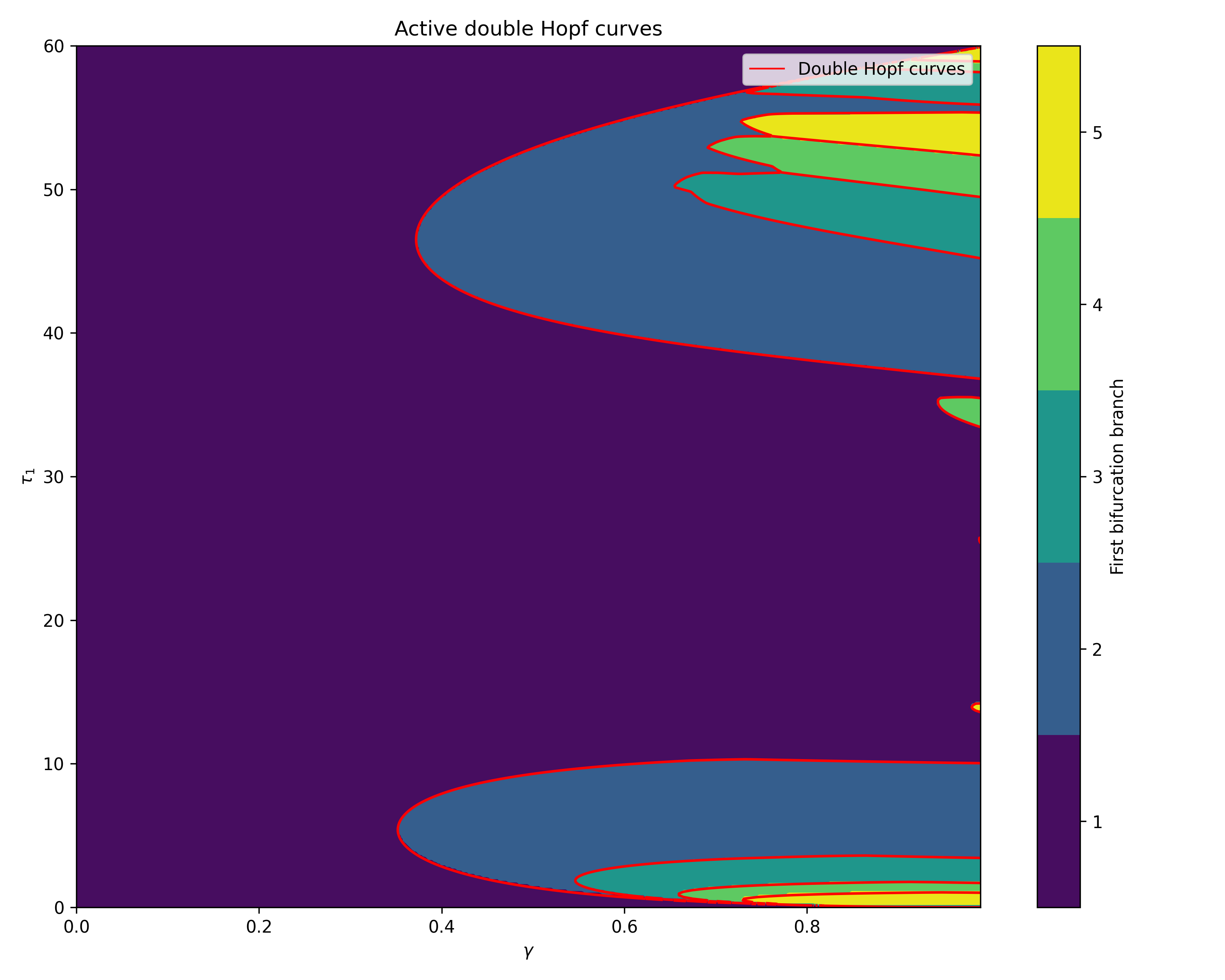}
    \caption{Minimal branch regions overlaid with red double Hopf curves at region boundaries. The red curve on the boundary of region $n$ and $m$ is the locus of points where $\alpha_{n,1}(\gamma,\tau_1) = \alpha_{m,1}(\gamma,\tau_1)$.}
    \label{fig:neutral:double-hopf1}
\end{figure}

We will refer to the branch with limit frequency $\beta_{n,1}$ as ``branch $n$,'' and we will refer to the branch which bifurcates first (in the sense that the associated $\alpha_{n,1}$ is minimal) as the ``first branch.'' This means that the double Hopf curves are also associated with ``branch switching,'' where the limit frequency of the first bifurcating branch jumps discontinuously. This lets us segment the $(\gamma,\tau_1)$ plane into regions corresponding to which branch index bifurcates first. Although one could search for points or curves where $\alpha_{n,1} = \alpha_{m,1}$ for any $m$ and $n$, clearly the ones of greatest interest are the curves on the boundaries of these regions, as these correspond to double Hopf points which bifurcate first. We can see these regions and the double Hopf curves along their boundaries in Figure \ref{fig:neutral:double-hopf1}. We will call these curves in the $(\gamma,\tau_1)$ plane of double Hopf points which bifurcate first the ``primary'' double Hopf points.
\vs
\subsection{\texorpdfstring{$\beta$}{β} resonances at primary double Hopf points}
\begin{figure}[tbp]
    \centering
        \centering
        \includegraphics[width=\linewidth]{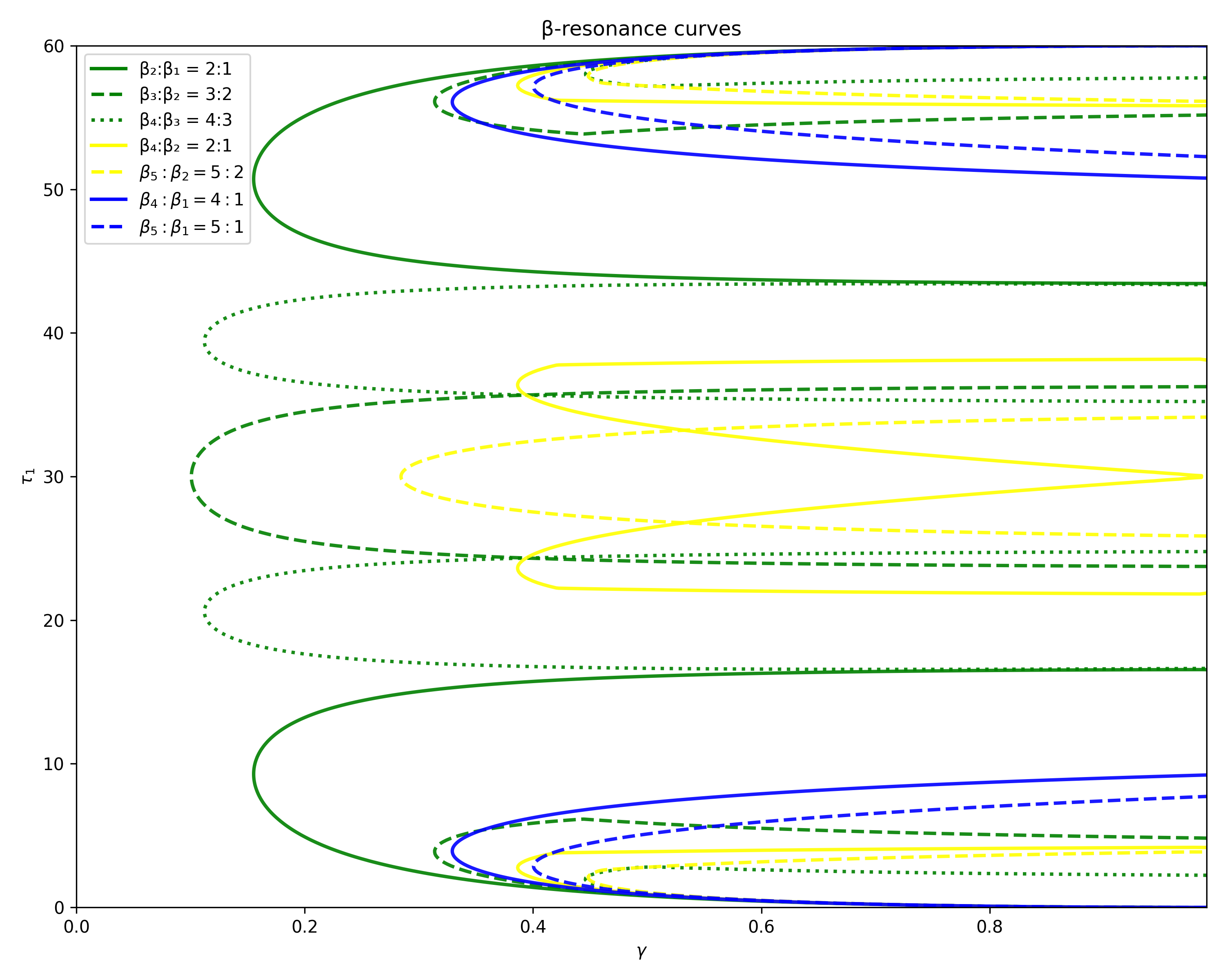}
    \caption{Curves of resonances between $\beta$ values which can coincide with branch switching regions. Most of these curves, however, do not intersect their corresponding double Hopf curve in \ref{fig:neutral:double-hopf1}.}
    \label{fig:neutral:beta-resonances}
\end{figure}

It can easily be seen numerically that $\beta_{n,1} \approx n\beta_{1,1}$, which also follows from the fact that solutions to \eqref{eq:neutral:beta-relation} occur near vertical asymptotes of the right hand side. Since each $\beta_{n,1}(\gamma,\tau_1)$ varies smoothly and independently with respect to $\gamma$ and $\tau_1$, this suggests that there could exist certain \emph{$\beta$-resonance curves} where, for example, $\beta_{2,1}:\beta_{1,1} = 2:1$. Since each $\beta_{n,1}(\gamma,\tau_1)$ can move only a limited amount for $(\gamma,\tau_1)\in [0,1)\times [0,\tau_2]$, and maintain their monotonic relationship throughout (see supplementary figures in Appendix \ref{appendix:supplementary-figs} showing contour plots and surfaces for $\beta_{n,1}(\gamma,\tau_1)$ for $n=1,2,3,4,5$), this implies that the interesting possible resonances at the primary double Hopf points are of the form $\beta_{n,1}:\beta_{m,1} = m:n$ at the boundary between the branch $n$ and branch $m$ regions. In other words, we are interested in points in the $(\gamma,\tau_1)$ plane where branch $m$ bifurcates first, where $\alpha_{m,1} = \alpha_{n,1}$, and where $\beta_{n,1}:\beta_{m,1} = m:n$.  (cf. Fig. \ref{fig:neutral:beta-resonances}). 
\vs

We can use these figures as a ``map'' for exploration of parameter regions with particularly interesting dynamics. Of course, the $(\gamma,\tau_1)$ values where these resonances or double Hopf points occur are very exact, so this map is best used to identify interesting approximate points which can then be refined to whatever level of precision is desired, by applying numerical root finding methods to find exact double Hopf points and/or $\beta$ resonances. For the remainder of this section, we will explore some of the particularly interesting dynamics discovered numerically along the primary double Hopf curves, particularly quasiperiodic multiconsensus, chaotic multiconsensus, and transient chaos leading to steady-state bifurcation. 

\begin{figure}[tbp]
    \centering
        \centering
        \includegraphics[width=\linewidth]{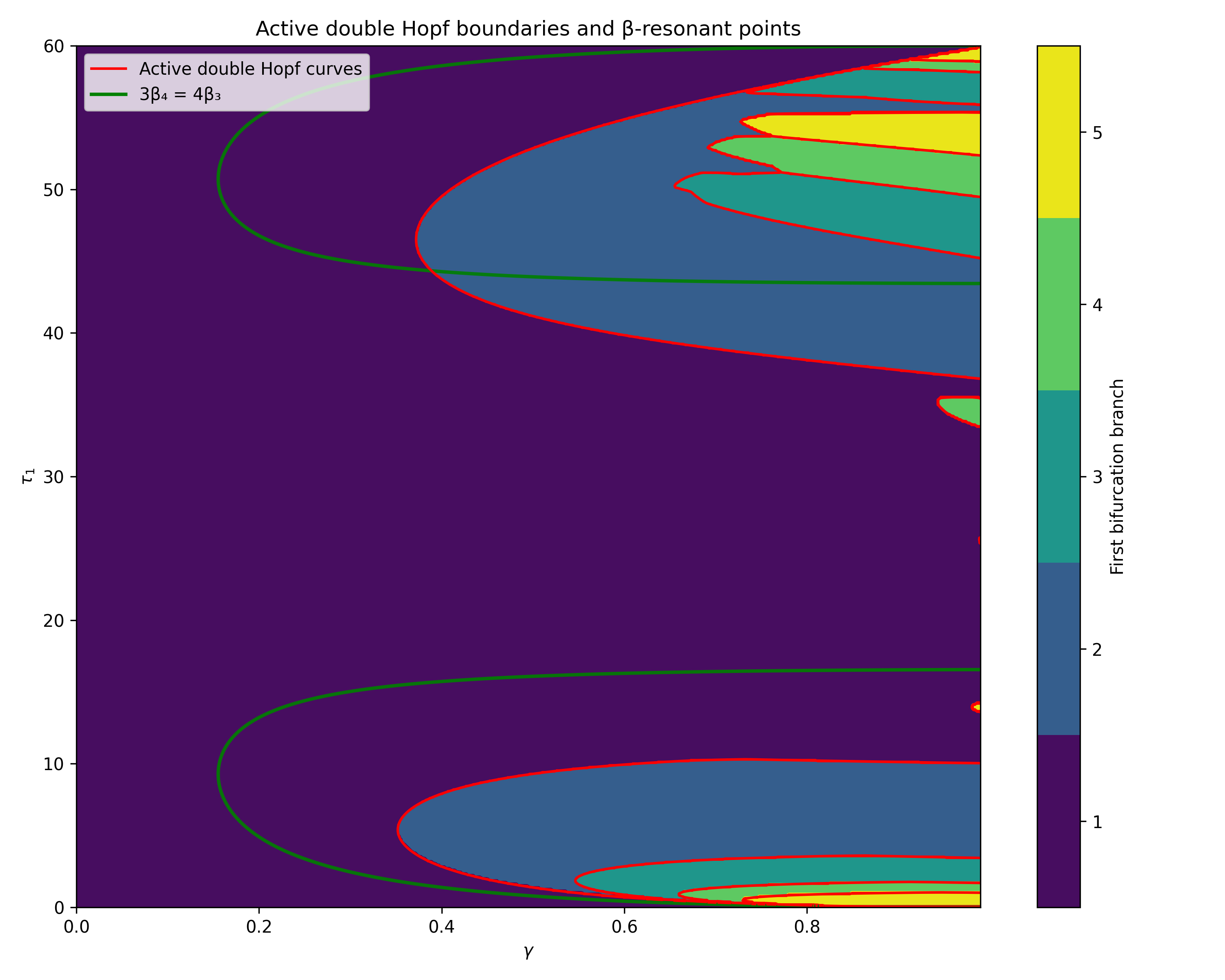}
    \caption{Branch regions with primary double Hopf curves and the $\beta_{2,1}:\beta_{1,1} = 2:1$ resonance curve. The point of intersection is a source of potentially chaotic dynamics. Due to nonlinear phase coupling effects, it is not necessary for values of $\beta$ to be exactly resonant on the double Hopf curve, and even nearby double Hopf points can exhibit these effects.}
    \label{fig:neutral:double-hopf-resonance}
\end{figure}

\subsection{Double Hopf bifurcation and the equivariant twisted degree}
First of all, we note that it is always possible for a limit cycle which emerges via Hopf bifurcation to undergo secondary bifurcation. These can include period doubling bifurcations as well as secondary Hopf bifurcations, also called torus bifurcations\footnote{This can also be viewed as a Neimark-Sacker bifurcation on the Poincar\`e map of a limit cycle.}. However, such situations cannot really be studied using the equivariant degree without linearizing along the bifurcating branch. A situation which can be studied using the degree, however, is \emph{double Hopf bifurcation}, where two distinct conjugate pairs of complex eigenvalues cross the imaginary axis simultaneously. This is a similar situation to secondary Hopf bifurcation of a limit cycle, which results in a 2-torus, but is also capable of producing many other rich types of dynamics.

\begin{remark}\rm
The double Hopf situation illustrates one of the great advantages of equivariant degree methods. In symmetric systems, it is quite ordinary and expected for multiple pairs of eigenvalues to cross the imaginary axis simultaneously. Indeed, coupling symmetries often force this situation. This results in symmetric limit cycles, whose symmetries are related to the maximal orbit types on the $\Gamma$-isotypic component where these eigenvalues crossed simultaneously. On the other hand, one could have pairs of eigenvalues crossing simultaneously whose frequencies are distinct and are either commensurate or incommensurate. If they are commensurate, they can be viewed as crossings on different $S^1$ modes for a single critical point $(\alpha_0,\beta_0,0)\in\Lambda$. If they are incommensurate, they can be viewed as crossings on the first mode on two critical points $(\alpha_0,\beta_0,0),(\alpha_0,\beta_1,0)\in\Lambda$. However, all of these situations are considered degenerate from the point of view of classical Hopf bifurcation theory, which does not distinguish between a high-dimensional center manifold caused by symmetrically linked eigenvalues versus crossings at distinct frequencies which can elicit more complicated multi-mode interactions. On the other hand, none of these situations are a problem from the point of view of the equivariant twisted degree, which gives us bifurcation invariants and therefore the existence of symmetric branches of solutions in any of these cases.   
\end{remark}

\subsection{Quasiperiodic multiconsensus and the Ruelle-Takens-Newhouse route to chaos}
\begin{figure}[tbp]
    \centering

    \includegraphics[width=\textwidth]{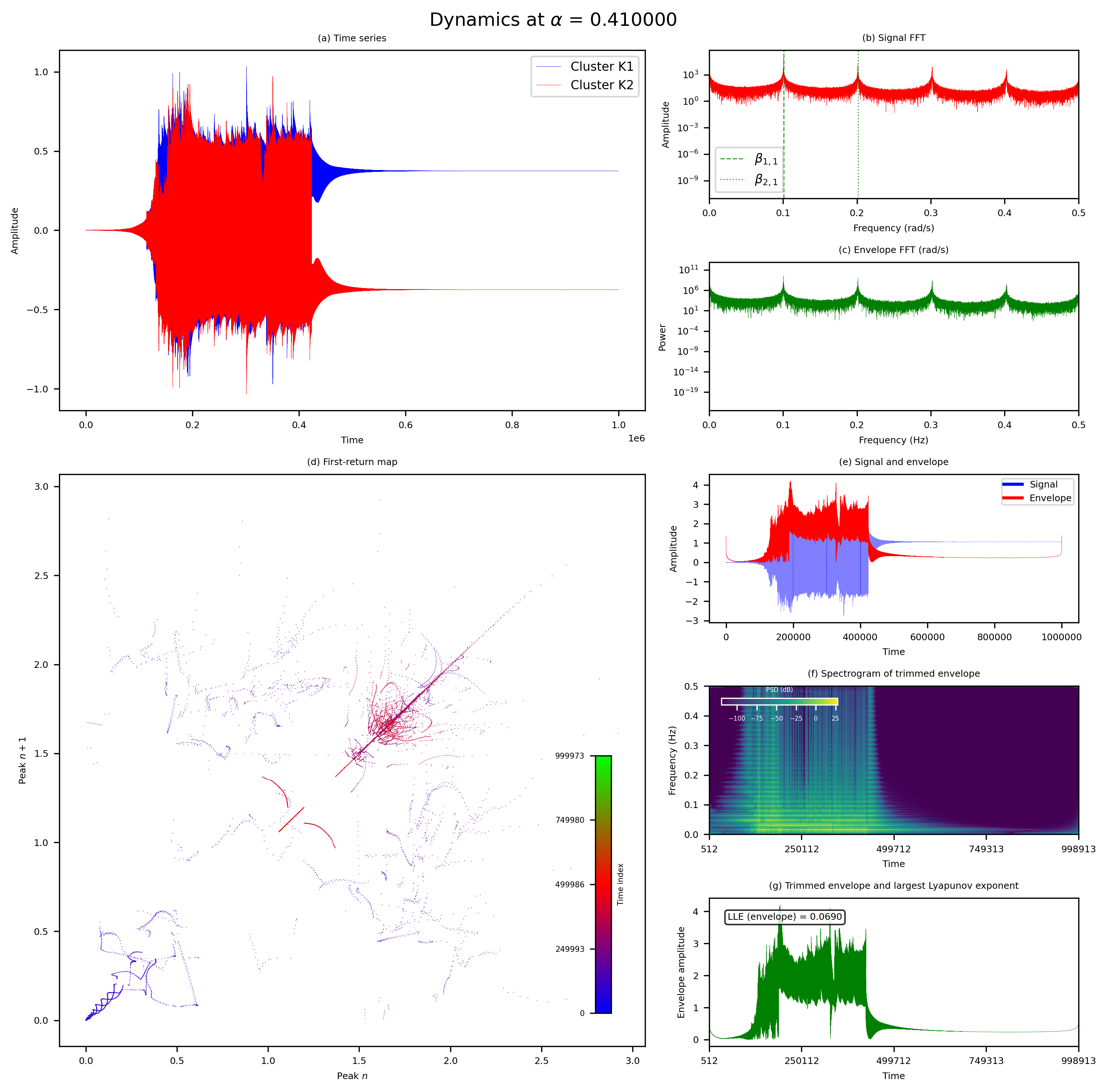}

    \caption{Chaotic saddle dynamics at $\alpha = 0.41$.}
    \label{fig:neutral:trajectory-0.41}
\end{figure}

When double Hopf bifurcation occurs for an incommensurate pair of eigenvalues, the two frequency modes tend to be decoupled, and can result in a quasiperiodic solution whose trajectory densely fills the 2-torus. When the frequency pair is commensurate with a resonant ratio (e.g. $\beta_{2,1}/\beta_{1,1} = 2:1$), this results in nonlinear phase-coupling terms which can yield extremely complex dynamics, including further bifurcation to chaos along several routes. One such route is the Ruelle-Takens-Newhouse (RTN) route \cite{ruelle1971nature}. In this scenario, a 2-torus bifurcates along a third frequency and then breaks down entirely. By checking where the $\beta$-resonant curves intersect the primary double Hopf curves, we obtain a list of candidate points for this scenario (cf. Fig. \ref{fig:neutral:double-hopf-resonance}). We will explore some of these points numerically and illustrate chaos in the form of chaotically amplitude-modulated periodic waves. One remarkable fact is that these chaotic multiconsensus solutions exhibit the same symmetries predicted by the equivariant degree theory for this isotypic component (i.e. they form two anti-phase clusters as before). We explore the intersection of the $\beta$-resonance curve and primary double Hopf curve at 
\begin{align*}
\tau_1 &= 0.1678667007,\\ 
\gamma &= 0.7760643397,\\
\alpha_1 = \alpha_2 &= 0.4029160833\\
\beta_1 &= 0.1007507146,\\
\beta_2 &= 0.2014401509,\\ 
\beta_2/\beta_1 &= 1.99939178,
\end{align*}
where all other parameters are the same as in the application in the previous section. Note that $\beta_1$ and $\beta_2$ are not in exact resonance. As we will see, they are close enough to the resonant regime to exhibit chaos. Since it is necessary to simulate trajectories out to very long time intervals in order to observe the low-frequency chaotic dynamics modulating the amplitude, we will give long time series plots showing the overall evolution of the signal (discarding initial transients), a Fourier transform showing the spectrum of the time series signal, which can be related to the frequencies $\beta_2$ and $\beta_1$. Since the chaos occurs in the amplitude modulation of this periodic oscillation, we will take a numerical Hilbert transform of the time series signal to extract its envelope, compute the Fourier transform of the envelope, and estimate the largest Lyapunov exponent of both the envelope and the underlying signal itself using the nolds library \cite{scholzel_nolds}. Finally, we will also show a first-return map taken from successive amplitude peaks of the periodic signal, projected onto the $V_1$ isotypic component. Since our system on the $V_1$ component basically exhibits one-dimensional dynamics (since all agents are either in phase or anti-phase), this is our equivalent of a Poincare section, and is sufficient to capture the amplitude-modulating dynamics. We believe that all these data taken together offer convincing numerical evidence of chaos, and strongly point to chaos through the Ruelle-Taken-Sackhouse route.
\vs
Sampling different values of $\alpha$ near the bifurcation point, we see different regions of quasiperiodicity, transient chaos (i.e. chaotic saddles), and sustained chaos. These different regimes alternate unpredictably at all scales of changes in $\alpha$. This, combined with the presence of chaotic saddles, where trajectories linger in a chaotic basin of attraction for a very long time before going to a different stable attractor (e.g. a stable quasiperiodic solution or stable steady-state solution), strongly suggest the possibility of riddled basins (in the sense of Alexander et al. \cite{alexander1992riddled}). This means that the basins of attraction of the chaotic attractor and of other attractors are intermingled, in the sense that for any ball around any point in the basin of an attractor, there is a set of positive Lebesgue measure belonging to the basin of a different attractor. This means that basins are intermingled in a measure-theoretic sense at every scale. The concept of riddled basins is closely connected to symmetries and the existence of invariant smooth manifolds. The existence of riddled basins could be strengthened numerically by computing the uncertainty exponent, which could be a promising topic for a future study. If riddled basins do exist, it would have profound implications for the robustness of multi-agent protocols in this type of system, as even infinitesimal errors in initialization could lead to qualitatively different multiconsensus patterns, a form of symmetry-induced decision fragility.
\vs
We will now show a few values of $\alpha$ near the bifurcation point and characterize the system's behavior. First, we choose $\alpha=0.405$. Here, we see an interesting transition between several regimes. As the solution evolves, it first forms a slowly growing quasiperiodic trajectory, then abruptly switches to a chaotic amplitude-modulated regime, then settles down to a periodic solution.  
\vs
The plots of chaotic dynamics in this section are made up of several subplots. In Figure \ref{fig:neutral:trajectory-0.41}, for example, subfigure \textbf{(a)} shows the time series solution over a long time scale. At this scale, we mainly see the envelope function, but this solution is actually a periodic function amplitude-modulated by a chaotic signal. Subfigure \textbf{(b)} shows the Fourier transform of the signal, which demonstrates clear peaks at the two Hopf frequencies and their higher harmonics, as well as broadband noise from the chaotic modulation. This can be seen more clearly in subfigure \textbf{(c)}, which shows the Fourier transform of the envelope function and shows broadband noise characteristic of chaos. This envelope function was extracted via Hilbert transform and is shown alongside the original signal in subfigure \textbf{(e)}. Subfigure \textbf{(f)} shows a spectrograph of the extracted and trimmed envelope, which is shown on its own in subfigure \textbf{(g)}. Note the transition to broadband noise which marks the clear transition both to and from chaotic dynamics. Finally, subfigure \textbf{(d)} shows a first-return map, obtained by calculating the peak amplitude of each peak and plotting successive peaks against each other, and filling the same role as a Poincare map. The first-return map is color-coded based on the position of the peak in the overall time series. The helical spirals in the bottom left of the first-return plot correspond to the quasiperiodic evolution of the solution before the abrupt transition to chaos, where these helical trajectories break apart into chaos. This observation is consistent with the RTN route to chaos. Note that these helical spirals would be closed curves, corresponding to quasiperiodic solutions, for a fully saturated quasiperiodic solution. They correspond to growing helices here because the solution transitions to chaos before becoming fully saturated.

% \begin{figure}[tbp]
%     \centering
%     % Top panel: Hilbert plot (wide)

%     % Bottom row: two panels
%     \begin{minipage}[t]{0.48\textwidth}
%         \centering
%         \includegraphics[width=\linewidth]{figures/chaos/alpha_sweep_plot_neutral_0.40400000.png}
%     \end{minipage}
%     \hfill
%     \begin{minipage}[t]{0.48\textwidth}
%         \centering
%         \includegraphics[width=\linewidth]{figures/chaos/frequency_diagram_neutral_0.40400000.png}
%     \end{minipage}
%     \\[1ex]   % a little vertical space
%     \includegraphics[width=\textwidth]{figures/chaos/envelope_alpha_0.40400000.png}

%     \caption{Chaotic trajectory diagnostics at $\alpha = 0.404$ . \textbf{Top:} Hilbert envelope analysis. \textbf{Bottom left:} time series and FFT. \textbf{Bottom right:} first‑return map.}
%     \label{fig:chaos-example}
% \end{figure}

It is interesting to note that the regions exhibiting chaotic dynamics still respect the symmetries of the bifurcation as predicted by the theoretical results. We observe a separation of dynamics in the form of chaotic amplitude-modulated oscillations, where the underlying fast oscillations are on the frequency scale of the critical frequencies $\beta_{1,1}$ and $\beta_{2,1}$, and some of their higher harmonics, while the slow dynamics modulating the amplitude of these oscillations are quasiperiodic or chaotic. We explore such solutions numerically, and infer chaos through numerical computation of Lyapunov coefficients of the envelope function, extracted from the time series data via numerical Hilbert transform. We strongly speculate that this chaos is through the RTN route, but cannot offer a conclusive proof. A detailed normal form exploration of double Hopf bifurcation in this system would be an interesting topic for future work.
\vs

\chapter{Conclusions and Future Work}\label{chapter:conclusion}
In this dissertation, we showed the local Hopf bifurcation, global unbounded continuation, and spatio-temporal classification of branches of non-constant periodic multiconsensus solutions for three multi-agent systems characterized by three different types of memory. In particular, we defined \emph{continuous memory} as a distributed delay of the form
\[
\int_0^\tau \phi(s) x(t-s)ds,
\]
where $\phi(s)$ is the delay/memory kernel which indicates how the agent continuously assigns weight to more recent versus more distant state histories, and $\tau$ is the time horizon (which can go to infinity if $\int_0^\infty\phi(s)ds <\infty$). In the three models we studied, we considered an equally weighted memory over a finite time horizon of the form 
\[
\int_0^\tau x(t-s)ds.
\]
We also considered a dependence on the instantaneous rate of change of a continuous memory term, 
\[
\frac{d}{dt}{\left[\int_0^\tau \phi(s)x(t-s)ds\right]},
\]
where we assume the memory kernel $\phi(s)\in L^1$, i.e. it is free of any Dirac atoms. This represents an agent's perception of recent trends, and we accordingly defined this as \emph{trend memory} or \emph{pseudoneutral memory}, as the closed-loop dynamics equations governing agents whose protocols include these types of terms become \emph{pseudoneutral equations}. We called these equations ``pseudoneutral'' because they are closely related to neutral equations, but are of purely retarded type. In particular, they can be viewed as an approximation of neutral delay terms. If one constructs a sequence of delay kernels which converges weakly to a Dirac delta function, the associated sequence of pseudoneutral equations will converge to a neutral equation. On the other hand, if we view distributed delays as giving the expected value of a variable delay which follows a probability distribution in a mean-field approximation, then pseudoneutral delays are the mean-field approximation of a distribution of neutral delays, which show how the pathological properties of neutral delays are ``smoothed out'' on a large scale if they are sufficiently distributed. 
\vs
For our systems, we considered a simplification of the above in the form of a trend memory of a sum of a nonlinearly transformed state history,
\[
\int_0^\tau g(x(t-s))ds.
\]
This is useful for implementing a saturating nonlinearity for $g$. such as $g(x) = \tanh(x)$, which allows us to give the trend memory a saturable maximum level of influence.
\vs
Finally, we considered a memory of past control actions, of the form
\[
g(x(t-\tau)).
\]
We defined this to be \emph{momentum memory} or \emph{neutral memory}. For agents with single-integrator dynamics, where they control their velocity $\dot x(t)$ via a protocol which depends on their state $x(t)$ (and possibly on memory), this term can represent the influence of remembered control actions or decisions.
\vs
For modeling networks of intelligent agents, both trend memory and momentum memory can be useful in modeling ``decision inertia,'' as both correspond to different ways that an agent's previous decisions can influence their current decision, and can model agents who may be reluctant to change course, who fall victim to sunken cost fallacies, or who possess certain kinds of physical or mental inertia. The key difference is that momentum memory involves a memory of every single control action taken, no matter how minute (or erroneous), whereas trend memory has natural smoothing effects and captures the generalities or ``directions'' of previous decisions. The other key difference is that momentum memory leads to closed-loop dynamics equations of neutral type, which have many difficult and sometimes pathological properties. 
\vs
We studied a system with just continuous memory, a system with both continuous memory and trend memory, and a system with continuous memory and momentum memory. For each system, we showed conditions for the local asymptotic stability of the trivial consensus equilibrium, derived the transcendental equations which describe the critical parameter values and limit frequencies, showed that the critical set is isolated, and obtained transversality conditions which determine the sign of the crossing number. Based on this, we then obtained conditions for the sign of the crossing numbers to be constant for all critical points, or for the sum of the signs of the crossing numbers taken over the critical set to be nonzero. 
\vs
Following this, we period normalized the system, extracting the frequency $\beta$ as a second parameter, reformulated the problem as a functional equation (as a compact perturbation of identity in Chapter \ref{chapter:distributed}, and as a condensing perturbation of identity in Chapters \ref{chapter:pseudo} and \ref{chapter:neutral}), and analyzed this as a two-parameter bifurcation problem. This allowed us to use the twisted equivariant degree to compute local bifurcation invariants around each critical point, which allowed us to prove local bifurcation by applying the equivariant Krasnosel'skii theorem. Global bifurcation was then achieved by applying the equivariant Rabinowitz alternative on a reduced problem in a fixed-point subspace which guaranteed that the global branches consisted entirely of non-constant periodic solutions throughout their full extent.
\vs
We then explored particular applications for each of the three systems. For the system with continuous memory in Chapter \ref{chapter:distributed}, we applied it as a formation control system for an octahedral formation of UAVs. In this application, periodic multiconsensus represents configurable patterns of periodic motions with zero mean displacement. We identified 27 distinct conjugacy classes of rhythmic motions, up to conjugacy, and identified how the desired pattern can be selected through a single transient forcing ``kick,'' after which further stable oscillations continue autonomously. In this context, multiconsensus can be viewed as a higher order of coordination and may be desirable, based on the application. This is comparable to the biological neural circuits known as Central Pattern Generators (CPGs). These circuits are closely involved in locomotion in animals ranging from invertebrates to mammals, and spontaneously generate periodic rhythmic excitations without any external forcing.
\vs
A particularly interesting property of this application is that these complicated synchronized choreographies were obtained using only autonomous controls within each UAV, without any external coordination or forcing (other than the small perturbation needed to select a particular pattern of interest, which could easily be programmed into the UAVs themselves). This means that as long as the UAVs have an onboard way of measuring each other's displacement from formation set points or relative distance from one another (which could be achieved in many ways, including via passive onboard sensing), these patterns are highly resilient to outside interference and cannot be ``jammed,'' because they are entirely self-guided and are not dependent on any outside instructions or control commands.
\vs
For the system with continuous memory and trend memory in Chapter \ref{chapter:pseudo}, we considered an economic application. In this application, inspired by other homogeneous agent models (HAM) which study markets with two populations of traders using different strategies, we consider each market itself as an agent whose goal is price discovery. We view $x_i(t)$ as representing the difference of the market price of asset $i$ from its fixed underlying value. An efficient market in equilibrium, therefore, should have $x_i(t)\equiv 0$, as the asset's price exactly matches its underlying value due to an equilibrium of supply and demand. We considered a protocol for each agent based on three actors trading within that market. There is the market maker or liquidity provider, which reacts to price changes instantly to close the bid-ask spread, represented by a constant term $-ax_i$. There are fundamentalist traders, who know the asset's underlying value, and track average price deviations averaged over a time horizon $\tau_2$, and represented by a continuous memory term. Lastly, there are chartists, who look at trends over a shorter time horizon $\tau_1$, and follow general market trends to buy when prices seem to be increasing (on average), and sell when they seem to be decreasing. Finally, we model the mutual influence of markets which share certain properties as part of a larger economy. 
\vs
In this context, the bifurcation parameter $\alpha$ represents the relative proportion (or aggressiveness) of fundamentalist traders. Once this exceeds a certain level, the market equilibrium becomes unstable, and we see emerging cycles of pricing bubbles and crashes, which grow in amplitude as the aggressiveness of traders grows. We also see arbitrage relationships emerging between different markets, where a bubble in one market corresponds to a crash in another. Whether these periodic multiconsensus cycles are advantageous in this application depends on one's perspective. They are probably worrying for the market maker and exciting for the arbitrageur. Interestingly, these inversions occurred in markets that had more in common. We found the period of the boom-bust cycles to be most strongly influenced by the length of memory of the fundamentalist traders. This ``forgetfulness'' could be interpreted as a certain threshold beyond which traders consider information to be of little value, or to not be representative of the current economy or current trends. The fact that relatively low-frequency boom-bust cycles correlate with the length of this memory matches informal observations that real life boom-bust cycles seem to occur on the time scale of institutional investors deciding that the lessons and mistakes of the past no longer apply to today's market conditions, but we will leave such conclusions for trained economists. 
\vs
In Chapter \ref{chapter:neutral} we considered a similar system, but with momentum memory in place of trend memory. In other words, where agents remember exact past decisions instead of a smoothed trend of past decisions. The closed-loop dynamics of this system form an actual neutral delay equation. We used the framework for such equations which was developed in Chapter \ref{chapter:pseudo} (though it was not yet strictly necessary for the equation in that chapter), which uses the Nussbaum-Sadovskii degree to analyze these equations, whose operator formulations fail to be compact. 
\vs
The application in this chapter was similar in spirit to that of Chapter \ref{chapter:neutral}, but with a different population of traders. In place of chartists, we consider momentum traders, who could also be viewed as algorithmic traders or high-frequency traders. These traders use market information as soon as they get it (and process it) to make very high frequency trades. The momentum memory in this model represents the inherent delay in getting market information from market makers, analyzing it to compute or estimate derivatives, and then choosing and executing their trades. In contrast to the chartists, who may examine trends over a period of days, the delay for these traders could be on the order of hours, minutes, or even seconds. We chose a delay of half a day as a reasonable value. We again found that the period of cycles on the component where bifurcation first occurred was largely driven by the fundamentalist traders, even when momentum traders in this application were given a substantially higher weight than the chartists in Chapter \ref{chapter:pseudo} ($\gamma = 0.2$ instead of $\gamma = 0.04$).
\vs
We numerically showed the stability of the periodic multiconsensus solutions emerging through Hopf bifurcation, corresponding to alternating pricing bubbles and crashes. We conducted our primary simulation with the same parameters as in Chapter \ref{chapter:pseudo}, showing that chartist traders and momentum traders such as algorithmic traders can produce similar dynamical effects in markets, at least under certain parameter ranges.
\vs
We also conducted a wider study of the $(\gamma,\tau_1)$ parameter space, focusing on the situation where $\tau_1\leq \tau_2$. We numerically computed the critical points where the first Hopf bifurcation occurs over a fine mesh of $\gamma,\tau_1$ values with $0\leq \gamma<1$ and $0\leq \tau_1 \leq \tau_2$. This allowed us to plot the branch switching regions, where the limit frequency of the first Hopf bifurcation switches from $\beta_{1,1}$ to $\beta_{n,1}$ (where $\beta_{n,1} \approx n\beta_{1,1}$, for $n=2,3,4,5$, and the associated smooth curves of double Hopf points, where $\alpha_{1,1}=\alpha_{n,1}$ and two complex eigenvalues $\beta_{1,1}$ and $\beta_{n,1}$ cross the imaginary axis simultaneously. We were also able to show the $\beta$-resonance curves where $\beta_{n,1}:\beta_{m,1} = n:m$. We showed that extremely rich and complex dynamics can emerge at the intersection of these curves, at resonant double Hopf points. 
\vs
We used numerical simulations and computations to show the existence of quasiperiodic multiconsensus and chaotic multiconsensus following double Hopf bifurcation at these critical points. These include chaotic saddles, where trajectories follow chaotic dynamics for extremely long timescales before returning to quasiperiodic or steady-state dynamics. This, along with the observed switching of dynamics at every scale of change in $\alpha$ past the critical value, is strongly suggestive of riddled basins. The chaotic multiconsensus dynamics were characterized by a fast carrier wave at the frequencies predicted by the imaginary eigenvalues at the double Hopf bifurcation, amplitude-modulated by a chaotic envelope function. We extracted the envelope function using a numerical Hilbert transform and took its Fourier transform, showing a broadband power spectrum consistent with chaos. We also numerically computed the largest Lyapunov exponent of both the extracted envelope (trimmed to remove distortion) and the underlying signal itself, using the nolds library for Python, showing the positive largest Lyapunov exponent also characteristic of chaotic trajectories. We also computed the first-return maps by comparing successive amplitude peaks. These show single fixed points for periodic solutions, closed curves (or projections of the 2-torus) for quasiperiodic solutions, and complex fractal structures for chaotic dynamics. Certain trajectories were observed to grow as quasiperiodic solutions, transition to chaos, and then transition to a quasiperiodic stable end-state, and their first return maps clearly illustrate the ``wrinkling'' and disintegration of a quasiperiodic torus into chaotic trajectories. This is consistent with the RTN route to chaos, where a 2-torus bifurcates to a 3-torus via a third incommensurate frequency, which is structurally unstable and results in chaotic solutions. 
\vs
The existence of chaotic dynamics in this type of system has powerful implications in protocol design for multi-agent systems. It illustrates that for protocols with memory, especially momentum/neutral memory, there are parameter regimes where, although the multiconsensus relationships (i.e. the spatio-temporal symmetries guaranteed by the degree) hold throughout, the trajectories can nonetheless evolve in complex and unpredictable ways. In particular, the likely existence of riddled basins has serious implications, as it implies the existence of parameter regimes that are not only chaotic, but within which any choice of any parameters, or any approximation, could eventually wind up in the basin of a totally different attractor, and even miniscule floating point rounding errors in precise numerical simulations are enough to not only initiate different chaotic trajectories, but to send solutions into completely different basins of attraction over long timescales. This is compounded upon the inherent unpredictability of chaotic dynamics, as it additionally becomes impossible to even predict if dynamics will be chaotic or quasiperiodic for a given set of parameter values, and solutions may dwell along the boundary of the basin of attraction for one of these attractors for a very long time (on the order of $t=10^6$ in our case) before eventually finding their way to a stable non-chaotic basin (which often includes steady-state solutions). This suggests an extreme fragility of decision-making at these thresholds.
\vs
In each chapter, we confirmed the theoretical results using numerical simulations, and further enhanced these results through numerical stability computations. In this way, we demonstrated a comprehensive set of tools for analyzing multiconsensus and consensus-breaking bifurcation in multi-agent systems in general, and ones featuring memory in particular. The combination of the strength of topological results (e.g. global continuation of branches, spatio-temporal symmetric classification and guarantees that solutions corresponding to maximal orbit types must exist) with robust numerical simulations is a powerful paradigm, where each method enhances the other and covers its limitations and deficiencies. This is especially the case in multi-agent systems with homogeneous agents and highly symmetric interaction topologies, which are natural modeling choices for many systems. 
\vs
\section{Future work}
\subsection{Alternative bifurcation parameters}
We consider the gain or weight of continuous memory as the bifurcation parameter. There were a number of motivations behind this decision. One reason is that this lets us explore how the influence of memory changes agent dynamics, as opposed to the length of the memory itself. Another is because we view memory as potentially a design feature of these protocols, or a useful component for modeling the behavior of real intelligent agents with memory such as humans, animals, or artificial intelligences. However, it is very common in the literature when considering time lags in multi-agent systems to take the delay length as the bifurcation parameter, and there is of course no obstruction to doing so with these models. One would only need to recompute the characteristic equations, re-derive the transcendental equations for the limit frequencies and critical points, and obtain new transversality conditions to compute the signs of crossing numbers. However, this could be somewhat more complicated for $\tau_1$ or $\tau_2$ as bifurcation parameters, since $\tau_1$ and $\tau_2$ appear inside trigonometric functions.

\subsection{Complexity collapse in pseudoneutral equations}
Recent work by Cassidy \cite{Cassidy2025} gives convergence results and error estimation for such quadrature of distributed delays, and partially addresses the observed phenomenon, known sometimes as ``complexity collapse,'' where delay equations with multiple discrete delays exhibit more complex dynamics as the number of delays is increased up to a certain point where this complexity sharply falls off. In some cases, this can be explained by multiple discrete delays approximating a single distributed delay sufficiently well. Given what we have observed in this dissertation regarding pseudoneutral equations, it would be interesting to study this same phenomenon in that context. In other words, is there a type of ``neutral complexity collapse'' that occurs when adding more neutral delay terms to an equation, such that the neutral terms are a sufficiently good approximation of an integral and thus the dynamics are sufficiently close to those of a pseudoneutral equation of retarded type? This is an interesting question which, to our knowledge at the time of publication, has not yet been studied. Such a phenomenon would show that pseudo‑neutral equations are not merely a mathematical curiosity but the natural limit of physical neutral networks as the number of independent lags grows.

\subsection{Deepening the UAV formation control application}
Our UAV formation control application made a number of assumptions which facilitated a clearer and more direct analysis, particularly regarding the symmetries and coupling. We believe these assumptions are entirely reasonable under reasonable conditions, particularly if the formation is moving sufficiently slowly and each UAV can sense its neighbors positions and deviations from formation set points. However, these requirements are only related to the isotypic decomposition and the specific forms of the axis-decoupled control functions, and are not related to the model itself. It would be entirely possible to consider protocols of the exact same form applied to the pitch, yaw, and roll axes of a UAV to control its motion at a much lower level, and to take the coupling functions in terms of something more physically reasonable (but mathematically more complex), such as a euclidean distance between a UAV and its neighbors. Both of these cases would result in more complex coupling, a larger state space, more $S_4\times \mathbb Z_2$-isotypic components, and richer spatio-temporal symmetries in the multiconsensus patterns. This could also be connected to other types of central pattern generator in other robotics applications. 

\subsection{Further analysis of chaotic multiconsensus in the neutral system}
The chaotic dynamics in the neutral system could be investigated on a much deeper level. To more definitively determine the existence of riddled basins and explore the shape of the various attractors near the resonant double Hopf points, we could conduct a numerical study to determine the boundaries of various basins, by computing the fraction of initial conditions in a fine mesh which reach each attractor. Ideally, a rigorous normal-form study could also be performed to show that the symmetry action satisfies the transversality conditions required for riddling. Additionally, we only explored one very small region of the parameter space in detail, and it is entirely possible that other very rich and complex forms of dynamics could prevail in other regions. 
\vs
Another potential enhancement would be to explore the double Hopf curves and formation of chaotic dynamics on other isotypic components, in particular the three-dimensional components $V_3$ and $V_4$. These components admit many more maximal orbit types, which could lead to many interesting dynamics and interactions. For example, Chossat and Golubitsky \cite{chossat1988symmetry} showed that when symmetric systems bifurcate to chaotic attractors, this can lead to symmetry \emph{increasing} bifurcations, where the bifurcating chaotic solution (in this interpretation, a chaotic multiconsensus) actually has more symmetries than the branch it bifurcated from. This is not really possible to explore on the $V_1$ component as it is one-dimensional and only admits a single maximal orbit type. Extending this to higher-dimensional isotypic components with more maximal orbit types could greatly enhance the study of chaotic multiconsensus structures.

\subsection{Analysis of double Hopf curves, branch switching, and chaos in the pseudoneutral system}
We did not perform the same type of study of the $(\gamma,\tau_1)$ parameter space for the pseudoneutral system. One reason for this is that, contrary to what one might expect, the corresponding equations in the pseudoneutral system are much less well-conditioned for this type of numerical analysis. However, we have found strong numerical and circumstantial evidence that double Hopf points and $\beta$-resonance points also occur in this system, though whether they form analytic closed curves as in the neutral case remains to be studied. In this case, it is also reasonable to expect that chaos could emerge after these double Hopf points, especially at or near $\beta$-resonances. 
Studying this phenomenon and finding the shape and structure of these double Hopf sets, and addressing the question of chaotic attractors and riddling in this system could also further illuminate the connections between pseudoneutral systems and neutral systems.

\subsection{Adding hysteretic memory protocols}
We have considered three types of memory in this dissertation with three distinct interpretations. One natural next step would be to explore hysteresis memory terms, which could model trust and betrayal between agents. For example, if the state-dependent coupling is governed by a Preisach operator, this could model agents whose response to their neighbors behavior could permanently change in response to extreme events, representing a gradual erosion of trust and the resulting changes in dynamics. This is a natural counterpart to time-based memory schemes such as the continuous memory, trend memory, and momentum memory protocols studied in this dissertation. A hysteretic Preisach operator would instead represent a kind of state-based memory, ideal for modeling trauma or long-term memory formation, while time-based memory is a better model for short-term memory. 
\vs
The Preisach operator has advantages in this sense because it is Lipshcitz continuous and can give a protocol which is Fr\`echet differentiable, if designed carefully (e.g. with a dead zone around the $\bm x \equiv 0$ consensus equilibrium). Such a system could then be formulated as a condensing perturbation of identity, similar to the framework introduced in Chapter \ref{chapter:pseudo} and used in Chapter \ref{chapter:neutral}. This allows us to use the twisted equivariant Nussbaum-Sadovskii degree to show Hopf bifurcation, and study the consensus-breaking bifurcation problem in these types of systems as well.

\subsection{Studying higher order systems with memory}
In this dissertation, we only considered single-integrator dynamics which yielded first order delay differential equations. There is no fundamental theoretical obstacle to applying these methods to higher order systems, and it would introduce a rich new typology of delays. For example, the trend and momentum memories as defined in this dissertation would no longer be pseudoneutral or neutral in nature. Rather, higher order terms involving delays in the highest order derivative could occur. There are already promising applications for second order systems with the types of memory discussed in this dissertation for single-lane car following traffic models, where drivers are agents whose protocol is based on maintaining an equilibrium of speed up to a speed-dependent safe stopping distance with the car in front of them. As the density of cars increases, the available space also decreases, and the maximum target speed diminishes as a result, yielding a good model for congestion. 
\vs
In this case, the distributed delay terms represent the expected value of reaction time delays to changes in the relative velocity of the nearest car, taken in a mean field approximation. Neutral memory, in this case, would correspond to the ``internal'' reaction delay of a driver reacting to changes in their own acceleration due to their own previous accelerator/brake inputs, modeling the fact that driver acceleration is not based on a known target acceleration level, but rather follows a heuristic pattern where acceleration and braking is done relative to a feedback loop of external information such as felt acceleration, visual stimuli such as speedometer readings, and other outside factors. In this case, taking a mean field approximation across a distribution of such driver-internal delay times would yield a pseudoneutral equation. If neutral delays are indeed a reasonable model of these types of inertial effects in decision-making, then the pseudoneutral mean field approximation can explain why associated real world systems often do not exhibit the same types of instabilities associated with neutral equations.

\subsection{Robustness of memory protocols}
An immediate practical question is whether the distributed delay protocols studied here are robust, in the sense that the stability of consensus or periodic multiconsensus solutions survive in the presence of signal noise, communication drop-outs, or parameter mismatches. The continuous memory term, mathematically a convolution with a step function in our models, acts as a low-pass filter, attenuating high-frequency fluctuations while preserving low-frequency dynamics. This suggests a resilience to noise which could be proven using the framework of robust consensus and $H_\infty$ analysis. This could be extended to the pseudoneutral and neutral memory types and reveal whether these memory components improve or degrade robustness. Numerical simulations already hint that the smoothing effect of continuous memory dominates (at least for reasonably small values of $\gamma$), but a formal proof would firmly establish continuous memory as a design principle for resilient multi-agent coordination.

\appendix

%MULTIPLE DISTRIBUTED
\chapter{Notation for Subdirect Products} \label{appendix:notation}
\section{Motivation}
Spatio-temporal symmetries are naturally characterized by product groups, where one of the factors describes symmetries in time (e.g. $S^1$) and the other describes symmetries in space (e.g. $S_4$). We often also wish to guarantee that solutions are odd, usually as a way of guaranteeing that a periodic solution is non-constant. Oddness can be naturally viewed as $\mathbb Z_2$-equivariance, but in fact, the situation more relevant to our purposes is showing that the subspace of solutions \emph{fixed} by some group can contain only odd solutions. In other words, for a map $x:\mathbb R\to V$, we need
\[
x(t) = -x(-t).
\]
This requires two independent actions: a group acting on $x\in V$ and a group acting on $x(t)\in V$. Consider the Banach space of continuous $2\pi$-periodic functions $\mathscr E:=C_{2\pi}(\mathbb R;\mathbb R)$. Then $\mathscr E$ is a natural isometric $S^1$-representation, where $\theta\in S^1$ acts on $x(t)\in \mathscr E$ as $x(t+\theta)$. It is also a natural $\mathbb Z_2$-representation, where $-1 \in \mathbb Z_2$ maps $x(t)$ to $-x(t)$. Notice that $(\pi,-1)\in S^1 \times \mathbb Z_2$ fixes $x\in \mathscr E$ if and only if $x(t)=-x(t+\pi)$, which implies that either $x$ is odd and non-constant, or $x \equiv 0$. If we consider $H \leq S^1$ to be the cyclic subgroup generated by $\pi \in S^1$, then clearly $H \cong \mathbb Z_2$. However, if $x$ is fixed by all of $H \times \mathbb Z_2 \leq S^1\times \mathbb Z_2$, then it must be fixed by $(0,-1)$ and $(\pi,1)$, which force $x\equiv 0$. 
\vs
This means that neither $S^1$ nor $\mathbb Z_2$ are sufficient to create fixed point subspaces of non-constant periodic solutions on their own. Nor is it enough to simply take the direct product of the subgroups of $S^1$ and $\mathbb Z_2$ which are generated by the elements needed to guarantee oddness. Rather, we need a subgroup composed of specific elements from $S^1$ paired with specific elements of $\mathbb Z_2$. Indeed, the subgroup of $S^1 \times \mathbb Z_2$ generated by $(\pi,-1)$ is of order 2 and is isomorphic to $\mathbb Z_2$, but it is related to the separate actions of $S^1$ and $\mathbb Z_2$ on the space $\mathscr E$. 
\vs
On the other hand, if there are more sophisticated spatial symmetries at play (e.g. a coupled system where $\bm x \in C_{2\pi}(\mathbb R;\mathbb R^n)$, then we may also have a group $\Gamma_0\leq S_n$ which acts on $\bm x(t) := (x_1(t),\dots, x_n(t))^T$ by permuting the indices of the components of $\bm x$. In this case, non-constant periodic solutions must either be locked in identical phase, or they must have certain symmetric phase relationships described by subgroups of $S^1 \times \Gamma_0$ (but not by direct products of subgroups of $S^1$ and $\Gamma_0$).
\vs
In short, we need a way to describe subgroups of product groups which are not products of subgroups. Such subgroups are known as \emph{subdirect products}, and also sometimes as \emph{twisted subgroups} (cf. Section \ref{sec:prelim:twisted-subgroups}). As it turns out, such subdirect products are almost always the most interesting spatio-temporal symmetries and, for a variety of reasons, maximal orbit types almost always take the form of subdirect products, as opposed to direct products of orbit types of the factor groups. In this appendix, we will describe the notational conventions used for such subdirect products throughout this dissertation, called \emph{amalgamated notation}, and we will describe in detail the structure of the subgroup lattices of the groups which are used in the applications across Chapters \ref{chapter:distributed}, \ref{chapter:pseudo}, and \ref{chapter:neutral}, namely $S_4 \times \mathbb Z_2$ and $S^1 \times S_4 \times \mathbb Z_2$.
\section{Amalgamated notation} \label{app:amalgamated}
The amalgamated notation, rigorously defined by Balanov, Krawcewicz, and Steinlein \cite{BalanovKrawcewiczSteinlein2006}, can be used to describe all possible subgroups of a direct product of two groups. Its basis is Goursat's lemma, which can be stated (in one form) as follows:
\begin{lemma}
    Let $G_1, G_2$ be two groups. Then for any $K \leq G_1 \times G_2$, there exist subgroups $H_1 \leq G_1, H_2 \leq G_2$, a group $L$, and epimorpshisms (surjective homomorphisms) $\varphi_1:H_1 \to L$, $\varphi_2:H_2\to L$, such that
    \[
    K = \{(h_1,h_2)\in H_1 \times H_2:\varphi_1(h_1) = \varphi_2(h_2)\}.
    \]
\end{lemma}
In other words, subdirect products can be described as fiber products. Moreover, by the first isomorphism theorem, this also means that $H_1' := \ker \varphi_1$ and $H_2' := \ker \varphi_2$ are normal subgroups of $H_1$ and $H_2$ respectively, and that $H_1/H_1' \cong L \cong  H_2/H_2'$. Therefore, we can view this construction as pairing cosets of $H_1'$ in $H_1$ with cosets of $H_2'$ in $H_2$. This means that elements in $H_1'$ can mix freely with those in $K_1'$, but elements in each coset in $H_1/H_1'$ can only be paired with elements in some coset $H_2/H_2'$ if they have the same image under $\varphi_1$ and $\varphi_2$. This pairing of cosets also induces an isomorphism $\psi :H_1/H_1'\to H_2/H_2'$. This means that a subdirect product can be uniquely specified by a quintuple $(H_1,H_2,\varphi_1,\varphi_2,L)$, which can be written in fiber product notation as
\[
{H_1}^{\varphi_1} \times_L {}^{\varphi_2}H_2.
\]
\vs
Notice that in general, one could pick a different epimorphism $\varphi_1'$ such that $\ker \varphi_1 = \ker \varphi_1'$, but $\varphi_1(h_1) \neq \varphi_1'(h_1)$ for some $h_1\in H_1$. However, a crucial result of \cite{BalanovKrawcewiczSteinlein2006} is that such choices yield conjugate subgroups. This means that from the point of view of orbit types, it is enough to specify the quadruple $(H_1, H_1', H_2,H_2')$, where $H_1'\trianglelefteq H_1, H_2'\trianglelefteq H_2$, and $L$ is implicitly understood as any representative of the isomorphism class of $H_1/H_1'$ (or equivalently $H_2/H_2'$). We then write
\[
{H_1}^{H_1'} \times_L {}^{H_2'}H_2.
\]
Since $L$ can be derived implicitly from the quadruple, it is often omitted, and we simply write
\[
{H_1}^{H_1'} \times {}^{H_2'}H_2.
\]
If either $H_1'=H_1$ or $H_2' = H_2$, (i.e. the corresponding epimorphism is trivial), we omit the corresponding superscript. For example, if both $H_1'=H_1$ and $H_2' = H_2$, then this is simply a product group, and we have
\[
{H_1}^{H_1'} \times_{\mathbb Z_1} {}^{H_2'}H_2 = H_1 \times H_2.
\]
\section{Notation for subgroups of \texorpdfstring{$S_4 \times \mathbb Z_2$}{S₄ × ℤ₂}}
As already described, our search for spatio-temporal symmetries leads us to consider products of three groups, generically of the form $S^1 \times \Gamma_0 \times \mathbb Z_2$, where $\Gamma_0$ is a finite group. While the amalgamated notation (and Goursat's lemma itself) is only defined for a direct product of two groups, it is of course possible to apply it to larger products inductively. However, this quickly creates an extreme notational burden even for simple groups. The solution advanced by \cite{BalanovKrawcewiczSteinlein2006} is to use an abbreviated notation for subgroups of $\Gamma_0 \times \mathbb Z_2$, on the basis that there are a limited number of ways that such amalgamated subgroups can be formed. Here we will provide the notational conventions used for subgroups of $S_4 \times \mathbb Z_2$:

\begin{itemize}
    \item For any $H\leq S_4$, $H^p := H \times \mathbb Z_2$ ($p$ stands for \emph{product}).
    \item For $D_n\leq S_4$, $D_n^z := {D_n}^{\mathbb Z_n} \times {}^{\mathbb Z_1}\mathbb Z_2$ ($z$ stands for \emph{cyclic}, as in $\mathbb Z_n$).
    \item For $D_{2n} \leq S_4$, ${D_{2n}^d}:={D_{2n}}^{D_n} \times {}^{\mathbb Z_1}\mathbb Z_2$ ($d$ stands for \emph{dihedral}, as in $D_n$).
    \item For $D_{2n} \leq S_4$, ${D_{2n}^{\tilde d}}:={D_{2n}}^{\widetilde D_n} \times {}^{\mathbb Z_1}\mathbb Z_2$, where $\widetilde D_n$ is the other subgroup of index 2 in $D_{2n}$ which is isomorphic but not conjugate to $D_n$.
    \item For $\mathbb Z_{2n}\leq S_4$, $\mathbb Z_{2n}^- := {\mathbb Z_{2n}}^{\mathbb Z_n} \times {}^{\mathbb Z_1}\mathbb Z_2$, (where $-$ stands generically for taking the subgroup of index 2 when there is only one such choice).
    \item $S_4^- := {S_4}^{A_4} \times {}^{\mathbb Z_1}\mathbb Z_2$.
\end{itemize}

One can easily confirm that these comprise all the index 2 quotient groups of orbit types of $S_4$. We can use these to list all subgroups of $S_4 \times \mathbb Z_2$:

\begin{align*}
\bz_2^-&=\{((1),1),((12)(34),-1)\},\\
%\bz_3^t&=\{((1),1),((123),\gamma),((132),\gamma^2)\},\\
%\bz_4^c&=\{((1),1),((1324),i),((12)(34),-1),((1423),-i)\},\\
\bz_4^-&=\{((1),1),((1324),-1),((12)(34),1),((1423),-1)\},\\
D_1^z&=\{((1),1),((12),-1)\},\\
V_4^-&=\{((1),1),((12)(34),1),((13)(24),-1),((14)(23),-1)\},\\
D_2^d&=\{((1),1),((12)(34),-1),((12),1),((34),-1)\},\\
D_2^z&=\{((1),1),((12)(34),1),((12),-1),((34),-1)\},\\
D_3^z&=\{((1),1),((123),1),((132),1),((12),-1),((23),-1),((13),-1)\}, \\
D_4^d&=\{((1),1),((1324),-1),((12)(34),1),((1423),-1),((34),1),\\
&~~~~~((14)(23),-1),((12),1),((13)(24),-1)\},\\
D_4^{\hat d}&=\{((1),1),((1324),-1),((12)(34),1),((1423),-1),((34),-1),\\
&~~~~~((14)(23),1),((12),-1),((13)(24),1)\},\\
D_4^z&=\{((1),1),((1324),1),((12)(34),1),((1423),1),((34),-1),\\
&~~~~~((14)(23),-1),((12),-1),((13)(24),-1)\},\\
S_4^-&=\{((1),1),((12),-1),((12)(34),1),((123),1),((1234),-1),((13),-1),\\
&~~~~~((13)(24),1),((132),1),((1342),-1),((14),-1),((14)(23),1),((142),1),\\
&~~~~~((1324),-1),((23),-1),((124),1),((1243),-1),((24),-1),((134),1),\\
&~~~~~((1423),-1),((34),-1),((143),1),((1432),-1),((243),1),((234),1)\}.
\end{align*}

\chapter{Software and Numerical Methods}\label{appendix:software}
Here we will go into more detail about the software packages and numerical methods used in this dissertation. For the numerical simulation of the DDEs, the JiTCDDE library in Python \cite{jitcxde} was used, and the DifferentialEquations.jl library for Julia \cite{rackauckas2017differentialequations} was used with the Rosenbrock23 solver. The Julia solver was secondary and was only used to confirm and validate the results obtained in Python, especially for the neutral equation, in situations where stiffness was suspected. Mathematica was used to analyze characteristic equations and to assist in their symbolic differentiation, and to obtain numerical solutions to the transcendental coincidence problems with arbitrary precision. Finally, GAP was used to compute basic degrees and maximal twisted orbit types. 

\section{Overview of the numerical method}
The JiTCDDE library for Python uses the the Bogacki-Shampine embedded Runge-Kutta pair (3(2)) with adaptive step-size control. This method is an explicit Runge-Kutta method with local error control. Delays are handled using the Shampine-Thompson method. The state and its derivative are stored at each successful integration step. The required past values $x(t-\tau)$ are obtained via piecewise cubic Hermite spline interpolation between stored anchor points consisting of a time, state, derivative triple. The JiTCDDE package is a just-in-time C compiled DDE solver. It takes SymPy expressions, symbolically differentiates them, then automatically translates this into C code which is compiled and linked dynamically. This gives a major speed boost over calculations done entirely in python which is quite significant for DDEs.
\vs
We chose a long enough time interval to obtain solutions near the stable limit cycle, the exact length of which differed slightly between each chapter. This length in time was dividied into an equally spaced mesh of points used as the time scale for the solver (although the solver has its own error-based adaptive step-size control within this mesh).

\subsection{Integration parameters}
The absolute and relative error tolerances were chosen as 1e-12 to ensure high accuracy. To resolve the fastest dynamics in the neutral equation, it was also necessary to set the maximum time step to 0.01. 

\section{Choosing state histories for numerical analysis of neutral equations}
When numerically simulating delay equations, one must provide the numerical solver with a state history in order for the problem to be well-posed. This presents an obvious paradoxical problem: in principle, the state history must satisfy the DDE for all $t \in [-r,0]$ (where $r = \max\{\tau_1,\tau_2\}$ in this case), but this requires knowledge of the very thing we are trying to simulate! In practice, we provide the solver with an invented history which agrees with the DDE at $t=0$ for all derivatives up to the highest order as is practicable. This is frequently just the first order, as determining these conditions requires one to solve systems of nonlinear equations which can be quite difficult. However, this means there are inevitable discontinuities around $t=0$ for all derivatives of a high enough order. 
\vs
For delay equations of retarded type, if we are interested in the long-term behavior of solutions, then this is not a problem, as the delay term (as well as the distributed delay, separately) act in a ``smoothing'' fashion and, after the delay term has elapsed a suitable number of times, the fabricated state history for $t\in [-r,0]$ will have been replaced with a real state history which is smooth for all derivatives up to as high an order as desired (provided the right-hand side is sufficiently smooth, of course). This means that, for all intents and purposes, we do not have to be much concerned with the validity of the state history and can set it to be identically zero, and can simply consider a long enough simulation time to allow this discontinuity to be smoothed. Such long simulation times are often necessary anyway to allow slowly growing quasiperiodic solutions to approach the periodic limit cycles of interest.
\vs
However, neutral equations pose a major problem in this regard. Due to the delay present in the derivative, any discontinuities in the supplied state history will propagate indefinitely through derivatives of all orders. This means that neutral equations \emph{do} have a strong and sensitive dependence on the state history given to the solver, and different state histories can provide dramatically different results in the solver at all time scales. This can result in the shifting of bifurcation points, the altering of critical frequencies, and other types of obstacles and uncertainties in the trustworthiness of numerical results for these types of problems.
\vs
\subsection{Choosing a consistent state history}
Here we will show a method for generating a relatively consistent state history for the neutral equation. We will then explain some of the complications this state history induces in the numerical problem. By a ``consistent state history,'' we mean $x_t\in C^1([-r,0];\mathbb R^n)$ where $r := \max\{\tau_1,\tau_2\}$ which agrees with the DDE (or at least its linearization) as much as is possible. At the very least, we desire that $x(0)$ and $\dot x(0)$ are equal to the starting conditions for the solver. As before, our numerical method involves choosing a perturbation near $x\equiv0$ in the form of a constant function $x(t) \equiv v_0$, where $v_0$ is a vector chosen on the $\Gamma$-isotypic component on which bifurcation occurs. This means that the state history $x_t \equiv 0$. We must also supply the derivative history $\dot x_t$. Of course, in theory, supplying a function segment $x_t \in C^1([-r,0];\mathbb R^8)$ would give both the state and derivative history at once. However, for practical reasons, it is usually easier to calculate the state and derivative at discrete convenient points in the past and interpolate a $C^1$ function between them. Setting $t=0$ in the linearized system, we obtain
\[
(\dot x(0) - \gamma \dot x(-\tau_1))v_0 = -av_0 -\alpha b \tau_2v_0 -Cv_0.
\]
Therefore, we can obtain the proper derivative at $t=0$ by choosing an appropriate value for $\gamma \dot x (-\tau_1)$. If we assume that $\dot x(-\tau_1) = 0$, then this implies $\dot x(0)v_0 = (-a_j-\alpha b \tau_2)v_0$, where $a_j := a + \mu_j$ and $\mu_j$ is the eigenvector of $C$ corresponding to $v_0\in E(\mu_j)$, as before. On the other hand, we can repeat this process recursively to solve for $\dot x(-\tau_1)$. If we consider $t=-N\tau_1$, where $N\in \mathbb N$, we get
\[
\dot x(-N\tau_1) = \gamma\dot x(-(N+1)\tau_1) -a_j -\alpha b\tau_2.
\]
Therefore,
\[
\dot x(-N\tau_1) = \frac{1}{\gamma}\Big(\dot x\big(-(N-1)\tau_1\big)-a_j-\alpha b\tau_2\Big), \quad \text{ for $N\geq1$}.
\]
If we assume that $\tau_2 = c\tau_1$ for some $c\in \mathbb N$, as will be the case in this example, then we obtain a sequence of $c+1$ points, $-c\tau_1,-(c-1)\tau_1,\dots,-\tau_1,0$ at which we can specify a derivative history which agrees with the DDE, with values in each intermediate interval $t\in(-m\tau_1),-(m-1)\tau_1)$ for $m=1,2,\dots,c$\footnote{The usage of $m$ in this section is local and self-contained, and is not to be confused with any other usage of this letter elsewhere} interpolated via cubic Hermite splines. By repeated back-substitution into the recurrence formula, we obtain a closed form expression for these points, for $m=0,1,\dots,c$:
\[
\dot x(-m\tau_1) = \gamma^{c-m} \dot x(-\tau_2)-(a_j+\alpha b\tau_2)(\gamma^{c-m-1}+\gamma^{c-m-2}+\dots+\gamma+1).
\]
Then we only need to decide on a convenient value of $\dot x (-\tau_2)$ at the beginning of the state history in order to obtain a set of consistent anchor points for the state history. We choose $\dot x(-\tau_2) = 0$, and thereby obtain
\[
\dot x(-m\tau_1) = -(a_j+\alpha b\tau_2)(\gamma^{c-m-1}+\gamma^{c-m-2}+\dots+\gamma+1).
\]
Since $(\gamma^{c-m-1}+\gamma^{c-m-2}+\dots+\gamma+1)$ is a partial sum of a geometric series, this can be simplified to
\[
\dot x(-m\tau_1) = -(a_j+\alpha b\tau_2)\frac{1-\gamma^{^{c-m}}}{1-\gamma}.
\]
\vs
Note that there is no guarantee at all that the interpolated cubic Hermite splines will be consistent with the DDE, and generically of course they will not be. However, for relatively small values of $\gamma$, any effect from these discrepancies in the state history is very small (and grows asymptotically smaller with each $\tau_1$ period), and has no qualitative effect on the results from the point of view of Hopf bifurcation. For larger values of $\gamma$, however, these effect could be significant. 

\subsection{Implementation of the distributed delay}
There are basically two ways to implement the distributed delay term in each of the three systems studied numerically in this dissertation. One method is quadrature (e.g. Gauss-Legendre quadrature), over a suitably fine mesh of discrete delays. The other is the ``linear chain trick''. This method involves writing the distributed delay as an auxiliary state variable
\[
D_i(t) := \int_0^{\tau_2}x(t-s)ds,\quad \dot D_i(t) = x(t) - x(t-\tau_2),
\]
and transforming the $n$-dimensional system of coupled DDEs into a $2N$-dimensional system made up of pairs
\begin{align*}
    \frac{d}{dt}\left[x_i(t)-\int_0^{\tau_1}g(x_i(t-s))ds\right]&=-ax_i(t) - f(D_i(t)) - h_i(x_i(t))\\
    \dot D_i(t) &= x_i(t) - x_i(t-\tau_2).
\end{align*}
The two methods are numerically equivalent if performed with sufficient precision, and the choice depends on the solver. We used both methods across both solvers for all numerical results in this dissertation. For JiTCDDE, which symbolically differentiates expressions and then compiles performant C code for integration, the quadrature method greatly lengthened compilation and simulation times without any enhancement in quality of the results, and so we primarily used the linear chain trick.
\vs

\section{NOLDS}
The NOnLinear measures for Dynamical Systems (NOLDS) library \cite{scholzel_nolds} for Python was used to numerically compute largest Lyapunov exponents as part of our numerical evidence for chaos in Chapter \ref{chapter:neutral}. In particular, we used NOLDS' lyap\_r method, which uses the algorithm developed by Rosenstein et al. \cite{rosenstein1993} to estimate largest Lyapunov exponents (LLEs) from time series data. The embedding dimension was chosen to be 4\footnote{This was chosen heuristically. Several choices of embedding dimension were tested, from 3 to 10, for representative chaotic signals, and did not produce any meaningful change in the LLE.}, and let NOLDS determine the optimal time separation by taking the Fourier transform of the signal, in order to ensure that the Lyapunov exponent calculations are based on the overall dynamics and not spurious relationships due to the periodic nature of the underlying carrier signal.

\subsection{Subsampling the chaotic signals}
The NOLDS lyap\_r method is $\mathcal O(N^2)$, and since the chaotic dynamics of the neutral system evolve over very long timescales, and we use an equally large number of points to ensure those dynamics are captured accurately, this means that it was not at all computationally feasible to apply the lyap\_r method to the entire signal. As an exmaple, many of the simulations for chaotic dynamics collected $10^6$ data points. Applying lyap\_r to this many points requires the initialization of arrays of at least $10^{12}$ 64 bit floats, which far outstripped the resources of the modest desktop computer on which these simulations were run. Accordingly, it was necessary to subsample our data to obtain a more reasonable number of points on which the lyap\_r algorithm could be run. We generally took around 15,000 points, which we found was within the capabilities of our hardware. The LLE is a statistical measure if the long-term divergence rate, and is known to be robust to downsampling.

\section{Mathematica}
Mathematica was used to generate some plots, particularly the plots of the transcendental equations for $\alpha$ and $\beta$, and the plots of complex roots of the characteristic equations. It was also used as a convenient way to numerically solve for roots of these characteristic equations with infinite precision, before truncating these to a level of precision suitable for use in Python. We also implemented a root solver in the Python code using the root finding algorithm from SciPy's optimization toolkit for convenience, which matched the results obtained in Mathematica. Mathematica was also used for symbolic differentiation of the characteristic equations, trigonometric reduction, and algebraic simplification of equations at various points.

\section{GAP}
GAP \cite{GAP4} was used in conjunction with the Equideg package developed by Haopin Wu \cite{equideg} to compute basic twisted degrees and maximal twisted orbit types. This was done by initializing the symmetric group $S_4$ and the cyclic group $\mathbb Z_2$ as permutation groups, taking their direct product, and then taking the direct product of this group with $S^!\cong SO(2)$. GAP was also used to compute the irreducible representations of $S_4 \times \mathbb Z_2$, which were used to determine the $S_4 \times \mathbb Z_2$-isotypic decompositions of $\mathbb R^6$ in Chapter \ref{chapter:distributed}, and $\mathbb R^8$ in Chapter \ref{chapter:pseudo} and Chapter \ref{chapter:neutral}. These were identified by their character tables in GAP, and the conjugacy classes associated with these characters were then identified with their corresponding classes in $S_6$ and $S_8$, respectively. The recent book of Balanov, Krawcewicz, Rachinskiy, Yu, and Wu \cite{BalanovEtAl2025} provides many examples and explicit code snippets for similar computations in GAP. 

\chapter{Supplementary plots, figures, and data}\label{appendix:supplementary-figs}
This appendix contains several additional plots and figures related to the $(\gamma,\tau_1)$ parameter space exploration in Chapter \ref{chapter:neutral}. These figures are included because they may be of use to other researchers who wish to study this system. For every time series plot and first-return map, we will provide exact parameters, although the strongly suspected presence of riddled basins in that parameter region will make it very difficult to exactly reproduce these results in any solver. Some of these first-return maps have been included because they very nicely illustrate the breakdown of a quasiperiodic torus, represented by expanding helical curves, into chaos; and also because we feel subjectively that they possess an eerie beauty.

\section{Critical surfaces in the parameter plane}
This section contains both contour maps and 3D plots of the critical $\alpha_{n,1}$ and $\beta_{n,1}$ surfaces for $n=1,2,\dots,5$ in the neutral system on the $V_1$ isotypic component, with $a=0.5, b=0.2, \tau_1= 60, c_1 = 0.15, c_2 = 0.05, c_3 = 0.01$, over the ranges $0\leq\gamma < 1$ and $0\leq \tau_1 \leq \tau_2$.

\begin{figure}[p] 
    \centering
    \subcaptionbox{Branch 1\label{fig:beta-branch1}}{
        \begin{minipage}[b]{0.43\textwidth}
            \centering
            \includegraphics[width=\linewidth]{figures/neutral/beta_branch_1_contour_smooth.png}
            \vspace{2pt}
            \includegraphics[width=\linewidth]{figures/neutral/beta_branch_1_surface_rotated.png}
        \end{minipage}
    }\hfill
    \subcaptionbox{Branch 2\label{fig:beta-branch2}}{
        \begin{minipage}[b]{0.43\textwidth}
            \centering
            \includegraphics[width=\linewidth]{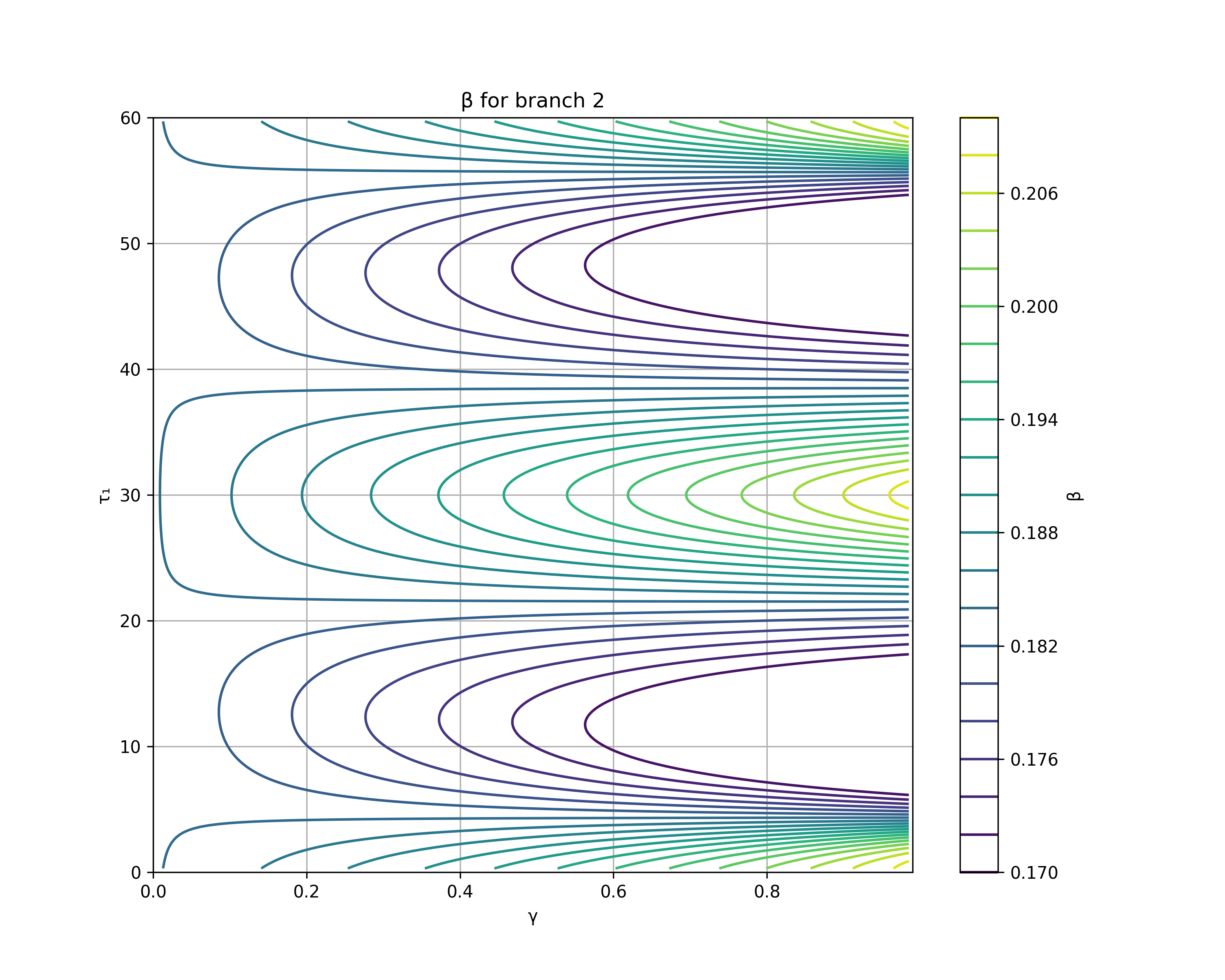}
            \vspace{2pt}
            \includegraphics[width=\linewidth]{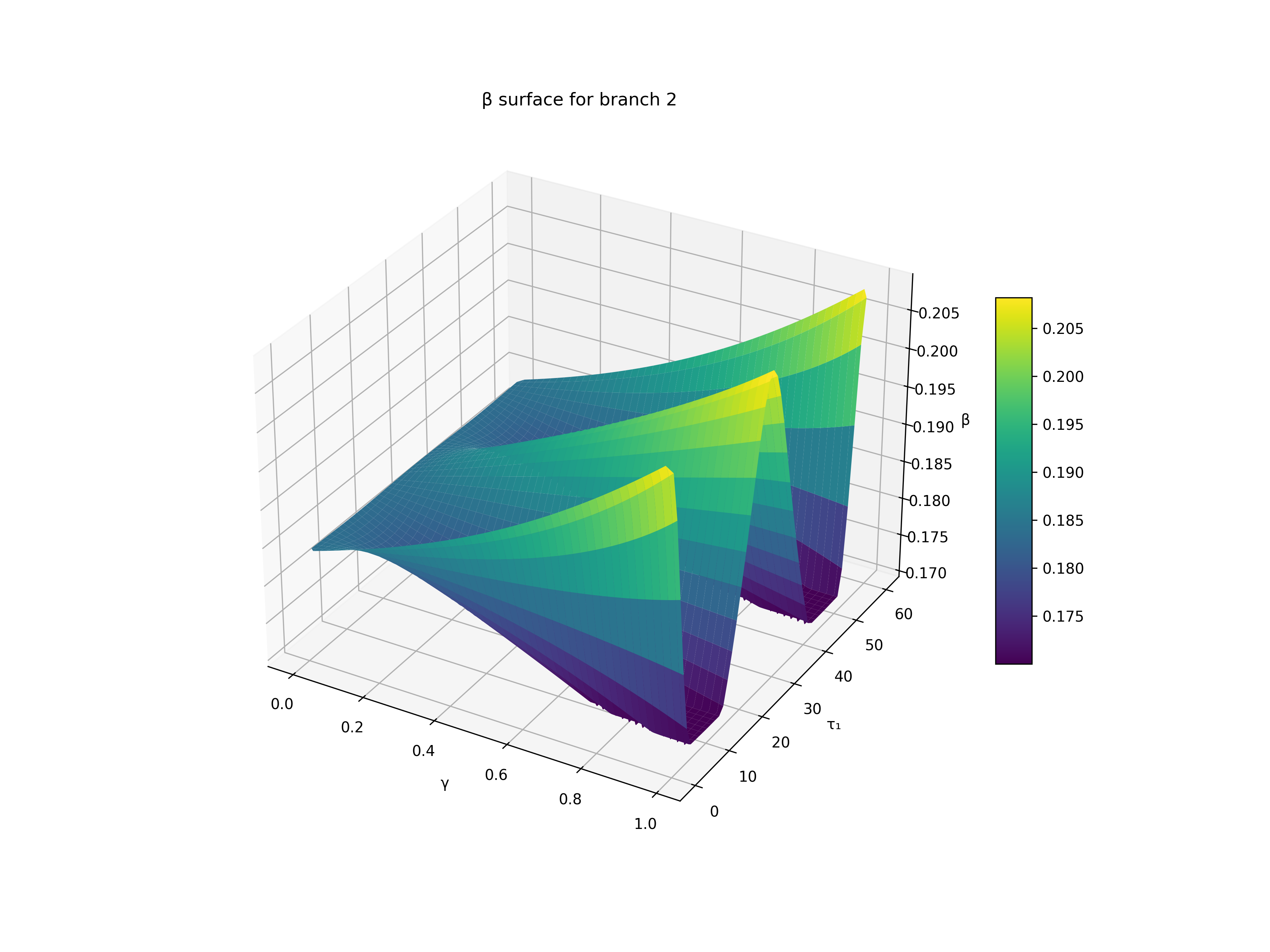}
        \end{minipage}
    }

    \vspace{\floatsep}  % equivalent spacing between rows of subfigures

    % Row 2: branches 3 and 4
    \subcaptionbox{Branch 3\label{fig:beta-branch3}}{
        \begin{minipage}[b]{0.43\textwidth}
            \centering
            \includegraphics[width=\linewidth]{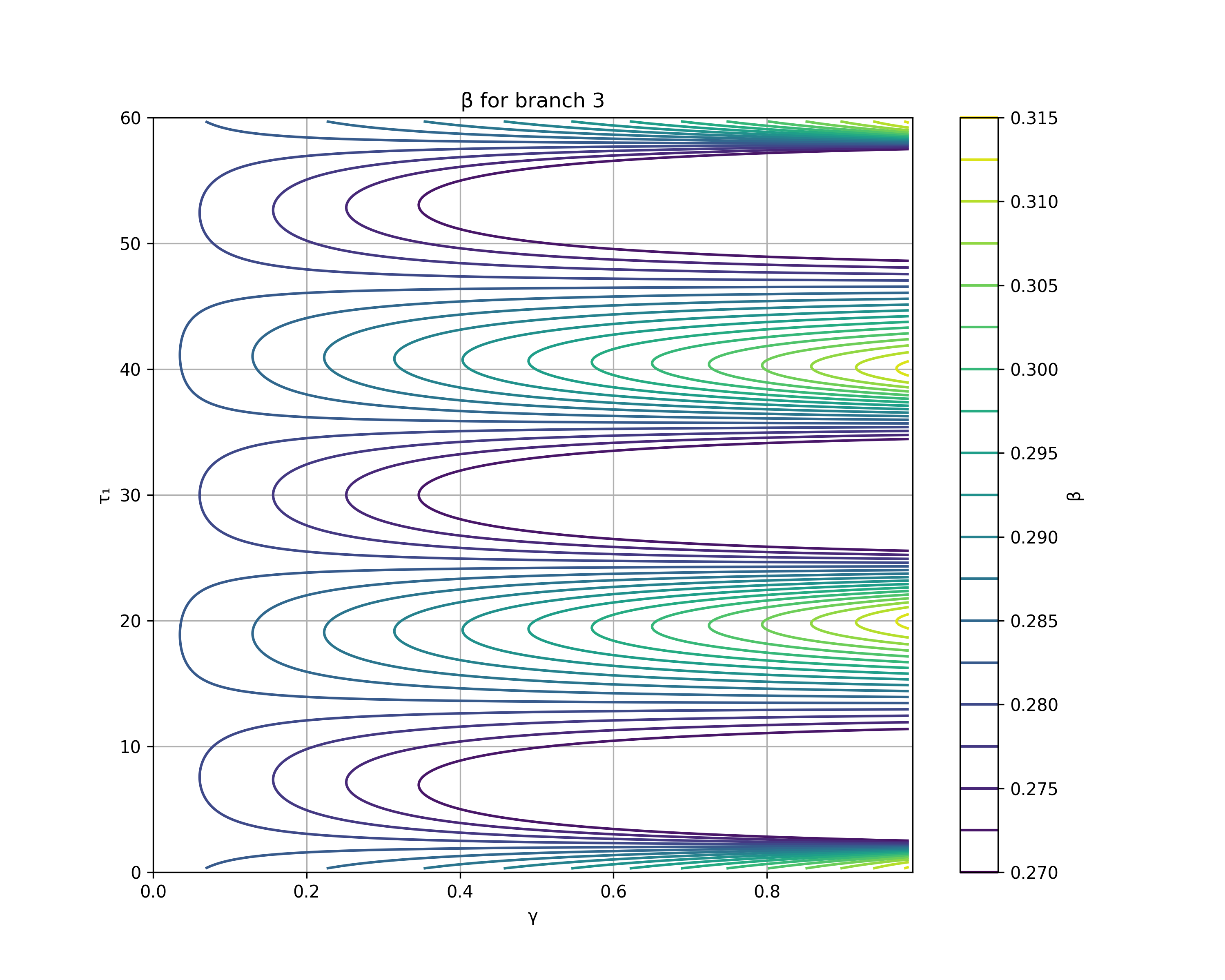}
            \vspace{2pt}
            \includegraphics[width=\linewidth]{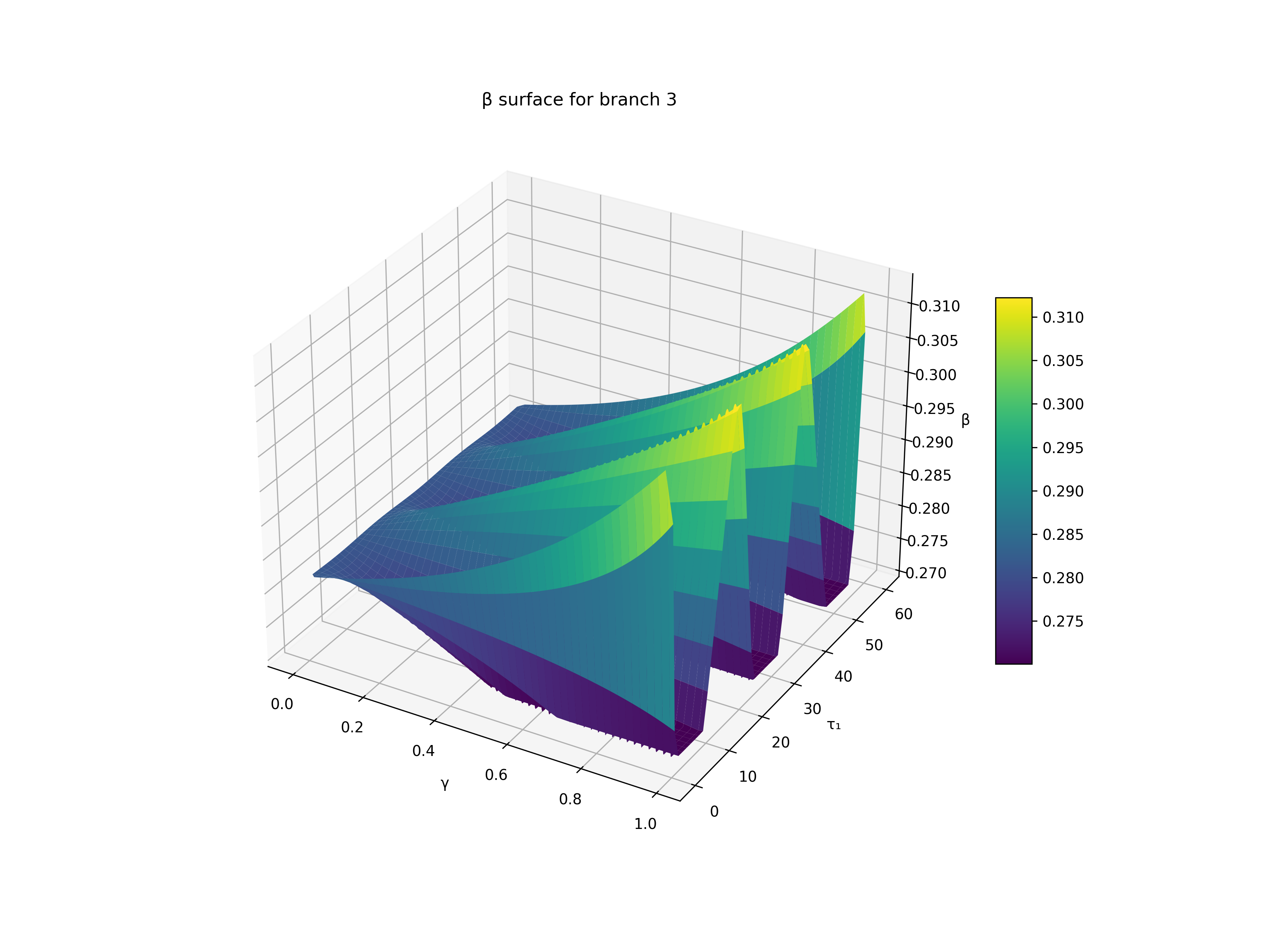}
        \end{minipage}
    }\hfill
    \subcaptionbox{Branch 4\label{fig:beta-branch4}}{
        \begin{minipage}[b]{0.43\textwidth}
            \centering
            \includegraphics[width=\linewidth]{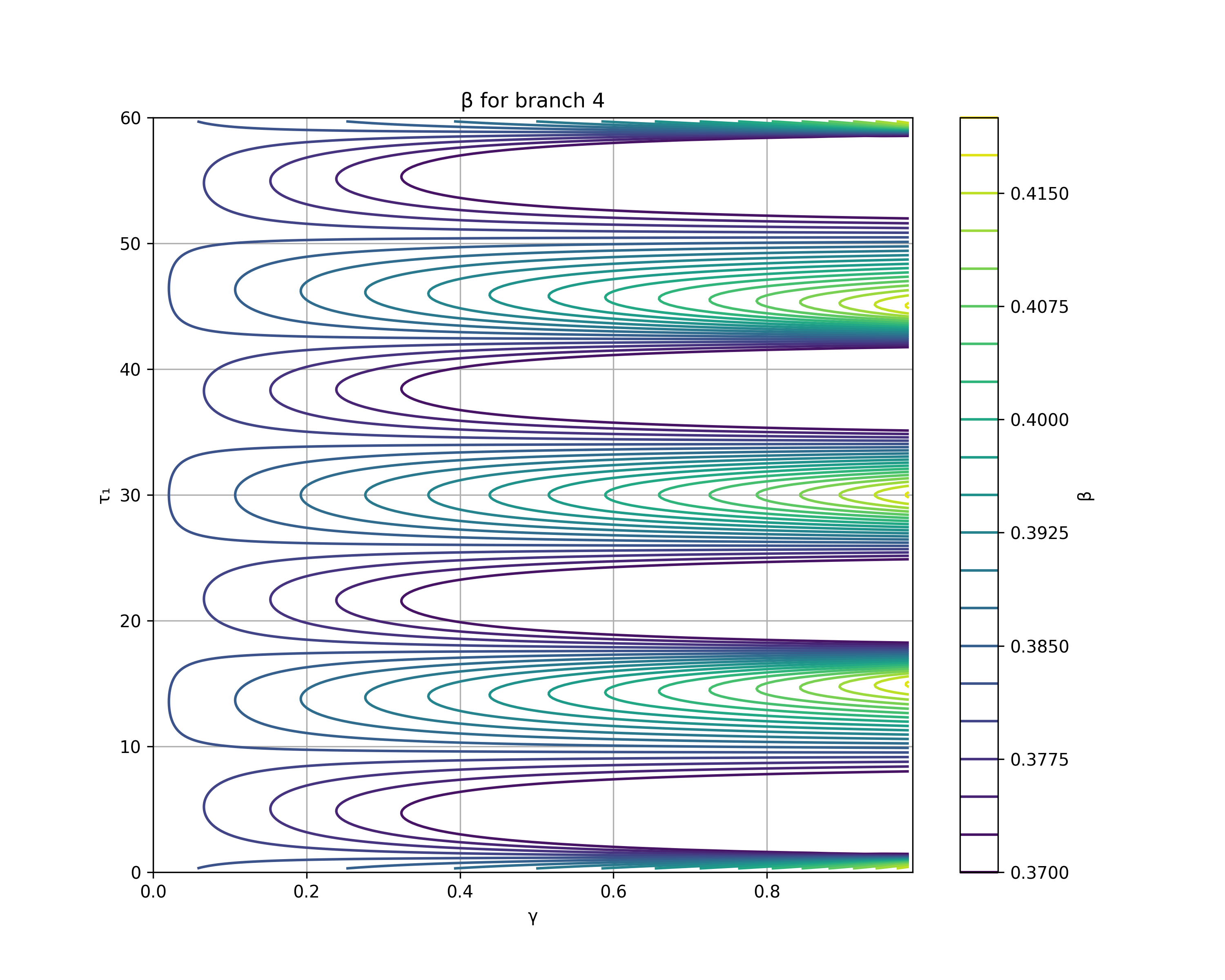}
            \vspace{2pt}
            \includegraphics[width=\linewidth]{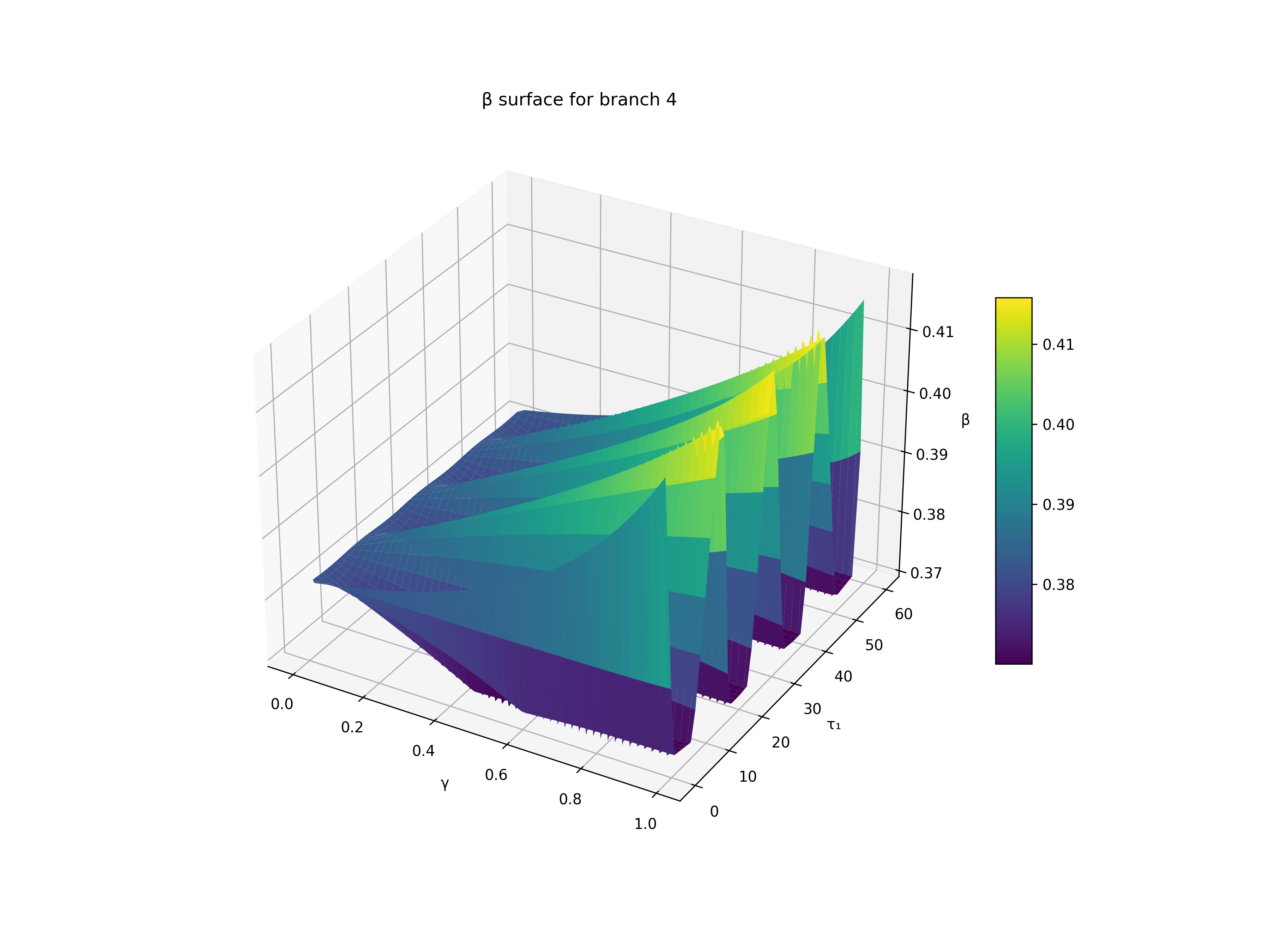}
        \end{minipage}
    }

    \vspace{\floatsep}

    % Row 3: branch 5, centered
    \subcaptionbox{Branch 5\label{fig:beta-branch5}}{
        \begin{minipage}[b]{0.43\textwidth}
            \centering
            \includegraphics[width=\linewidth]{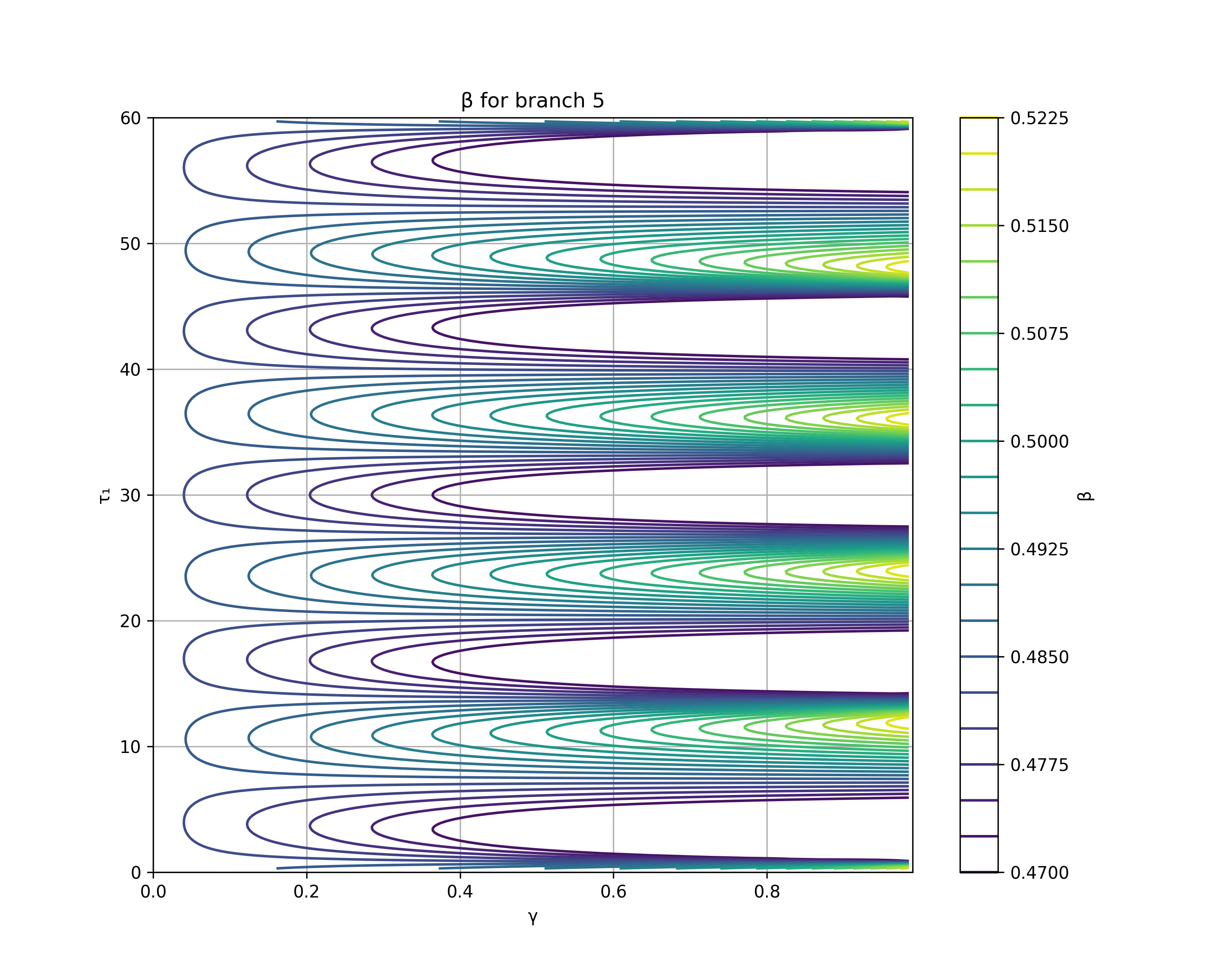}
            \vspace{2pt}
            \includegraphics[width=\linewidth]{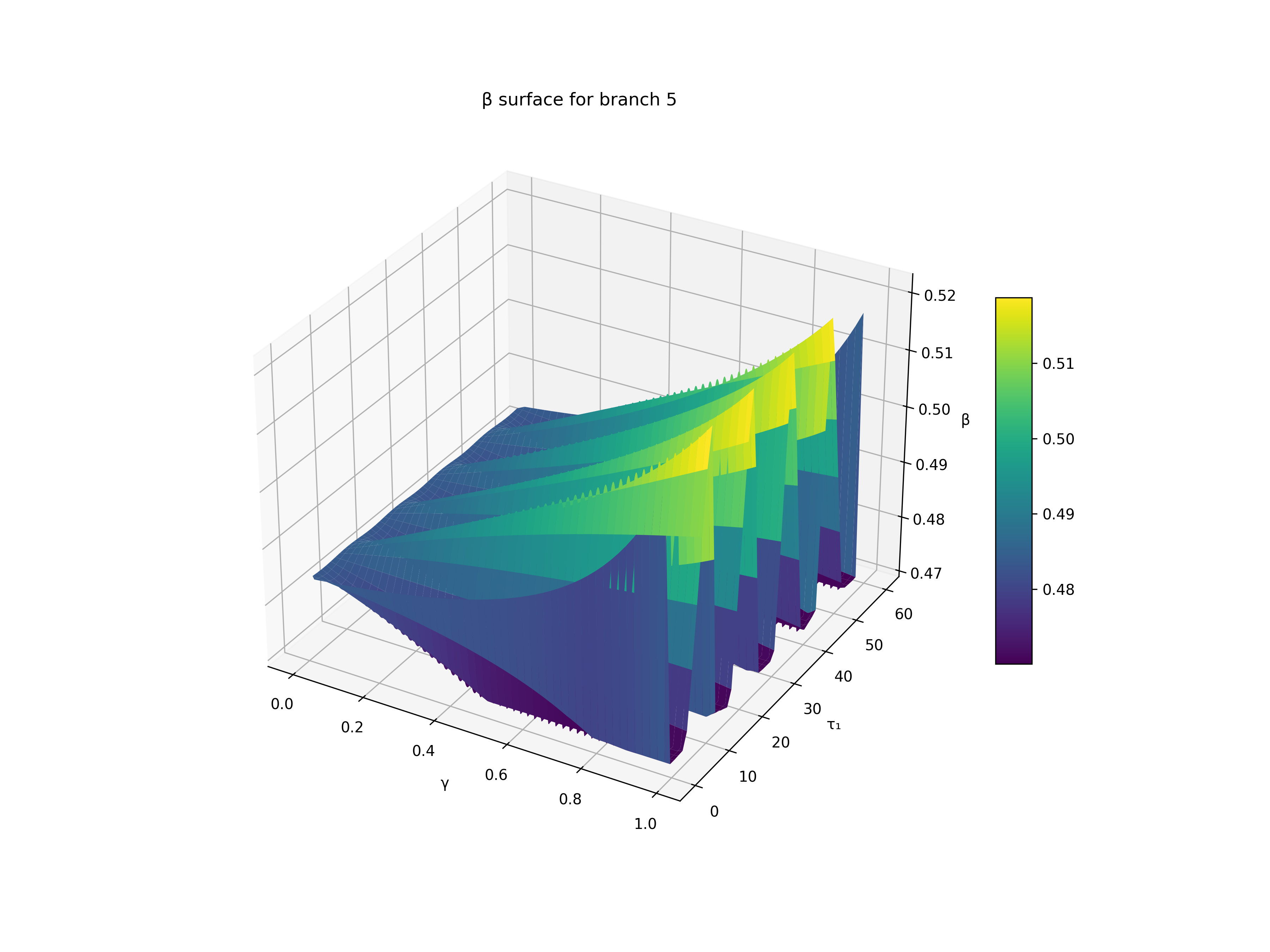}
        \end{minipage}
    }

    \caption{Contour and surface plots of $\beta_{n,1}(\gamma,\tau_1)$ for branches $n=1,\dots,5$. In each subfigure, the left image is the contour map and the right image is the 3D surface.}
    \label{fig:beta-all}
\end{figure}

\begin{figure}[tbp]
    \centering
    \begin{minipage}{0.43\textwidth}
        \centering
        \includegraphics[width=\linewidth]{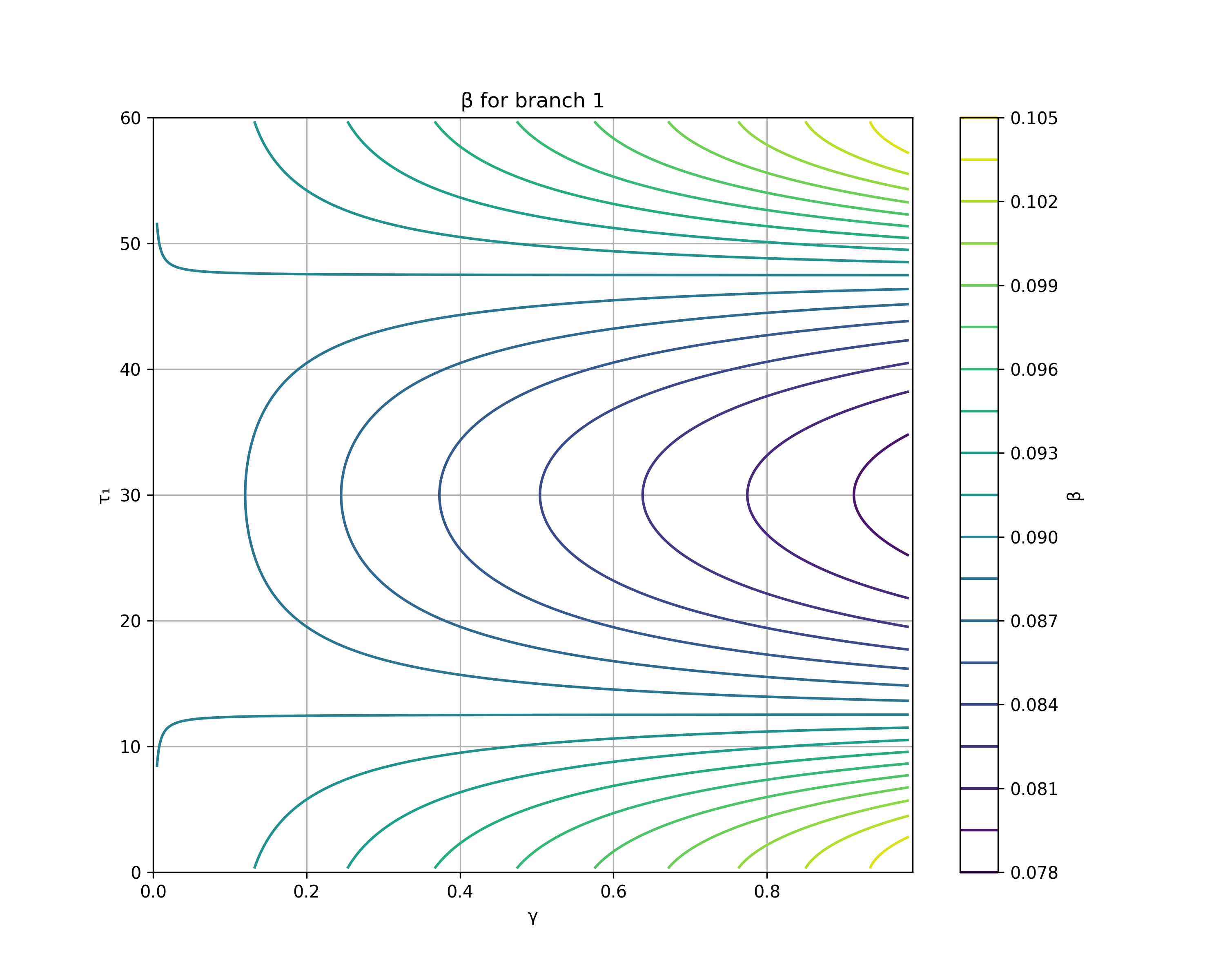}
    \end{minipage}\hfill
    \begin{minipage}{0.43\textwidth}
        \centering
        \includegraphics[width=\linewidth]{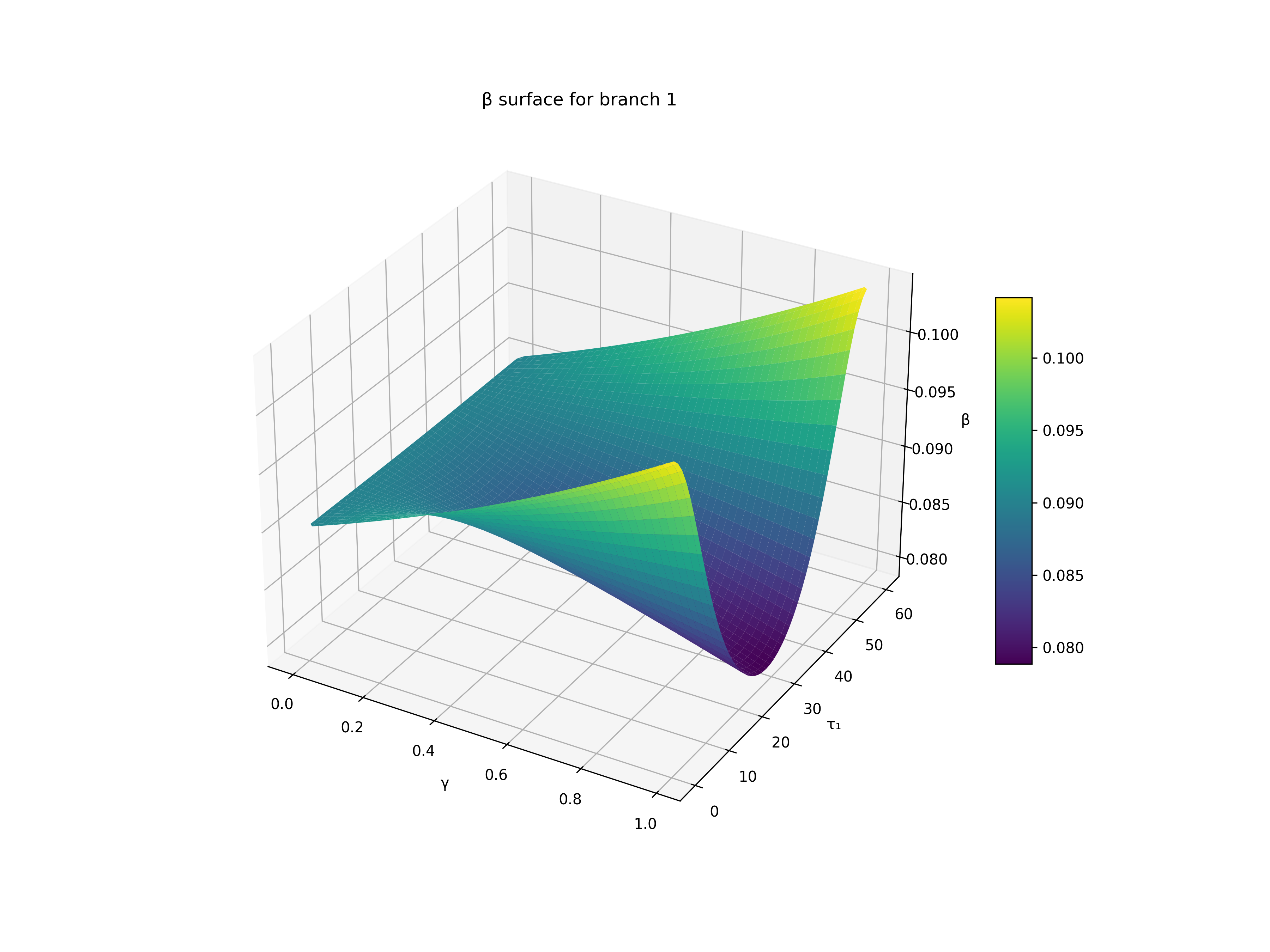}
    \end{minipage}
    \caption{Contour plot and surface for $\beta_{1,1}(\gamma,\tau_1)$ on branch 1. This corresponds to the first solution to \eqref{eq:neutral:beta-relation} for the given values of $\gamma,\tau_1$.}
    \label{fig:neutral:beta-branch-1}
\end{figure}

\begin{figure}[tbp]
    \centering
    \begin{minipage}{0.43\textwidth}
        \centering
        \includegraphics[width=\linewidth]{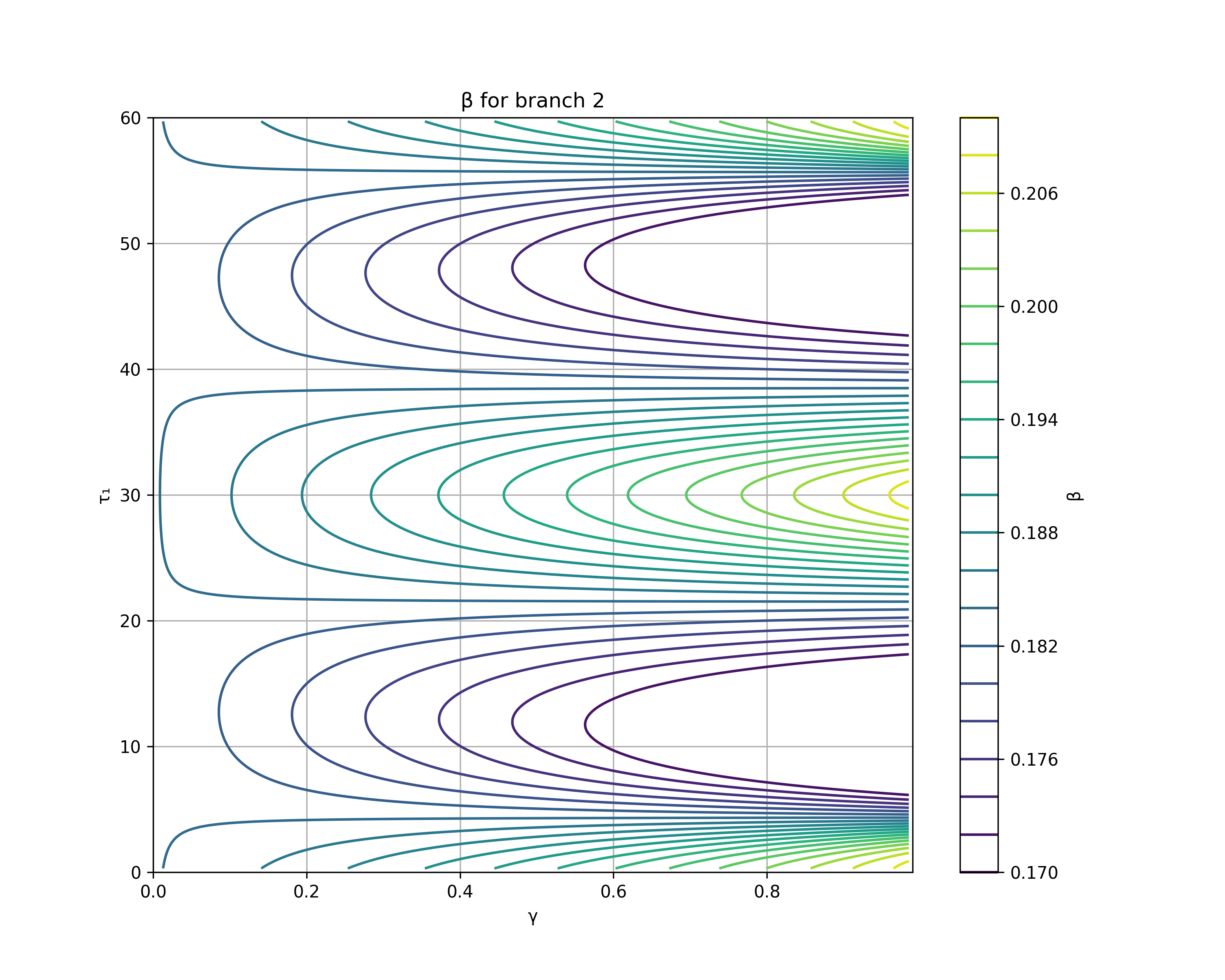}
    \end{minipage}\hfill
    \begin{minipage}{0.43\textwidth}
        \centering
        \includegraphics[width=\linewidth]{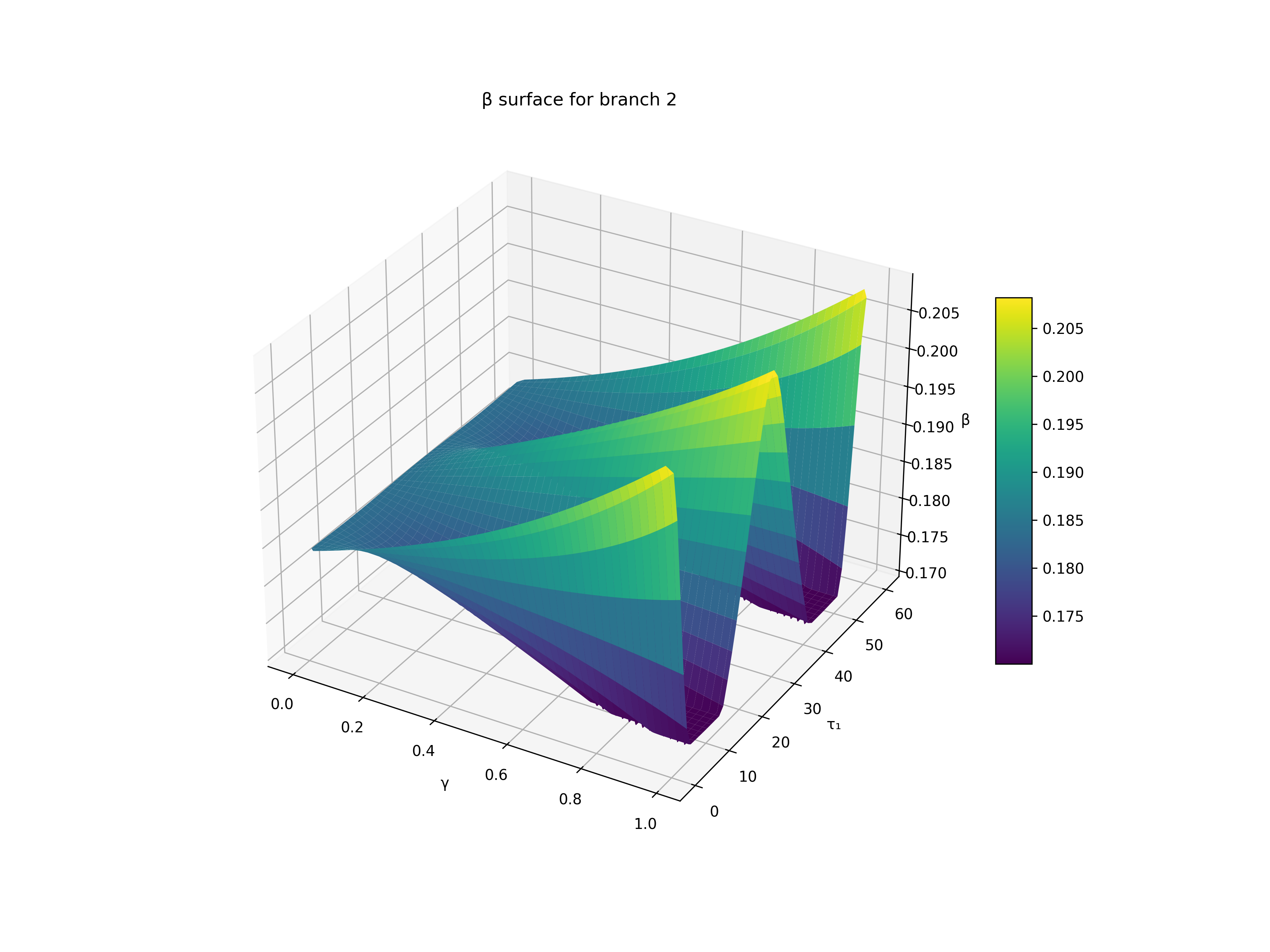}
    \end{minipage}
    \caption{Contour plot and surface for $\beta_{1,1}(\gamma,\tau_1)$ on branch 2. This corresponds to the second solution to \eqref{eq:neutral:beta-relation} for the given values of $\gamma,\tau_1$.}
    \label{fig:neutral:beta-branch-2}
\end{figure}

\begin{figure}[tbp]
    \centering
    \begin{minipage}{0.43\textwidth}
        \centering
        \includegraphics[width=\linewidth]{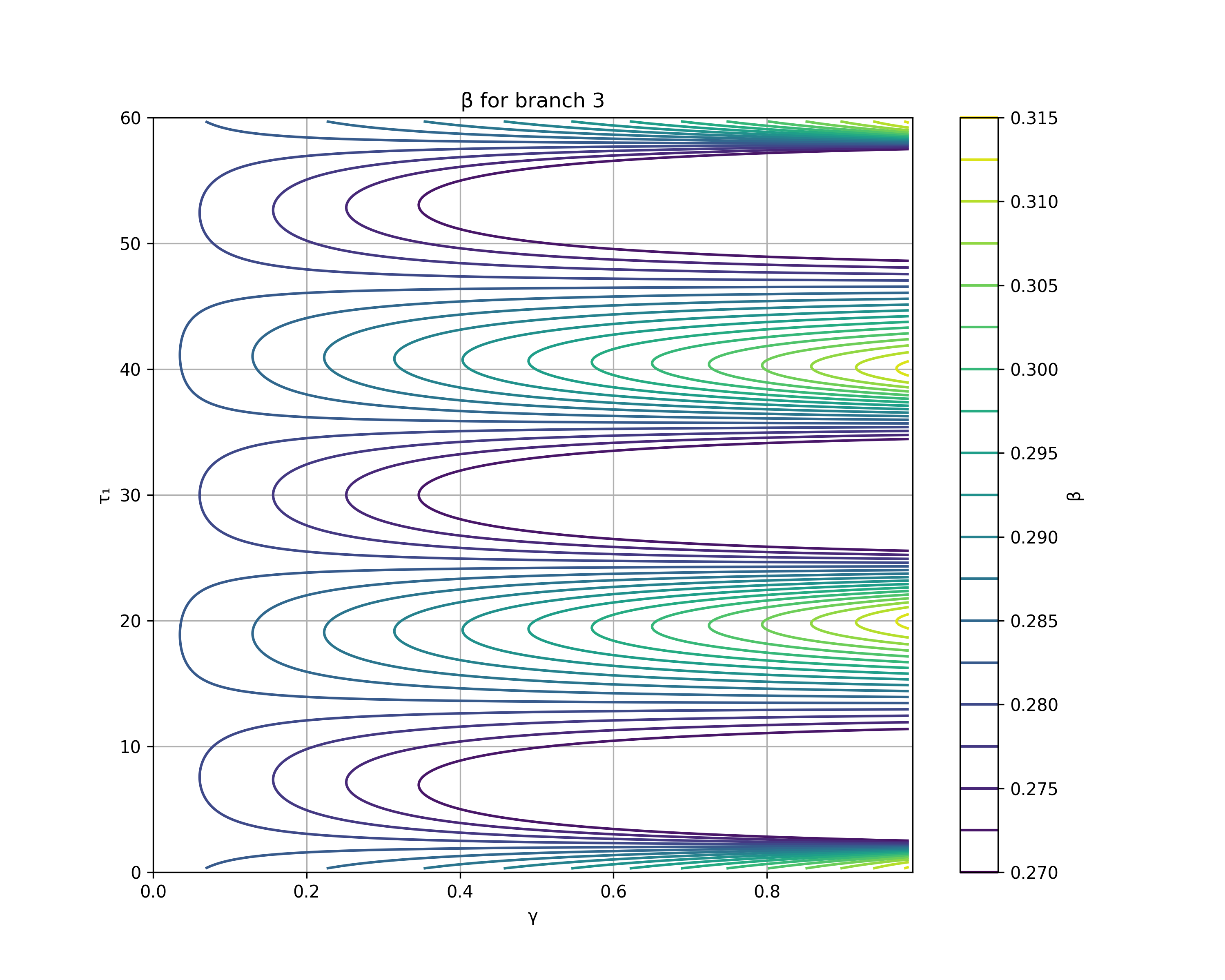}
    \end{minipage}\hfill
    \begin{minipage}{0.43\textwidth}
        \centering
        \includegraphics[width=\linewidth]{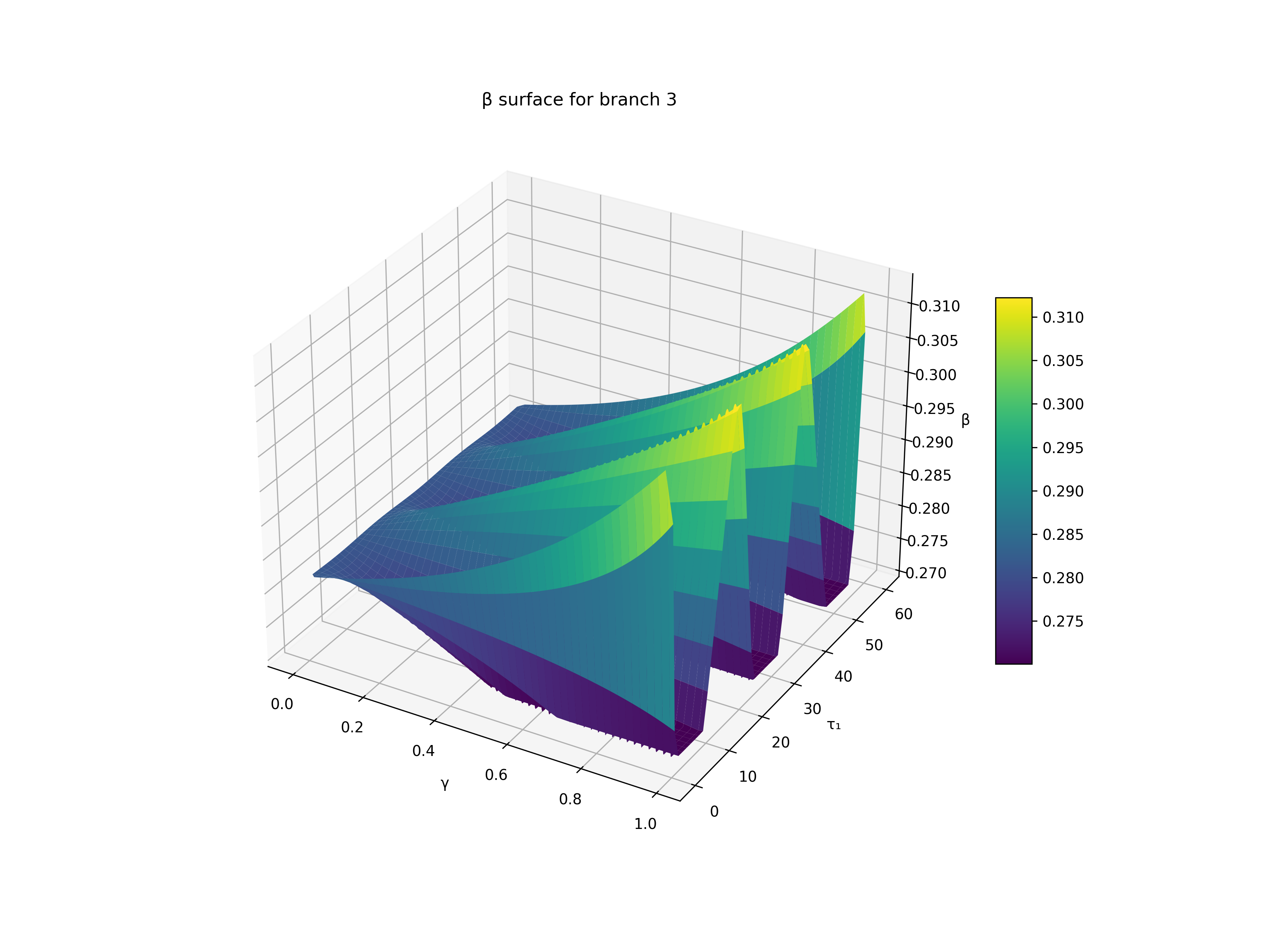}
    \end{minipage}
    \caption{Contour plot and surface for $\beta_{1,1}(\gamma,\tau_1)$ on branch 3. This corresponds to the third solution to \eqref{eq:neutral:beta-relation} for the given values of $\gamma,\tau_1$.}
    \label{fig:neutral:beta-branch-3}
\end{figure}

\begin{figure}[tbp]
    \centering
    \begin{minipage}{0.43\textwidth}
        \centering
        \includegraphics[width=\linewidth]{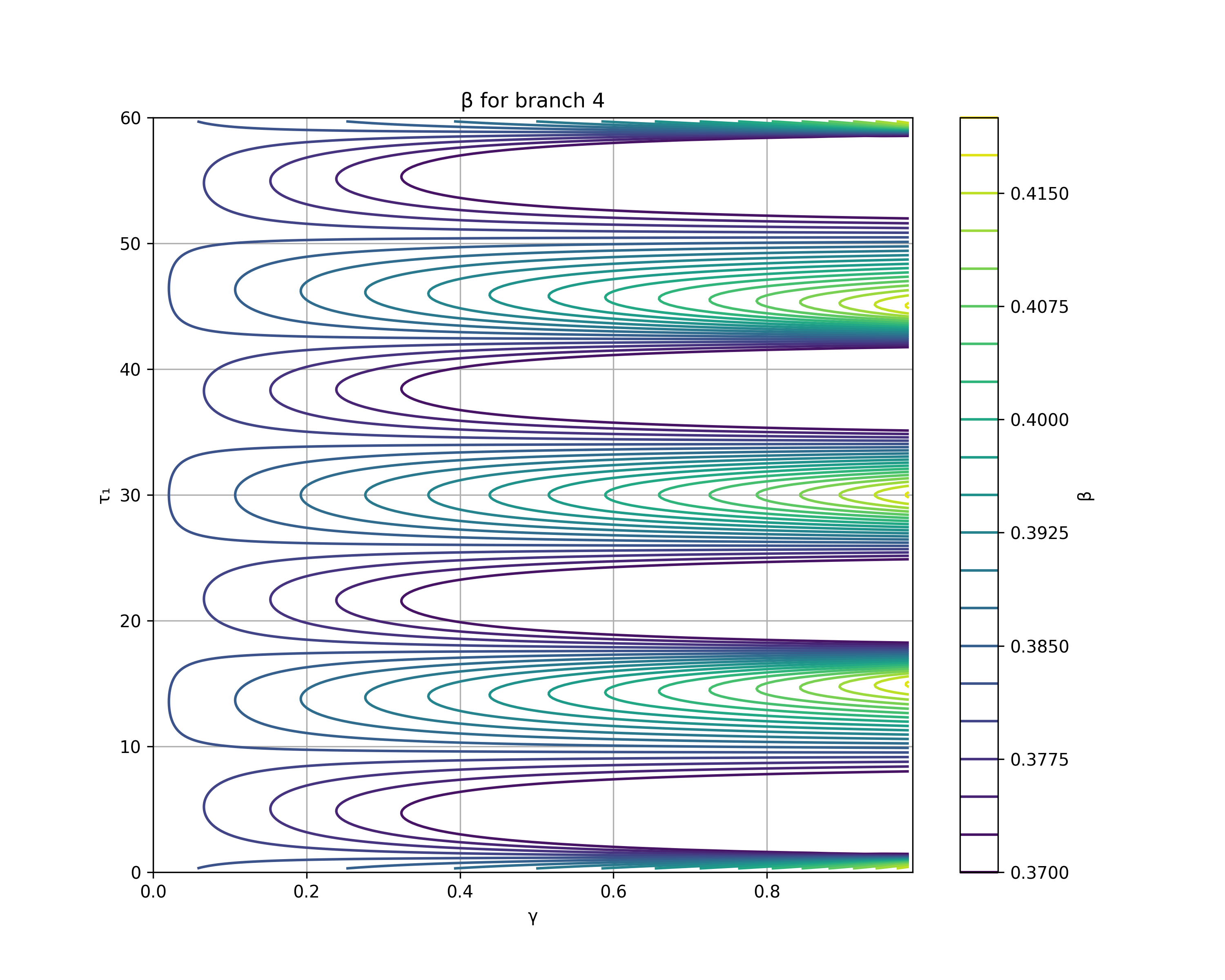}
    \end{minipage}\hfill
    \begin{minipage}{0.43\textwidth}
        \centering
        \includegraphics[width=\linewidth]{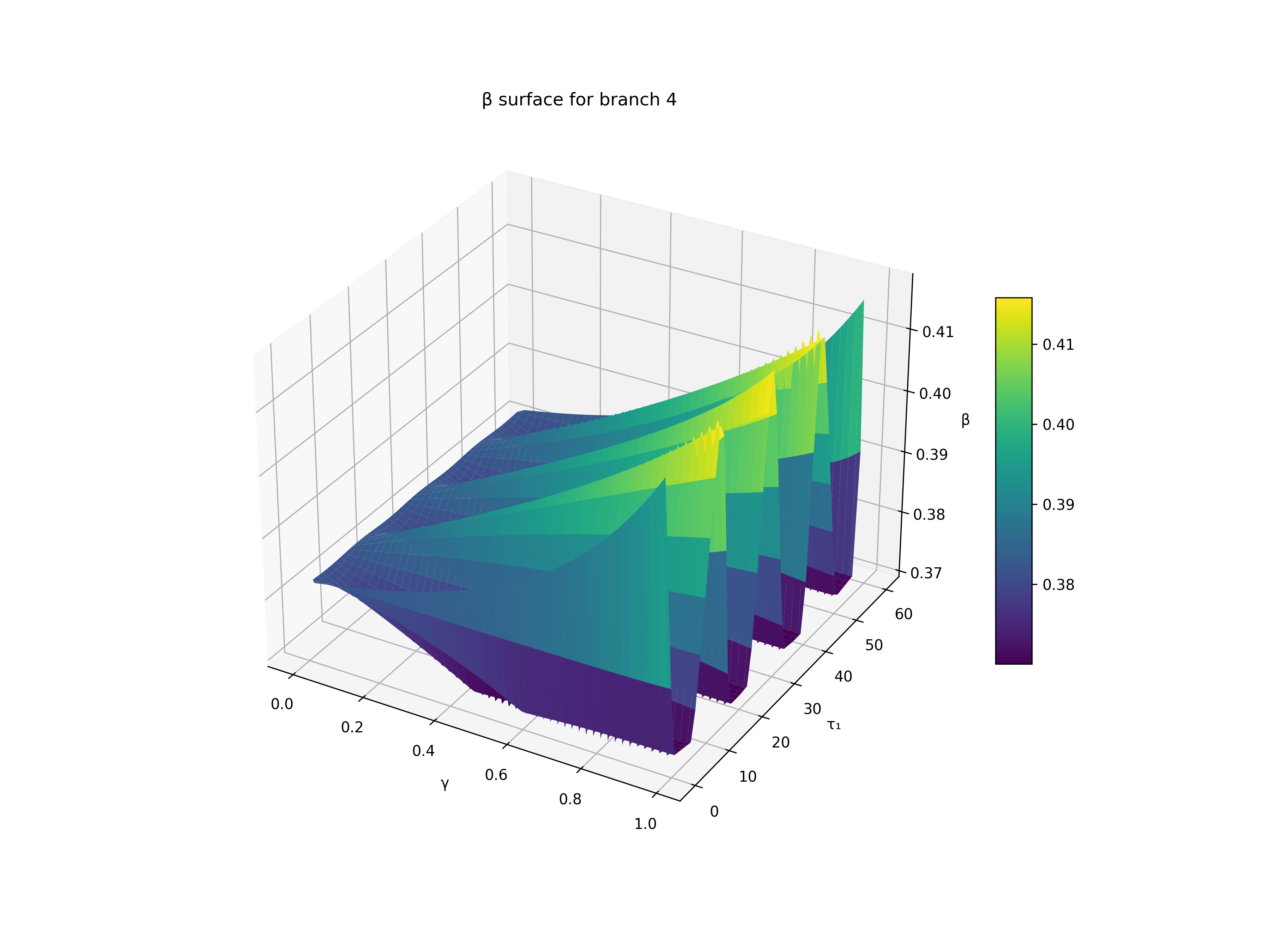}
    \end{minipage}
    \caption{Contour plot and surface for $\beta_{1,1}(\gamma,\tau_1)$ on branch 4. This corresponds to the fourth solution to \eqref{eq:neutral:beta-relation} for the given values of $\gamma,\tau_1$.}
    \label{fig:neutral:beta-branch-4}
\end{figure}

\begin{figure}[tbp]
    \centering
    \begin{minipage}{0.43\textwidth}
        \centering
        \includegraphics[width=\linewidth]{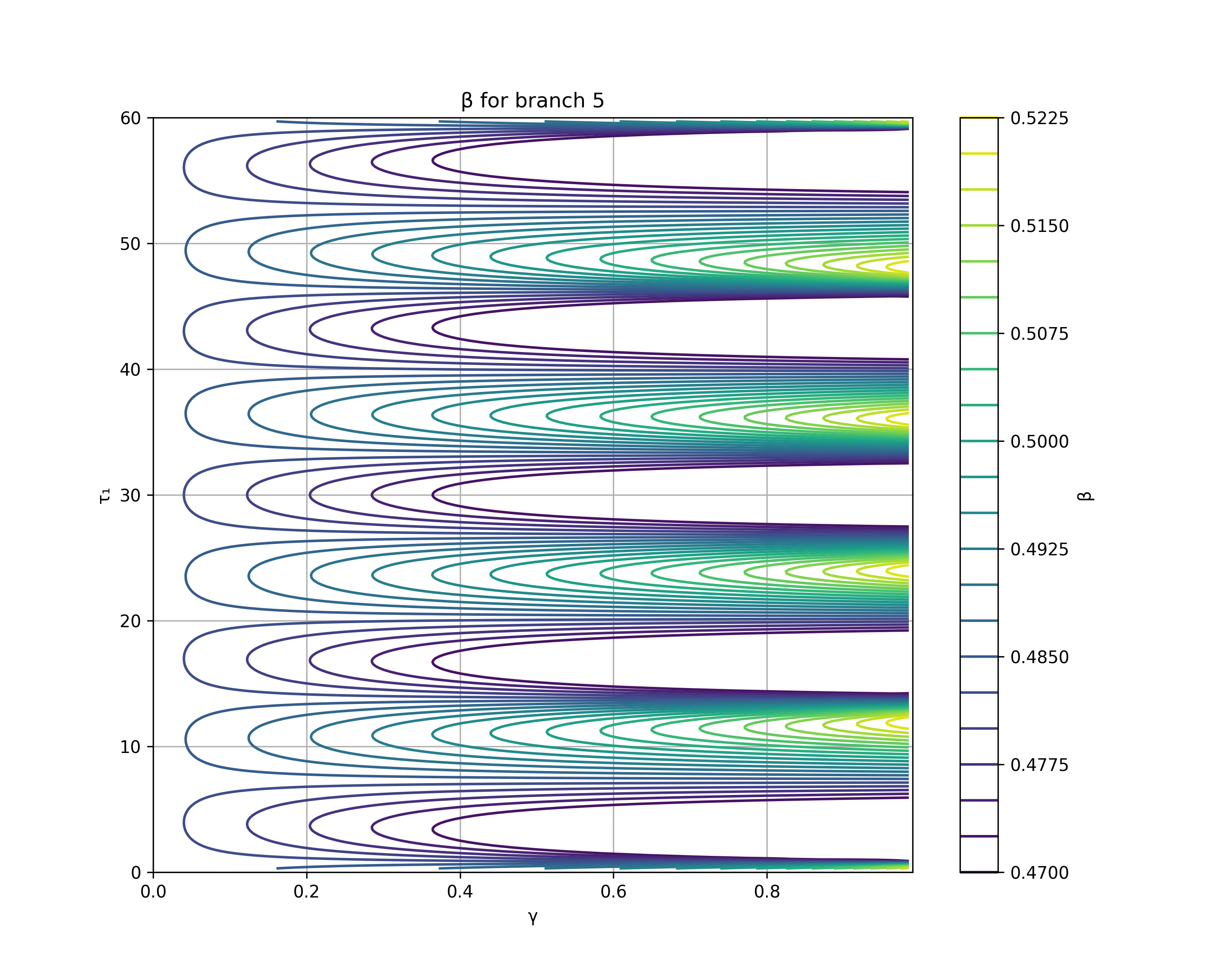}
    \end{minipage}\hfill
    \begin{minipage}{0.43\textwidth}
        \centering
        \includegraphics[width=\linewidth]{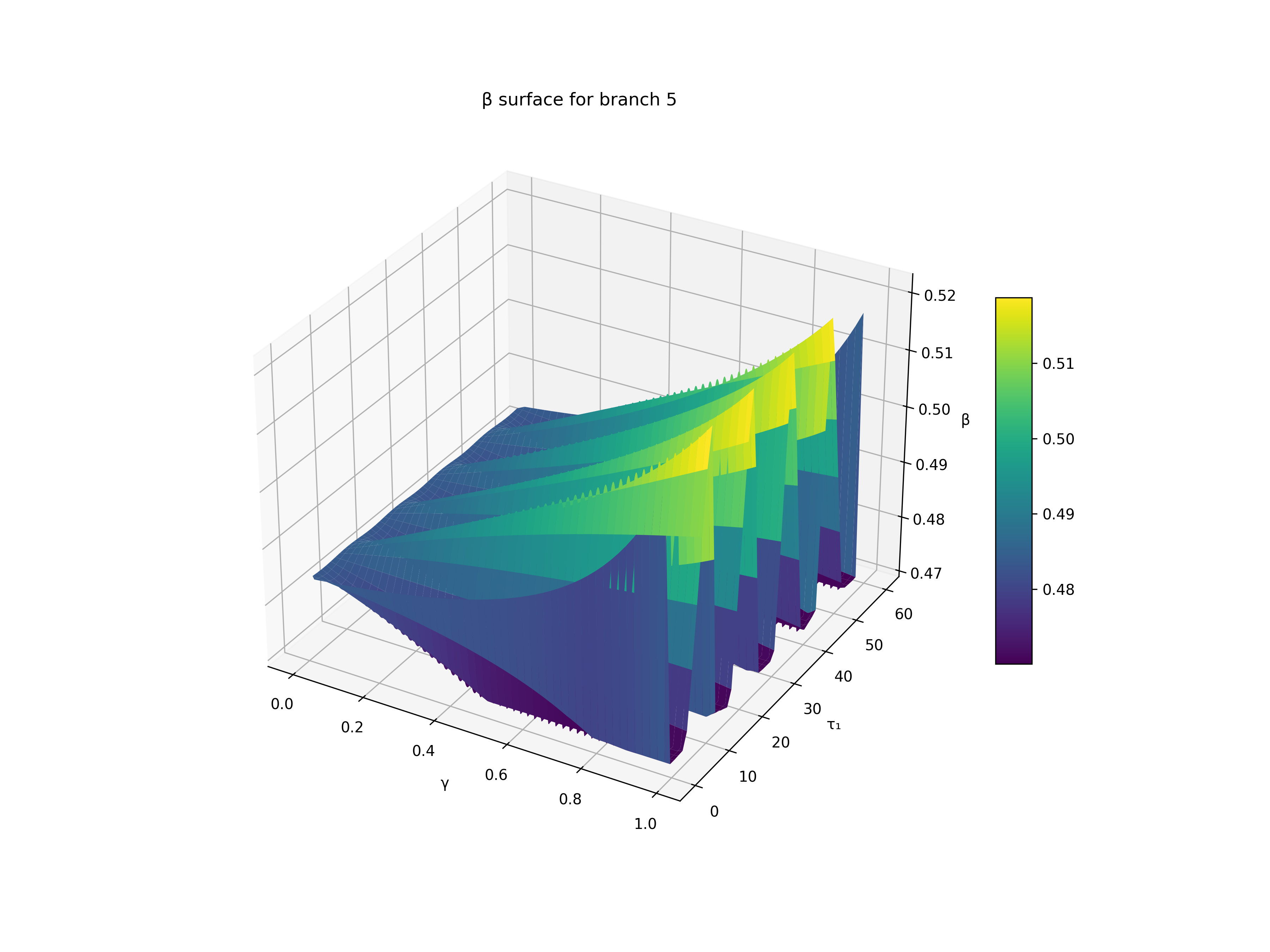}
    \end{minipage}
    \caption{Contour plot and surface for $\beta_{1,1}(\gamma,\tau_1)$ on branch 5. This corresponds to the fifth solution to \eqref{eq:neutral:beta-relation} for the given values of $\gamma,\tau_1$.}
    \label{fig:neutral:beta-branch-5}
\end{figure}

\begin{figure}[tbp]
    \centering
    \begin{minipage}{0.43\textwidth}
        \centering
        \includegraphics[width=\linewidth]{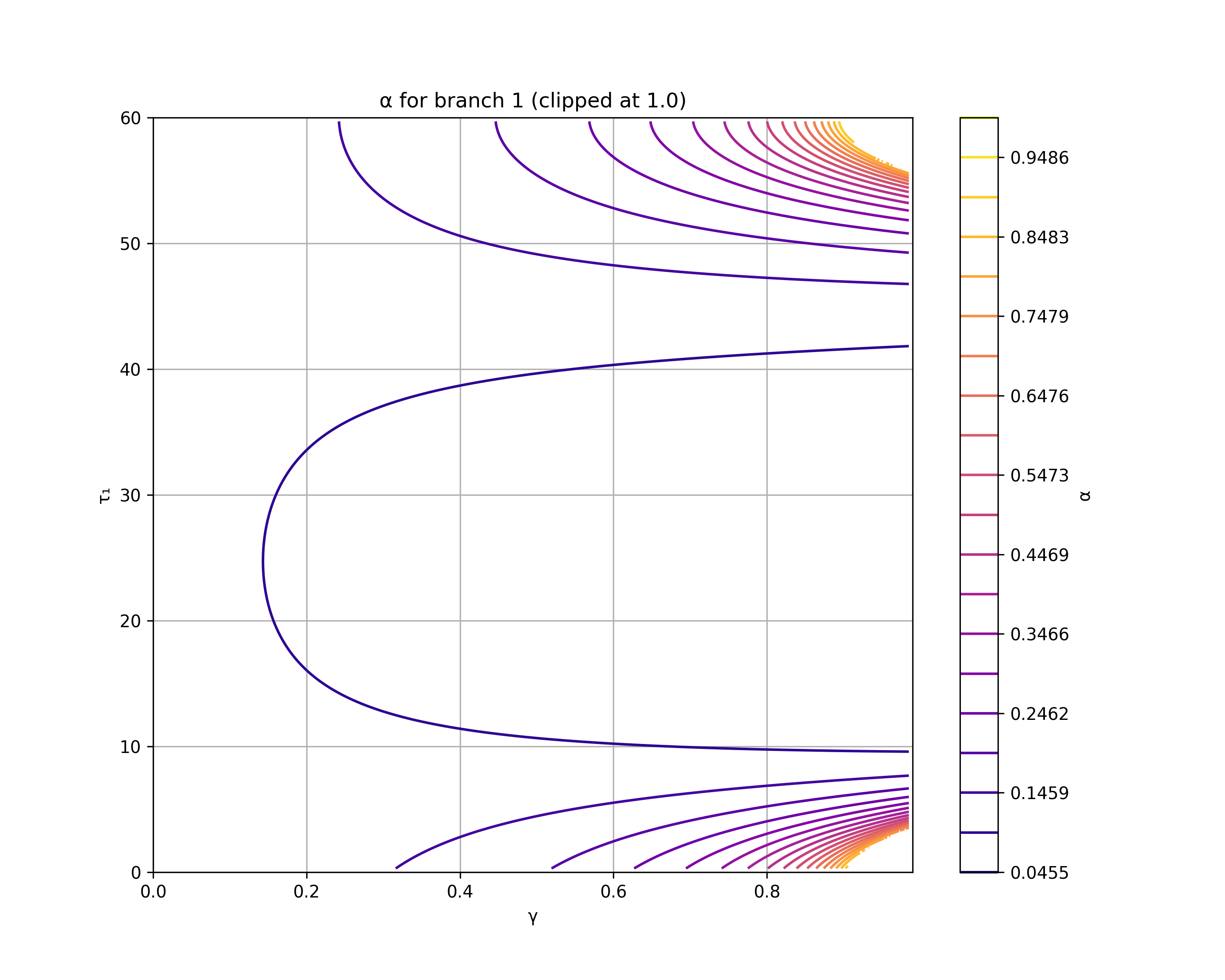}
    \end{minipage}\hfill
    \begin{minipage}{0.43\textwidth}
        \centering
        \includegraphics[width=\linewidth]{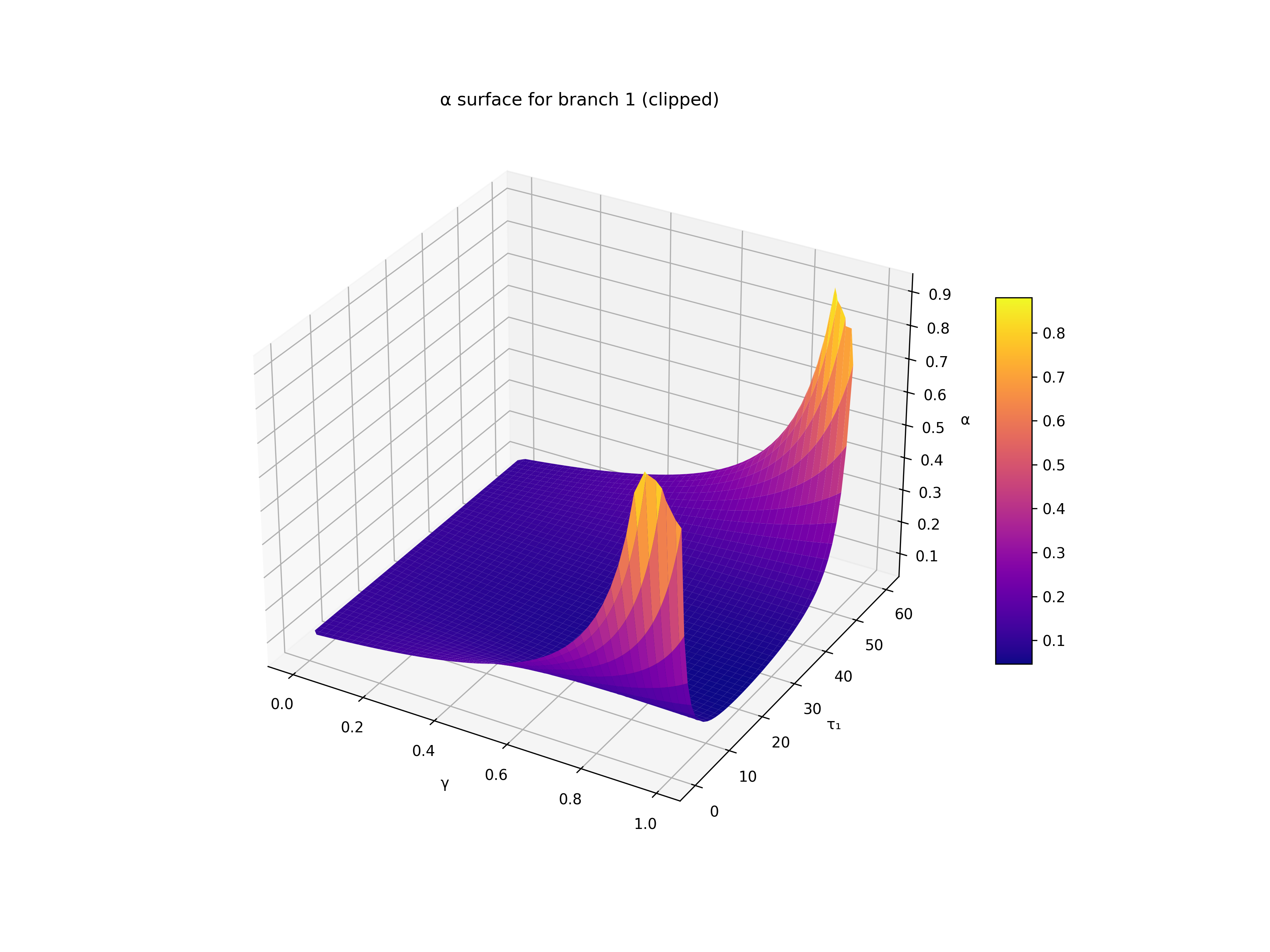}
    \end{minipage}
    \caption{Contour plot and critical surface for $\alpha_{1,1}(\gamma,\tau_1)$ on branch 1 (i.e. corresponding to $\beta_{1,1}(\gamma,\tau_1)$. Values clipped to $\alpha<1$.}
    \label{fig:neutral:alpha-branch-1}
\end{figure}

\begin{figure}[tbp]
    \centering
    \begin{minipage}{0.43\textwidth}
        \centering
        \includegraphics[width=\linewidth]{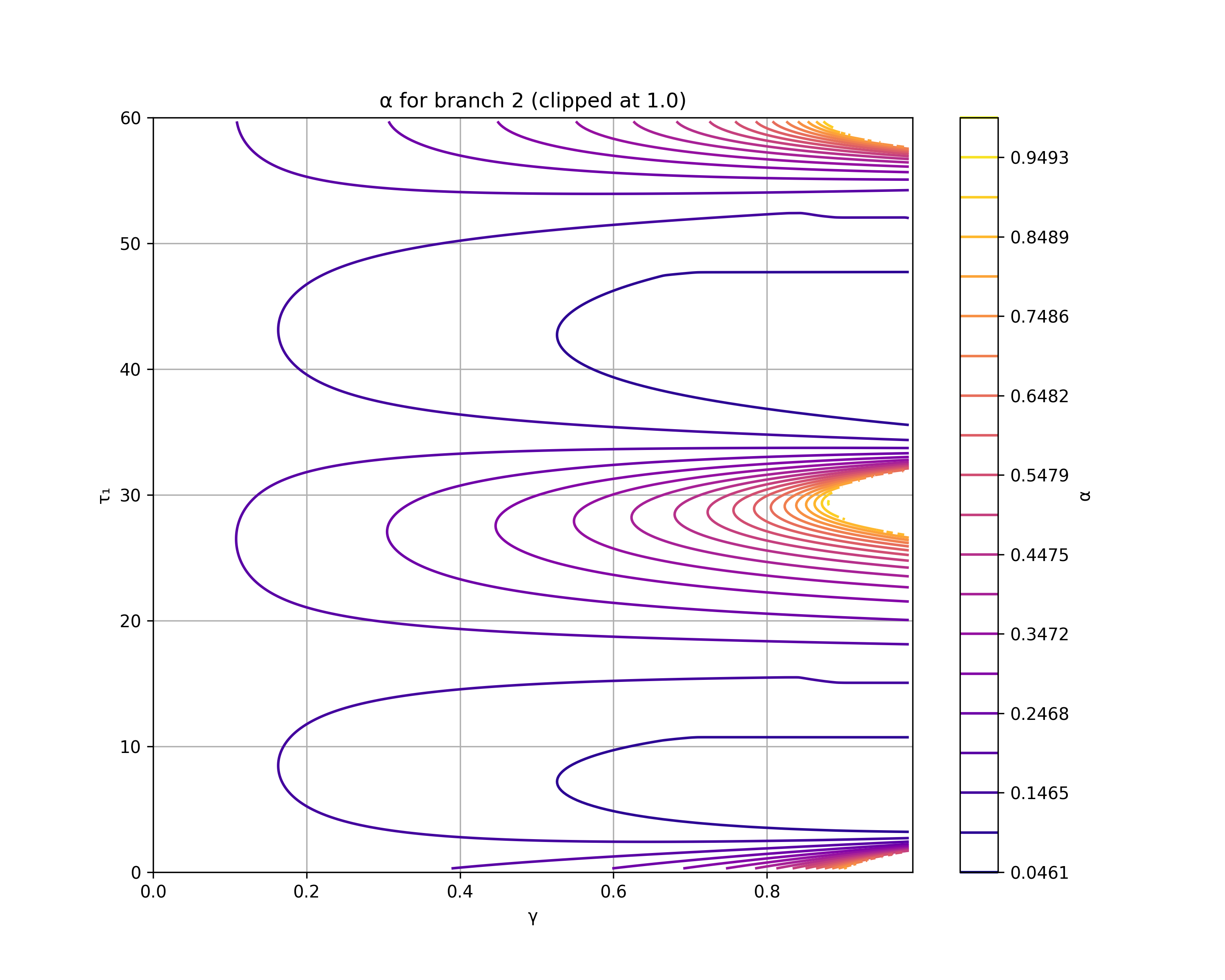}
    \end{minipage}\hfill
    \begin{minipage}{0.43\textwidth}
        \centering
        \includegraphics[width=\linewidth]{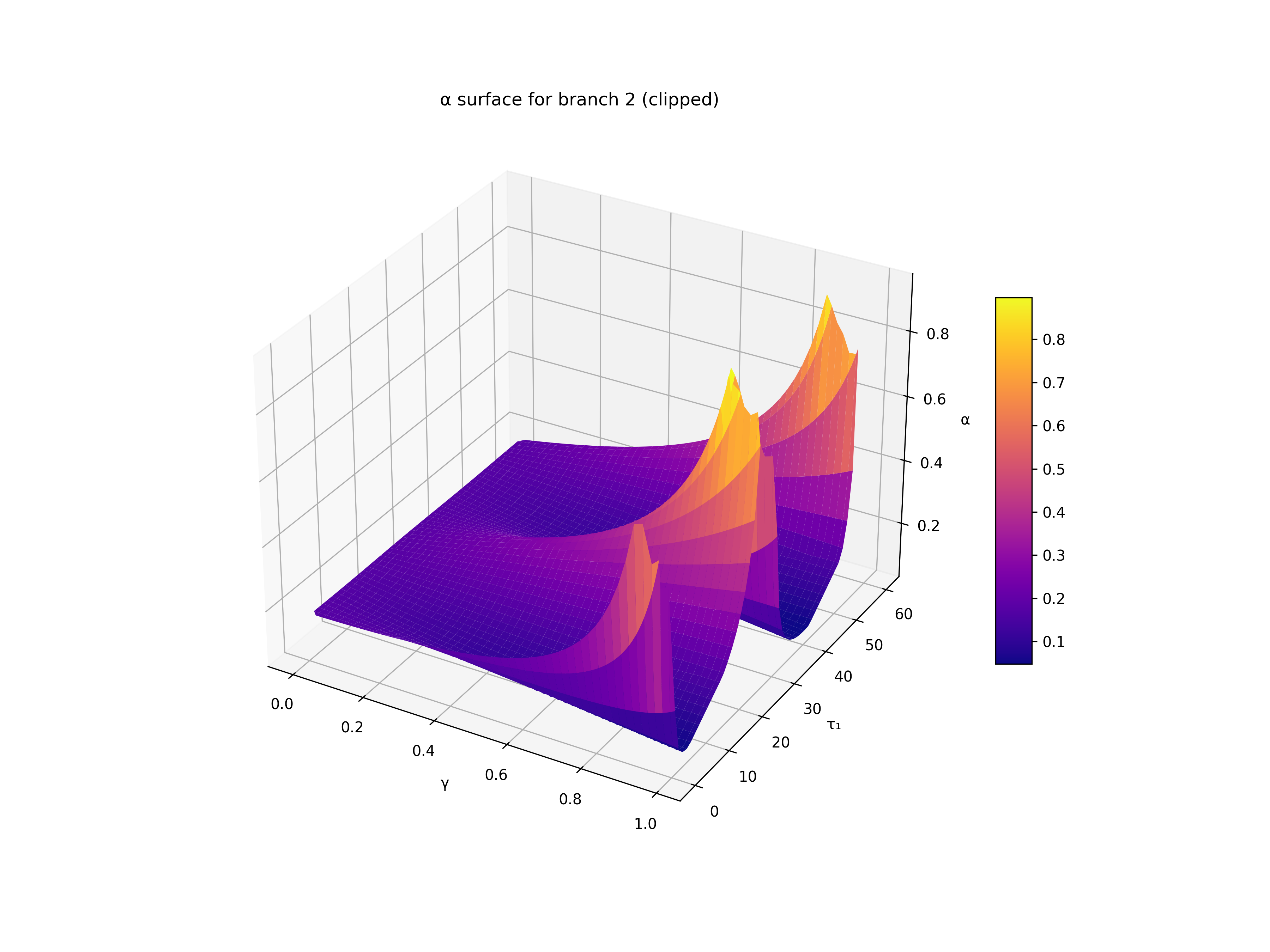}
    \end{minipage}
    \caption{Contour plot and critical surface for $\alpha_{2,1}(\gamma,\tau_1)$ on branch 2 (i.e. corresponding to $\beta_{2,1}(\gamma,\tau_1)$. Values clipped to $\alpha<1$.}
    \label{fig:neutral:alpha-branch-2}
\end{figure}

\begin{figure}[tbp]
    \centering
    \begin{minipage}{0.43\textwidth}
        \centering
        \includegraphics[width=\linewidth]{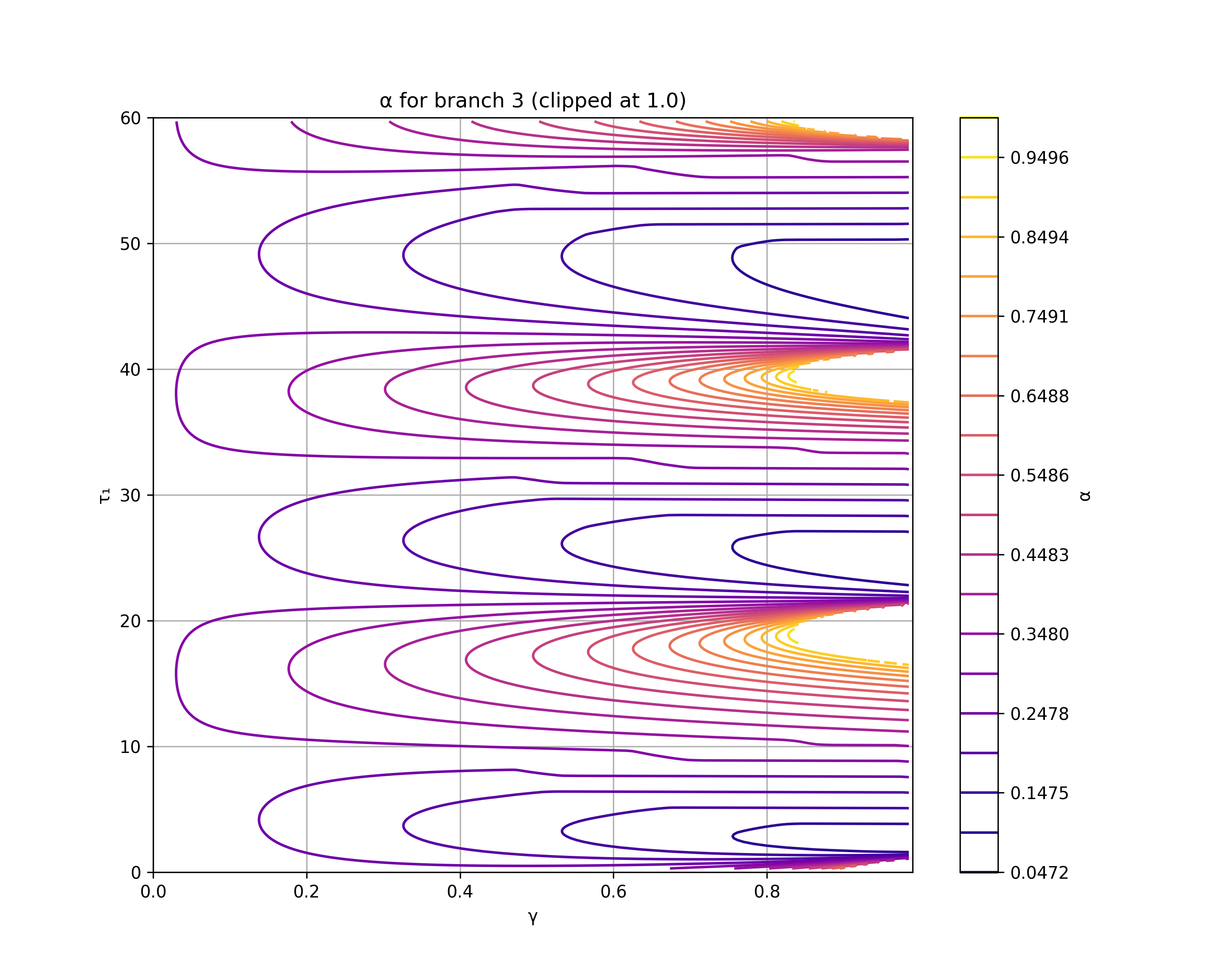}
    \end{minipage}\hfill
    \begin{minipage}{0.43\textwidth}
        \centering
        \includegraphics[width=\linewidth]{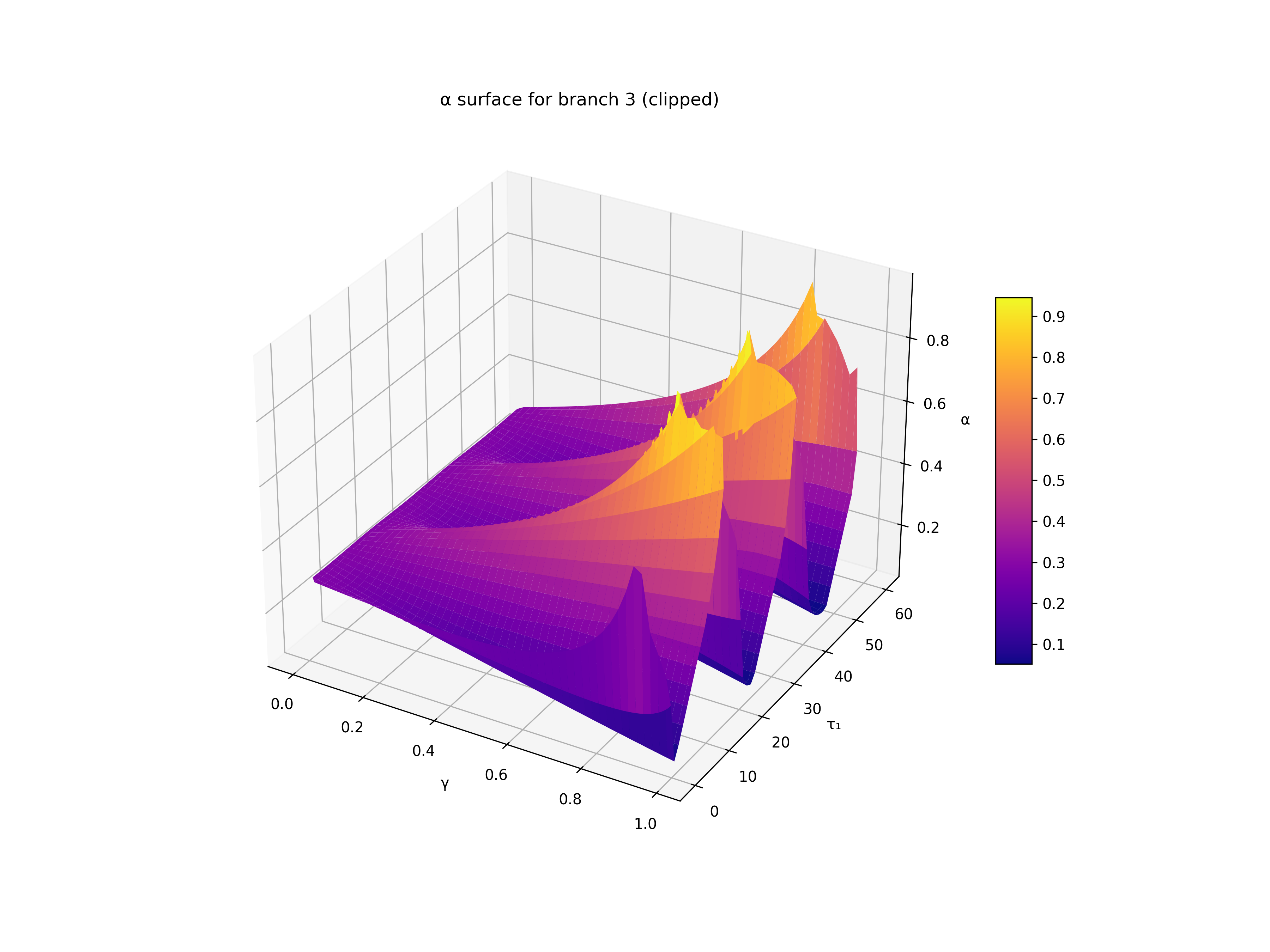}
    \end{minipage}
    \caption{Contour plot and critical surface for $\alpha_{3,1}(\gamma,\tau_1)$ on branch 3 (i.e. corresponding to $\beta_{3,1}(\gamma,\tau_1)$. Values clipped to $\alpha<1$.}
    \label{fig:neutral:alpha-branch-3}
\end{figure}

\begin{figure}[tbp]
    \centering
    \begin{minipage}{0.43\textwidth}
        \centering
        \includegraphics[width=\linewidth]{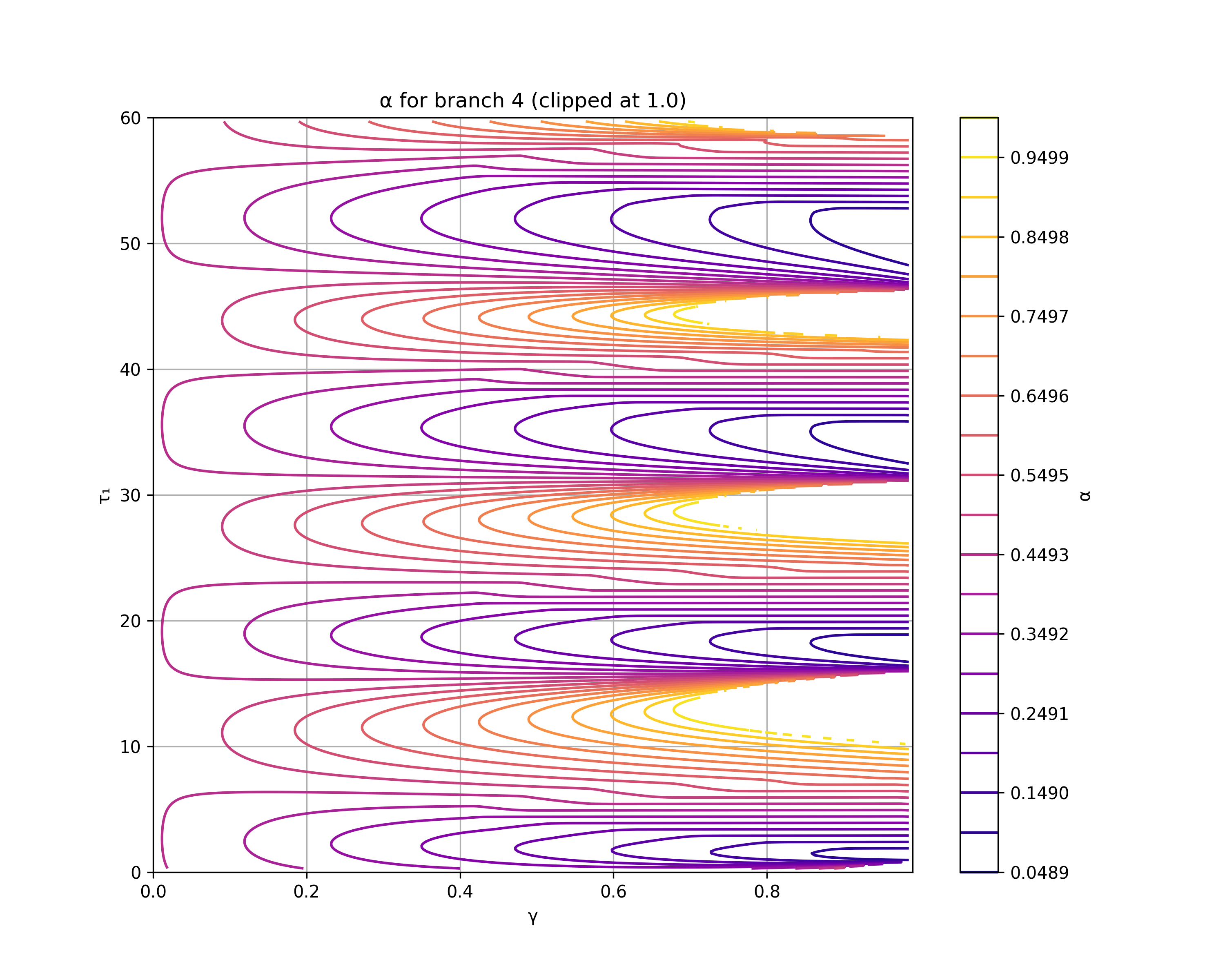}
    \end{minipage}\hfill
    \begin{minipage}{0.43\textwidth}
        \centering
        \includegraphics[width=\linewidth]{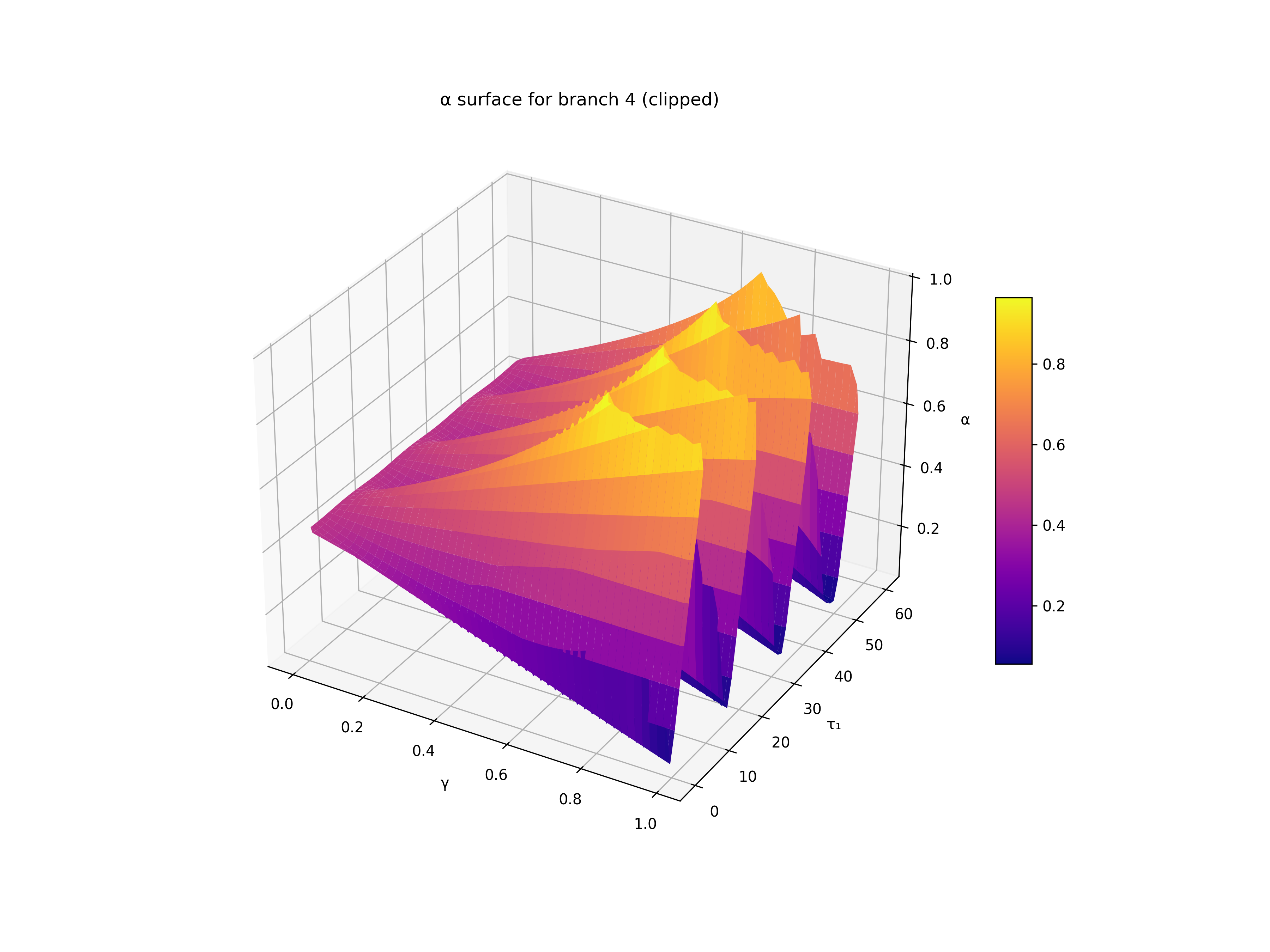}
    \end{minipage}
    \caption{Contour plot and critical surface for $\alpha_{4,1}(\gamma,\tau_1)$ on branch 4 (i.e. corresponding to $\beta_{4,1}(\gamma,\tau_1)$. Values clipped to $\alpha<1$.}
    \label{fig:neutral:alpha-branch-4}
\end{figure}

\begin{figure}[tbp]
    \centering
    \begin{minipage}{0.43\textwidth}
        \centering
        \includegraphics[width=\linewidth]{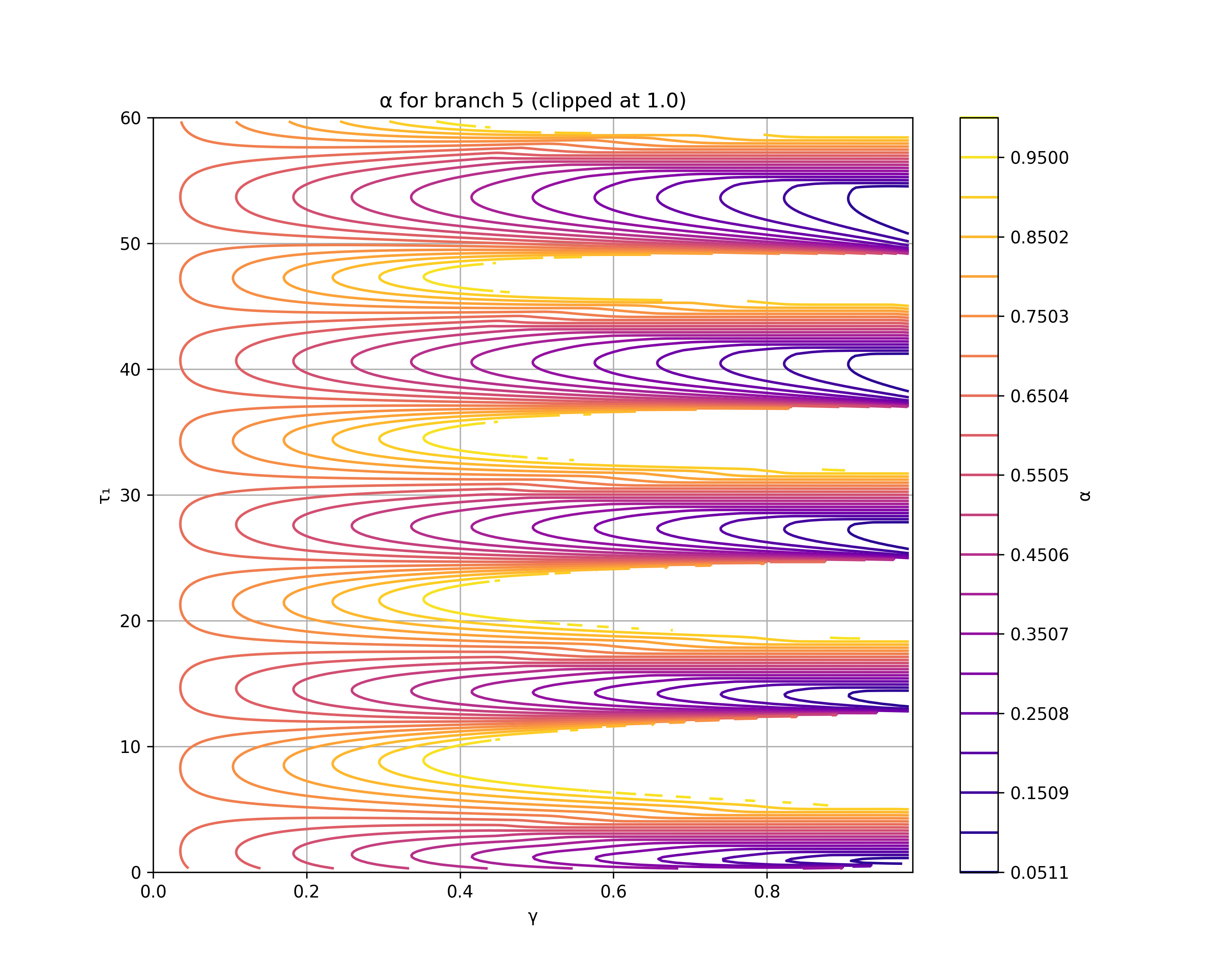}
    \end{minipage}\hfill
    \begin{minipage}{0.43\textwidth}
        \centering
        \includegraphics[width=\linewidth]{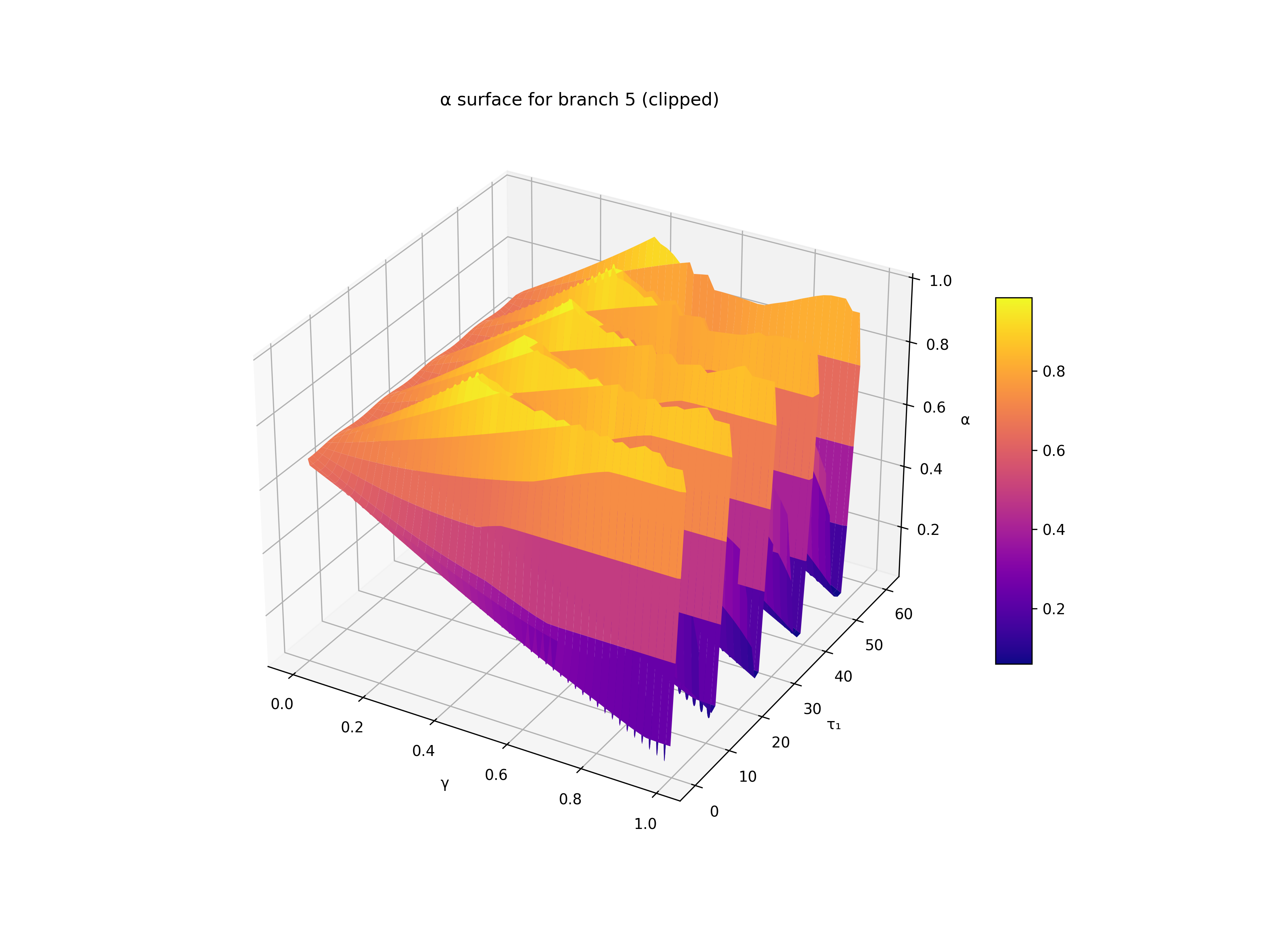}
    \end{minipage}
    \caption{Contour plot and critical surface for $\alpha_{5,1}(\gamma,\tau_1)$ on branch 5 (i.e. corresponding to $\beta_{5,1}(\gamma,\tau_1)$. Values clipped to $\alpha<1$.}
    \label{fig:neutral:alpha-branch-5}
\end{figure}

\begin{figure}[tbp]
    \centering
    \begin{minipage}{0.43\textwidth}
        \centering
        \includegraphics[width=\linewidth]{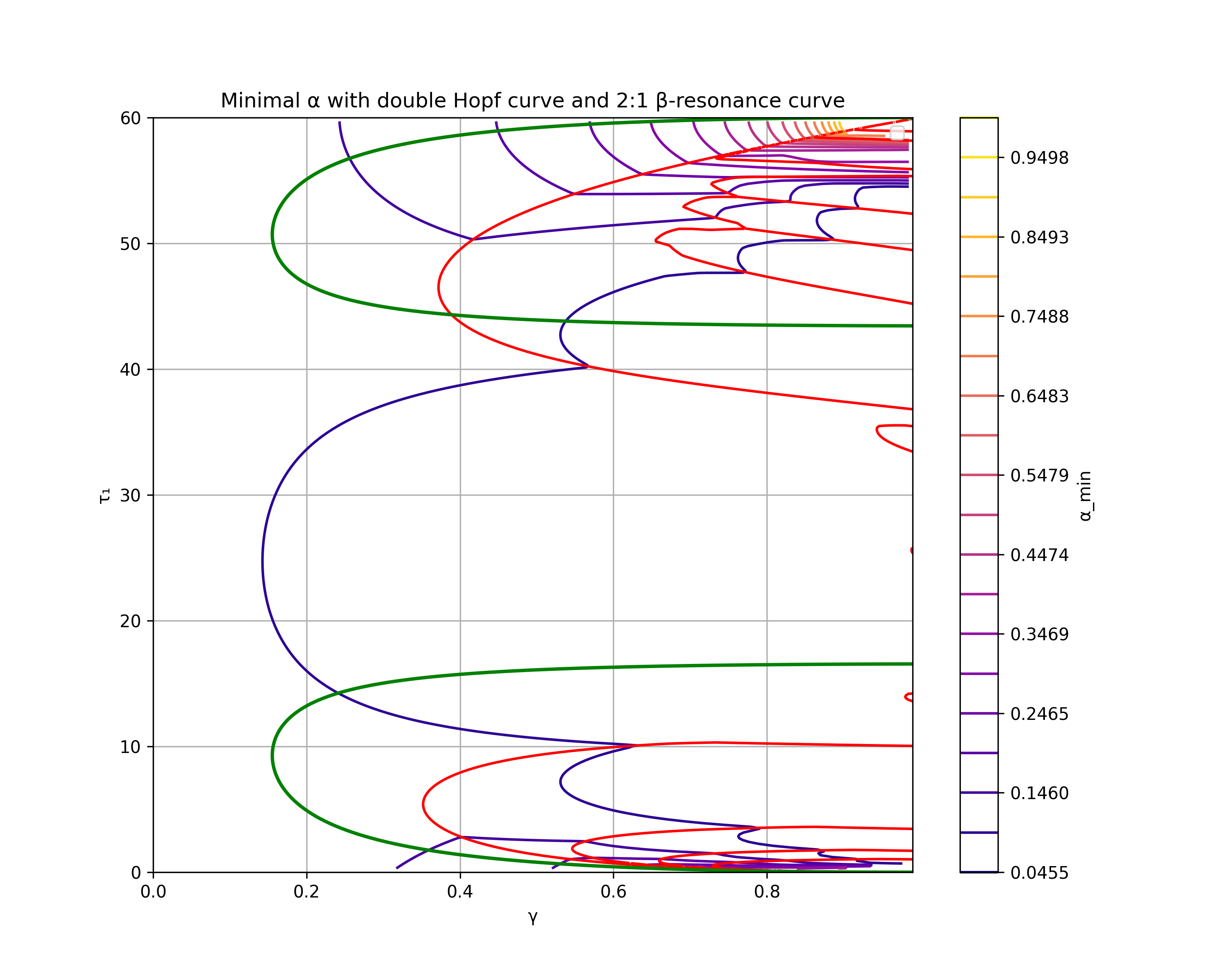}
    \end{minipage}\hfill
    \begin{minipage}{0.43\textwidth}
        \centering
        \includegraphics[width=\linewidth]{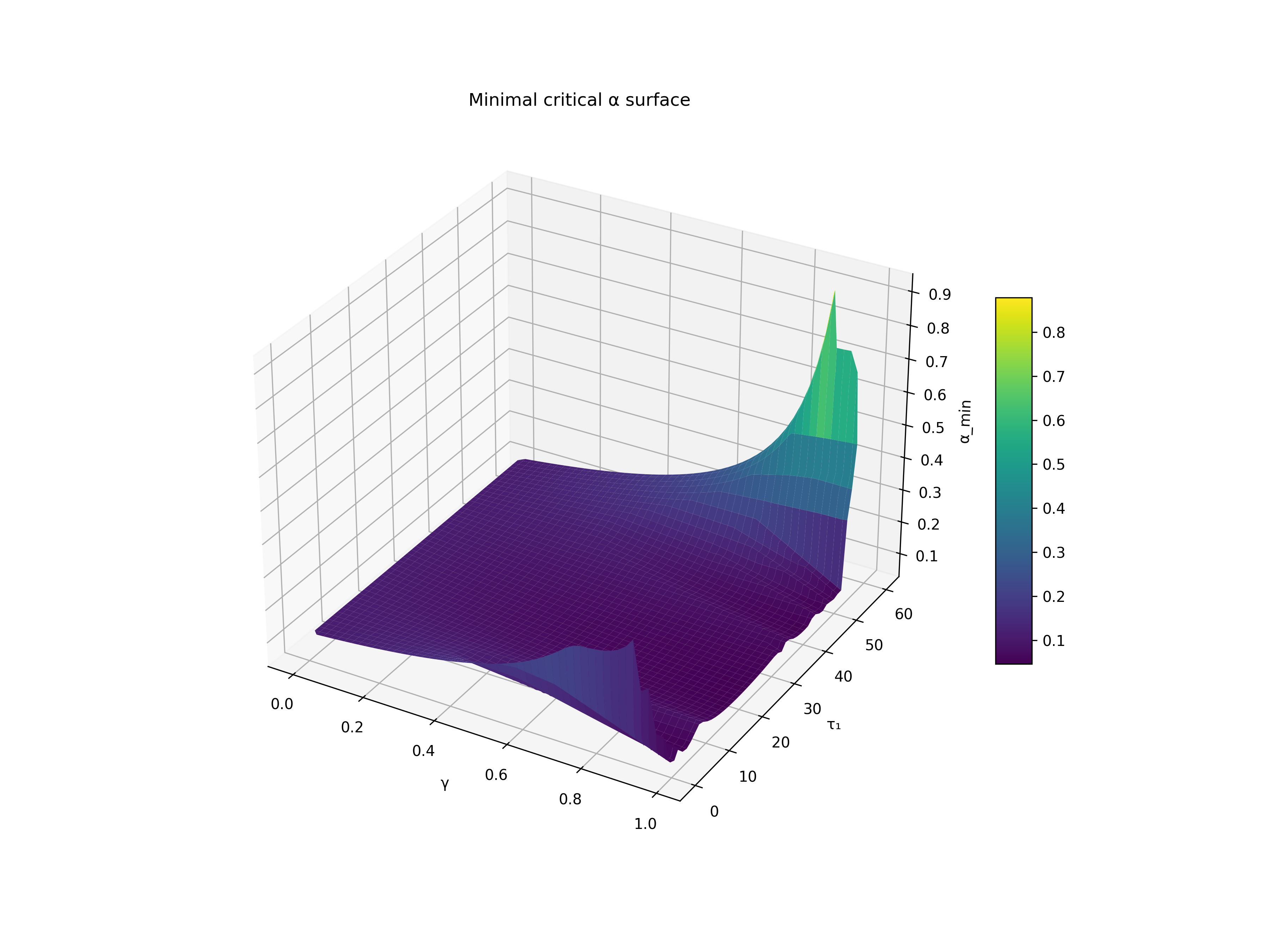}
    \end{minipage}
    \caption{Contour plot and critical surface for the first Hopf bifurcation, taken as the minimum over the previous surfaces. The contour plot has been overlaid with the double Hopf curves.}
    \label{fig:neutral:alpha-min-with-resonance}
\end{figure}

\FloatBarrier

\section{Sample plots of chaotic and quasiperiodic multiconsensus solutions}
On the following pages we provide some additional time series taken at various values of $\alpha$ near the double Hopf point at $\gamma =0.7760643397, \tau_1=0.1678667007, \alpha_{1,1}=\alpha_{2,1}=0.4029160833, \beta_{1,1}=0.1007507146,\beta_{2,1}=0.2014401509$. All other parameter values are taken to be the same as in Chapter \ref{chapter:neutral}. The purpose of these plots is to give some examples of the types of quasiperiodic solutions which exist near the bifurcation point, and give a (imprecise) view of the boundary of the chaotic attractor as $\alpha$ moves past the critical value at $\alpha_{1,1}$. A more detailed exploration of the structure and boundary of the attractor is planned for future work, and could allow for more rigorous conclusions on, for example, the presence of riddled basins.

\begin{figure}[t]
    \centering

    \includegraphics[width=\textwidth]{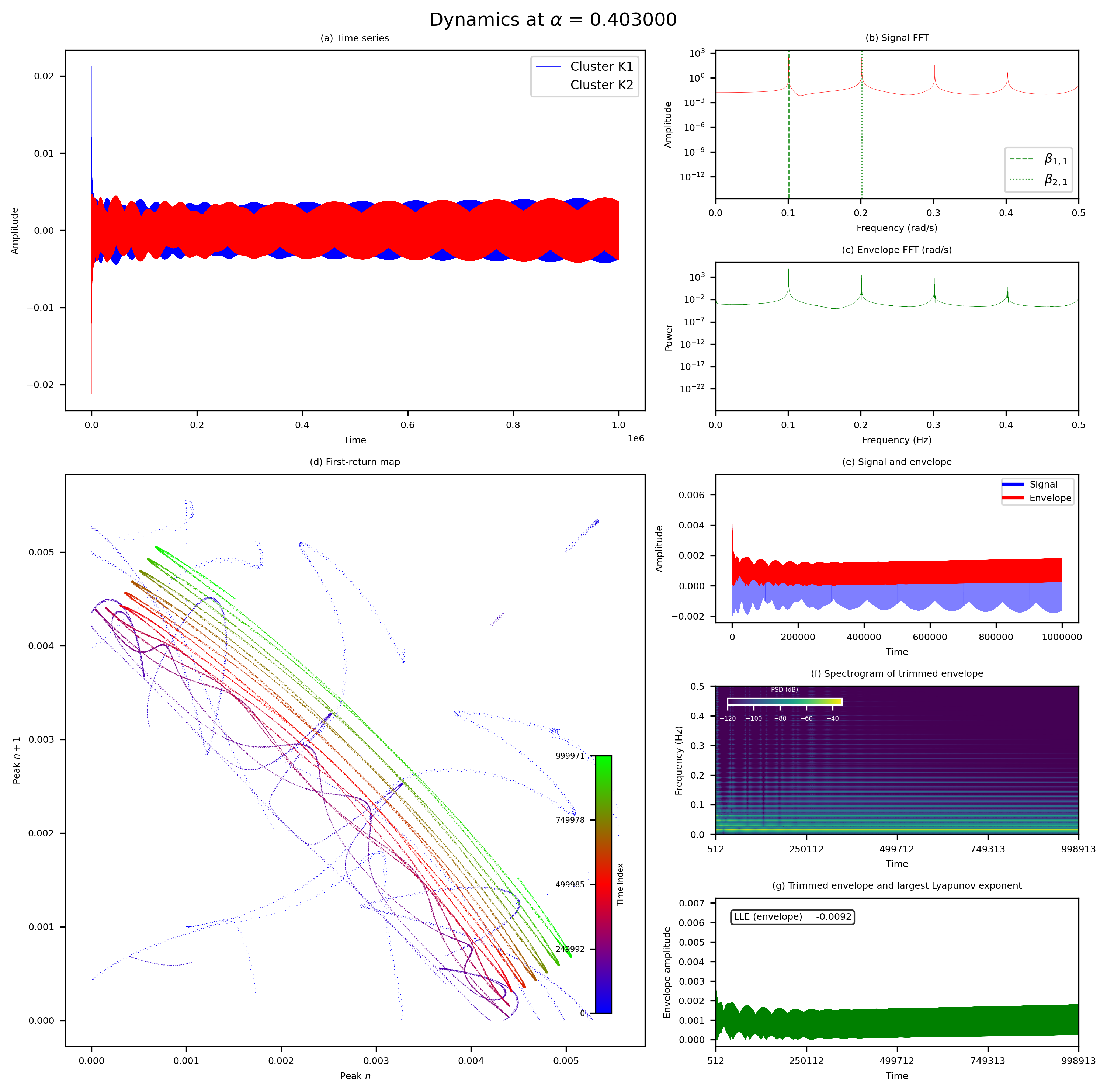}

    \caption{Dynamics at $\alpha = 0.403$, very close to the double Hopf point. We see an initial transient, possibly chaotic but very brief, which then evolves towards a quasiperiodic solution.}
    \label{fig:neutral:trajectory-0.403}
\end{figure}

\begin{figure}[tbp]
    \centering

    \includegraphics[width=\textwidth]{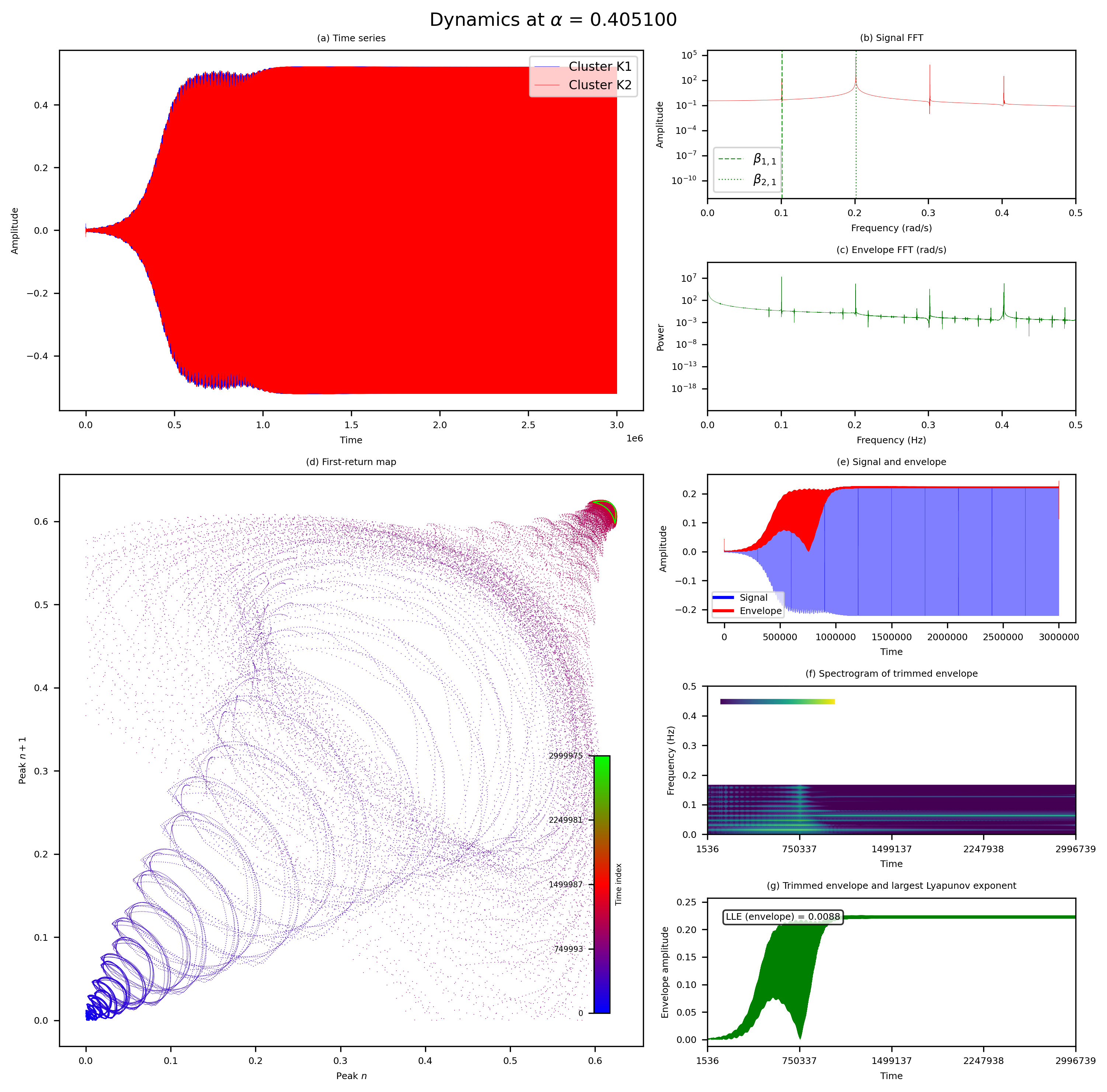}

    \caption{Dynamics at $\alpha = 0.4051$, showing a solution which grows quasiperiodically, hitting a possibly chaotic peak at around $t=7.5\times10^5$ before converging to a quasiperiodic solution, shown as the green closed curve in the first-return plot.}
    \label{fig:neutral:trajectory-0.4051}
\end{figure}

\begin{figure}[tbp]
    \centering

    \includegraphics[width=\textwidth]{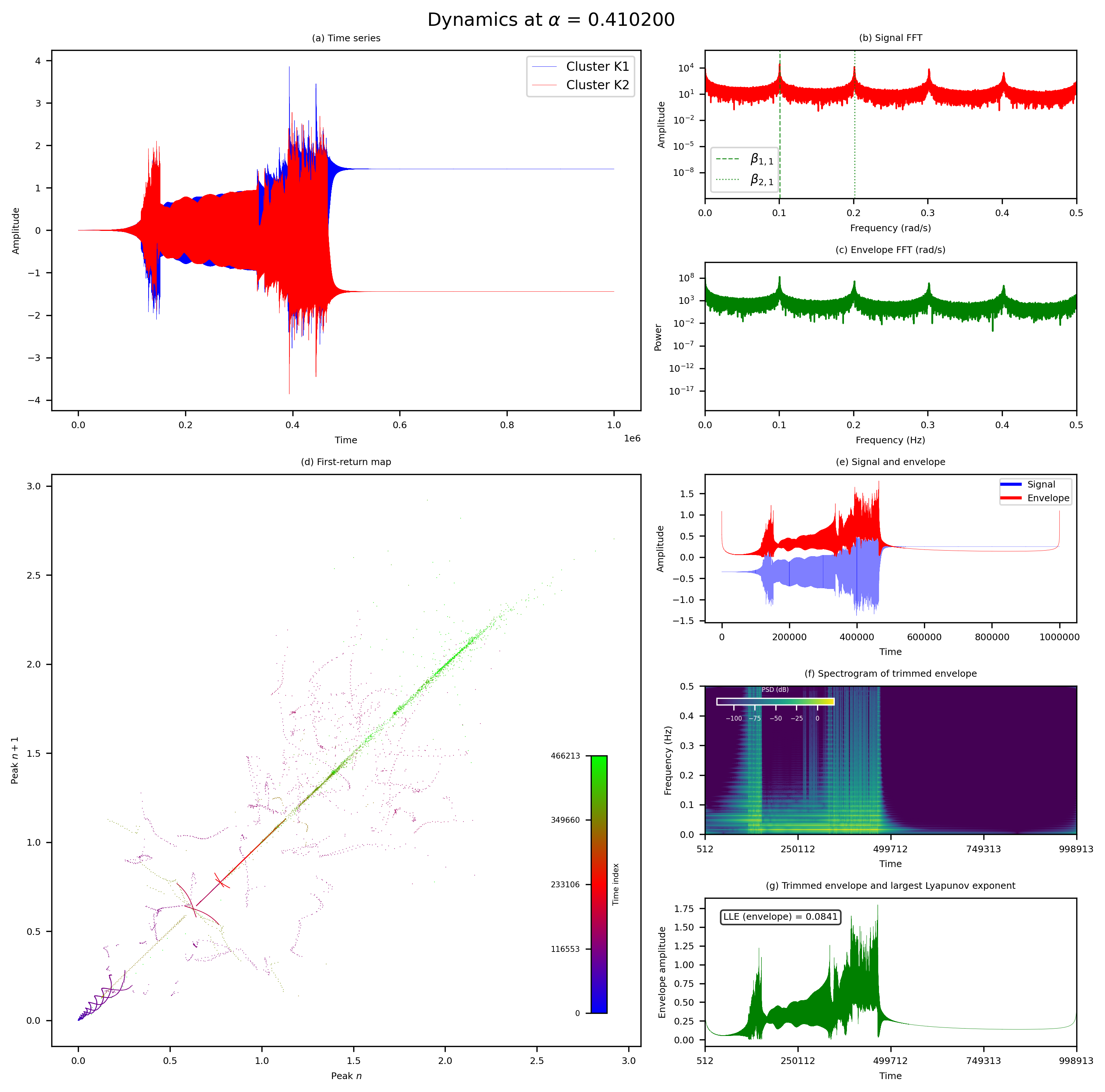}

    \caption{Dynamics at $\alpha = 0.4102$, showing a chaotic transient leading to quasiperiodic dynamics, which return to chaos before leading to steady-state bifurcation.}
    \label{fig:neutral:trajectory-0.4102}
\end{figure}

\begin{figure}[tbp]
    \centering

    \includegraphics[width=\textwidth]{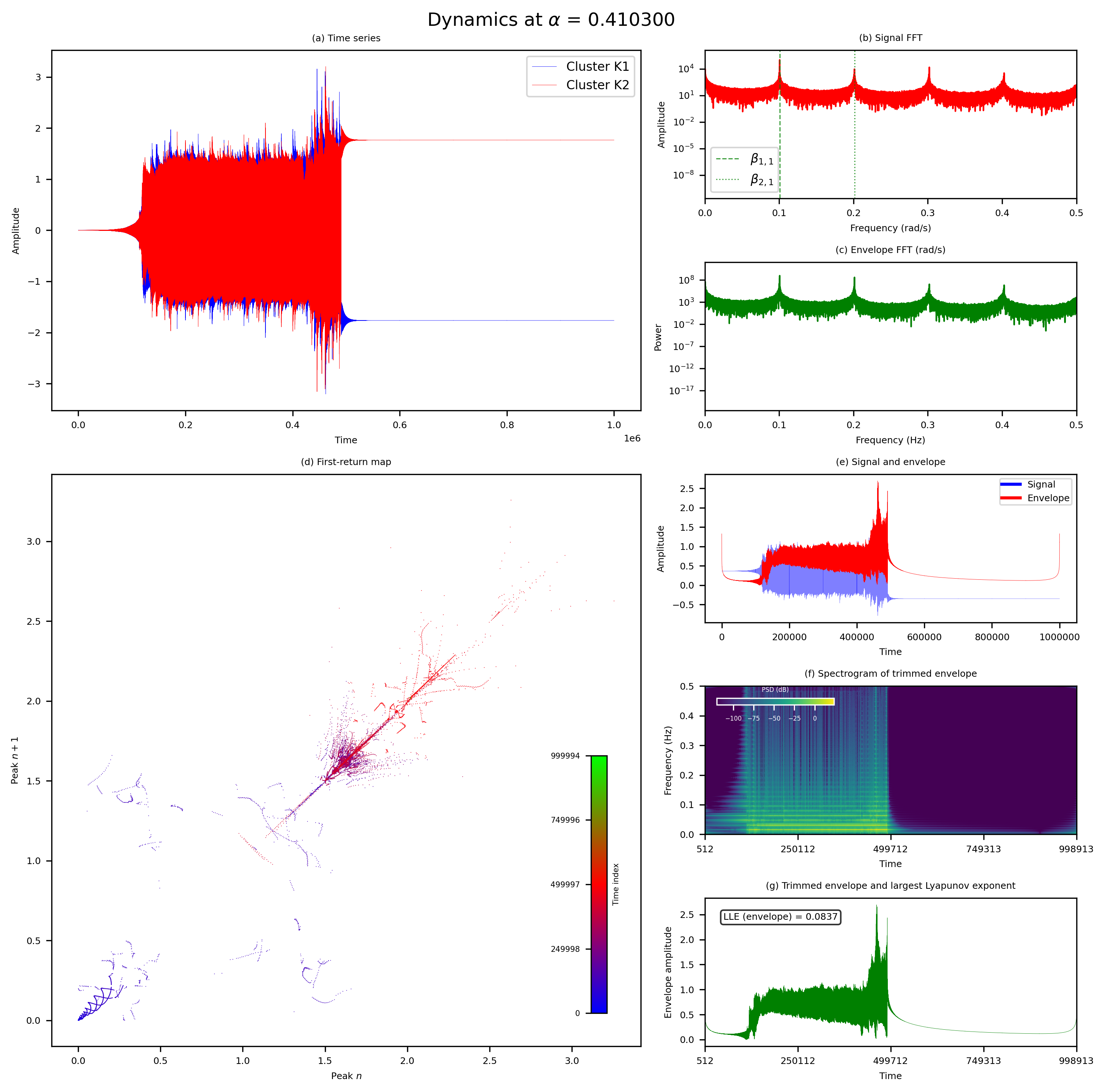}

    \caption{Dynamics at $\alpha = 0.4103$, showing a chaotic transient leading to steady-state bifurcation.}
    \label{fig:neutral:trajectory-0.4103}
\end{figure}

\begin{figure}[tbp]
    \centering

    \includegraphics[width=\textwidth]{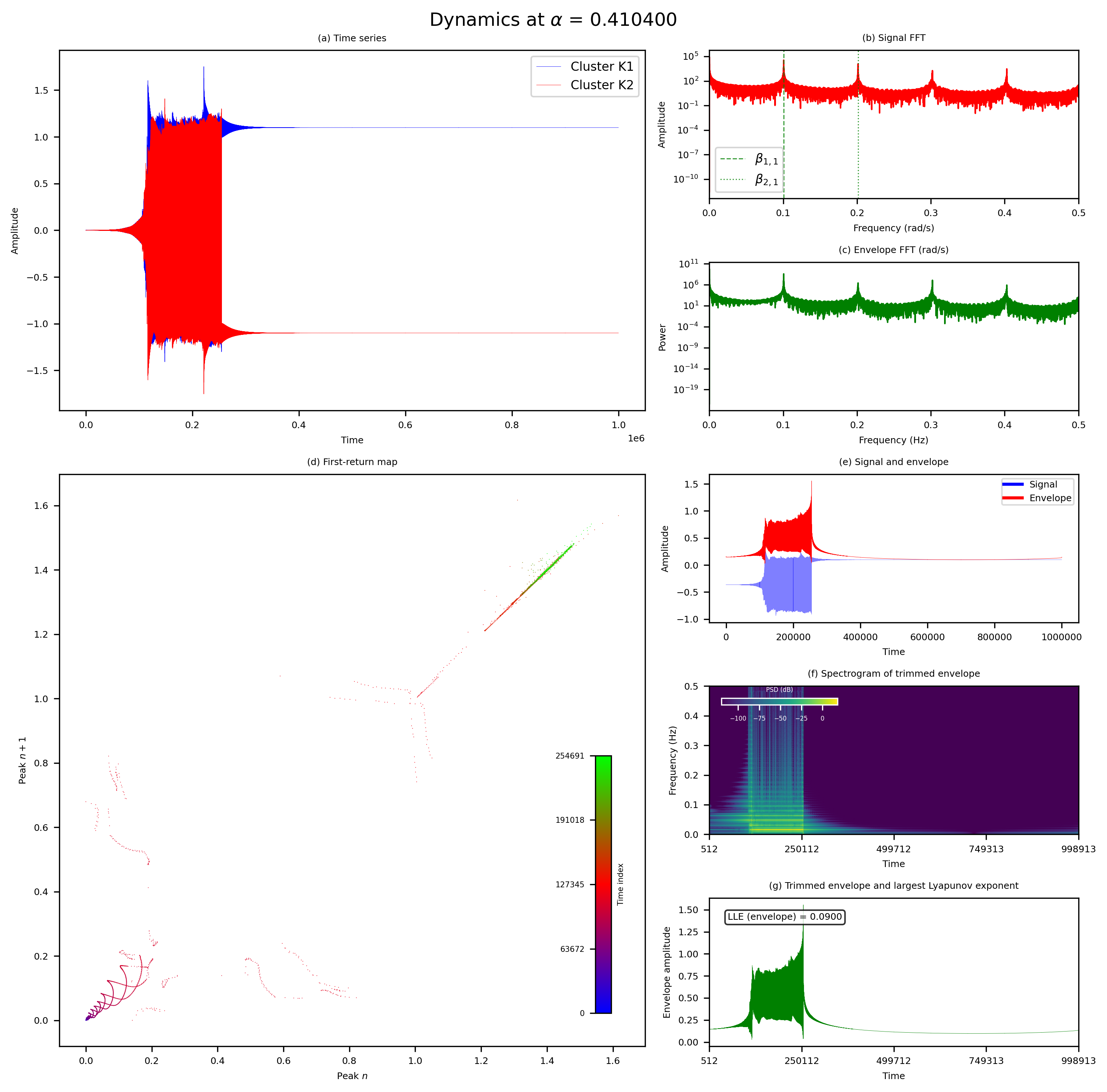}

    \caption{Dynamics at $\alpha = 0.4104$, showing a chaotic transient leading to steady-state bifurcation.}
    \label{fig:neutral:trajectory-0.4104}
\end{figure}
\backmatter
\bibliography{bibliography}

@book{Hale1977,
  author    = {Hale, J. K.},
  title     = {Theory of Functional Differential Equations},
  series    = {Applied Mathematical Sciences},
  volume    = {3},
  publisher = {Springer-Verlag},
  address   = {New York},
  year      = {1977}
}

@article{Hutchinson1948,
  author  = {Hutchinson, G. E.},
  title   = {Circular causal systems in ecology},
  journal = {Annals of the New York Academy of Sciences},
  volume  = {50},
  pages   = {221--246},
  year    = {1948}
}

@book{Krasnoselskii1956,
  author    = {Krasnosel'skii, M. A.},
  title     = {Topological Methods in the Theory of Nonlinear Integral Equations},
  year      = {1956},
  publisher = {Gosudarstv. Izdat. Tehn.-Teor. Lit.},
  address   = {Moscow},
  language  = {Russian},
  note      = {English translation: Pergamon Press, 1964}
}

@article{Rabinowitz1971,
  author  = {Rabinowitz, Paul H.},
  title   = {Some global results for nonlinear eigenvalue problems},
  journal = {Journal of Functional Analysis},
  volume  = {7},
  number  = {3},
  pages   = {487--513},
  year    = {1971}
}

@incollection{Rabinowitz1985,
  author    = {Rabinowitz, Paul H.},
  title     = {Global aspects of bifurcation},
  booktitle = {Topological Methods in Bifurcation Theory},
  series    = {S\'{e}minaire de Math\'{e}matiques Sup\'{e}rieures},
  volume    = {91},
  pages     = {63--112},
  publisher = {Presses de l'Universit\'{e} de Montr\'{e}al},
  address   = {Montr\'{e}al},
  year      = {1985}
}

@book{BalanovEtAl2025,
  author    = {Balanov, Z. and Krawcewicz, W. and Rachinskii, D. and Yu, J. and Wu, H.},
  title     = {Degree Theory and Symmetric Equations Assisted by GAP System},
  publisher = {American Mathematical Society},
  year      = {2025}
}

@book{KrawcewiczWu1997,
  author    = {Krawcewicz, W. and Wu, J.},
  title     = {Theory of Degrees with Applications to Bifurcations and Differential Equations},
  publisher = {John Wiley \& Sons},
  year      = {1997}
}

@article{Erbe1993,
  author  = {Erbe, L. H. and Krawcewicz, W. and Wu, J. H.},
  title   = {A composite coincidence degree with applications to boundary value problems of neutral equations},
  journal = {Transactions of the American Mathematical Society},
  volume  = {335},
  number  = {2},
  pages   = {459--478},
  year    = {1993}
}

@article{Ize1989,
  author  = {Ize, J. and Massabò, I. and Vignoli, A.},
  title   = {Degree theory for equivariant maps, I},
  journal = {Transactions of the American Mathematical Society},
  volume  = {315},
  number  = {2},
  pages   = {433--510},
  year    = {1989}
}

@article{Ize1992,
  author  = {Ize, J. and Massabò, I. and Vignoli, A.},
  title   = {Degree theory for equivariant maps, the general $S^1$-action},
  journal = {Memoirs of the American Mathematical Society},
  volume  = {100},
  number  = {481},
  year    = {1992}
}

@article{BalanovKrawcewiczRuan2006a,
  author    = {Balanov, Z. and Krawcewicz, W. and Ruan, H.},
  title     = {Applied equivariant degree. {I}: An axiomatic approach to primary degree},
  journal   = {Discrete and Continuous Dynamical Systems},
  volume    = {15},
  number    = {3},
  pages     = {983--1016},
  year      = {2006},
  zbl       = {1115.58013},
  doi       = {10.3934/dcds.2006.15.983},
  url       = {https://doi.org/10.3934/dcds.2006.15.983}
}

@book{BalanovKrawcewiczSteinlein2006,
  author    = {Balanov, Z. and Krawcewicz, W. and Steinlein, H.},
  title     = {Applied Equivariant Degree},
  series    = {AIMS Series on Differential Equations \& Dynamical Systems},
  volume    = {1},
  publisher = {American Institute of Mathematical Sciences},
  year      = {2006}
}

@article{Wang2024,
  author  = {Wang, H. and Han, Q.-L.},
  title   = {Distribution of roots of quasi-polynomials of neutral type and its application. II: Consensus protocol design of multi-agent systems using delayed state information},
  journal = {IEEE Transactions on Automatic Control},
  volume  = {69},
  number  = {6},
  pages   = {4058--4065},
  year    = {2024}
}

@book{Keller2011,
  author    = {Keller, A. A.},
  title     = {Time-Delay Systems: With Applications to Economic Dynamics and Control},
  publisher = {Lambert Academic Publishing},
  year      = {2011},
}

@incollection{Hale2006History,
  author    = {Hale, J. K.},
  title     = {History of Delay Equations},
  booktitle = {Delay Differential Equations and Applications},
  editor    = {Arino, O. and Hbid, M. L. and Ait Dads, E.},
  series    = {NATO Science Series II: Mathematics, Physics and Chemistry},
  volume    = {205},
  pages     = {1--28},
  publisher = {Springer},
  address   = {Dordrecht},
  year      = {2006},
}

@article{FariaMagalhaes1995,
  author  = {Faria, T. and Magalh{\~a}es, L. T.},
  title   = {Normal forms for retarded functional differential equations with parameters and applications to Hopf bifurcation},
  journal = {Journal of Differential Equations},
  volume  = {122},
  number  = {2},
  pages   = {181--200},
  year    = {1995},
}

@article{Wearing2005,
  author  = {Wearing, Helen J. and Rohani, Pejman and Keeling, Matt J.},
  title   = {Appropriate models for the management of infectious diseases},
  journal = {PLoS Medicine},
  volume  = {2},
  number  = {7},
  pages   = {e174},
  year    = {2005},
}

@article{Hutt2005,
  author  = {Hutt, Axel and Atay, Fatihcan M.},
  title   = {Analysis of nonlocal neural fields for both general and gamma-distributed connectivities},
  journal = {Physica D: Nonlinear Phenomena},
  volume  = {203},
  number  = {1–2},
  pages   = {30--54},
  year    = {2005},
}

@article{Hespanha2007,
  author  = {Hespanha, Jo{\~a}o P. and Naghshtabrizi, Payam and Xu, Yonggang},
  title   = {A survey of recent results in networked control systems},
  journal = {Proceedings of the IEEE},
  volume  = {95},
  number  = {1},
  pages   = {138--162},
  year    = {2007},
}

@article{Sharma2013,
  author    = {Sharma, S. K. and Karmeshu},
  title     = {Neuronal model with distributed delay: {E}mergence of unimodal and bimodal {ISI} distributions},
  journal   = {IEEE Transactions on NanoBioscience},
  volume    = {12},
  number    = {1},
  pages     = {1--12},
  year      = {2013},
}

@article{Yu2016,
  author    = {Yu, Jinchen and Peng, Mingshu},
  title     = {Stability and bifurcation analysis for the {K}aldor-{K}alecki model with a discrete delay and a distributed delay},
  journal   = {Physica A: Statistical Mechanics and its Applications},
  volume    = {460},
  pages     = {66--75},
  year      = {2016},
}

@book{Tarasov2021,
  author    = {Tarasov, Vasily E. and Tarasova, Valentina V.},
  title     = {Economic Dynamics with Memory: Fractional Calculus Approach},
  publisher = {Walter de Gruyter GmbH},
  year      = {2021},
  address   = {Berlin/Boston},
  isbn      = {9783110624601}
}

@inproceedings{Sipahi2015,
  author    = {Sipahi, Rifat and Atay, Fatihcan M. and Niculescu, Silviu-Iulian},
  title     = {Stability analysis of a constant time-headway driving strategy with driver memory effects modeled by distributed delays},
  booktitle = {IFAC-PapersOnLine},
  volume    = {28},
  number    = {12},
  pages     = {376--381},
  year      = {2015},
  note      = {12th IFAC Workshop on Time Delay Systems (TDS 2015), Ann Arbor, MI, USA}
}

@article{Cassidy2025,
author = {Cassidy, Tyler},
title = {Using Multidelay Discrete Delay Differential Equations to Accurately Simulate Models With Distributed Delays},
journal = {Studies in Applied Mathematics},
volume = {154},
number = {6},
pages = {e70069},
year = {2025}
}

@article{tavakoli2020multi,
  title={Multi-delay complexity collapse},
  author={Tavakoli, S Kamyar and Longtin, Andr{\'e}},
  journal={Physical Review Research},
  volume={2},
  number={3},
  pages={033485},
  year={2020},
  publisher={APS}
}

@article{Brouwer1912,
  author    = {L. E. J. Brouwer},
  title     = {\"Uber Abbildung von Mannigfaltigkeiten},
  journal   = {Mathematische Annalen},
  volume    = {71},
  pages     = {97--115},
  year      = {1912},
}

@article{LeraySchauder1934,
  author    = {Jean Leray and Jules Schauder},
  title     = {Topologie et \'equations fonctionnelles},
  journal   = {Annales Scientifiques de l'\'Ecole Normale Sup\'erieure},
  series    = {3},
  volume    = {51},
  pages     = {45--78},
  year      = {1934},
}

@incollection{Dylawerski1988,
  author    = {Grzegorz Dylawerski},
  title     = {An {$S^1$}-degree and {$S^1$}-maps between representation spheres},
  booktitle = {Algebraic Topology and Transformation Groups},
  editor    = {Tammo tom Dieck},
  series    = {Lecture Notes in Mathematics},
  volume    = {1361},
  pages     = {14--28},
  publisher = {Springer},
  address   = {Berlin},
  year      = {1988},
}

@article{munz2011consensus,
  author  = {Ulrich M{\"u}nz and Antonis Papachristodoulou and Frank Allg{\"o}wer},
  title   = {Consensus in multi-agent systems with communication delays: A survey},
  journal = {Annual Reviews in Control},
  volume  = {35},
  number  = {1},
  pages   = {77--97},
  year    = {2011}
}

@inbook{2013guo,
  title = {Introduction to Functional Differential Equations},
  booktitle = {Bifurcation Theory of Functional Differential Equations},
  publisher = {Springer, New York},
  author = {Guo, S. and Wu, J.},
  year = {2013},
  pages = {41-60},
}

@article{jitcxde,
	author = {Ansmann, G.},
	title = {Efficiently and easily integrating differential equations with {JiTCODE}, {JiTCDDE}, and {JiTCSDE}},
	year = {2018},
	journal = {Chaos},
	volume = {28},
	number = {4},
	pages = {043116},
	doi = {10.1063/1.5019320},
}

@article{vicsek1995novel,
  author  = {Vicsek, T. and Czirók, A. and Ben-Jacob, E. and Cohen, I. and Shochet, O.},
  title   = {Novel type of phase transition in a system of self-driven particles},
  journal = {Physical Review Letters},
  volume  = {75},
  pages   = {1226--1229},
  year    = {1995}
}

@article{reynolds1987boid,
  author  = {Reynolds, C. W.},
  title   = {Flocks, herds and schools: A distributed behavioral model},
  journal = {SIGGRAPH Computer Graphics},
  volume  = {21},
  pages   = {25--34},
  year    = {1987}
}

@article{saber2004consensus,
  author  = {Olfati-Saber, R. and Murray, R. M.},
  title   = {Consensus problems in networks of agents with switching topology and time-delays},
  journal = {IEEE Transactions on Automatic Control},
  volume  = {49},
  pages   = {1520--1533},
  year    = {2004}
}

@article{jadbabaie2003coordination,
  author  = {Jadbabaie, A. and Lin, J. and Morse, A. S.},
  title   = {Coordination of groups of mobile autonomous agents using nearest neighbor rules},
  journal = {IEEE Transactions on Automatic Control},
  volume  = {48},
  pages   = {988--1001},
  year    = {2003}
}

@article{ren2005consensus,
  author  = {Ren, W. and Beard, R. W.},
  title   = {Consensus seeking in multiagent systems under dynamically changing interaction topologies},
  journal = {IEEE Transactions on Automatic Control},
  volume  = {50},
  pages   = {655--661},
  year    = {2005}
}

@article{chen2014multiconsensus,
  author  = {Chen, J. and Chi, M. and Guan, Z.-H. and Liao, R.-Q. and Zhang, D.-X.},
  title   = {Multiconsensus of second-order multiagent systems with input delays},
  journal = {Mathematical Problems in Engineering},
  volume  = {2014},
  pages   = {424537},
  year    = {2014}
}

@article{xie2015second,
  author  = {Xie, D. and Xu, S. and Li, Z. and Zou, Y.},
  title   = {Second-order group consensus for multi-agent systems with time delays},
  journal = {Neurocomputing},
  volume  = {153},
  pages   = {133--139},
  year    = {2015}
}

@article{wang2021global,
  author  = {Wang, Y. and Cao, J. and Lu, B. and Cheng, Z.},
  title   = {Global asymptotic consensus of multi-agent internet congestion control system},
  journal = {Neurocomputing},
  volume  = {446},
  pages   = {50--64},
  year    = {2021}
}

@article{krawcewiczwuxia1993,
  author  = {Krawcewicz, W. and Wu, J. and Xia, H.},
  title   = {Global Hopf bifurcation theory for condensing fields and neutral equations with applications to lossless transmission problems},
  journal = {Canadian Applied Mathematics Quarterly},
  volume  = {1},
  pages   = {167--220},
  year    = {1993}
}

@article{sadovskii1971,
  author  = {Sadovskii, B. N.},
  title   = {Applications of topological methods to the theory of periodic solutions of nonlinear differential-operator equations of neutral type},
  journal = {Soviet Mathematics Doklady},
  volume  = {12},
  pages   = {1543--1547},
  year    = {1971}
}

@article{darbo1955,
  author  = {Darbo, G.},
  title   = {Punti uniti in trasformazioni a codominio non compatto},
  journal = {Rendiconti del Seminario Matematico della Università di Padova},
  volume  = {24},
  pages   = {84--92},
  year    = {1955}
}

@article{nussbaum1972,
  author  = {Nussbaum, R. D.},
  title   = {Degree theory for local condensing maps},
  journal = {Journal of Mathematical Analysis and Applications},
  volume  = {37},
  pages   = {741--766},
  year    = {1972}
}

@book{izevignoli2003,
  author    = {Ize, J. and Vignoli, A.},
  title     = {Equivariant Degree Theory},
  publisher = {de Gruyter},
  address   = {Berlin},
  year      = {2003}
}

@incollection{balanov2008symmetric,
  author    = {Balanov, Z. and Krawcewicz, W.},
  title     = {Symmetric Hopf bifurcation: Twisted degree approach},
  booktitle = {Handbook of Differential Equations: Ordinary Differential Equations},
  editor    = {Battelli, F. and Fečkan, M.},
  volume    = {4},
  pages     = {1--131},
  publisher = {Elsevier/North Holland},
  address   = {Amsterdam},
  year      = {2008}
}

@article{geba1994equivariant,
  author  = {Gęba, K. and Krawcewicz, W. and Wu, J.},
  title   = {An equivariant degree with applications to symmetric bifurcation problems. Part 1: Construction of the degree},
  journal = {Proceedings of the London Mathematical Society},
  volume  = {69},
  number  = {2},
  pages   = {377--398},
  year    = {1994}
}

@article{dylawerski1991s1,
  author  = {Dylawerski, G. and Gęba, K. and Jodel, J. and Marzantowicz, W.},
  title   = {An $S^1$-equivariant degree and the Fuller index},
  journal = {Annales Polonici Mathematici},
  volume  = {52},
  number  = {3},
  pages   = {243--280},
  year    = {1991}
}

@article{crane2026,
  author  = {Chen, C. and Crane, C. and Hensley, T.},
  title   = {Global Hopf bifurcation in symmetric configurations of distributed delay differential equations},
  journal = {Communications on Pure and Applied Analysis},
  volume  = {27},
  pages   = {170--194},
  year    = {2026}
}

@software{scholzel_nolds,
  author       = {Schölzel, Christopher},
  title        = {Nonlinear measures for dynamical systems},
  month        = jun,
  year         = 2019,
  publisher    = {Zenodo},
  version      = {0.5.2},
  doi          = {10.5281/zenodo.3814723},
  url          = {https://doi.org/10.5281/zenodo.3814723},
}

@article{alexander1992riddled,
  author  = {J. C. Alexander and Celso Grebogi and Edward Ott and James A. Yorke},
  title   = {Riddled basins},
  journal = {International Journal of Bifurcation and Chaos},
  volume  = {2},
  number  = {4},
  pages   = {795--812},
  year    = {1992}
}

@article{beja1980price,
  title={On the dynamic behavior of prices in disequilibrium},
  author={Beja, A. and Goldman, M. B.},
  journal={Journal of Finance},
  volume={35},
  number={2},
  pages={235--248},
  year={1980}
}

@article{chiarella1992dynamics,
  title={The dynamics of speculative behaviour},
  author={Chiarella, C.},
  journal={Annals of Operations Research},
  volume={37},
  pages={253--287},
  year={1992}
}

@article{dibeh2007dynamics,
  title={Dynamics of a speculative market with time delays},
  author={Dibeh, G.},
  journal={Chaos, Solitons \& Fractals},
  volume={31},
  number={2},
  pages={511--518},
  year={2007}
}

@article{he2010dynamics,
  title={Dynamics of moving average rules in a continuous-time financial market model},
  author={He, X. Z. and Zheng, M.},
  journal={Journal of Economic Behavior \& Organization},
  volume={76},
  number={3},
  pages={615--634},
  year={2010}
}

@article{matsumoto2016heterogeneous,
  title={A heterogeneous agent model of asset price with three time delays},
  author={Matsumoto, A. and Szidarovszky, F.},
  journal={Frontiers in Applied Mathematics and Statistics},
  volume={2},
  pages={15},
  year={2016}
}

@article{dobrescu2016asset,
  title={Asset price dynamics in a chartist-fundamentalist model with time delays: A bifurcation analysis},
  author={Dobrescu, L. I. and Neamtu, M. and Mircea, G.},
  journal={Discrete Dynamics in Nature and Society},
  volume={2016},
  article={4907468},
  year={2016}
}

@article{guerrini2014heterogeneous,
  title={Heterogeneous fundamentalists in a continuous time model with delays},
  author={Guerrini, L. and Sodini, M. and Szidarovszky, F.},
  journal={Discrete Dynamics in Nature and Society},
  volume={2014},
  article={959514},
  year={2014}
}

@article{chossat1988symmetry,
  author  = {Pascal Chossat and Martin Golubitsky},
  title   = {Symmetry-increasing bifurcation of chaotic attractors},
  journal = {Physica D},
  volume  = {32},
  number  = {3},
  pages   = {423--436},
  year    = {1988}
}

@article{ruelle1971nature,
  author  = {David Ruelle and Floris Takens},
  title   = {On the nature of turbulence},
  journal = {Communications in Mathematical Physics},
  volume  = {20},
  pages   = {167--192},
  year    = {1971}
}

@article{rackauckas2017differentialequations,
  title={DifferentialEquations.jl – A Performant and Feature-Rich Ecosystem for Solving Differential Equations in Julia},
  author={Rackauckas, Christopher and Nie, Qing},
  journal={Journal of Open Research Software},
  volume={5},
  number={1},
  pages={15},
  year={2017},
  publisher={Ubiquity Press},
  doi={10.5334/jors.151}
}

@manual{GAP4,
  organization = {The GAP Group},
  title = {{GAP} -- Groups, Algorithms, and Programming, Version 4.13.1},
  year = {2024},
  url = {https://www.gap-system.org},
}

@manual{equideg,
	author = {Wu, H.},
	title = {GAP Equivariant Degree Library},
	year = {2024},
	organization = {Github},
    url = {https://github.com/psistwu/equideg},
    note = {GAP package},
}

@article{rosenstein1993,
title = {A practical method for calculating largest Lyapunov exponents from small data sets},
journal = {Physica D: Nonlinear Phenomena},
volume = {65},
number = {1},
pages = {117-134},
year = {1993},
issn = {0167-2789},
author = {Michael T. Rosenstein and James J. Collins and Carlo J. {De Luca}}
}

\end{document}